\documentclass{acta_2022}

\usepackage[longnamesfirst]{natbib} \usepackage{har2nat}
\setcitestyle{comma,aysep={}}
\shortcites{samplepreprint} 

\usepackage{amssymb,amsmath}
\usepackage{newtxtext}
\usepackage{newtxmath}
\DeclareSymbolFont{CMlargesymbols}{OMX}{cmex}{m}{n}
\DeclareMathDelimiter{(}{\mathopen} {operators}{"28}{CMlargesymbols}{"00}
\DeclareMathDelimiter{)}{\mathclose}{operators}{"29}{CMlargesymbols}{"01}
\DeclareMathAlphabet\mathcal{OMS}{cmsy}{m}{n}
\SetMathAlphabet\mathcal{bold}{OMS}{cmsy}{b}{n}
\pagerange{\pageref{firstpage}--\pageref{lastpage}}
\allowdisplaybreaks
\doi{XXXXXXXX}  
\usepackage[UKenglish]{babel} 
\newcommand{\ignore}[1]{}
\usepackage[twoside,
	paperheight=247mm,
	paperwidth=174mm,
	marginratio={3:5,2:3},
	left=18mm,
	top=20mm]{geometry}
\usepackage[hang]{subfigure}
\numberwithin{figure}{section}
\numberwithin{table}{section}
\usepackage{graphics,graphicx}
\usepackage[cmyk]{xcolor}
\usepackage{enumitem} 
\usepackage{color}
\usepackage{shuffle}

\newtheorem{theorem}{Theorem}[section]

\newtheorem{algorithm}[theorem]{Algorithm}

\newcommand{\Z}{\mathbb{Z}}
\newcommand{\C}{\mathbb{C}}

\newcommand{\e}{\ensuremath{\mathrm{e}}}
\newcommand{\tr}{\mathrm{tr}}

\newcommand{\W}{\mathcal{W}}
\newcommand{\M}{\mathcal{M}}
\newcommand{\cinfM}{C^{\infty}(\M,\mathbb{R})}

\newcommand{\blue}[1]{{\color{blue} #1}}

\begin{document}

\title[Splitting methods for differential equations]{Splitting methods for differential equations}

\author[S. Blanes, F. Casas, and A. Murua]{Sergio Blanes\\
Instituto de Matem\'atica Multidisciplinar, \\ Universitat Polit\`ecnica de Val\`encia,\\ 46022-Valencia, Spain\\
E-mail: {\sf serblaza@imm.upv.es}\\
\and
Fernando Casas\\
Departament de Matem\`atiques and IMAC, Universitat Jaume I,\\ 12071-Castell\'on, Spain\\
E-mail: {\sf Fernando.Casas@uji.es}\\
\and
Ander Murua\\
Konputazio Zientziak eta A.A. Saila, Informatika Fakultatea, EHU/UPV,\\ Donostia/San Sebasti\'an, Spain\\
E-mail: {\sf Ander.Murua@ehu.es}
}

\maketitle
\label{firstpage}

\vspace*{2cm}

\begin{center}
\today
\end{center}
\begin{abstract}

This overview is devoted to splitting methods, a class of numerical integrators intended for differential equations that can be subdivided into
different problems easier to solve than the original system. Closely connected with this class of integrators are composition methods,
in which one or several low-order schemes are composed to construct higher-order numerical approximations to the exact solution. We analyze
in detail the order conditions that have to be satisfied by these classes of methods to achieve a given order, and provide some insight about their qualitative
properties in connection with geometric numerical integration and the treatment of highly oscillatory problems. 
Since splitting methods have received considerable attention in the realm of partial differential
equations, we also cover this subject in the present survey, with special attention to parabolic equations and their problems. An exhaustive list of methods
of different orders is collected and tested on simple examples. Finally, some applications of splitting methods in different
areas, ranging from celestial mechanics to statistics, are also provided.

\end{abstract}

\tableofcontents 

\vspace{2mm}

\section{Introduction}
\label{sect1}

\subsection{Lie--Trotter and Strang methods}
\label{subsec.1.1}

If, as has sometimes been argued, there are only ten big ideas in numerical analysis and all the rest are merely variations on those themes, \emph{splitting}
is undoubtedly one of them \cite{macnamara16os}. In fact, one could say that, at least since Descartes stated his four rules of logic in the \emph{Discours de la m\'ethode}\footnote{\emph{Le second, de diviser chacune des difficult\'es que j'examinarais, en autant de parcelles qu'il se pourrait, et qu'il serait requis pour les mieux r\'esoudre. Le troisi\`eme, de conduire par ordre mes pens\'ees, en commen\c{c}ant par les objets les plus simples et les plus ais\'es \`a conna\^{i}tre, pour monter
peu \`a peu, comme par degr\'es jusques \`a la connaissance des plus compos\'es; et supposant m\^{e}me de l'ordre entre ceux qui ne se pr\'ec\`edent point 
naturellement les uns les autres.} (Ren\'e Descartes, \emph{Discours de la m\'ethode}, Seconde partie).}, 
the notion of sub-dividing a complicated problem into its simpler constituent parts, solving each one of them separately
and combining those separated solutions in a controlled way to get a solution to the original overall problem constitutes a guiding principle in all areas of
science and philosophy. 

In the realm of numerical analysis of differential equations, this basic principle can be stated as follows. Suppose one has the abstract initial value problem 
\begin{equation} \label{ivp.1}
  x' \equiv \frac{d x}{dt} = f(x), \qquad x(0) = x_0
\end{equation}
associated to an ordinary differential equation (ODE) or a partial differential equation (PDE),
in which case $f$ is a certain spatial partial differential operator. Furthermore, suppose that $f$ does not depend explicitly on time and can be decomposed as 
\begin{equation} \label{ivp.2}
  f = f_1 + \cdots + f_m, \qquad m \ge 2,
\end{equation}
so that each initial value problem $x^{\prime} = f_j(x)$, $x(0) = x_0$ is easier to
solve than (\ref{ivp.1}). Most commonly, they can be integrated exactly in closed form. Then, it is possible to take advantage of the decomposition (\ref{ivp.2}) to get accurate approximations to the solution of (\ref{ivp.1})
by means of \emph{splitting methods}.

The best example to start with is perhaps the case of a linear differential equation defined in $\mathbb{R}^D$ and $m=2$, namely
\begin{equation} \label{ivp.3}
   x' = f_1(x) + f_2(x) = F_1 \, x+ F_2 \, x, \qquad x(0) = x_0,
\end{equation}
where $F_1$ and $F_2$ are $D \times D$ matrices and $x \in \mathbb{R}^D$. 
The solution reads 
\[
   x(t) = \e^{t (F_1 + F_2)} x_0,
\]     
so that by computing this matrix exponential directly we would have solved (\ref{ivp.3}) without requiring splitting methods. Associated with (\ref{ivp.3})
one has the matrix differential equation
\begin{equation} \label{ivp.3b}
   \frac{d X}{dt} = (F_1 + F_2) X, \qquad X(0) = I,
\end{equation}
in the sense that 
\[
  x(t) = X(t) \, x_0 \qquad \mbox{ and } \qquad X(t) = \e^{t (F_1 + F_2)}.
\]  
This is useful (mainly) for theoretical purposes, since one usually tries to compute $\e^{t (F_1 + F_2)} x_0$ directly instead of first computing the matrix
exponential and then multiplying it by $x_0$.

It often happens
that evaluating the action of $X(t)$ on $x_0$
is difficult and/or computationally expensive. If, however, this is not the case for each $\e^{t F_j}$ separately, then 
one may use the well known Lie product formula \cite[p. 295]{reed80momI}
\begin{equation} \label{ivp.4}
  \e^{t(F_1 + F_2)} = \lim_{n \rightarrow \infty} \left( \e^{\frac{t}{n} F_2}  \, \e^{\frac{t}{n} F_1 } \right)^n.
\end{equation}
To get an approximate solution of (\ref{ivp.3}) at the final time $t = t_f$, one subdivides the interval $[0, t_f]$ into $N$ steps of length $h$, with $N h = t_f$, and
computes the sequence
\begin{equation} \label{ivp.5}
   x_{n+1} = \e^{h F_2} \, \e^{h F_1} \, x_n, \qquad n \ge 0,
\end{equation}  
so that $x_{n+1} \approx x(t_{n+1} = (n+1) h)$. This is the so-called \emph{Lie--Trotter scheme}.
When the matrices commute, the sequence produces the exact solution. To put it another way, if the commutator $[F_1, F_2] \equiv F_1 F_2 - F_2 F_1 = 0$,
then $\exp(h (F_1 + F_2)) = \exp(h F_2) \exp(h F_1)$. Otherwise, a direct calculation shows that
\[
  \e^{h (F_1 + F_2)} - \e^{h F_2} \e^{h F_1} = \frac{1}{2} h^2 [F_1, F_2] + \mathcal{O}(h^3)
\]
as $h \rightarrow 0$, and hence the previous approximation is only of first order of accuracy. 
Another version of the method is possible of course by reversing the order of $F_1$ and $F_2$, namely
\begin{equation} \label{ivp.5b}
   x_{n+1} = \e^{h F_1} \, \e^{h F_2} \, x_n, \qquad n \ge 0
\end{equation}  
has the same order of accuracy and properties as (\ref{ivp.5}). Needless to say, for any $m > 2$, it results in
\[
   x_{n+1} = \e^{h F_m} \,  \cdots \, \e^{h F_2} \,  \e^{h F_1} \, x_n
\]
(or any other permutation of the matrices $F_j$).

A higher order approximation can be achieved by considering a symmetrized version of (\ref{ivp.5}), 
\begin{equation} \label{bo2}
  \e^{h (F_1 + F_2)} - \e^{\frac{1}{2} h F_1} \e^{h F_2} \e^{\frac{1}{2} h F_1} = C h^3+ \mathcal{O}(h^4),
\end{equation}
where the constant $C$ can be obtained either by comparing Taylor series or by applying the Baker--Campbell--Hausdorff (BCH) formula as
\cite{varadarajan84lgl}
$C = \frac{1}{24} ([F_1,[F_1,F_2]] + 2 [F_2,[F_1,F_2]])$. Therefore, the sequence
\begin{equation} \label{ivp.6}
   x_{n+1} = \e^{\frac{1}{2} h F_1} \e^{h F_2} \e^{\frac{1}{2} h F_1} \, x_n, \qquad n \ge 0
\end{equation}  
produces a second-order approximation 
for the solution of (\ref{ivp.3}). 
This corresponds to the \emph{Strang} (splitting) \emph{scheme}.
Again, if the role of $F_1$ and $F_2$ is interchanged, one has another version of the Strang splitting scheme.

Simple generalizations to
the case $m > 2$ include in particular the product
\[
   \e^{\frac{1}{2} h F_1} \,  \e^{\frac{1}{2} h F_2} \, \cdots  \, \e^{ h F_m} \, \cdots \, \e^{\frac{1}{2} h F_2} \, \e^{\frac{1}{2} h F_1}.
\]
Higher order splitting methods could in principle be constructed by including more exponentials with their corresponding coefficients in a time step, namely
\begin{equation} \label{ho.1}
 \Psi(h) = \e^{a_{s+1} h F_1} \,  \e^{b_s h F_2} \, \e^{a_s h F_1} \, \cdots \, \e^{b_1 h F_2} \, \e^{a_1 h F_1}.
\end{equation}
The number $s$ as well as the
coefficients $a_j$, $b_j$ are chosen so that 
\[
  \Psi(h) = \e^{h (F_1 + F_2)} + \mathcal{O}(h^{r+1})
\]
as $h \rightarrow 0$ for a given order $r$. In (\ref{ho.1}), the first and last exponentials correspond to $F_1$. This format is convenient for implementation, since the
last exponential in one step can be concatenated with the first one at the next step, thus reducing the number of evaluations by one. This corresponds to the well known
FSAL (First Same As Last) property. 
The situation when one has an exponential of $F_2$
as the first and last term is recovered by taking $a_1= a_{s+1} =0$.

\subsection{Flows and differential operators}
\label{subsec.1.2}

The Lie--Trotter scheme can be easily generalized to any system (\ref{ivp.1})-(\ref{ivp.2}) when the solution is no longer given by
exponentials, as in the linear case. If $m=2$, it is equivalent to the following.
\begin{algorithm} (Lie--Trotter). 
\label{alg-LT}
Starting from $x_{0} = x(0)$, for $n \ge 0$
\begin{itemize}
  \item solve $y_{1}^{\prime} = f_1(y_1)$, $\ y_{1}(t_n) = x_n$, in $[t_n, t_{n+1}]$;
  \item set $y_{n+1/2} = y_{1}(t_{n+1})$;
  \item solve $y_{2}^{\prime} = f_2(y_{2})$, $\ y_{2}(t_n) = y_{n+1/2}$, in $[t_n, t_{n+1}]$;
  \item finally, set $x_{n+1} = y_{2}(t_{n+1})$.
\end{itemize}  
\end{algorithm}
Alternatively, if we denote the solution of equation (\ref{ivp.1}) for each $t \in \mathbb{R}$ as
\[
  x(t) = \varphi_t^{[f]}(x_0),
\]
then Algorithm \ref{alg-LT} can be formally expressed as
\begin{equation} \label{LT.2}
     x_{n+1} = \chi_h(x_n) \equiv \varphi_h^{[2]} \big( \varphi_h^{[1]} (x_n) \big) = \big(\varphi_h^{[2]} \circ \varphi_h^{[1]}\big) (x_n),
\end{equation}
where we have used the simplified notation $\varphi_t^{[j]}(x_0)$ for the solutions $y(t)= \varphi_t^{[f_j]}(x_0)$
 of the subproblems $y^{\prime} = f_j(y)$, $y(0) = x_0$. Analogously, the Strang splitting
(\ref{ivp.6}) is generalized as  follows.
\begin{algorithm} (Strang). From $x_{0} = x(0)$, for $n \ge 0$,
\label{alg-S}
\begin{itemize}
  \item solve $y_{1}^{\prime} = f_1(y_{1})$, $\ y_{1}(t_n) = x_n$, in $[t_n, t_{n+ \frac{1}{2}}]$, with $t_{n+ \frac{1}{2}} = (n + \frac{1}{2}) h$;
  \item set $y_{n+\frac{1}{2}} = y_{1}(t_{n+\frac{1}{2}})$;
  \item solve $y_{2}^{\prime} = f_2(y_{2})$, $\ y_{2}(t_n) = y_{n+\frac{1}{2}}$, in $[t_n, t_{n+1}]$;
  \item set $\hat{y}_{n+\frac{1}{2}} = y_{2}(t_{n+1})$;
  \item solve $y_{1}^{\prime} = f_1(y_1)$, $\ y_1(t_{n + \frac{1}{2}}) = \hat{y}_{n+\frac{1}{2}}$, in $[t_{n+ \frac{1}{2}}, t_{n+1}]$;
  \item finally, set $x_{n+1} = y_1(t_{n+1})$  
\end{itemize}
\end{algorithm}
or, in short,
\begin{equation} \label{S.2}
  x_{n+1} = S_h^{[2]} (x_n) \equiv  \big( \varphi_{h/2}^{[1]} \circ \varphi_h^{[2]} \circ \varphi_{h/2}^{[1]} \big)(x_n).
\end{equation}
If equation (\ref{ivp.1}) corresponds to a (nonlinear) ordinary differential equation evolving in $\mathbb{R}^D$,
\begin{equation} \label{ode.1}
  x' = f(x), \qquad x(0) = x_0 \in \mathbb{R}^D,
\end{equation}
$f$ is called the \emph{vector field}. If (\ref{ode.1}) admits for each $x_0 \in \mathbb{R}^D$ a unique solution $x(t)$ defined for all $t \in \mathbb{R}$, the map
\[
\begin{aligned}
  \varphi_t^{[f]}: \; \; & \mathbb{R}^D \longrightarrow \mathbb{R}^D \\
                         & x_0 \longmapsto x(t) = \varphi_t^{[f]}(x_0)
\end{aligned}
\]
is referred to as the $t$-\emph{flow} \cite{arnold89mmo}. Thus, for each value of the real parameter $t$, $\varphi_t^{[f]}$ maps
$\mathbb{R}^D$ in $\mathbb{R}^D$ in such a way that $\varphi_t^{[f]}(z)$ is the value at time $t$ of the solution of the system with initial
value $z$ at time $0$, whereas, for fixed $x_0$ and varying $t$, $\varphi_t^{[f]}(x_0)$ is the solution of the initial value problem  (\ref{ode.1}).  

It is worth mentioning that the solution $x(t)$ of (\ref{ivp.1}) is in general defined for a {\em maximal time-interval}  $(t_{\min}(x_0), t_{\max}(x_0))$ (with $-\infty \leq t_{\min}(x_0) \leq 0 \leq t_{\max}(x_0)\leq +\infty$). Furthermore, the vector field $f$ of many systems of ordinary differential equations is singular (or undefined) for some $x\in \mathbb{R}^D$. Thus, in general, $f$ is defined for some open set $\mathcal{U} \subset \mathbb{R}^D$. Hence, for a given $t \in \mathbb{R}$, the $t$-flow $\varphi_t^{[f]}$ is a map from 
\begin{equation*}
\mathcal{D}_t = \{x_0 \in \mathcal{U}\ :\ t \in (t_{\min}(x_0), t_{\max}(x_0))\}
\end{equation*}
to $\mathcal{U}$. In general, one may have different domains of definition $\mathcal{U}$ (resp. $\mathcal{D}_t$) for each $f_j$ (resp. $\varphi_t^{[f_j]}$).
In this general situation the compositions in (\ref{LT.2}) and (\ref{S.2}) are not well defined for all $x_0$, $h$ and $n$. In order to avoid these technicalities, we will assume in what follows that $\mathcal{U}=\mathbb{R}^D$, $t_{\min}=-\infty$, and $t_{\max}=+\infty$, for each vector field $f_j$.

Associated with the vector field $f$ is the Lie derivative or \emph{Lie operator} $F$ \cite{arnold89mmo}, mapping smooth functions  
$g:\mathbb{R}^D \rightarrow \mathbb{R}$ into the real-valued function $F  \, g: \mathbb{R}^D \rightarrow \mathbb{R}$ such that, for $x \in \mathbb{R}^D$,        
\[
   (F \, g) (x) = \left.\frac{d}{d t}\right|_{t=0} g \big(\varphi_{t}^{[f]}(x) \big), 
\]
that is,
\[
  (F \, g)(x) = f(x) \cdot \nabla g(x).
\]
Then, the flow of (\ref{ode.1}) verifies \cite{sanz-serna94nhp,hairer06gni}
\[
  g \big(\varphi_t^{[f]}(x) \big) = \big( \e^{t F} g \big)(x) \equiv \left( \sum_{k \ge 0} \frac{t^k}{k!} F^k g \right)(x).
\]
The operator $X(t) \equiv \e^{t F}$ is called \emph{Lie transformation}, and can be seen as the formal solution of the operator equation
\begin{equation} \label{ope.1}
  \frac{dX}{dt} = X \, F, \qquad X(0) = I.
\end{equation}
This can be seen as follows: on the one hand,
\[
  \frac{d }{dt} g \big(\varphi_t^{[f]}(x) \big) = \frac{d }{dt} \big(X(t) g \big)(x) = \left( \frac{d X(t)}{dt} g \right) (x),
\]
and on the other hand
\[
      \frac{d }{dt} g \big(\varphi_t^{[f]}(x)\big)  = (F g) \big(\varphi_t^{[f]}(x) \big) = X(t) (F g)(x).
\]      
Lie operators satisfy some remarkable properties \cite{arnold89mmo}. In particular, although they do not commute, their
commutator is nevertheless a first-order linear differential operator. Specifically, let $F$ and $G$ be the Lie operators
associated with $f$ and $g$, respectively, and $u: \mathbb{R}^D \rightarrow \mathbb{R}$ a given smooth function. Then,
\[
  [F, G] u = (F G - G F) \, u = \sum_{i,j=1}^D \left( f_j \frac{\partial g_i}{\partial x_j} - g_j \frac{\partial f_i}{\partial x_j} \right)
  \frac{\partial u}{\partial x_i}
\]
and it is possible to associate a new vector field to this differential operator, $w = (f,g)$, 
with components
\[
  w_i = (f,g)_i = \sum_{j=1}^D \left( f_j \frac{\partial g_i}{\partial x_j} - g_j \frac{\partial f_i}{\partial x_j} \right).
\]  
It is called the \emph{Lie--Poisson bracket} of $f$ and
$g$, and its Lie operator $W$ satisfies $W = [F, G]$.

Now suppose that $f(x) =  f_1(x) + f_2(x)$, so that each part $x^\prime = f_j(x)$   
is exactly solvable (or can be numerically solved up to round-off accuracy) with flow
$ x(t)=\varphi^{[j]}_{t}(x_0)$. Denoting by $F_1$ and $F_2$ the Lie operators associated with $f_1$ and $f_2$, respectively, it holds that
\[
  g \big(\varphi^{[1]}_{t}(x) \big) = \left(\e^{t F_1} g\right)(x), \qquad  g \big( \varphi^{[2]}_{t}(x) \big) = \left( \e^{t F_2} g \right)(x).
\]
Then, for the first-order approximation $\chi_h = \varphi^{[2]}_{h} \circ \varphi^{[1]}_{h}$ furnished by Algorithm \ref{alg-LT} 
one has $g \big( \chi_h(x) \big) = \big(\Psi(h) g \big)(x)$, 
where $\Psi(h)$ is a series of linear differential operators defined as
\begin{equation} \label{chi2}
  \Psi(h) = \e^{h F_1} \, \e^{h F_2}.
\end{equation}
Notice that the exponentials of Lie derivatives in (\ref{chi2}) appear in reverse order with respect to the maps in the integrator 
\cite[p. 88]{hairer06gni}. Of course, the same procedure can be applied to the Strang splitting, resulting in the product
\begin{equation} \label{slt1}
  \Psi(h) = \e^{\frac{h}{2} F_1} \, \e^{h F_2} \, \e^{\frac{h}{2} F_1}. 
\end{equation}

These considerations show that: (i) splitting methods for the problem (\ref{ivp.1})-(\ref{ivp.2}) defined in a certain function space
(as for instance with partial differential equations) can be formulated in terms of the solution of each subproblem (either exact or approximate)
by means of Algorithms \ref{alg-LT} and \ref{alg-S}, and (ii)
splitting methods applied to nonlinear ODEs evolving in $\mathbb{R}^D$ can also be formally
expressed as products of exponentials of differential operators, since it is possible to transform the original nonlinear problem into a linear one with the Lie
formalism.
This observation is very useful when analyzing the order conditions for a method to be of a
given order. In particular, one has the same order conditions for linear and nonlinear ODEs (see Section \ref{sect2}).

The previous integrators are sometimes called \emph{multiplicative operator-splitting methods}, especially in the literature concerning the
numerical treatment of partial
differential equations. In that area, one has still to specify how to solve each initial value sub-problem in Algorithms \ref{alg-LT}-\ref{alg-S} as well as the
boundary conditions. 
Moreover, one should take into account that, for a given differential equation, different ways to carry out the splitting in fact lead to different
integrators. 

\subsection{Adjoint method, conjugate method}
\label{sub1.2b}

 The flow $\varphi_t^{[f]}$ of  (\ref{ode.1}) verifies
$\big(\varphi_{-t}^{[f]}\big)^{-1} = \varphi_t^{[f]}$, but this property is not shared by many numerical integrators, and 
in particular by the map $\chi_h$ corresponding to the Lie--Trotter scheme. 

In general, if $\psi_h(x)$ represents a numerical method of order at least one, i.e., $\psi_h(x) = x + h \, f(x) +
     \mathcal{O}(h^2)$, then $(\psi_{-h})^{-1}(x) = x + h \, f(x) + \mathcal{O}(h^2)$, so that
   \[
        \psi_h^* \equiv   (\psi_{-h})^{-1}
   \]
   is also a numerical method of order at least one. It is called the \textit{adjoint method} of $\psi_h$ \cite{sanz-serna94nhp}. In words, 
  stepping forwards with the
   given method $\psi_h$ is the same as stepping backwards with the inverse of its adjoint $\psi_h^*$. If 
   $\psi_h \equiv \chi_h = \varphi_h^{[2]} \circ \varphi_h^{[1]}$, then clearly $\chi_h^* = \varphi_h^{[1]} \circ \varphi_h^{[2]}$. Additional examples are the
   explicit and implicit Euler methods:
 \[
   x_{n+1} = \psi_h^{e} (x_n) = x_n + h f(x_n), \qquad   x_{n+1} = \psi_h^{i} (x_n) = x_n + h f(x_{n+1}), 
 \]
 since $\psi_h^{i} = (\psi_h^{e})^*$.   

Whenever an integrator satisfies
   \[
       \psi_h = \psi_h^* = (\psi_{-h})^{-1},
   \]
   it is called a \textit{time-symmetric} or \textit{self-adjoint} method. Alternatively, $x_{n+1} = \psi_h(x_n)$
    is time-symmetric if and only if exchanging $h \leftrightarrow -h$ and $x_n \leftrightarrow x_{n+1}$ one gets the
    same expression, i.e., $\psi_{-h}(x_{n+1}) = x_n$. The Strang scheme (\ref{S.2}) is an example of a time-symmetric method.
    
 It is in fact straightforward to construct time-symmetric methods using the adjoint: given an arbitrary method $\psi_h$
   of order $r \ge 1$, then the compositions
\begin{equation} \label{eq.MethodAdjoint}
  \psi_{h/2} \circ \psi_{h/2}^*  \quad \mbox{ and } \quad  \psi_{h/2}^* \circ \psi_{h/2}
\end{equation}
 are time-symmetric methods of order $r \ge 2$ \cite{sanz-serna94nhp}. Moreover, 
symmetric methods are necessarily of even order, as we will show in Section \ref{sect2}. Notice that the Strang method (\ref{S.2}) is simply
\[
  S_h^{[2]} = \chi_{h/2}^* \circ \chi_{h/2},
\]
where $\chi_h$ is given by eq. (\ref{LT.2}).
Additional examples are the trapezoidal rule $\psi_h^{t} = \psi_{h/2}^{i} \circ  \psi_{h/2}^{e}$ and the midpoint rule $\psi_h^{m} = \psi_{h/2}^{e} \circ  \psi_{h/2}^{i}$.

The Strang scheme can also be expressed as
\begin{equation} \label{conju.1}
\begin{aligned}
  S_h^{[2]} & =  \varphi_{-h/2}^{[1]} \circ  (\varphi_{h}^{[1]} \circ  \varphi_{h}^{[2]}) \circ  \varphi_{h/2}^{[1]} \\
  & =  \pi_h^{-1} \circ \chi_h^* \circ \pi_h, 
\end{aligned}  
\end{equation}
with $\pi_h=  \varphi_{h/2}^{[1]}$. 
In the terminology of dynamical systems, the Strang and Lie--Trotter schemes are said to be \textit{conjugate} to each other 
by the ($\mathcal{O}(h)$-near to identity) map 
$\pi_h =  \varphi_{h/2}^{[1]}$, which can be considered as a change of coordinates. Furthermore,
the result of $n$ applications of the Strang scheme can
be recovered from $n$ applications of Lie--Trotter by carrying out just an initial transformation at the initial step and its \emph{inverse} 
at the final step.

Since many dynamical properties are invariant under changes of coordinates,
conjugate methods provide the same characterization of these properties. Other examples of conjugate methods are the trapezoidal and midpoint rules, the map
$\pi_h$ being in this case the implicit Euler method \cite{hairer06gni}. We will treat conjugate methods in detail in Sections \ref{sect4} and \ref{sect5}.

\subsection{The mathematical pendulum}
\label{subsec.1.3}

The simple mathematical pendulum constitutes a standard example of a nonlinear Hamiltonian system. In appropriate units, the corresponding
Hamiltonian function reads
\begin{equation} \label{pendulum1}
  H(q,p) = \frac{1}{2} p^2 + (1- \cos q),
\end{equation}
where $q$ denotes the angle from the vertical suspension point and $p$ is the associated momentum. 

As is well known, the equations of motion of a generic Hamiltonian system with Hamiltonian $H(q,p)$, and $q,p \in \mathbb{R}^d$,
are given by \cite{goldstein80cme}
\begin{equation} \label{eq-mo}
   \frac{dq}{dt} = \nabla_{p} H, \qquad\quad
   \frac{dp}{dt} = -\nabla_{q} H,
 \end{equation}
the function $H(q,p)$ remains constant along the evolution, and the corresponding $t$-flow, denoted as $\varphi_t^{[H]}$, 
is a \emph{symplectic} transformation \cite{arnold89mmo}: its Jacobian matrix $\varphi_t^{\prime [H]}$ verifies the identity
\[
   (\varphi_t^{\prime [H]})^T \, J \, \varphi_t^{\prime [H]} = J \qquad \mbox{ for } \; t \ge 0,
\]
where $J$ is the basic canonical matrix
\begin{equation} \label{canonJ}
   J = \left( \begin{array}{rr}
   		0_d & I_d \\
		- I_d & 0_d
		\end{array} \right).
\end{equation}
In the particular case of (\ref{pendulum1}), $d=1$, and the equations of motion are
\[
  \frac{dq}{dt} = p, \qquad\qquad \frac{dp}{dt} = - \sin q.
\]
Given a Hamiltonian $H(q,p)$ that can be decomposed as  
\begin{equation} \label{Hsplit}
  H(q,p) = H_1(q,p) + H_2(q,p),
\end{equation}
it makes sense to split the equations of motion (\ref{eq-mo}) as
\[
\frac{d}{dt} 
\left( \begin{matrix} q \\ p \end{matrix} \right) = \left( \begin{matrix} \nabla_p H_1(q,p) \\ -\nabla_q H_1(q,p)\end{matrix} \right)  + 
\left( \begin{matrix}\nabla_p H_2(q,p) \\ -\nabla_q H_2(q,p)\end{matrix}\right), 
\end{equation*}
so that each subsystem is itself Hamiltonian. In that case,
we can then apply Algorithm \ref{alg-LT} by composing both maps and form the 
first-order scheme
\begin{equation}   \label{LieTrotterHam}
 x_{n+1} = \chi_h(x_n) \equiv  \big( \varphi_{h}^{[H_2]} \circ \varphi_h^{[H_1]}\big)(x_n), \quad
  n=0,1,2,\ldots,
\end{equation}
where $x_n = (q_n, p_n)^T$. 
Similarly, Algorithm \ref{alg-S} gives the second-order scheme
\begin{equation} \label{StrangHam}
  x_{n+1} = S_h^{[2]} (x_n) \equiv  \big( \varphi_{h/2}^{[H_1]} \circ \varphi_h^{[H_2]} \circ \varphi_{h/2}^{[H_1]}\big)(x_n), \quad
  n=0,1,2,\ldots.
\end{equation}
Notice that, since both $\chi_h$ and  $S_h^{[2]}$  are defined as compositions of flows of Hamiltonian systems and the composition of symplectic maps is also symplectic
\cite{arnold89mmo}, then, both  (\ref{LieTrotterHam}) and (\ref{StrangHam}) are \emph{symplectic integrators} \cite{sanz-serna94nhp}.

The fact that schemes (\ref{LieTrotterHam}) and (\ref{StrangHam}) share the symplectic property with the exact flow  has
 remarkable consequences in practice
concerning the preservation of properties and the error propagation on long time integrations, as we will shortly illustrate.

For the Hamiltonian (\ref{pendulum1}) describing the pendulum (and in fact, for many other mechanical systems), 
one can separate the contributions of the kinetic energy
$T(p) = \frac{1}{2} p^2$ and the potential energy $V(q) = 1- \cos q$, so that a natural splitting of the form (\ref{Hsplit}) is then
\begin{equation} \label{eq:Hsep}
  H(q,p) = T(p) + V(q).
\end{equation}
 This corresponds to splitting the equations of motion (\ref{eq-mo}) into the subsystems
\begin{equation}   \label{eq.fsp0}
    \left( \begin{array}{c}
       q'  \\
       p' \end{array} \right) = 
      \left( \begin{array}{c}
         \nabla_p T(p) \\
         0 \end{array} \right)
    \qquad\quad   \mbox{ and } \qquad\quad
    \left( \begin{array}{c}
       q'  \\
       p' \end{array} \right) = 
      \left( \begin{array}{c}
         0 \\
         -\nabla_q V(q) \end{array} \right),
 \end{equation}   
which in turn implies that
\begin{equation}   \label{eq.fsp1}
  \varphi_t^{[T]}:  \; \left( \begin{matrix} q_0 \\ p_0 \end{matrix} \right) \longmapsto 
  \left( \begin{matrix} q_0+ t \, \nabla T_p(p_0) \\ p_0 \end{matrix} \right) 
  \end{equation}
and 
\begin{equation}   \label{eq.fsp2} 
  \varphi_t^{[V]}:  \; \left( \begin{matrix} q_0 \\ p_0 \end{matrix} \right) \longmapsto 
  \left( \begin{matrix} q_0 \\ p_0 - t \, \nabla V_q(q_0)\end{matrix} \right).
  \end{equation}
Then, the first order scheme (\ref{LieTrotterHam}) reduces to 
\begin{align}   \label{Euler-sympl2}
      q_{n+1} & =  q_n + h \, \nabla_p T(p_n), \quad
      p_{n+1} =  p_n - h \, \nabla_q V(q_{n+1}), \qquad n=0,1,2\ldots
\end{align}
Compared to the explicit Euler method 
\[
 q_{n+1} = q_n + h \, \nabla_p T(p_n), \qquad p_{n+1} =  p_n - h \, \nabla_q V(q_{n}),
\] 
it only differs in that $\nabla_q V$ is evaluated at the updated value $q_{n+1}$ instead of $q_n$. It makes sense, then, to call
scheme (\ref{Euler-sympl2}) the symplectic Euler-VT method: one first computes the gradient of the kinetic energy $T$ and then computes the gradient of the
potential energy $V$.  

In accordance with our treatment in subsection \ref{sub1.2b}, the adjoint of (\ref{Euler-sympl2}) corresponds 
to composing the maps $\varphi_t^{[T]}$ and $\varphi_t^{[V]}$ in reverse order,
\begin{align}   \label{Euler-sympl}
      p_{n+1} & =  p_n - h \, \nabla_q V(q_{n}), \qquad q_{n+1}  =  q_n + h \, \nabla T_p(p_{n+1}),
    \end{align}
so that we call it the symplectic Euler-TV method.

Obviously, our discussion of schemes (\ref{Euler-sympl2}) and  (\ref{Euler-sympl}) above applies for any Hamiltonian system whose Hamiltonian function can be written in the (so-called separable) form (\ref{eq:Hsep}). 
That is, the Lie--Trotter scheme leads to the two variants of the symplectic Euler method when it is applied to separable Hamiltonian systems.

As for the Strang splitting scheme, described in Algorithm~\ref{alg-S} in general, and in (\ref{StrangHam}) for Hamiltonian problems,
when $H(q,p) = T(p) + V(q)$  it reduces to the much celebrated \emph{St\"ormer--Verlet method} \cite{hairer03gni}. Specifically, depending on the order in which
both parts are evaluated, one has the following two variants:
\begin{algorithm} (St\"ormer--Verlet-VTV). From $(q_0, p_0) = (q(0), p(0))$, for $n \ge 0$,
\label{alg-SV-VTV}
\begin{itemize}
  \item $p_{n+1/2} = p_n - \frac{h}{2} \nabla_q V(q_n)$ 
  \item $q_{n+1} = q_n + h \nabla_p T(p_{n+1/2})$
  \item $p_{n+1} = p_{n + 1/2} - \frac{h}{2} \nabla_q V(q_{n+1})$,
\end{itemize}
\end{algorithm}  
and
\begin{algorithm} (St\"ormer--Verlet-TVT). From $(q_0, p_0) = (q(0), p(0))$, for $n \ge 0$,
\label{alg-SV-TVT}
\begin{itemize}
  \item $q_{n+1/2} = q_n + \frac{h}{2} \nabla_p T(p_n)$ 
  \item $p_{n+1} = p_n - h \nabla_q V(q_{n+1/2})$
  \item $q_{n+1} = q_{n + 1/2} + \frac{h}{2} \nabla_p T(p_{n+1})$,
\end{itemize}
\end{algorithm}  
respectively. Clearly, Algorithms~\ref{alg-SV-VTV} and \ref{alg-SV-TVT} correspond to  \emph{time-symmetric}  methods and can be obtained 
by composing the Euler-TV method and its adjoint. Specifically, if $\chi_{h}$ corresponds to method (\ref{Euler-sympl}), then
\begin{equation}  \label{leapfrog}
 S_h^{[2]} \equiv \chi_{h/2}^* \circ \chi_{h/2} = \varphi_{h/2}^{[V]} \circ \varphi_h^{[T]} \circ  \varphi_{h/2}^{[V]},
\end{equation}
recovers St\"ormer--Verlet-VTV, whereas the TVT version corresponds to $\varphi_{h/2}^{[T]} \circ \varphi_h^{[V]} \circ  \varphi_{h/2}^{[T]}$.

\begin{figure}[htb]
\centering

\includegraphics[scale=0.5]{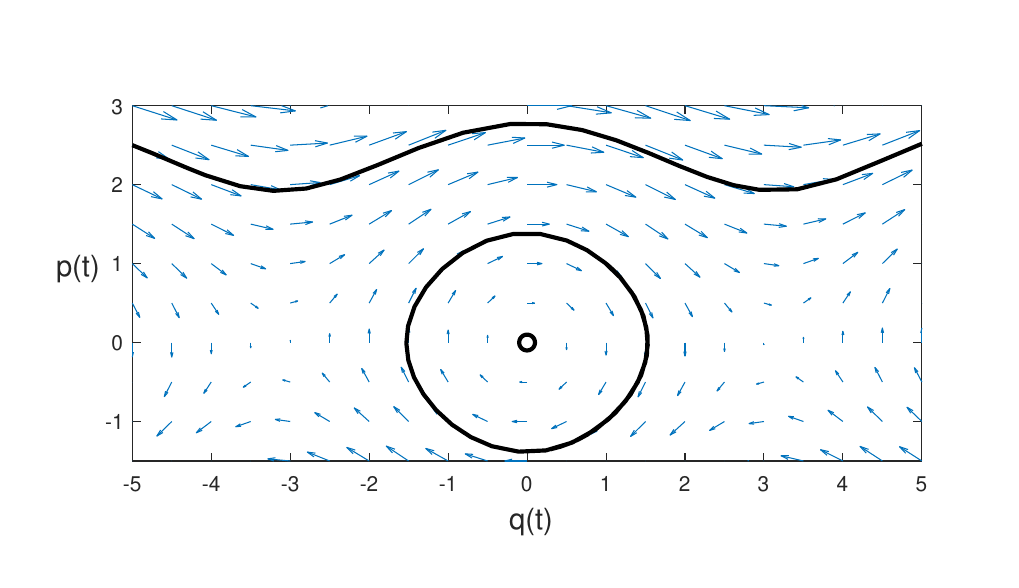} \\
\includegraphics[scale=0.45]{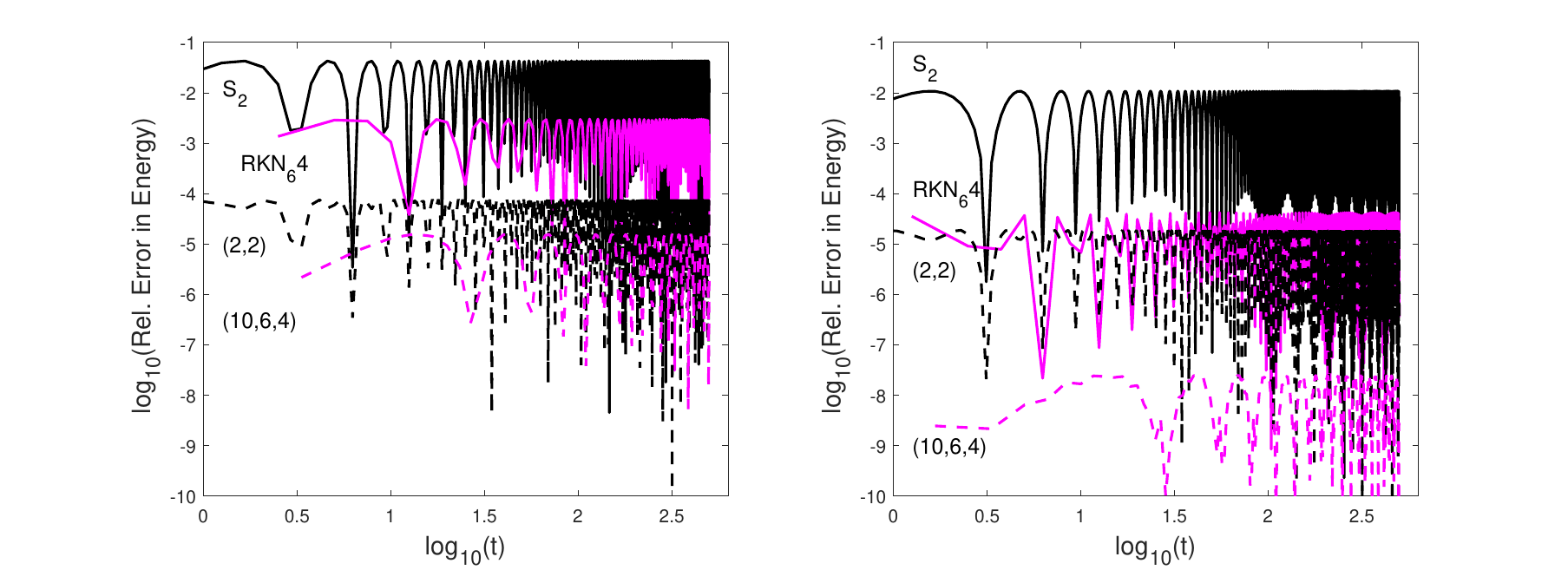}
\caption{Simple pendulum. Top: phase space and three trajectories with initial conditions $(q_0,p_0)=(-5,\frac52),(1,1),(\frac1{10},0)$. Bottom: relative
error in energy
committed by different splitting methods along the solution with initial condition $(q_0,p_0)=(\frac1{10},0)$ in the interval
$t\in[0,500]$ with (left) 1200 evaluations  and (right) 2400 evaluations of the potential.}
\label{fig_Pendulo}
\end{figure}

The top panel of Figure~\ref{fig_Pendulo} shows trajectories of the pendulum (\ref{pendulum1}), starting from three different initial conditions:
$(q_0,p_0)=(-5,\frac52),(1,1),(\frac1{10},0)$, corresponding to different regions of the phase space: for $(q_0,p_0)=(-5,\frac52)$ it holds that $T(p(t))>V(q(t))$,
whereas for $(q_0,p_0)=(1,1)$ one has $T(p(t))\simeq V(q(t))$ on average. Finally, for $(q_0,p_0)=(\frac1{10},0)$, the system can be seen as a slightly perturbed harmonic
oscillator. In this case one could consider the following decomposition:
\begin{equation} \label{pendulum2}
  H= \frac12 (p^2 + q^2) + (1-\frac12q^2-\cos q) \equiv H_1(q,p)+H_2(q),
\end{equation}
where $H_1(q,p) = \frac{1}{2}(p^2 + q^2)$ corresponds to the harmonic oscillator, whose exact solution is known (a rotation in phase space). Moreover,
$|H_2(q)| = \varepsilon |H_1(q,p)|$, with $\varepsilon \approx 10^{-3}$ along the orbit originated in $(\frac1{10},0)$. With splitting (\ref{pendulum2}), the map
$\varphi_h^{[H_1]} \circ \varphi_h^{[H_2]}$ reads
\begin{equation} \label{LTpert}
 \left(  \begin{array}{c}
      q_{n+1} \\
      p_{n+1} \end{array} \right) = R(h) 
       \left(  \begin{array}{c}
      q_{n} \\
      p_{n} + h (q_n - \sin q_n) \end{array} \right), \quad \mbox{with} \quad 
      R(h) = \left( \begin{array}{cc}
      		\cos h & \sin h \\
		-\sin h & \cos h
		  \end{array} \right), 
\end{equation}
whereas the second-order scheme $\varphi_{h/2}^{[H_1]} \circ \varphi_h^{[H_2]} \circ  \varphi_{h/2}^{[H_1]}$ can be formulated as
\begin{equation} \label{soint}
 \begin{aligned}
   &   
     \left(\begin{array}{c}
          q_{n+1/2} \\
	  p_{n+1/2}
           \end{array} \right)
					 =  R(h/2) \left(
     \begin{array}{c}
           q_n\\
            p_n
           \end{array} \right) \\
  &    p_{n+1/2}^*  =  p_{n+1/2} + h (q_{n+1/2} - \sin q_{n+1/2}) \\
     &  \left(\begin{array}{c}
          q_{n+1} \\
	  p_{n+1}
           \end{array} \right)
					 =  R(h/2)
\left(
     \begin{array}{c}
           q_{n+1/2}\\
            p_{n+1/2}^* 
           \end{array} \right).
 \end{aligned}		
\end{equation}
To illustrate the different splittings, we take $(\frac1{10},0)$ as the initial condition, integrate until the final time $t_f=500$
and measure the relative error in energy,  $|H(q_n,p_n)-H(q_0,p_0)|/|H(q_0,p_0)|$, along the trajectory. The step size is taken so that all the methods
tested require the same number of evaluations of the potential (and thus essentially involve the same computational cost): 1200 evaluations (bottom left panel) and
2400 evaluations (bottom right panel). The schemes we test are the following. On the one hand, 
the St\"ormer--Verlet algorithm  (\ref{leapfrog}) (denoted as $S_2$ in the graphs) and the 4th-order Runge--Kutta--Nystr\"om splitting method proposed in \cite{blanes02psp}
(RKN$_64$), as representatives of the $T+V$ splitting. On the other hand, the specially adapted schemes (\ref{soint}), denoted $(2,2)$ and 
the $(10,6,4)$ integrator presented in \cite{blanes13nfo}, both for the case where $H = H_1 + H_2$, with $H_2$ small compared with $H_1$.
The notation  $(10,6,4)$ refers to the fact that the local error of the method 
is of order $\mathcal{O}(\varepsilon h^{11} + \varepsilon^2 h^7 + \varepsilon^3 h^ 5)$ if $H_2$ is $\varepsilon$ times smaller than $H_1$. For analogy, we label Strang method applied to the perturbed harmonic oscillator  (\ref{pendulum2}) as $(2,2)$.
Notice that the splitting (\ref{pendulum2}) is more advantageous 
for this initial condition (the simple
method   (\ref{soint}) behaves better than the 4th-order scheme), and that the improvement with respect to the $T+V$ splitting is approximately of the size 
of $\varepsilon$. We therefore see that for this type
of problem it is possible to construct integrators providing much more accurate results
with the same computational effort.

\begin{figure}[htb]
\centering
\includegraphics[scale=0.4]{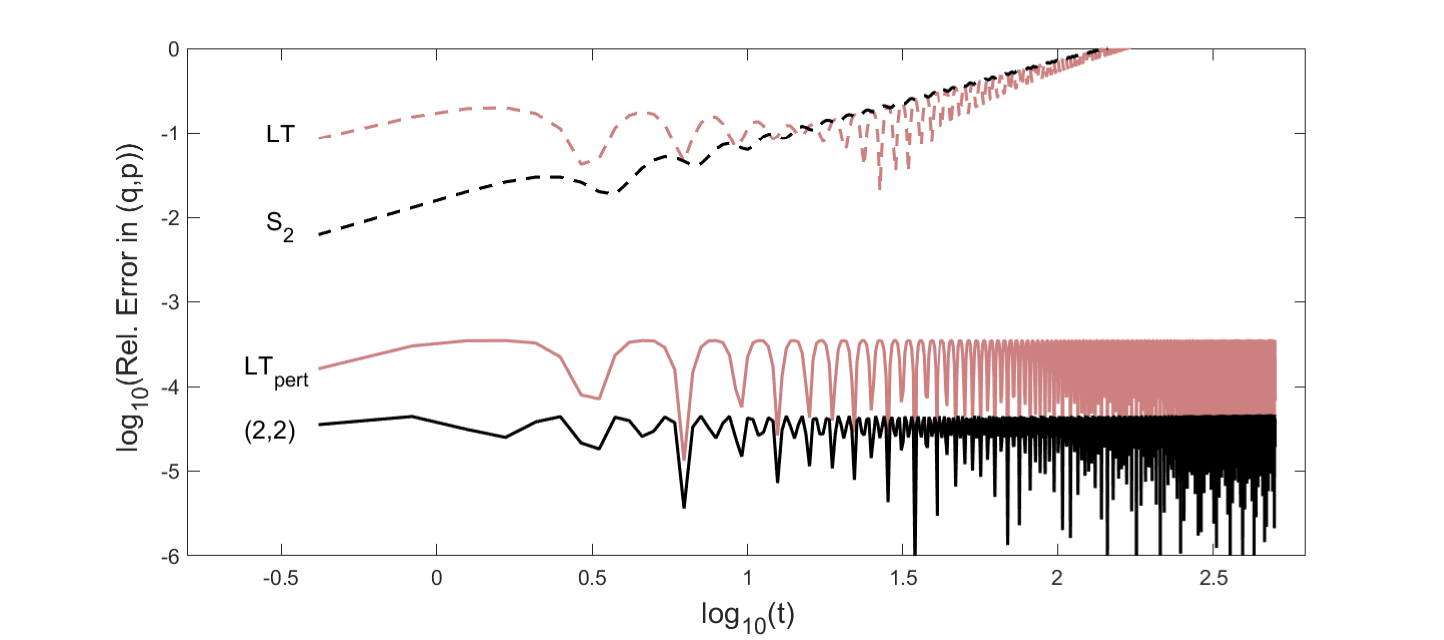}
\caption{Pendulum. Relative error in phase space
for different splitting methods along the solution with initial condition $(q_0,p_0)=(\frac1{10},0)$ in the interval
$t\in[0,500]$ and time step $h=\frac5{12}$.}
\label{fig_PenduloEG}
\end{figure}

Since all the schemes are symplectic integrators and the evolution is taking place in a compact domain, the error in energy is bounded (in contrast with
standard non-symplectic methods, whose error in energy usually grows linearly with $t$), whereas the error
in phase space $(q,p)$ grows linearly when applied to near-integrable Hamiltonian systems  \cite{hairer06gni}.
One should notice, however, that if the scheme is conjugate to another more accurate
integrator, then the global error in phase space will remain bounded for some time interval before it starts growing linearly (see Subsection
\ref{sub_processing1} for more details).
This feature is illustrated in
Figure~\ref{fig_PenduloEG}: we integrate the system starting with the same initial condition and final time $t_f=500$ with $h=\frac5{12}$ and compute the relative error
\[
   \|(q(t_n),p(t_n))-(q_n,p_n)\|/\|(q(t_n),p(t_n))\|,
\]
with different schemes. The reference solution $(q(t_n),p(t_n))$ is computed numerically with very high accuracy. Specifically, we test the following integrators:
the Lie--Trotter method, eq. \eqref{LieTrotterHam} (denoted LT in the graph), and the St\"ormer--Verlet ($S_2$) method for the splitting $H = T+V$, and
scheme $(2,2)$ and  \eqref{LTpert} (LT$_{\rm pert}$), which corresponds
to the Lie--Trotter method applied to the perturbed harmonic oscillator \eqref{pendulum2}.
We observe that, since LT and $S_2$ are conjugate to each other (see eq. \eqref{conju.1}), their errors are quite similar after some time interval. On the other hand,
$(2,2)$ and LT$_{\rm{pert}}$ (which are also conjugate to each other) show a different behavior: no linear growth is visible, and the error of LT$_{\rm{pert}}$ is larger than that of (2,2) by approximately the same factor for the whole time interval considered in Figure~\ref{fig_PenduloEG}.  This is related to the fact that $(2,2)$ and LT$_{\rm{pert}}$ are conjugate to another method with a local error essentially $\mathcal{O}(h^3\varepsilon^2)$ 
(see Subsections \ref{sub_processing1} and \ref{apss}).
More comments on these observations along with additional explanations  will be given in Subsection \ref{sub_processing1}.

\subsection{The gravitational N-body problem}
\label{subsec.1.4}

Another popular example to illustrate the behavior and performance of splitting methods corresponds to the important problem in 
Classical Mechanics of a planetary system modeled as 
$N$ bodies (a massive star and $N-1$ planets) under mutual gravitational Newtonian interaction. This is also a Hamiltonian system with  
\begin{equation}   \label{eq:HamNbody} 
   H(q,p) = \sum_{i=0}^{N-1}  \frac{1}{2m_i}p_i^Tp_i -
    G \sum_{i=1}^{N-1} \sum_{j=0}^{i-1}
    \frac{m_im_j}{\|q_i-q_j\|}.
 \end{equation}
Here $(q,p)$ denote the `supervectors' composed by the positions $q_i \in \mathbb{R}^3$ and momenta $p_i \in \mathbb{R}^3$ 
of the `Sun' ($i=0$) and the $N-1$ planets ($i=1,\ldots, N-1$) in some Cartesian coordinate system:
$q = (q_0, q_1, \ldots, q_{N-1})^T$, $p = (p_0, p_1, \ldots, p_{N-1})^T$.
In (\ref{eq:HamNbody}), $m_i$ is the mass of the $i$-th body, and $G$ is the universal gravitational constant. 

Now the equations of motion (\ref{eq-mo}) read
\begin{equation} \label{eq:Nbody}
\begin{split}
   \frac{dq_i}{dt} &= \nabla_{p_i} H = \frac{1}{m_i} p_i, \qquad i=0,\ldots, N-1, \\
   \frac{dp_i}{dt} &= -\nabla_{q_i} H = - G \sum_{j \neq i}  \frac{m_i m_j}{\|q_i-q_j\|^3}(q_i-q_j), \qquad i,j=0,\ldots, N-1.
   \end{split}
\end{equation}
Since the kinetic energy $T(p)$ and the potential energy $V(q)$ are in this case
\begin{equation}
\label{eq:TV-Nbody}
  T(p) = \sum_{i=0}^{N-1}  \frac{1}{2m_i}p_i^Tp_i, \qquad\quad  V(q) = -G \sum_{i=1}^{N-1} \sum_{j=0}^{i-1} \frac{m_im_j}{\|q_i-q_j\|},
\end{equation}
it also makes sense to separate the Hamiltonian (\ref{eq:HamNbody}) as $H(q,p) = T(p) + V(q)$, so that the symplectic Euler  and the
St\"omer--Verlet schemes can be applied to arbitrary configurations of the bodies. However,  this choice is suboptimal for planetary systems,
where planets describe near-Keplerian orbits around the central star. An alternative procedure taking advantage of the hierarchical nature of the
motion of the planets around the central massive body was first proposed in \cite{wisdom91smf}, and is known in the literature as the 
Wisdom--Holman integration map. It essentially consists in changing coordinates so that the transformed $H$ can be written
as an integrable part $H_1$ (corresponding to the 
Keplerian motion of the planets), and a small perturbation $H_2$ (that accounts for the gravitational interaction of the planets among themselves), and then 
applying the second-order scheme (\ref{StrangHam}) to this new Hamiltonian.

Specifically, \cite{wisdom91smf} consider a linear canonical change of variables to rewrite (\ref{eq:HamNbody}) in the so-called Jacobi coordinates 
$(\hat q_i, \hat p_i)$, $i=0,1,\ldots,N-1$. Here $\hat q_0$ is the position of the center of mass of the system, $\hat q_1$ is the relative position 
of the first planet with respect to the central star, and for $i=2,\ldots,N-1$, $\hat q_i$ is the position of the $i$-th planet relative to the center of mass 
of the central star and the planets with lower indices.  That is,   
\begin{equation}
\label{eq:qhatq}
\hat q_0 =\frac{1}{M_{N-1}}\sum_{j=0}^{N-1} m_j \, q_j, \qquad 
\hat q_i = q_i - \frac{1}{M_i}\sum_{j=0}^{i-1} m_j \, q_j, \quad i=1,\ldots,N-1,
\end{equation}
where $M_i = \sum_{j=0}^{i} m_j$ for $i=0,\ldots,N-1$. This can be written in a more compact way as $\hat q = A q$, where $A$ is an $N\times N$ invertible matrix with mass-dependent entries, $\hat q = (\hat{q}_0, \ldots, \hat{q}_{N-1})^T$ and $q = (q_0, \ldots, q_{N-1})^T$.

The conjugate momenta $\hat p_i$, $i=0,1,\ldots,N-1$, are uniquely determined 
by the requirement that the change of variables be canonical (i.e., $p = A^T \hat p$), so that 
the transformed Hamiltonian function is obtained by rewriting (\ref{eq:HamNbody}) in the new variables:
$\hat H(\hat{q},\hat{p})  \equiv H( A^{-1} \hat{q}, A^T \hat{p})$.
 
It is straightforward to check that  the kinetic energy, expressed in terms  of $\hat{p}$, has the same diagonal structure as in $p$:
\begin{equation*}
T =  \sum_{i=0}^{N-1}  \frac{1}{2m_i}p_i^Tp_i =  \sum_{i=0}^{N-1}  \frac{1}{2 \hat m_i} \hat p_i^T \hat p_i,
\end{equation*}
where
\begin{equation*}
\hat m_0 = M_{N-1}, \qquad \hat m_i = \frac{M_{i-1}}{M_{i}} m_i, \qquad i=1,\ldots,N-1.
\end{equation*}
 The Hamiltonian in the new variables $\hat H(\hat{q},\hat{p})$ can
be split as $\hat H(\hat{q},\hat{p}) = H_1(\hat q, \hat p) + H_2(\hat q)$, where
\begin{equation} \label{jacobi.1}
\begin{aligned}
H_1(\hat q, \hat p)  &=  \frac{1}{2 \hat m_0} \hat p_0^T \hat p_0  + \sum_{i=1}^{N-1} \big( \frac{1}{2 \hat m_i} \hat p_i^T \hat p_i - \frac{m_0 m_i}{\|\hat q_i\|} \big), \\
H_2(\hat q) &= V(q) + G \sum_{i=1}^{N-1} \frac{m_0 m_i}{\|\hat q_i\|} \\
&= G m_0 \sum_{i=1}^{N-1} m_i \big( \frac{1}{\|\hat q_i\|}- \frac{1}{\|q_i-q_0\|} \big)
-G \sum_{i=2}^{N-1} \sum_{j=1}^{i-1} \frac{m_i m_j}{\|q_i-q_j\|}
\end{aligned}
\end{equation}
and $q$ has to be expressed in terms of $\hat{q}$ according to $q = A^{-1} \hat q$.

Observe that the potential energy does not depend on $\hat q_0$, so $\hat p_0$ is constant (in fact, it is the linear momentum of the system)
and therefore we can remove it from $\hat H$ if we assume that the center of the mass is at rest.

Clearly, for fixed $\hat q$ and varying mass ratios, 
\[
\hat q_i = q_i - q_0 +  \mathcal{O}(\varepsilon), \quad H_2(\hat q) = \mathcal{O}(\varepsilon) \quad \mbox{as} \quad
 \varepsilon\equiv \frac{1}{m_0} \max_{1\leq i \leq N-1} m_i \to 0.
\]
Hamiltonian $H_1$ can then be considered as a collection of $N-1$ two-body problems, and $H_2$ as a perturbation.
It turns out that the flow $\varphi_h^{[H_1]}$ can be computed with the algorithm
proposed e.g. in \cite[p. 165]{danby88foc}, whilst $H_2$ only depends on $\hat q$, and thus its flow $\varphi_h^{[H_2]}$ can be explicitly 
evaluated in an efficient way. Notice that the number of terms in $H_1$ grows linearly with the number of bodies $N$, whereas the number of
terms in $H_2$ grows quadratically.

\begin{figure}[htb]
\centering

\includegraphics[scale=0.55]{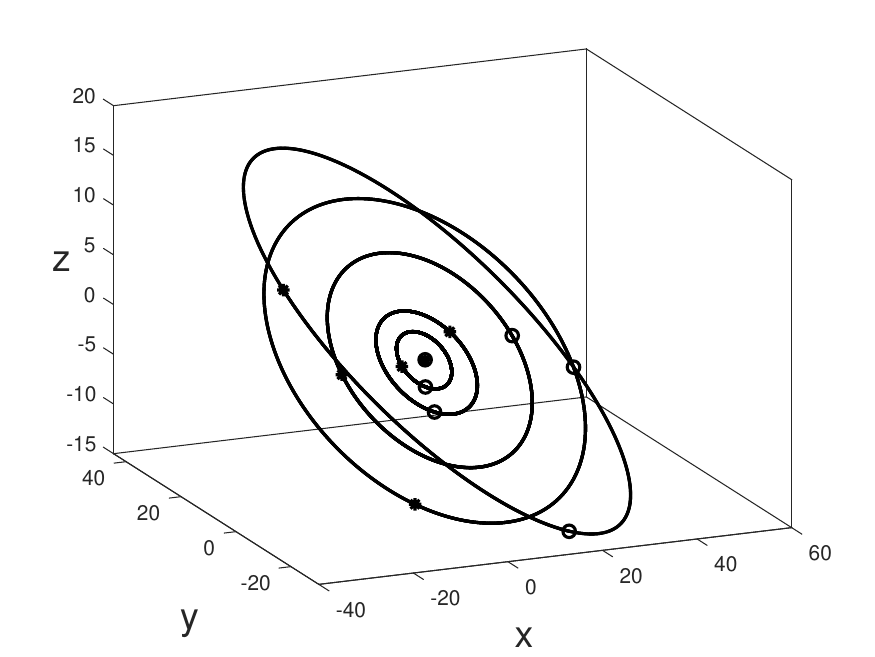} \\
\includegraphics[scale=0.45]{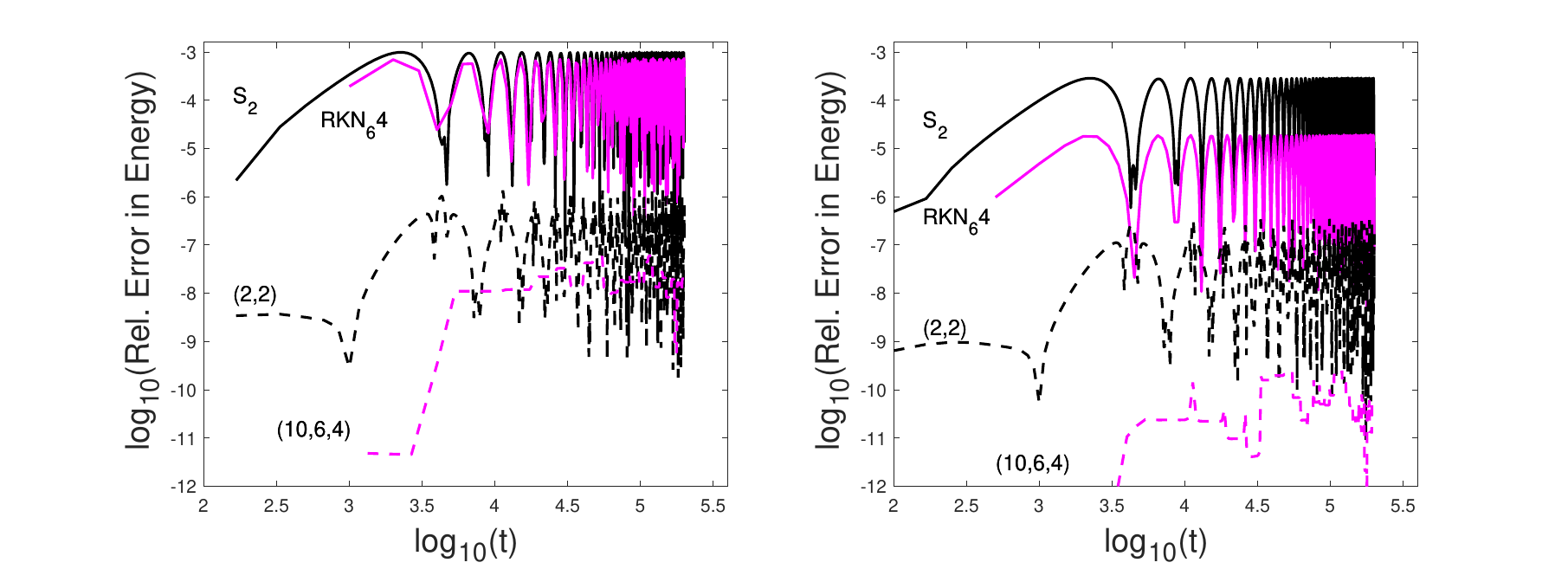}
\caption{Top:  Trajectories of the six-body system modeling the outer Solar System. Bottom: relative error in energy as a function of time for an interval of
200000 days obtained with different splitting methods with (left) 1200 evaluations  and (right) 2400 evaluations of the force.}
\label{fig_NBody1}
\end{figure}

As an illustrative example, we next consider the outer Solar System modeled as a six-body system with the inner Solar System ($i=0$), the four giant planets ($i=1,2,3,4$), and Pluto ($i=5$), all considered as point masses. The initial conditions for each planet are taken at Julian time (TDB)
2440400.5 (28 June 1969), obtained from the DE430 ephemerides \cite{folkner14tpa} and normalized so that the center of mass of the system is at rest.
A schematic diagram of the trajectories is shown in Figure \ref{fig_NBody1} (top) with the initial (circles) and final (stars) positions of each object after 
200000 days. For this problem we
test the same methods as for the pendulum in Figure \ref{fig_Pendulo}. On the one hand,  St\"ormer--Verlet (\ref{leapfrog}), $S_2$, 
and the 4th-order Runge--Kutta--Nystr\"om splitting method of \cite{blanes02psp}, RKN$_64$ when the Hamiltonian (\ref{eq:HamNbody}) is separated into kinetic and potential energy (\ref{eq:TV-Nbody}). On the other hand, 
 the specially adapted schemes (\ref{soint}) (called in this setting the
Wisdom-Holman integrator, and denoted (2,2) as before) and 
the $(10,6,4)$ integrator presented in \cite{blanes13nfo} when $H$ is expressed as $H(\hat{q},\hat{p}) = H_1(\hat q, \hat p) + H_2(\hat q)$, with
(\ref{jacobi.1}) in Jacobi coordinates. We integrate for a relatively short time interval, $t_f = 200000$ days (or approximately 46 periods of Jupyter and two periods of Pluto)
and compute the relative error in the energy with each integrator with a step size so that all of them require the same number of 
force evaluations. The results are displayed in Figure \ref{fig_NBody1} with 1200  (bottom left) and 2400 evaluations (bottom right). As in the previous example of the
pendulum, the error in energy remains bounded in all cases, and scheme $(2,2)$ provides an error energy almost $1000$ times smaller than 
$S_2$.
This illustrates the fact that,  by taking a splitting adapted
to the structure of the problem and designing integrators taking these specific features into account, it is possible to greatly improve the efficiency. 

\begin{figure}[htb]
\centering

\includegraphics[scale=0.5]{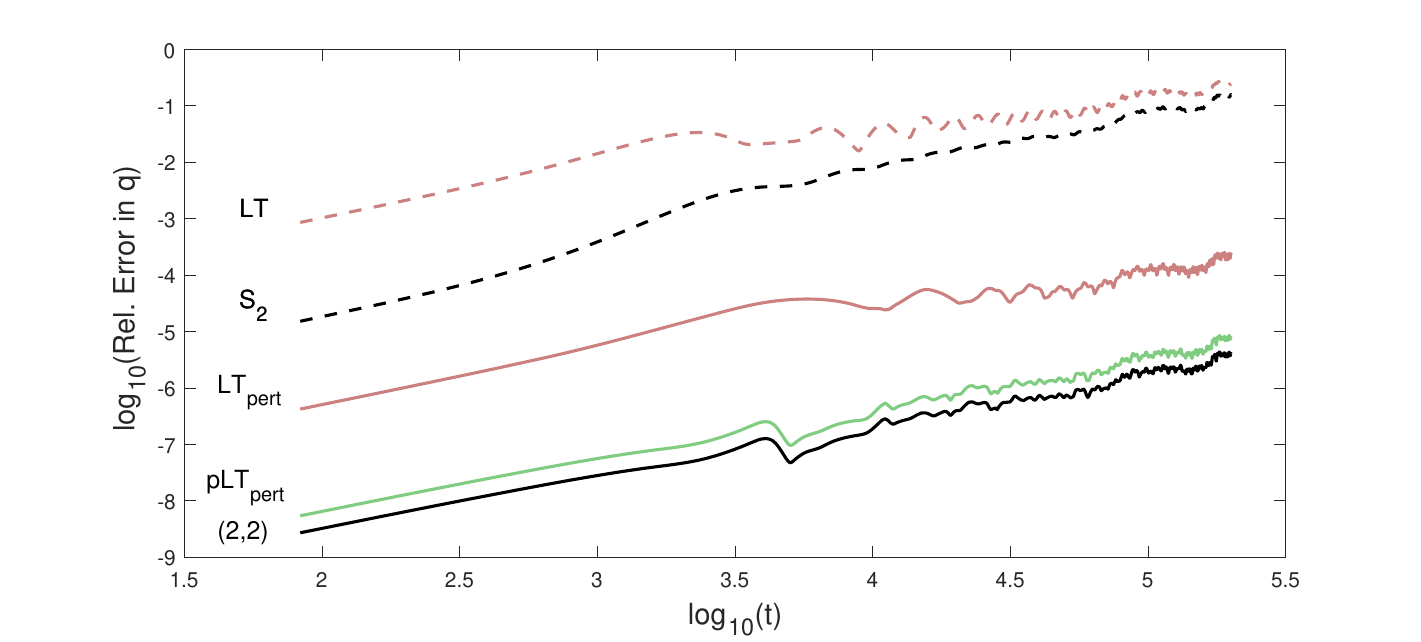}
\caption{Outer Solar System. 
Relative error in position as a function of $t$ for a time interval $[0,t_f=200000]$ days with $h=t_f/1200$ obtained with different splitting methods.}
\label{fig_NBody1EG}
\end{figure}

We display in Figure~\ref{fig_NBody1EG} the relative errors in positions $q=(q_0,\ldots,q_5) \in \mathbb{R}^{18}$
along the time integration for the same methods as in Figure~\ref{fig_PenduloEG}. As in the pendulum problem,  the error for LT is larger than for $S_2$ at the beginning of the integration interval, but they become similar after some time. Here $(2,2)$ is also more accurate than $S_2$ by a factor  significantly smaller than $\varepsilon$. Compared to the pendulum problem,  
the error grows linearly right from the beginning for schemes LT$_{\rm{pert}}$ and $(2,2)$ (although with a smaller slope than LT and $S_2$).
The curve labelled by pLT$_{\rm{pert}}$ corresponds to the relative error obtained by applying LT$_{\rm{pert}}$ with initial conditions $(\bar q_0,\bar p_0) = \varphi_{h/2}^{[H_1]}(q_0,p_0)$. Notice that the phase errors of pLT$_{\rm{pert}}$ are very similar to those of $(2,2)$. All these observations will be accounted for in Subsection \ref{sub_processing1}.

\subsection{The time-dependent Schr\"odinger equation}
\label{subsec.1.3b}

The basic object to study the time evolution of a system in Quantum Mechanics is the (time-dependent) Schr\"odinger equation. In the case
of one particle of unit mass in a potential $V(x)$, it reads \cite{messiah99qme}
\begin{equation}   \label{Schr0}
  i \hbar \frac{\partial}{\partial t} \psi (x,t)  = -\frac{1}{2}
  \Delta \psi (x,t) + V(x) \psi (x,t), \qquad \psi(x,0)=\psi_0(x).
\end{equation}
Here
$\psi:\Omega \subset \mathbb{R}^3\times\mathbb{R} \longrightarrow \mathbb{C}$ is the wave function, 
representing the state of the system and $\hbar$ is the reduced Planck constant. The quantity $|\psi(x,t)|^2$ represents a probability density for the
position of the particle, in the sense that the probability of the particle to be located in $S \subset \Omega$ at time $t$ is 
$\int_{S} |\psi(x,t)|^2 dx$. The equation is then defined in the Hilbert space $L^2(\Omega, \mathbb{C})$.

If we introduce the self-adjoint operators $\hat{T}$, $\hat{V}$ acting on $\psi \in L^2(\Omega, \mathbb{C})$ as
\[
  \hat{T} \psi = -\frac{1}{2} \Delta \psi, \qquad\quad  \hat{V} \psi = V(x) \psi,
\]
then,  a straightforward calculation shows that $[\hat{V}, [\hat{T},\hat{V}]] \psi = |\nabla V|^2 \psi$, and therefore
\begin{equation} \label{rknS}
  [\hat{V},[\hat{V}, [\hat{T},\hat{V}]]] \psi = 0.
\end{equation}  

A standard approach for applying splitting methods in this setting consists in discretizing first the equation in space. If we consider the
one-dimensional case for simplicity and if the wave function is negligible outside a space interval $[a,b]$ on the time interval of interest, then we can limit ourselves
to the study of the equation on that finite interval with periodic boundary conditions \cite{lubich08fqt}. 
After rescaling, the periodic interval can always be restricted to
$[-\pi, \pi]$. In this way, the original problem is transformed into ($\hbar = 1$)
\begin{equation} \label{Schr1}
  i  \frac{\partial}{\partial t} \psi (x,t) = -\frac{1}{2} \frac{\partial^2 \psi}{\partial x^2}(x,t) + V(x) \psi(x,t), \qquad x \in [-\pi,\pi],
\end{equation}
with $\psi(-\pi,t) = \psi(\pi,t)$ for all $t$.

The wave function is then approximated by a trigonometric polynomial $u(x,t)$ whose coefficients are obtained by requiring that the approximation satisfies   
(\ref{Schr1}) in a grid of $M$ equispaced points $x_j = -\pi + j \cdot 2 \pi/(M-1)$ on the interval $[-\pi, \pi]$. The vector $u = (u_{0}, \ldots,
u_{M-1})^T \in \mathbb{C}^M$ formed by the grid values $u_{j} \approx \psi(x_{j},t)$, $j=0,1,\ldots,M-1$, then verifies the $M$-dimensional linear ODE
\begin{equation}   \label{lin1}
  i \frac{d }{dt} u(t) = H \, u(t) = (T + V) u(t),  \qquad
    u(0) = u_{0} \in \mathbb{C}^M.
\end{equation}
Here $V = \mbox{diag}(V(x_j))$ and $T = -\frac{1}{2} D$, where $D$ is the second-order periodic spectral differentiation matrix \cite{trefethen00smi}.
As is well known, $T u = \mathcal{F}^{-1} D_T
\mathcal{F} u$, where $\mathcal{F}$ and $\mathcal{F}^{-1}$
are the forward and backward discrete Fourier transform, and $D_T$ is also diagonal. The transformations $\mathcal{F}$ and $\mathcal{F}^{-1}$
are computed with the fast Fourier transform
(FFT) algorithm, requiring $\mathcal{O}(M \log M)$ operations.

Notice that solving equations $i u' = T u$ and $i u' = V u$ is done trivially by using 
exponentials of diagonal matrices and FFTs, namely,
\[
\left(\e^{\tau V} u \right)_j =
\e^{\tau V(x_j)}u_j, \qquad\quad 
\e^{\tau T} u = \mathcal{F}^{-1} \e^{\tau D_T} \mathcal{F} \, u,
\] 
for a time step $h$, with $\tau = -i h$. Therefore, splitting methods constitute a valid alternative to approximate the 
solution $u(t) = \e^{\tau H} u_0 = \e^{ \tau (T+V)} u_0$, which may be prohibitively expensive to evaluate for large values of $M$. 
Thus, the Lie--Trotter scheme  reads
\begin{equation} \label{LT_Schro}
  \e^{\tau (T+V)} = \e^{\tau T} \, \e^{\tau V} + \mathcal{O}(\tau^2),
\end{equation}
whereas the 2nd-order Strang splitting constructs the numerical approximation $u_{n+1}$ at time $t_{n+1} = t_n + \Delta t$ by
\begin{equation} \label{ssTV}
  u_{n+1} =   \e^{\tau/2 V} \, \e^{\tau T} \, \e^{\tau/2 V} \, u_n.
\end{equation}
The resulting scheme is called the \emph{split-step Fourier} method in the chemical literature \cite{feit82sot}, 
and has some remarkable properties. In particular, it is unitary and 
time-symmetric \cite{lubich08fqt}, as is the exact solution $\e^{\tau H}$. 

Relation (\ref{rknS}) still holds for the matrices $T$ and $V$ if the number of points $M$ in the space discretization is sufficiently large, and in
fact $[V,[T,V]]$ is diagonal if the derivatives of the potential are computed first and then evaluated at the corresponding space grid.

\begin{figure}[!htb]
\centering
\includegraphics[scale=0.6]{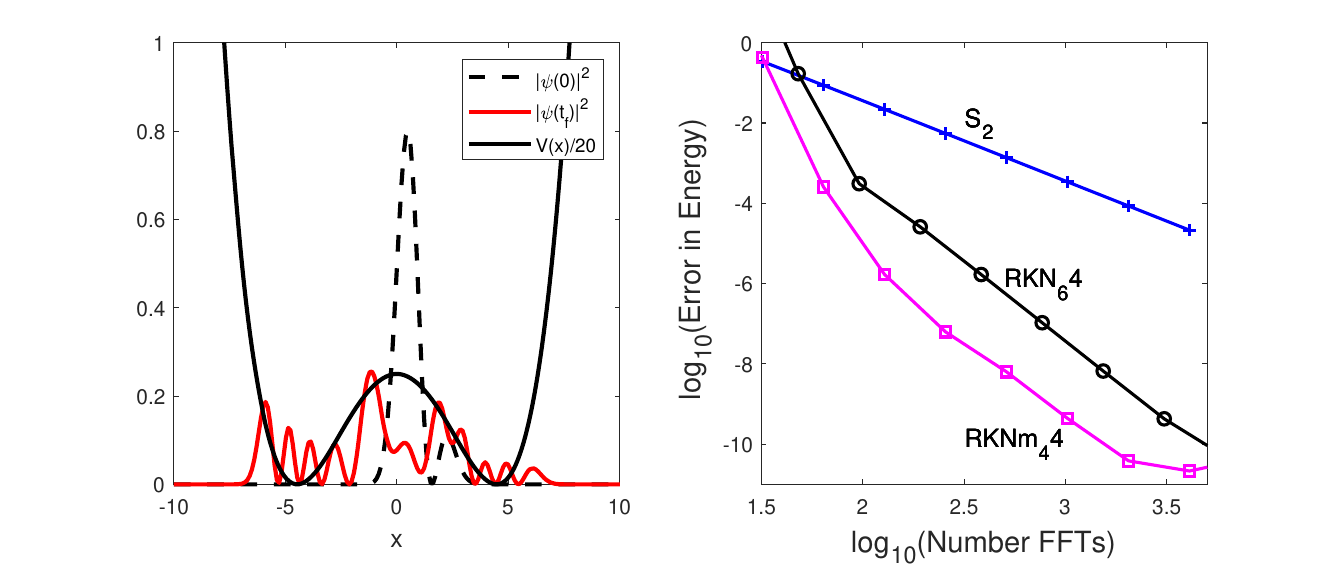}
\caption{Time-dependent Schr\"odinger equation with a double-well potential. Left: $V(x)$, initial and final wave function with $M=256$ discretization points. Right: relative error in energy at the final time vs. number of FFTs for different values of the time step obtained with the Strang method $S_2$, and two 4th-order splitting schemes: one involving 6 evaluations
of $V$ (RKN$_64$) and another with 4 evaluations of $V$, and incorporating in addition the double commutator $[V,[T,V]]$ (RKNm$_44$).}
\label{fig_Schrod1}
\end{figure}

We next illustrate the procedure with the one-dimensional double-well potential 
\begin{equation}  \label{Schr1.Potential}
  V(x)=\frac1{80} (x^2 - 20)^2,
\end{equation}
and the initial wave function $\psi (x,0)=\psi_0(x)=\sigma \cos^2(x)\e^{-\frac12 (x-1)^2}$, where $\sigma$ is an appropriate normalizing constant. We take
$M=256$ discretization points on the interval $x \in [-13,13]$
and integrate the resulting linear ODE (\ref{lin1}) in the interval $0 \le t \le t_f = N h =10$. Figure~\ref{fig_Schrod1} shows in the left panel 
$|\psi_0(x)|^2$, $|\psi (x,t_f)|^2$ and the potential $V(x)$, whereas the right panel shows an efficiency diagram. Specifically, we display
the error in energy measured at the final time, $|u_N^THu_N-u_0^THu_0|$, 
as a function of the number of FFT calls (and its inverse) as an estimate of the computational effort of each method. The lines correspond to the 
Strang splitting (\ref{ssTV}), $S_2$, with time steps 
$h = 10/2^k$, $k=1,2,\ldots,12$, the 4th-order RKN splitting method RKN$_64$ from \cite{blanes02psp}, already illustrated in Figures \ref{fig_Pendulo} 
and \ref{fig_NBody1},  and another 4th-order scheme including the double commutator
$[V,[T,V]]$ into its formulation (denoted RKNm$_44$). In this diagram
the slope of each line for sufficiently small $h$ indicates the order of the scheme.
As in the previous examples, by taking into account the specific features of the problem at hand, 
it is possible
to construct more accurate and efficient numerical approximations.

\begin{figure}[!htb]
\centering
\includegraphics[scale=0.45]{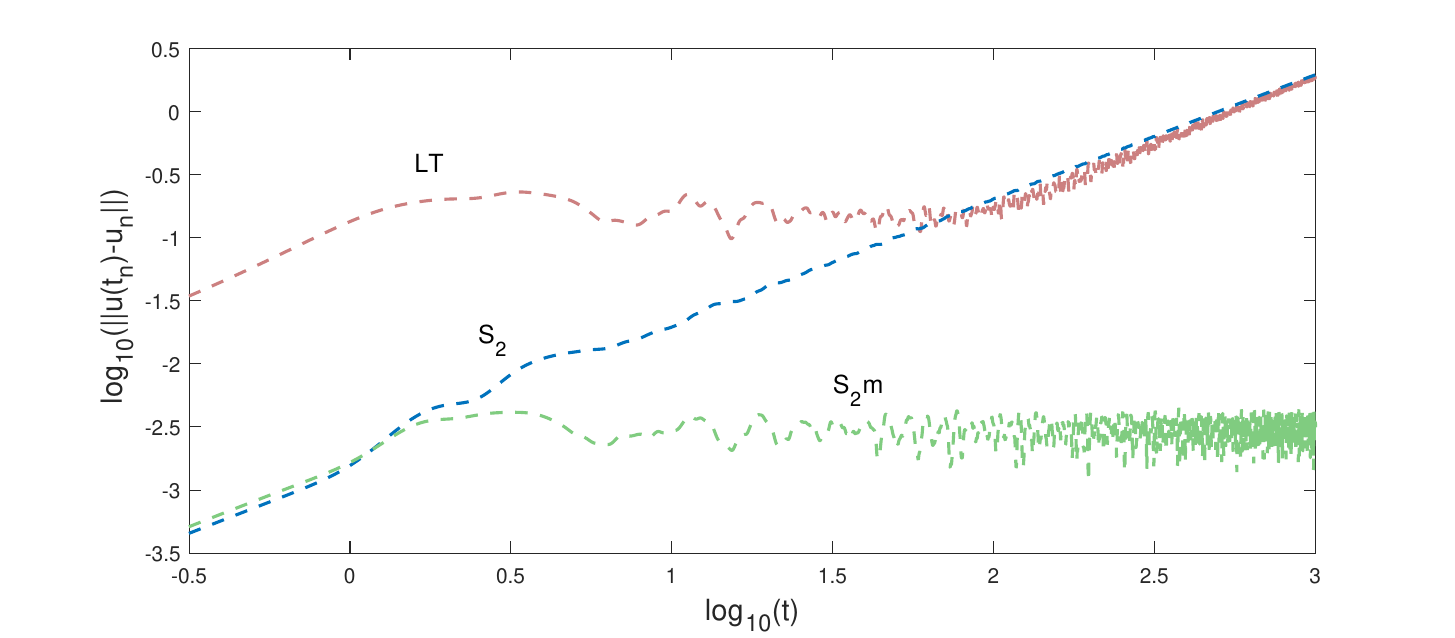}
\caption{Time-dependent Schr\"odinger equation with a double-well potential: relative error in the wave function as a function of $t$ with $h=1/20$ for Strang ($S_2$), 
 Lie-Trotter (LT) and scheme $S_2m$ involving a double commutator and conjugate to a method of order 4.}
\label{fig_Schrod1EG}
\end{figure}

The error in energy also remains bounded for these unitary integrators, as in the previous examples involving classical Hamiltonian problems, whereas the error in the
wave function grows linearly with $t$, unless the scheme is conjugate to another more accurate one, in which case it is bounded for some time before linear growth
takes place. To illustrate this feature, in Figure \ref{fig_Schrod1EG} we depict how the error in the solution $\|u(t_n)-u_n\|$ evolves with time for a longer integration interval
$t\in[0,1000]$, with step size $h=\frac1{20}$. As usual, the reference solution is computed numerically with sufficiently high accuracy, and the tested schemes are
as follows: Lie--Trotter, eq. \eqref{LT_Schro} (LT), Strang, eq. (\ref{ssTV}), and a variant of Strang involving a double commutator, namely
\begin{equation} \label{TI_act}
    S_2m: \qquad u_{n+1} =   \e^{\tau/2 V + (\tau^3/48) [V,[T,V]]} \, \e^{\tau T} \, \e^{\tau/2 V + (\tau^3/48) [V,[T,V]]}  \, u_n.
\end{equation}    
This scheme does not require any additional FFTs and in addition is actually conjugate to a method of order 4. We observe that LT approaches $S_2$ 
after a transition time, whereas the error for $S_2m$ remains bounded for the whole interval. In fact, its linear error growth only appears when the time interval is doubled.

\subsection{Splitting methods as geometric numerical integrators}

As the examples of subsections \ref{subsec.1.2}-\ref{subsec.1.3b}
illustrate, splitting schemes, even of low order of consistency such as Lie--Trotter and Strang  methods, preserve
by construction
structural properties of the exact solution, such as symplecticity (in classical Hamiltonian dynamics) and unitarity (in quantum evolution problems). This
feature gives them qualitative superiority with respect to other standard integrators in practice, especially when long time intervals are concerned. 
In this sense, splitting methods constitute an important class of \emph{geometric numerical integrators}.

Although the idea that numerical integrators applied to an ordinary differential equation should preserve as many properties of the system as possible has been 
implicitly assumed  since the early days of numerical analysis, it is fair to say that in the classical consistency/stability approach the emphasis has been
on other issues. In particular, the goal has mainly been to compute the solution of (\ref{ivp.1}) at time $t_f = N h$ with a global error 
$\|x_N - x(t_f)\|$ smaller than a prescribed tolerance and as efficiently
as possible. To do that one chooses the class of method (one-step,
multistep, extrapolation, etc.), the order (fixed or adaptive) and the time step
(constant or variable). This approach has proved to be very fruitful, giving rise to highly tuned and thoroughly tested software packages generally available
to solve a great variety of problems.

On the other hand, there are special types of problems arising in many fields of science and applied mathematics that 
possess an underlying geometric structure which influences the qualitative character of their solutions, 
and so one aims naturally to construct numerical approximations that preserve this structure. Classical Hamiltonian systems such as those illustrated previously
constitute a case in point.
It turns out, however, that many numerical integrators included in standard software packages do not take into account these distinctive 
features of the equations to be solved, and the question is whether it is possible to design, analyze and apply new schemes providing approximate solutions 
that share one or several geometric properties with the exact solution. This is precisely the realm of Geometric Numerical Integration, 
a terminology introduced in  \cite{sanz-serna97gi}.

According to \cite{mclachlan06gif}, \emph{`Geometric integration' is the term used to describe numerical methods for computing the solution
of differential equations, while preserving one or more physical/mathematical properties of the system} exactly \emph{(i.e., up to round-off error)}. 
Thus, rather than primarily taking into account prerequisites such as consistency and stability, the aim is to reproduce the qualitative
features of the solution of the differential equation being discretized, in particular its geometric properties, such as the symplectic character (for Hamiltonian
systems) and unitarity (quantum mechanics), but also the phase-space volume (for divergence-free vector fields), time-reversal symmetries, first integrals of motion (energy, linear and angular
momentum), Casimirs, Lyapunov functions, etc. In these structure-preserving methods one tries to incorporate as many of these properties as possible
and, as a result, they exhibit improved qualitative behavior. In addition, they typically
allows for a significantly more accurate
integration for long-time intervals than with general-purpose methods \cite{hairer06gni,blanes16aci}.

Although splitting methods have a long history in numerical mathematics and have been applied, sometimes with different names, 
in many different contexts (partial differential equations, quantum statistical mechanics, chemical physics, molecular dynamics, etc.), 
it is fair to say that  the interest in splitting has revived with the advent of Geometric Numerical Integration,  and new and very efficient schemes have been 
put to use to solve a wide variety of problems. The reason is clear: if the problem (\ref{ivp.3}) has some property that is deemed to preserve 
(symplectic, unitary, volume-preserving, etc.) and each subproblem $x' = f_j(x)$, $x(0) = x_0$ can be integrated exactly or by means of a numerical method
preserving these properties, then the splitting
method constructed by composing the solution of the subproblems is also symplectic, unitary, volume-preserving, etc. 
In other words, splitting methods
provide by construction approximations lying in the same group of diffeomorphisms as the system $x' = f(x)$ \cite{mclachlan02sm}. Here we are 
assuming of course that each subproblem $x' = f_j(x)$ possesses the same characteristic feature as the total problem considered.

\subsection{Relevance of splitting methods}

Given a certain differential equation $x' = f(x)$, the application of splitting methods to solve the corresponding initial value problem
involves three main steps \cite{mclachlan02sm}:
\begin{itemize}
 \item[(1)] Choosing the set of terms $f_j$ such that $f = f_1 + \cdots + f_m$. Different ways of decomposing $f$ may give rise to integrators
 with different qualitative behavior and efficiency, as we have seen in the previous examples.
 \item[(2)] Solving each subproblem $x' = f_j(x)$ either exactly or approximately.
 \item[(3)] Combining these solutions to get an approximation for the original overall problem.
\end{itemize} 

Being such a simple idea, it is hardly surprising that the splitting principle can be used in so many different settings. In particular, one may 
\begin{itemize}
  \item split  the differential equation into linear and nonlinear parts;
  \item in an ODE describing a Hamiltonian system with an additional small dissipation, separate the Hamiltonian part and the dissipation;
  \item decompose into parts describing different physical processes, as for example diffusion and reaction in partial differential equations;
  \item get approximations to the original problem by solving one space direction at a time (dimensional splitting in PDEs).
\end{itemize}
In addition, splitting methods possess some advantages concerning their implementation, in particular the following.
\begin{itemize}
 \item They are typically explicit.
 \item Their storage requirements are quite small. The algorithms are sequential and the solution at intermediate stages can be stored in
 the solution vectors.
 \item Programming higher order schemes is not more difficult than Lie--Trotter and Strang splitting methods, at least in the context of ODEs.
 Usually, a few more lines of code is all what is required to deal with the additional stages.
 \item As stated earlier, they can preserve a wide variety of structures possessed by the differential equation.
\end{itemize} 
They also present some disadvantages, of course. Among them, we can mention the following.
\begin{itemize}
 \item Splitting schemes of order three or higher necessarily involve negative coefficients. In other words, they require sub-steps that go backwards
 in time, and this has severe repercussions when applying them to, for instance, reaction-diffusion equations (see Section \ref{sect6}).
 \item Although it is possible to construct accurate high order splitting methods, stability can be an issue, in the sense that their stability interval
might be reduced to render them useless in practice. This aspect has to be taken seriously into account when designing  new methods.
 \item Ordinary splitting does not capture the correct steady state solutions (where these solutions exist) \cite{macnamara16os}, in the sense that the 
 numerical solutions obtained converge to limits that are not steady-state solutions, but just approximations of them. This can lead to unacceptable errors, for instance
 in the simulation of combustion. In this context, balanced splitting techniques have been introduced to correct this flaw \cite{speth13bsa}.
\end{itemize} 

Splitting methods constitute an important tool in different areas of science where the evolution of systems is governed by differential equations.
In addition to Hamiltonian systems,
they can be successfully applied in the numerical study of Poisson systems,
systems possessing integrals
of motion (such as energy and angular momentum) and systems
with (continuous, discrete and time-reversal)
symmetries. In fact, splitting methods have been designed (often independently) and extensively used
in fields as distant as molecular dynamics, simulation of storage rings
in particle accelerators, celestial
mechanics, astronomy, quantum (statistical) mechanics, plasma physics, hydrodynamics, and Markov Chain Monte Carlo methods.

Operator-splitting methods also appear outside the realm of differential equations, and in particular in optimization, in a variety of different
special forms and different denominations (gradient-projection, proximal-gradient, alternating direction method of multipliers or ADMM, 
split Bregman, etc.). All of them are 
related with special types of splitting methods, such as Douglas--Rachford and Peaceman--Rachford schemes. More details can be found in several
contributions collected in the comprehensive book \cite{glowinski16smi}.

\subsection{Some historical remarks}

There is ample consensus that the beginning of splitting is related to the product formula (\ref{ivp.4}). What is not so clear is the origin of the formula
itself. Thus, \cite[p. 295]{reed80momI} establish it as `the classical theorem of Lie', but give no exact source, whereas \cite{chorin78pfa} call it 
`the 1875 formula of S. Lie', citing the classical treatise \cite{lie88tdt}, and, based on this reference, \cite{glowinski16smi} even ascribe to Lie himself `the
first operator-splitting scheme recorded in history'. The problem is that the reference \cite{lie88tdt} is clearly not from 1875, and it is not evident (at least to us)
that this formula appears there explicitly. 

On the other hand, as pointed out in \cite{cohen82eif}, the result (\ref{ivp.4}) can be found in several references published during the 1950s,
namely \cite{butler55pff,golden57sto}, whereas it was Trotter who generalized it to self-adjoint linear operators \cite{trotter59otp} without mentioning Lie or these previous references.
Subsequently, formula (\ref{ivp.4}), even in the matrix case, has been attributed to Trotter \cite[p. 181]{bellman70itm}. In view of the situation, we believe
we are not committing an act of historical injustice by referring to 
the approximation (\ref{ivp.5}) and Algorithm \ref{alg-LT} as the \emph{Lie--Trotter} scheme.

With respect to the splitting method (\ref{ivp.6}), it first appeared in print in \cite{strang68otc} as an alternative way to solve multidimensional problems 
with one-dimensional operators. 
We have already seen that, when applied to Hamiltonian systems of the form $H(q,p) = T(p) + V(q)$, it leads to the 
\blue{Algorithm~\ref{alg-SV-TVT}}
when composing the flows associated to $T(p)$ and $V(q)$. It is called the \emph{St\"ormer--Verlet method} since it was used by the astronomer
C. St\"ormer in 1907 in his computations of the motion of ionized particles in the Earth's magnetic field \cite{stormer07slt}, 
and by L. Verlet in molecular dynamics \cite{verlet67ceo}.
It is also referred to as the \emph{leapfrog method} in the context of PDEs describing wave propagation and as the \emph{Wisdom--Holman method} when applied to
the splitting (\ref{jacobi.1}) 
\cite{wisdom91smf}
In fact, it can be found in several classical references, the oldest one being perhaps Newton's \emph{Principia}\footnote{\emph{Philosophiae Naturalis Principia Mathematica},
Book 1, Section 2, Proposition 1.}. For a detailed account the reader is referred
to the enlightening review \cite{hairer03gni}. 

As well as symplecticity when applied to Hamiltonian systems, the St\"ormer--Verlet method preserves many other geometric properties of the exact
flow associated with an ordinary differential equation. This includes the preservation of all linear first integrals (such as the linear momentum), and quadratic
first integrals of the form $I(q,p) = p^T C q$ for Hamiltonian systems, where $C$ is a symmetric matrix. In other words, $I(q,p)$, computed along the numerical trajectory, 
is constant. A classical example is 
the angular momentum in
$N$-body problems if the forces only depend on the distances of the particles. 

All these favorable properties, in addition to its optimal stability property and reversibility, help to understand why this method is probably the most
used splitting scheme and geometric integrator, especially in molecular dynamics \cite{schlick10mma}, condensed matter simulations \cite{ceperley95pii}, 
 sampling with the hybrid Monte Carlo method \cite{neal11muh}, etc.

 The convenience of designing numerical integration methods that, by construction, preserve the symplectic structure when applied to Hamiltonian systems
 was duly recognized during the 1950s in the field of accelerator physics. Thus, in the words of an early pioneer, \emph{if one wishes to examine solutions to differential equations, adoption of a `Hamiltonian' or `canonical'
integration algorithm would be reassuring} \cite{laslett86nda}. This was the point of view adopted by 
De Vogelaere in his pioneering paper \cite{devogelaere56moi}\footnote{See also \cite{skeel20otu} 
for the context of the work and the preprint itself, typeset in \LaTeX.}, 
where he devoted himself to the task of designing \emph{a method of integration
which, if there was no round-off error, would give a solution with the contact transformation property}. Here ``contact transformation''
has to be understood as ``symplectic transformation''. The first order schemes proposed in \cite{devogelaere56moi}, although implicit in general, turn out
to be explicit when $H(q,p) = T(p) + V(q)$, in which case they reproduce the symplectic Euler schemes (\ref{Euler-sympl2}) and (\ref{Euler-sympl}).

 It was another accelerator physicist, R. Ruth, who presented in 1983 what is probably the first splitting method of order three
 \cite{ruth83aci}. This paper can be considered as the actual starting point in the
systematic exploration of symplectic integrators along several parallel avenues:  
(i) the use of generating functions in the context of Hamiltonian mechanics to produce appropriate canonical transformations 
approximating the exact flow in each integration step \cite{feng87tsm,channell90sio}; (ii) the conditions that 
Runge--Kutta methods have to satisfy to be symplectic \cite{sanzserna88rks,lasagni88crk,suris88pos}; (iii) the design of explicit symplectic methods of order 4 and higher for
Hamiltonian systems that can be split into two pieces which can be solved exactly when considered as independent systems \cite{neri88laa,forest92sol}, with
the help of the Lie formalism. This approach was further elaborated in \cite{yoshida90coh}. Working in the context of the hybrid Monte Carlo algorithm for dynamical fermions, a splitting method of order 4 was also independently proposed around the same time in \cite{campostrini90aco}. In parallel developments,
what is now called the Suzuki--Yoshida composition technique for increasing the order of numerical integrators appeared in \cite{creutz89hoh}
for Monte Carlo simulations and in \cite{suzuki90fdo} and 
\cite{yoshida90coh}.

We should also mention the papers \cite{raedt83aot} and \cite{takahashi84mco}, who pioneered the use of double commutators to 
get approximations of higher order than those
obtained by the St\"ormer--Verlet method in path integral Monte Carlo simulations: 
in the first case by constructing a 4th-order splitting scheme, and in the second, a method that it is also of order four by conjugation. In fact, scheme (\ref{TI_act})
is conjugate to the one proposed in \cite{takahashi84mco}.
The paper 
\cite{ruth83aci} also presents a third-order method using double commutators.

That splitting and composition methods could be used to construct integrators for problems evolving in groups other than the symplectic group was emphasized
in \cite{forest90fos} and further developed  in \cite{feng92fps}, with the aim of constructing schemes able to preserve different structures.

 The 1990s saw a dramatic increase in the interest and applications of splitting integrators in several fields, often with spectacular results. We should mention in particular 
 those achieved in \cite{wisdom91smf}, revealing the existence of chaotic phenomena in the Solar System by numerically integrating
 the planetary equations of motion over very large time intervals.

The state of the art of splitting methods in the context of geometric numerical integration
was masterfully summarized in 2002 by McLachlan and Quispel in their review paper \cite{mclachlan02sm}, which has 
greatly influenced subsequent investigation in the field, as testified by its growing number of citations over the years, in many different areas.

Among the huge number of published works on splitting methods, the following surveys are worth highlighting.
\begin{itemize}
 \item The monograph \cite{yanenko71tmo} (an English translation of the Russian edition published in 1967) is perhaps the first to be devoted to the method
 of splitting (or \emph{method of fractional steps}) for solving ``complicated problems of mathematical physics in several variables''. Those include the
 numerical treatment of parabolic and hyperbolic equations, as well as boundary value problems for the Laplace and Poisson equations, and several
 applications in elasticity theory and hydrodynamics. It is based on the early contributions of Peaceman, Rachford, Douglas, and several authors of
 the Soviet school (Dyakonov, Marchuk, Samarskii, Yanenko, and others).
 \item The exhaustive review article \cite{marchuk90saa} 
 included in the Volume I of the Handbook of Numerical Analysis \cite{ciarlet90hon}. It can be seen as an update of the previous work,
with a systematic study of operator splitting and alternating direction methods for solving linear and nonlinear partial differential equations. It includes
convergence analyses and new applications to problems in hydrodynamics, meteorology and oceanography.
\item The review paper \cite{mclachlan02sm}, mainly focused on the application of splitting methods as geometric numerical integrators for various classes
of ordinary differential equations. In that context, a classification of ODEs and their integration methods into different categories is carried out, 
and the question of how to decompose a given vector field into much simpler vector fields is also examined, as well as the composition of these elementary
flows.
\item Books and monographs dealing with geometric numerical integration and structure-preserving algorithms contain plenty of material on splitting
methods. Among them, we can cite the influential work \cite{sanz-serna94nhp}, the canonical reference \cite{hairer06gni}, as well as \cite{leimkuhler04shd,leimkuhler15md} and \cite{blanes16aci}. 
\item The multi-author book \cite{glowinski16smi} constitutes an excellent illustration of the ample scope and wide range of applications 
that toay's operator-splitting methods are able to deal with. These include the numerical solution of problems modeled by linear
and nonlinear partial differential equations and inequalities, problems in information sciences and image processing, 
and large-scale optimization problems, among others.
\end{itemize}

\subsection{Plan of the paper}

In this paper we will focus on splitting methods applied to evolutionary problems, mostly described by ordinary differential equations. These can directly
model the problem one is interested in, or they can result from evolutionary PDEs previously discretized in space. Particular attention will be addressed to
problems possessing special properties, very often from a geometric origin, that are worth preserving by the numerical methods. In so doing, we will
follow a strategy similar to that  in \cite{blanes08sac}, trying to avoid any duplication of the material already collected in the classic references
cited above, and including new results, schemes and applications which have appeared in the literature during the last few years.

In particular, no general rule is provided here on how to split the defining operator in (\ref{ivp.1}). As mentioned earlier, this issue is further analyzed in
\cite{mclachlan02sm}, and in fact some of the open problems listed there are related with it. We have already seen in the examples provided in this section
that several splittings of the same problem are possible, often leading to methods with very different performances. Moreover, in certain cases,
 the original system has
several geometric properties that are simultaneously preserved along the evolution, whereas different splittings preserve different properties and it is
generally difficult to find one splitting that preserves most of them.

With these considerations in mind, the rest of the paper is organized as follows. In Section \ref{sect2} we first review the general composition technique and then
provide a detailed analysis of the order conditions required by splitting and composition methods to achieve a given order of accuracy. There are some
relevant problems, however, whose particular structure allows one to design specially adapted methods, and some of them are reviewed in Section 
\ref{sect3}, where we also show how to adapt existing splitting methods to non-autonomous problems.

We summarize in Section \ref{sect4} some of the qualitative properties possessed by splitting methods in the context of geometric numerical integration of ordinary
differential equations, with special attention to the idea of processing, whereas Section \ref{sect5} is devoted to the treatment of highly oscillatory problems. 

Splitting methods are particularly well adapted to deal with partial differential equations whose defining operator contains contributions coming from very different
physical sources, and so they have a long history in this area. Section \ref{sect6} contains a brief survey, with special emphasis on Schr\"odinger
equations and general parabolic evolution equations. The existence of negative coefficients in the methods, however, leads to an order barrier for parabolic
equations, and Section \ref{sect7} reviews splitting methods with complex coefficients as a possible way to overcome this order barrier.

In Section \ref{sect8} we present an extended list of existing methods, classifying them into different families and giving the appropriate references. Their
corresponding coefficients are also provided as supplementary material at the website

\begin{center}
  \url{https://www.gicas.uji.es/SplittingMethods.html}
\end{center}  

These methods are numerically tested on simple examples in the Appendix. Finally, some relevant applications of splitting
methods in different fields are discussed in Section \ref{sect9}.


\section{High order splitting and composition methods}
\label{sect2}

The Lie--Trotter and Strang splitting methods, in spite of their low order of accuracy, provide a fairly good description of the systems they are approximately
solving. In fact, for many problems, including molecular dynamics applications and reaction-diffusion equations, Verlet and Strang splitting are the most
popular integrators, perhaps an illustration that, according to 
 \cite{macnamara16os},  ``it is a meta-theorem of numerical analysis that 
second order methods often achieve the right balance between accuracy and complexity''. There are other areas, however, where
a higher degree of precision is required, in addition to the preservation of qualitative properties. A classical example is the long-term
numerical integration of the Solar System, both forwards (e.g to analyze the
existence of chaos \cite{laskar89ane,sussman92ceo}) and backwards in time (to study the insolation quantities of the Earth \cite{laskar04alt}).
 Thus, in this
section, after reviewing a general technique to get high order integrators by composing low order ones, 
we analyze from different perspectives the order conditions a splitting method has to satisfy to achieve a given order. Such analysis allows one to provide complementary information about the integrators:
number of order conditions, explicit expressions, remainders in the asymptotic expansions, etc.

\subsection{Raising the order by composition}
\label{subsec21}

\subsubsection{Composition of Strang maps}

Starting from the Strang splitting $S_h^{[2]}= \varphi_{h/2}^{[1]} \circ \varphi_h^{[2]}  \circ \varphi_{h/2}^{[1]}$, the composition
\begin{equation}  \label{compm.1}
   \psi_h = S_{\gamma_s h}^{[2]} \circ S^{[2]}_{\gamma_{s-1} h} \circ \cdots \circ S^{[2]}_{\gamma_1 h}
\end{equation}
is at least of order three if 
\begin{equation} \label{eq:ord3}
   \sum_{j=1}^s \gamma_j = 1, \qquad \mbox{ and } \qquad  \sum_{j=1}^s \gamma_j^{3} = 0.
\end{equation}
The smallest value of $s$ for which equations (\ref{eq:ord3}) admit real solutions is $s=3$. In that case, by imposing the symmetry $\gamma_1 = \gamma_3$,
we indeed get a method of order 4, sometimes called the  \emph{triple jump}:
\begin{equation} \label{suzu1}
   S_{h}^{[4]} = S_{\gamma_3 h}^{[2]} \circ S_{\gamma_2 h}^{[2]} \circ  S_{\gamma_1 h}^{[2]}, \qquad \mbox{ with } \qquad
   \gamma_1 = \gamma_3 = \frac{1}{2 - 2^{1/3}}, \qquad  \gamma_2 = 1 - 2 \gamma_1.
\end{equation}
In general, the recursion
\begin{equation} \label{eq:triple_jump}
   S_{h}^{[2k]} = S_{\gamma_1 h}^{[2k-2]} \circ S_{(1- 2 \gamma_1) h}^{[2k-2]} \circ  S_{\gamma_1 h}^{[2k-2]},  \qquad \mbox{ with } \qquad
     \gamma_1= \frac{1}{2 - 2^{1/(2k-1)}}
\end{equation}
can be used to get methods  of arbitrarily high order $2k$ ($k \ge 2$) \cite{creutz89hoh} starting from the Strang map $S^{[2]}_h$ (notice that such methods can be written in the form (\ref{compm.1})). The price to be paid is the existence of large positive and negative coefficients $\gamma_j$ and the great number of elementary
flows in (\ref{eq:triple_jump}) for high orders.
The alternative formed by the five maps composition  (\emph{quintuple jump})
\begin{equation} \label{eq:quintuple_jump}
   S_{h}^{[2k]} = S_{\gamma_1 h}^{[2k-2]} \circ S_{\gamma_1 h}^{[2k-2]} \circ S_{(1- 4 \gamma_1) h}^{[2k-2]} \circ  S_{\gamma_1 h}^{[2k-2]}\circ  S_{\gamma_1 h}^{[2k-2]},  \qquad 
     \gamma_1=  \frac{1}{4 - 4^{1/(2k-1)}}
\end{equation}
also gives methods $S_{h}^{[2k]}$ of order $2k$ of the form (\ref{compm.1}) with relatively smaller coefficients $\gamma_j$  but even larger numbers of elementary flows.

 In general, other choices for the coefficients $\gamma_j$ in (\ref{compm.1}) 
are more appropriate if one is interested in achieving orders $\ge 6$ with a lower number of elementary flows and relatively small coefficients.

Condition (\ref{eq:ord3}), using the approach based on linear differential operators discussed in Subsection~\ref{subsec.1.2}, can be derived as follows: the Lie transformation associated with the Strang map $S_h^{[2]}$, eq. (\ref{slt1}), can be written as
$\e^{h/2 F_1} \e^ {h F_2} \e^{h/2 F_1} = \e^{Y(h)}$, where
\[
\begin{aligned}
Y(h) &= \sum_{n=1}^\infty h^{2n-1} Y_{2n-1} \\
&:= \log(\e^{h/2 F_1} \e^ {h F_2} \e^{h/2 F_1})\\
&= h (F_1 + F_2) - \frac{h^3}{24} [F_1,[F_1,F_2]] - \frac{h^3}{12} [F_2,[F_1,F_2]]  + \cdots,
\end{aligned}
\]
that is, $Y_1=F_1 + F_2$, $Y_3 = -\frac{1}{24} [F_1,[F_1,F_2]] - \frac{1}{12} [F_2,[F_1,F_2]]$, and for each $n>2$, $Y_{2n-1}$ is certain linear combination of $(2n-1)$-fold commutators of $F_1$ and $F_2$. In consequence,
the Lie transformation $\Psi(h)$ of (\ref{compm.1}) formally satisfies
\begin{equation}
\label{eq:Psi(h)_0}
\Psi(h) =  \e^{Y(\gamma_1 h)} \cdots \, \e^{Y(\gamma_s h)},
\end{equation}
so that $g(\psi_h(x))=(\Psi(h)g)(x)$ for any $x \in \mathbb{R}^D$ and any smooth function $g:\mathbb{R}^D \to \mathbb{R}$. It is straightforward to check that
\begin{equation}
\label{eq:Psi(h)3}
\begin{split}
\Psi(h) &=  \e^{\gamma_1 h Y_1 + \gamma_1^3 h^ 3 Y_3 + \cdots} \cdots \, \e^{\gamma_s h Y_1 + \gamma_s^3 h^3 Y_3 + \cdots} \\
&= \e^{h \left(\sum_{j=1}^s \gamma_j \right) (F_1+F_2)} + h^3 \left(\sum_{j=1}^s \gamma_j^3 \right) Y_3 + \mathcal{O}(h^4). 
\end{split}
\end{equation}
This shows that the composition (\ref{compm.1}) of Strang maps is of order at least 3 if condition (\ref{eq:ord3}) holds. 

In fact, (\ref{eq:Psi(h)3}) is also true if in (\ref{compm.1}) the Strang map $S^{[2]}_h$ is replaced by any second order time-symmetric map. Furthermore, 
 the triple jump recursion (\ref{eq:triple_jump}) (resp. the quintuple jump recursion (\ref{eq:quintuple_jump})) also gives rise to $2k$th order maps $S_h^{[2k]}$
starting from an arbitrary time-symmetric second order map $S_h^ {[2]}$.  Indeed, this is a consequence of the following four statements:
\begin{itemize}
\item[(1)]  Given an arbitrary near-identity map $\chi_h:\mathbb{R}^D \to \mathbb{R}^D$ (i.e., $\chi_h(x) = x + \mathcal{O}(h)$ as $h\to 0$),  there exists  a series 
\begin{equation*}
Y(h) = \sum_{n\geq 1} h^n Y_n 
\end{equation*}
of (first order) differential operators acting on smooth functions 
such that formally, $g(\chi_h(x)) = (\e^{Y(h)} g)(x)$ for each $x \in \mathbb{R}^D$ and $g\in C^{\infty}(\mathbb{R}^D, \mathbb{R})$. Moreover, $\chi_h$ is an $r$th order integrator for the ODE system
$x' = f_1(x) + f_2(x)$
if and only if 
\begin{equation}
\label{eq:ocond_Y}
 Y_1 = F_1+F_2, \qquad Y_n = 0 \quad  \mbox{for} \quad 2 \leq n \leq r.
\end{equation}
This statement can be proved as follows. Given    
a basic integrator $\chi_h:\mathbb{R}^D \to \mathbb{R}^D$,  consider the linear differential operators $X_n$ ($n\geq 1$) acting on smooth functions 
$g\in C^{\infty}(\mathbb{R}^D, \mathbb{R})$
 as
\begin{equation}
\label{eq:X_n}
X_n g(y) =  \left. \frac{1}{n!}\frac{d^n}{dh^n}\right|_{h=0}  g(\chi_h(y)), \qquad y \in \mathbb{R}^D,
\end{equation}
so that formally $g(\chi_h(x)) = (X(h)g)(x)$, where
\begin{equation}
\label{eq:X(h)}
X(h) = I + \sum_{n\geq 1} h^n X_n,
\end{equation}
and $I$ denotes the identity operator. Each $X_n$ is an $n$th-order differential operator.
Thus, the integrator $\chi_h$ is of order $r$ if
\begin{equation*}
X_n = \frac{1}{n!} (F_1+F_2)^n, \qquad 1 \leq n \leq r.
\end{equation*}
Now consider the series of differential operators
\[
Y(h) = \sum_{n\geq 1} h^n Y_n := \log(X(h)) =
\sum_{m\geq1} \frac{(-1)^{m+1}}{m} \left( h X_1 + h^2 X_2+\cdots \right)^m,
\]
that is,
\begin{equation*}
Y_n = \sum_{m\geq1}^n \frac{(-1)^{m+1}}{m} \sum_{j_1+\cdots+j_m=n}
X_{j_1} \cdots X_{j_m},
\end{equation*}
so that $X(h)= \exp(Y(h))$, and formally, $g(\chi_h(x)) =
(\exp(Y(h))g(x)$. It can be shown that each $Y_n$ is a first order differential operator.
Clearly,  the basic integrator is of order $r$
if (\ref{eq:ocond_Y}) holds.

\item[(2)] The map $\chi_h$ is time-symmetric if and only if $Y(h) = h Y_1 + h^3 Y_3 + \cdots$, which implies that time-symmetric methods \emph{are necessarily of even order}. Indeed, for the adjoint integrator $\chi^*_h=\chi_{-h}^{-1}$, one obviously gets
 $g(\chi^*_h(x)) = \e^{-Y(-h)}g(x)$. Hence,  $\chi_h$ is time-symmetric if and only if  $-Y(-h)\equiv Y(h)$.
 
 \item[(3)] If  $S_{h}^{[2k-2]}$ is a time-symmetric integrator of order at least $2k-2$, then  the composition
\begin{equation}  \label{compm.k}
   \psi_h = S_{\gamma_s h}^{[2k-2]} \circ S^{[2k-2]}_{\gamma_{s-1} h} \circ \cdots \circ S^{[2k-2]}_{\gamma_1 h}
\end{equation}
is of order at least $2k-1$ if 
\begin{equation} \label{eq:ord=2k-1}
   \sum_{j=1}^s \gamma_j = 1, \qquad\quad  \sum_{j=1}^s \gamma_j^{2k-1} = 0.
\end{equation}
In fact, if $\e^{h(F_1 + F_2) + h^{2k-1} Y_{2k-1} + \cdots}$ is the Lie transform of $S_h^{[2k-2]}$, then  $g(\psi_h(x))=(\Psi(h)g)(x)$ for any $x \in \mathbb{R}^D$ and any smooth function $g:\mathbb{R}^D \to \mathbb{R}$, where 
\begin{align*}
\Psi(h) &=  \e^{\gamma_1 h (F_1+F_2) + \gamma_1^{2k-1} h^{2k-1} Y_{2k-1} + \mathcal{O}(h^{2k})} \cdots \, \e^{\gamma_s h (F_1+F_2) + \gamma_s^{2k-1} h^{2k-1} Y_{2k-1} + \mathcal{O}(h^{2k})}\\
& = 
 \e^{\left(\sum_{j=1}^s \gamma_j \right) (F_1+F_2)} + h^{2k-1} \left(\sum_{j=1}^s \gamma_j^{2k-1} \right) Y_{2k-1} + \mathcal{O}(h^{2k}). 
\end{align*}
\item[(4)]  If  $S_{h}^{[2k-2]}$ is a time-symmetric integrator and the sequence $(\gamma_1,\dots,\gamma_s)$ is palindromic, in the sense that for all $j$
 \begin{equation*}
\gamma_{s-j+1} = \gamma_j,
\end{equation*}
then clearly the composition (\ref{compm.k}) is time-symmetric, and hence of even order.
\end{itemize}

 \subsubsection{Composition of Lie-Trotter maps}
 \label{ss2cltm}

One could consider an analogous composition to (\ref{compm.1}), but this time with the Lie--Trotter scheme $\chi_h = \varphi_h^{[1]} \circ \varphi_h^{[2]}$
as the basic method. In that case, the Lie transform of $\chi_h$ is of the form $\e^{Y(h)}$ with 
\begin{equation}
\label{eq:Y(h)}
Y(h) =  h (F_1 + F_2) + h^2 Y_2 + h^3 Y_3 + \cdots, 
\end{equation}
so that  the Lie transform $\Psi(h)$ associated to the composition
$\chi_{\gamma_s h} \circ  \cdots \circ \chi_{\gamma_2 h} \circ \chi_{\gamma_1 h}$
is 
\begin{align*}
\Psi(h) &=  \e^{Y(\gamma_1 h)} \cdots \e^{Y(\gamma_s h)} \\
&= \e^{\left(\sum_{j=1}^s \gamma_j \right) (F_1+F_2)} + h^2 \left(\sum_{j=1}^s \gamma_j^2 \right) Y_2 + \mathcal{O}(h^3). 
\end{align*}
This shows that such a scheme is of order 2 if 
\begin{equation} \label{or2_lt}
\sum_{j=1}^s \gamma_j = 1, \qquad\quad  \sum_{j=1}^s \gamma_j^2 = 0.
\end{equation}
Obviously, such a system of equations does not admit real solutions. 

The situation is different, however, if one composes $\chi_h$ with its adjoint $\chi_h^* = \varphi_h^{[2]} \circ \varphi_h^{[1]}$, that is,
\begin{equation} \label{eq:compint}
\psi_h =  \chi_{\alpha_{2s} h} \circ \chi^*_{\alpha_{2s-1}h} \circ \cdots \circ \chi_{\alpha_{2}h} \circ \chi^*_{\alpha_{1}h}.
\end{equation}
If $\e^{Y(h)}$, with $Y(h)$ given by (\ref{eq:Y(h)}), is the Lie transform of $\chi_h$, then the Lie transform of its adjoint
$\chi^*_h$ is  $\e^{-Y(-h)}$ and the Lie transform $\Psi(h)$ of (\ref{eq:compint}) satisfies
\begin{equation} \label{lt_comp1}
\begin{aligned}
\Psi(h) &=  \e^{-Y(-\alpha_1 h)} \e^{Y(\alpha_2 h)} \cdots \, \e^{-Y(-\alpha_{2s-1} h)}\e^{Y(\alpha_{2s} h)} \\
&= \e^{\left( \sum_{j=1}^{2s} \alpha_j \right) (F_1+F_2)} + h^2 \left(\sum_{j=1}^s  \big(\alpha_{2j}^2 - \alpha_{2j-1}^2 \big) \right) Y_2 + \mathcal{O}(h^3),
\end{aligned}
\end{equation}
 whence composition (\ref{eq:compint}) is at least of order 2 if the coefficients satisfy
 \[
     \sum_{j=1}^{2s} \alpha_j = 1,  \qquad\quad  \sum_{j=1}^s  \big(\alpha_{2j}^2 - \alpha_{2j-1}^2 \big) = 0.
\]   
The argument above also holds when $\chi_h$ is a first order integrator other than Lie-Trotter.
   
As mentioned earlier, the simplest situation corresponds to $s=1$, in which case we recover the Strang splitting, $S_h^{[2]} = \chi_{h/2} \circ \chi_{h/2}^*$.
 In fact, the general scheme (\ref{eq:compint}) can be
 rewritten as the splitting method
 \begin{equation} \label{eq:splitting}
\psi_h = \varphi^{[1]}_{a_{s+1} h}\circ \varphi^{[2]}_{b_{s}h}\circ \varphi^{[1]}_{a_{s} h}\circ \cdots\circ
 \varphi^{[1]}_{a_{2}h}\circ
 \varphi^{[2]}_{b_{1}h} \circ \varphi^{[1]}_{a_{1}h},
\end{equation}
where $a_{1}=\alpha_{1}$, and for $j=1,\ldots,s$,
\begin{equation}
  \label{eq:abalpha}
a_{j+1} = \alpha_{2j} + \alpha_{2j+1}, \qquad\quad  b_j=\alpha_{2j-1}+\alpha_{2j}
\end{equation}
 (with $\alpha_{2s+1}=0$). Conversely, any integrator of the form (\ref{eq:splitting}) satisfying the condition
$\sum_{j=1}^{s+1} a_j = \sum_{j=1}^{s} b_j$
  can be expressed in the form (\ref{eq:compint}), as shown
  in \cite{mclachlan95otn}.

Clearly, any splitting scheme (\ref{eq:splitting}) with a palindromic sequence of coefficients, that is, satisfying 
\begin{equation*}
a_{s-j+2} = a_j, \qquad b_{s-j+1}=b_j \qquad \mbox{ for all } j,
\end{equation*}
is time-symmetric, and thus of even order.  Written in the composition format (\ref{eq:compint}),  it is time-symmetric if 
\begin{equation*}
\alpha_{2s-j+1} = \alpha_j.
\end{equation*}

\subsection{Order conditions I: splitting schemes with the BCH formula}
\label{subsection:BCH}

In the analysis of the order conditions for splitting methods, and without lost of generality, we will consider the linear case 
(\ref{ivp.3b}), so that the treatment is essentially based on matrices, or more generally, on linear operators. Thus, the integrator (\ref{eq:splitting}) corresponds 
in this setting to the product of exponentials 
\begin{equation} \label{eq:Psi(h)}
   \Psi(h) = \e^{a_{s+1} h F_1} \, \e^{b_{s} h F_2} \, \e^{a_{s} h F_1} \, \cdots \, \e^{a_2 h F_1} \, \e^{b_1 h F_2} \, \e^{a_{1} h F_1},
\end{equation}
intended to approximate $\e^{h(F_1+F_2)}$.

In the case of systems of ODEs $x' = f_1(x) + f_2(x)$ with (non necessarily linear) vector fields $f_1$ and $f_2$, one may
consider the Lie operators $F_j$ of $f_j$, and compare the Lie transformation $\e^{h (F_1+F_2)}$ with the operator
\begin{equation}
\label{eq:LieTransformationPsi}
\Psi(h) = \e^{a_{1} h F_1} \, \e^{b_{1} h F_2} \, \e^{a_{2} h F_1} \, \cdots \, \e^{a_s h F_1} \, \e^{b_s h F_2} \, \e^{a_{s+1} h F_1}
\end{equation}
associated with the map $\psi_h$ (in the sense that $g(\psi_h(x))=\Psi(h)g(x)$ for arbitrary smooth $g \in C^{\infty}(\mathbb{R}^D,\mathbb{R})$). 
Thus, all the formulas derived below for the linear case can be applied to the Lie transformation 
(\ref{eq:LieTransformationPsi}) associated with the map (\ref{eq:splitting}) just by reversing the order of the sequence $(a_1,b_1,a_2,\ldots,a_s,b_s,a_{s+1})$ of coefficients of the splitting scheme. Alternatively, we may reverse the order of the products of the operators in all the expressions involved.

Equivalent results could also be obtained for the ODE case using the concept of {\em word series} and related techniques introduced in~\cite{murua17wsf}.
 The Lie transform approach is more general, however, as it can be directly applied to the case of splitting methods for differential equations on manifolds, as discussed at the end of Subsection~\ref{subsec41}.

Generally speaking, the order conditions for a method of order $r$
are systems of polynomial equations in the coefficients obtained by requiring that the Taylor expansions in the
step size $h$ of both the exact and numerical solution agree up to terms in $h^r$. A standard approach to obtain the order conditions of scheme (\ref{eq:Psi(h)})
consists in formally using the Baker--Cambpell--Hausdorff (BCH) 
formula to express $\Psi(h)$ as one exponential of a series of operators in powers of
$h$, and finally to compare this series with $\e^{h(F_1 + F_2)}$. In this way, one gets
\begin{equation} \label{bch2}
\begin{aligned}
  &  \log(\Psi(h)) = h (w_{1} F_1 + w_{2} F_2)  + h^2 w_{12} F_{12} + h^3 (w_{122} F_{122}  + w_{112} F_{112}) \\
  & \qquad + h^4 ( w_{1222} F_{1222} + w_{1122} F_{1122} + w_{1112} F_{1112}) + \mathcal{O}(h^5),
\end{aligned}
\end{equation}
where
\[
\begin{aligned}
   & F_{12} = [F_1, F_2], \qquad F_{122} = [F_{12},F_2], \qquad F_{112} = [F_1,F_{12}] \\
   & F_{1222} = [F_{122},F_2], \qquad F_{1122} = [F_1, F_{122}], \qquad F_{1112} = [F_1,F_{112}]
 \end{aligned}
 \]
 and $w_1$, $w_2$, $w_{12}$, $w_{122}, \ldots$ are polynomials of homogeneous degree in the parameters $a_j$, $b_j$. In particular 
 \[
 \begin{aligned}
   & w_1 = \sum_{i=1}^{s+1} a_i, \qquad w_2 = \sum_{i=1}^{s} b_i, \qquad w_{12} = \frac{1}{2} w_1 w_2 - \sum_{1 \le i < j \le s+1} b_i a_j, \\
   &   w_{122} = \frac{1}{12} w_{1} w_{2}^2 -\frac{1}{2}  \sum_{1\leq i\leq j < k \leq s} b_i a_j b_k, \\
   &  w_{112} =  \frac{1}{12} w_{1}^2 w_{2} - \frac{1}{2} \sum_{1\leq i < j \leq k \leq s+1} a_i b_j a_k, 
\end{aligned}
 \]    
with $b_{s+1}=0$. From (\ref{bch2}), it is clear that a characterization of the
order of the splitting scheme (\ref{eq:splitting}) is obtained by requiring the consistency conditions $w_1=w_2=1$, that is, 
\begin{equation}
\label{eq:consistency_cond}
  \sum_{j=1}^{s+1} a_j = \sum_{j=1}^{s} b_j = 1
\end{equation}
(ensuring that the scheme (\ref{eq:splitting}) is at least of order 1)
 and $w_{12}=w_{122}=w_{112}=\cdots =0$ up
to the required order. The set of order
conditions thus obtained will be independent in the general case if
the operators $F_1$, $F_2$, $F_{12}$, $F_{122}$, $F_{112}, \ldots$
form a basis of the
free Lie algebra in the alphabet $\{1,2\}$. We have considered in (\ref{bch2}) the so-called \emph{Lyndon basis},  
associated to the set of Lyndon words in the alphabet $\{1,2\}$
\begin{equation}
\{1,2,12,122,112,1222,1122,1112,\ldots\}.
\end{equation}
They are defined as follows~\cite{reutenauer93fla}:
a word $\ell_1\cdots \ell_m$ is a Lyndon word if
$(\ell_1\cdots \ell_k) \prec (\ell_{k+1}\cdots \ell_m)$ for each $1 \leq k <
m$, where $\prec$ is the lexicographical order (i.e., the order used when
ordering words in the dictionary) on the set of words on the alphabet $\{1,2\}$. 
 For instance,  $112$ is a Lyndon word, while  neither $211$ nor $121$ are, as $2 \not \prec 11$ and as $12 \not \prec 1$, respectively.  
 
 The element of the basis associated to a Lyndon word $\ell_1\cdots \ell_m$ with $m\geq 2$ is given as $F_{\ell_1\cdots \ell_m} = [F_{\ell_1\cdots \ell_n}, F_{\ell_{n+1}\cdots \ell_m}]$, where $n$ is the smallest number such that both $\ell_1\cdots \ell_n$ and $\ell_{n+1}\cdots \ell_m$ are themselves Lyndon words.

   An efficient algorithm
(based on the results in~\cite{murua06tha}) for the BCH formula and
related calculations in the Lyndon basis (and some other basis)  that allows one to obtain
(\ref{bch2}) up to terms of arbitrarily high degree is
presented in~\cite{casas09aea}.

Of course, if another basis of the free Lie algebra in the alphabet $\{1,2\}$ is used to expand $\log(\Psi(h))$ in (\ref{bch2}), a different characterization of the order conditions will be obtained, with a different set of polynomial functions on $(a_1,b_1,\ldots,a_s,b_s,a_{s+1})$. In any case, the number of such independent
 conditions arising at each order $n$ can be obtained just by determining the dimension of $\mathcal{L}_n(F_1,F_2)$, the linear span of all commutators 
 containing $n$ operators $F_1$, $F_2$. This number, denoted $c_n$, is given in Table \ref{table_orcon} (see \cite{munthe-kaas99cia,mclachlan02sm}).

\begin{table}
\begin{center}
  \begin{tabular}{|c|ccccccccccc|} \hline
  Order $n$ & 1 & 2 & 3 & 4 & 5 & 6 & 7  &  8 & 9 & 10 & 11 \\ \hline\hline
$c_n$ & 2 & 1 & 2 & 3 & 6 & 9 & 18 & 30 & 56& 99 & 186 \\
$d_n$ & 2 & 1 & 2 & 2 & 4 & 5 & 10 & 14 & 25 & 39 & 69  \\ 
$m_n$ & 1 & 0 & 1 & 1 & 2 & 2 &  4 & 5  & 8 & 11 & 17 \\ \hline
\end{tabular}
\end{center}
  \caption{Number of independent order conditions for general splitting methods, $c_n$, and for RKN-type splitting methods, $d_n$.
  The number $m_n$ corresponds to the number of order conditions for compositions of a 2nd-order time-symmetric method.}
  \label{table_orcon}
\end{table}

\subsection{Order conditions II: splitting schemes with Lyndon words}

Whereas the previous characterization of the order of the splitting scheme (\ref{eq:Psi(h)}) allows one to easily get the number of order conditions,
obtaining explicit expressions for the polynomials $w_{\ell_1 \cdots \ell_n}$ is much more difficult when the considered order increases. The following
alternative characterization, based on the direct comparison of the power series expansions of $\e^{h\, (F_1 + F_2)}$ and (\ref{eq:Psi(h)}), tries to
ameliorate this difficulty.

\subsubsection{Basic expansions and necessary order conditions}

Both  $\e^{h\, (F_1 + F_2)}$ and (\ref{eq:Psi(h)}) admit an expansion in series indexed by the set
\begin{equation*}
\W = \{1, 2, 11, 12, 21, 22, 111, 112, 121, 211, 122, \cdots\}
\end{equation*}
of words in the alphabet $\{1,2\}$. More precisely,  $\e^{h\, (F_1 + F_2)}$ can be expanded as
\begin{equation}
\label{eq:alpha-series}
I + h \alpha_1  F_1 + h \alpha_2 F_2 +  h^2 \alpha_{11} F_1 F_1 + h^2 \alpha_{12} F_1 F_2 +
 h^2 \alpha_{21} F_2 F_1 + h^2 \alpha_{22} F_2 F_2 + \cdots,  
\end{equation}
with $\alpha_{\ell_1 \cdots \ell_n} = 1/n!$. As for (\ref{eq:Psi(h)}), it can be expanded, for arbitrary $s$, as (\ref{eq:alpha-series}), where for each word $\ell_1 \cdots \ell_n \in \W$ with $n$ letters, the corresponding coefficient
\[
\alpha_{\ell_1 \cdots \ell_n} = u_{\ell_1 \cdots \ell_n}(a_1,b_1,\ldots,a_s,b_s,a_{s+1})
\]
 is a homogeneous polynomial of degree $n$ in the variables $a_1,b_1,\ldots,a_s,b_s,a_{s+1}$.  

 It is straightforward to check that such polynomials satisfy the following  relations, which allows one to compute them recursively: 
  \begin{itemize}
 \item if $\ell_1=\cdots=\ell_{j} = 1$ and $\ell_{j+1}\neq 1$ with $j\geq 1$,
\begin{align*}
u_{\ell_1 \cdots \ell_n}(a_1,b_1,\ldots,a_s,b_s,a_{s+1})&= 
\sum_{k=0}^{j}  u_{\ell_{k+1} \cdots \ell_{n}}(a_1,b_1,\ldots,a_s,b_s,0)\, \frac{a_{s+1}^k}{k!},\\
u_{\ell_1 \cdots \ell_n}(a_1,b_1,\ldots,a_s,b_s,0)&=u_{\ell_1 \cdots \ell_n}(a_1,b_1,\ldots,a_s).
\end{align*}
\item if $\ell_1=\cdots=\ell_{j} = 2$ and $\ell_{j+1}\neq 2$ with $j\geq 1$,
\begin{equation*}
u_{\ell_1 \cdots \ell_n}(a_1,b_1,\ldots,a_s,b_s,a_{s+1})= 
\sum_{k=0}^{j}  u_{\ell_{k+1} \cdots \ell_{n}}(a_1,b_1,\ldots,a_s)\, \frac{b_{s}^k}{k!},
\end{equation*}
 \item if $\ell_n=\cdots=\ell_{1} = 1$,
\begin{equation*}
u_{\ell_1 \cdots \ell_n}(a_1)= \frac{a_{1}^n}{n!},
\end{equation*}
 \item if $\ell_{1} \neq 1$,
\begin{equation*}
u_{\ell_1 \cdots \ell_n}(a_1)=0.
\end{equation*}
\end{itemize}

In this way (\ref{eq:Psi(h)})  is at least of order $r$ if and only if the conditions
\begin{equation}
\label{eq:ocond_ab}
u_{\ell_1 \cdots \ell_n}(a_1,b_1,\ldots,a_s,b_s,a_{s+1})=\frac{1}{n!}
\end{equation}
 hold for each word $\ell_1\cdots\ell_n$ with $n\leq r$ letters in the alphabet $\{1,2\}$.   For illustration, in Table~\ref{tab.ocond_ab} we explicitly give these conditions  corresponding to words with up to two letters. 
 Notice that (\ref{eq:ocond_ab}) for the single-letter words $1$ and $2$  coincide with the  consistency conditions (\ref{eq:consistency_cond}).

\begin{table}[h!]
\begin{center}
\begin{tabular}{|c||c||c||c|} \hline
  Word  \hspace*{-0.2cm}&  \hspace*{0.4cm}  Condition & Word  \hspace*{-0.2cm} &  \hspace*{0.4cm}  Condition \\ \hline \hline
      $1$ & $\displaystyle \sum_{j=1}^{s+1} a_j=1$ & $2$   &  $\displaystyle  \sum_{j=1}^{s} b_j=1$ \\ \hline
      $11$   &  $\displaystyle  \frac{1}{2} \sum_{j=1}^{s+1} \frac{a_j^2}{2} + \sum_{1\leq i < j \leq s+1}a_i  a_j = \frac12$ & 
      $12$ & $\displaystyle \sum_{1 \le i < j \le s+1} b_i a_j=\frac12$ \\  \hline
      $21$ & $\displaystyle \sum_{1 \le j \le i \le s} b_i a_j=\frac12$ &
      $22$   &  $\displaystyle  \frac{1}{2} \sum_{j=1}^{s} \frac{b_j^2}{2} + \sum_{1\leq j < i \leq s}b_j  b_i = \frac12$ \\ \hline
\end{tabular}
\end{center}
\caption{Conditions (\ref{eq:ocond_ab}) corresponding to words with up to two indices.\label{tab.ocond_ab}}
\end{table}
However, such order conditions are not all independent.   For instance, from Table~\ref{tab.ocond_ab}, one can check that 
\begin{equation*}
u_{1}^2 = 2 u_{11}, \quad u_{2}^2 = 2 u_{22}, \quad u_{1} u_{2} = u_{12} + u_{21}.
\end{equation*}
For a consistent method, $u_{1} = u_{2} =1$, hence $u_{11}=\frac12$, $u_{22}=\frac12$, and
\begin{equation*}
 u_{12}-\frac12 =\frac12 - u_{21}, 
\end{equation*}
which implies that, if $u_{12}=\frac12$, then automatically $u_{21} = \frac12$.

A complete characterization of the relations among the order conditions (\ref{eq:ocond_ab}) will be obtained in paragraph~\ref{sss:shuffle} below. As a previous step, we obtain integral representations of both   $\e^{h\, (F_1 + F_2)}$ and (\ref{eq:Psi(h)}) , which in addition 
give useful expressions for the remainders of their truncated series expansions. 

\subsubsection{Integral representation and remainders}
\label{sss:remainders}

Consider the solution $Y(\tau)$  of the initial value problem
\begin{equation}
\label{eq:Yode}
\frac{d}{d\tau} Y(\tau) = h  A(\tau) Y(\tau), \qquad Y(0) = I,  
\end{equation}
with $A(\tau) = d_1(\tau) F_1 + d_2(\tau) F_2$ and 
\begin{equation}
\label{eq:def_d}
(d_1(\tau),d_2(\tau)) = 
\left\{
\begin{matrix}
(a_i, 0) & \mbox{ if } & \tau \in [2i-2,2i-1], \ i \in \{1,\ldots,s+1\}, \\
(0,b_i) & \mbox{ if } & \tau \in [2i-1,2i], \ i \in \{1,\ldots,s\},\\
(0,0) & \mbox{ if } & \tau >2s+1.
\end{matrix}
\right.
\end{equation}
It is straightforward to check that 
\begin{equation*}
Y(\tau) =
\left\{
\begin{matrix}
\e^{(\tau-(2i-2)) h a_{i+1} F_1} \, \e^{b_i h F_2}\, \e^{a_i h F_1}\, \cdots \, \e^{a_2 h F_1} \, \e^{b_1 h F_2}\,  \e^{a_1 h F_1}, 
& \tau \in [2i-2,2i-1], \\
\e^{(\tau-(2i-1)) h b_i F_2}\, \e^{a_i h F_1}\, \cdots\, \e^{a_2 h F_1}\, \e^{b_1 h F_2} \, \e^{a_1 h F_1},
&  \tau \in [2i-1,2i].
\end{matrix}
\right.
\end{equation*}
In particular, $\Psi(h) = Y(2s+1)$. The solution $Y(\tau)$ of (\ref{eq:Yode}) satisfies
\begin{equation}
\label{eq:Y(tau)}
Y(\tau) = I + h \int_0^\tau  A(\tau_1) Y(\tau_1) d\tau_1. 
\end{equation}
From that, one obtains 
\[
Y(\tau) = I + h \int_0^\tau A(\tau_1) d\tau_1 
+ h^2 \int_0^\tau \int_0^{\tau_1}  A(\tau_1) A(\tau_2)Y(\tau_2)  d\tau_2 d\tau_1
\]
and more generally,
\begin{equation}
\label{eq:Yodesol}
Y(\tau) =  I + \sum_{n=1}^{m} h^n \int_0^\tau \int_0^{\tau_1} \cdots \int_0^{\tau_{n-1}} 
A(\tau_1) \cdots A(\tau_n) d\tau_n \cdots d\tau_1
+ h^{m+1} \mathcal{R}_{m+1}(\tau,h),
\end{equation}
where for each $n$, the remainder $\mathcal{R}_{n}(\tau,h)$ satisfies 
\begin{equation}
\label{eq:Yoderem}
\mathcal{R}_{n}(\tau,h) = \int_0^\tau \int_0^{\tau_1} \cdots \int_0^{\tau_{n-1}} A(\tau_1)  \cdots A(\tau_n) Y(\tau_n) d\tau_n \cdots d\tau_1.
\end{equation}

By substituting $A(\tau_j) = d_1(\tau_j) F_1 + d_2(\tau_j) F_2$, then expanding all the products and taking constant linear operators out from integral signs, we obtain 
\[
Y(\tau) =  I + \sum_{n=1}^{m} h^n \sum_{\ell_1,\ldots,\ell_n \in \{1,2\}}
\alpha_{\ell_1 \cdots \ell_n}(\tau) \, F_{\ell_1} \cdots F_{\ell_n}
+ h^{m+1}\mathcal{R}_{m+1}(\tau,h),
\]
where 
\begin{equation}
\label{eq:alpha_iterated_integral}
\alpha_{\ell_1 \cdots \ell_n}(\tau) =  \int_0^\tau \int_0^{\tau_1} \cdots \int_0^{\tau_{n-1}} d_{\ell_1}(\tau_1)\cdots  d_{\ell_n}(\tau_n) d\tau_n \cdots d\tau_1.
\end{equation}
In particular, we have
\begin{equation*}
\begin{split}
\Psi(h) &=  I + \sum_{n=1}^{m} h^n \sum_{\ell_1,\ldots,\ell_n \in \{1,2\}}
u_{\ell_1 \cdots \ell_n}(a_1,b_1,\ldots,a_s,b_s,a_{s+1}) \, F_{\ell_1} \cdots F_{\ell_n}\\
&+ h^{m+1}\mathcal{R}_{m+1}(2s+1,h),
\end{split}
\end{equation*}
where $u_{\ell_1 \cdots \ell_n}(a_1,b_1,\ldots,a_s,b_s,a_{s+1}) = \alpha_{\ell_1 \cdots \ell_n}(2s+1)$, that is,
\[
\begin{aligned}
 & u_{\ell_1 \cdots \ell_n}(a_1,b_1,\ldots,a_s,b_s,a_{s+1}) = \\
  & \qquad\qquad   \int_0^{2s+1} \int_0^{\tau_1} \cdots \int_0^{\tau_{n-1}} d_{\ell_1}(\tau_1)\cdots  d_{\ell_n}(\tau_n) d\tau_n \cdots d\tau_1.
\end{aligned}
\]
We now consider (\ref{eq:Yode}) with $d_1(\tau) \equiv d_2(\tau) \equiv 1$. Clearly, in that case, $Y(\tau)=\e^{\tau h (F_1 + F_2)}$, and 
\[
\begin{aligned}
\e^{h (F_1 + F_2)} &=  I + \sum_{n=1}^{m} h^n \int_0^1 \int_0^{\tau_1} \cdots \int_0^{\tau_{n-1}} (F_1 + F_2)^n d\tau_n \cdots d\tau_1
+ h^{m+1} \overline{\mathcal{R}}_{m+1}(1,h),\\
&= I + \sum_{n=1}^{m} h^n \sum_{\ell_1,\ldots,\ell_n \in \{1,2\}}
\frac{1}{n!} \, F_{\ell_1} \cdots F_{\ell_n}
+ h^{m+1}\overline{\mathcal{R}}_{m+1}(1,h),
\end{aligned}
\]
where for each $n$, the remainder $\overline{\mathcal{R}}_{n}(\tau,h)$ is given by
\begin{equation*}
\overline{\mathcal{R}}_{n}(\tau,h) = \int_0^\tau \int_0^{\tau_1} \cdots \int_0^{\tau_{n-1}} (F_1+F_2)^n 
\e^{\tau_n h (F_1+F_2)}  d\tau_n \cdots d\tau_1.
\end{equation*}
We finally arrive at the following expression for the local error of the splitting scheme (\ref{eq:Yode}):
\[
\begin{aligned}
 & \Psi(h) - \e^{h (F_1 + F_2)} =  \\ 
 & \quad \sum_{n=1}^{m} h^n \sum_{\ell_1,\ldots,\ell_n \in \{1,2\}}
\left( u_{\ell_1 \cdots \ell_n}(a_1,b_1,\ldots,a_s,b_s,a_{s+1}) - \frac{1}{n!} \right)\, F_{\ell_1} \cdots F_{\ell_n}\\
& \quad + h^ {m+1} \left( \mathcal{R}_{m+1}(2s+1,h) - \overline{\mathcal{R}}_{m+1}(1,h) \right).
\end{aligned}
\]
Hence, if the scheme is of order $r$ (that is, (\ref{eq:ocond_ab})
 holds for each word $(\ell_1,\ldots,\ell_n)$ with $n\leq r$ letters in the alphabet $\{1,2\}$), then
\begin{equation*}\Psi(h) - \e^{h (F_1 + F_2)}= h^ {n+1} \left( \mathcal{R}_{n+1}(2s+1,h) - \overline{\mathcal{R}}_{n}(1,h) \right).
\end{equation*}

\subsubsection{Iterated integrals and shuffle relations}
\label{sss:shuffle}

Iterated integrals of the form (\ref{eq:alpha_iterated_integral}) were first considered and studied in \cite{chen57iop}.
It is well known that the integration-by-parts formula gives (for an arbitrary integrable path $(d_1(\tau),d_2(\tau))$) the relations
\begin{align*}
\alpha_{\ell_1}(\tau) \alpha_{\ell_2}(\tau) &= \alpha_{\ell_1\ell_2}(\tau) + \alpha_{\ell_2\ell_1}(\tau), \\
\alpha_{\ell_1}(\tau) \alpha_{\ell_2\ell_3}(\tau) &= \alpha_{\ell_1\ell_2\ell_3}(\tau) + \alpha_{\ell_2\ell_1\ell_3}(\tau) + \alpha_{\ell_2\ell_3\ell_1}(\tau),
\end{align*}
and more generally, 
\begin{equation}
\label{eq:shuffle_relations}
\alpha_{\ell_1 \cdots \ell_n}(\tau) \alpha_{\ell_{n+1}\cdots \ell_{n+m}}(\tau) = \sum_{\sigma \in \mathrm{Sh}(n,m)} \alpha_{\ell_{\sigma(1)} \cdots  \ell_{\sigma(n+m)}}(\tau),
\end{equation}
where $\mathrm{Sh}(n,m)$ is the set of the $(n+m)!/(n!m!)$ permutations $\sigma$ of $(1,\ldots,n+m)$ that are obtained by interleaving $(1,\ldots,n)$ and $(n+1,\ldots,n+m)$ while preserving their respective ordering.  

It will be useful to interpret the relations (\ref{eq:shuffle_relations}) in terms of the so-called shuffle product of words:  the shuffle product $\shuffle$ of two words $\ell_1\cdots \ell_n$ and $\ell_{n+1}\cdots \ell_{n+m}$ is defined as the following formal sum of words,
\begin{equation*}
\ell_1\cdots \ell_n \shuffle \ell_{n+1}\cdots \ell_{n+m} =  \sum_{\sigma \in \mathrm{Sh}(n,m)} \ell_{\sigma(1)} \cdots \ell_{\sigma(n+m)}.
\end{equation*}
By extending $\alpha_w(\tau)$ linearly to the case where $w$ is a linear combination of words, 
the relations (\ref{eq:shuffle_relations}) can be interpreted as
\begin{equation*}
\alpha_{w}(\tau) \alpha_{w'}(\tau) = \alpha_{w  \shuffle  w'}
\end{equation*}
for arbitrary words $w=\ell_1\cdots \ell_n$ and $w'=\ell_{n+1}\cdots \ell_{n+m}$ in the alphabet $\mathcal{A} = \{1,2\}$.

The shuffle product $\shuffle$ defines a commutative algebra (the so-called shuffle algebra) over the vector space of formal linear combinations of words in the alphabet $\mathcal{A}$. The shuffle algebra is freely generated by the set of {\em Lyndon words}~\cite{reutenauer93fla}. 

The fact that the coefficients of the series expansions of both $\e^{h(F_1+F_2)}$ and $\Psi(h)$  satisfy the shuffle relations, together with the fact that the set of Lyndon words freely generate the shuffle algebra, implies that
 a set of independent conditions for a consistent splitting scheme to attain order $r$ can be obtained by considering (\ref{eq:ocond_ab})  for each Lyndon word $(\ell_1,\ldots,\ell_n)$ of length $n\leq r$.

\subsection{Order conditions III: splitting methods with Lyndon multi-indices}
\label{subsec.2.4}

   In~\cite{blanes13nfo}, yet another characterization of the order conditions is obtained in terms of explicitly given polynomial equations. We next describe
   this alternative formulation. To do that, we always assume that the consistency conditions (\ref{eq:consistency_cond}) hold, so that the method is at least of first order. In this case, the polynomial equations are expressed in terms of the coefficients $b_1,\ldots,b_{s}$ and the coefficients $c_1,\ldots,c_s$ given by 
\begin{equation}
  \label{eq:ci}
c_i = \sum_{j=1}^{i} a_{j},  \qquad i=1,2,\ldots,s.
\end{equation}

We begin by rewriting (\ref{eq:Psi(h)}) as
\[
\begin{aligned} 
  &  \Psi(h) = \e^{a_{s+1} h F_1} \, \e^{b_s h F_2} \,  \e^{a_s h F_1} \cdots  \e^{b_1 h F_2} \,  \e^{a_1 h F_1}\\
  & \qquad =  \e^{h F_1} \, \left( \prod_{j=s}^1 \e^{-c_j h F_1} \, \e^{b_j h F_2} \, \e^{c_j h F_1}\right)
  = \e^{h F_1} \e^{b_s h C(c_s h)} \cdots \e^{b_1 h C(c_1 h)},
\end{aligned}
\]
where
\begin{equation}
\label{eq:C(h)}
C(h) = \e^{ -h F_1} \,  F_2\, \e^{h F_1}  = \sum_{n=1}^{\infty} h^{n-1} C_n = C_1 + h C_2 + h^2 C_3 + h^ 3 C_4 + \cdots,
\end{equation}
with $C_1 = F_2$, and 
\[ 
  C_k = \frac{1}{(k-1)} [C_{k-1}, F_1] \qquad \mbox{ for }  \, k > 1.
\]  
Now the order of the scheme (\ref{eq:splitting}) is established by comparing the expansion in powers of $h$ of $\e^{-h F_1} \Psi(h)$ with that of $\e^{-h F_1}\e^{h (F_1 + F_2)}$. 

Clearly,  $Y(\tau) := \e^{-h \tau F_1} \e^{\tau h (F_1 + F_2)}$ is the solution of  
 (\ref{eq:Yode}) with $A(\tau) = C(\tau h)$. Since the solution $Y(\tau)$ of (\ref{eq:Yode}) admits the representation (\ref{eq:Yodesol}) with remainder (\ref{eq:Yoderem}), and  $\e^{-h F_1} \e^{h (F_1 + F_2)} = Y(1)$ we conclude that 
\begin{equation*}
\begin{split}
\e^{-h F_1} \e^{h (F_1 + F_2)}
 &=  I + \sum_{k=1}^{n} h^k \int_0^1 \int_0^{\tau_k} \cdots \int_0^{\tau_2} 
C(\tau_k h) \cdots C(\tau_1 h) d\tau_1 \cdots d\tau_k\\
&+ h^{n+1} \mathcal{R}_{n+1}(1,h),
\end{split}
\end{equation*}
where for each $k$, 
\begin{equation*}
\mathcal{R}_k(1,h) = \int_0^1 \int_0^{\tau_k} \cdots \int_0^{\tau_2} C(\tau_k h) \cdots C(\tau_1 h) Y(\tau_1) d\tau_1 \cdots d\tau_k.
\end{equation*}
By substitution of (\ref{eq:C(h)}), we obtain that for each $k$,
\begin{gather*}
h^k \int_0^1 \int_0^{\tau_k} \cdots \int_0^{\tau_2} 
C(\tau_k h)  \cdots C(\tau_1 h) d\tau_1 \cdots d\tau_k \\
=  \sum_{i_1,\ldots,i_k \geq 1}h^{i_1+\cdots +i_k}  
 \left( \int_0^1 \int_0^{\tau_k} \cdots \int_0^{\tau_2}
 \tau_1^{i_1-1} \cdots \tau_k^{i_k-1}
  d\tau_1 \cdots d\tau_k
 \right)
   \,C_{i_k} \cdots C_{i_1} \\
   =  \sum_{i_1,\ldots,i_k \geq 1}  
    \frac{h^{i_1+\cdots +i_k}}{(i_1+\cdots + i_k)\cdots (i_1+i_2) i_1} \, C_{i_k} \cdots C_{i_1}.
\end{gather*}

As for $\e^{-h F_1}\, \Psi(h)$,
\begin{equation}
\begin{aligned} 
   &  \e^{-h F_1}\, \Psi(h)=   \prod_{j=s}^1 \e^{b_j h C(c_j h)} = 
   \prod_{j=s}^1 \left(I + \sum_{k\geq 1} \frac{h^k b_j^k}{k!} C(c_j h)^k \right) \\
   & \quad = I + \sum_{k\geq 1} \frac{h^k}{k!} \sum_{j=1}^s b_j^k C(c_j h)^k 
   + \sum_{k_1, k_2\geq 1} \frac{h^{k_1+k_2}}{k_1! k_2!} \sum_{1\leq j_1 < j_2 \leq s} b_{j_1}^{k_1} b_{j_2}^{k_2} C(c_{j_2} h)^{k_2}  C(c_{j_1} h)^{k_1}\\
      & \qquad + \sum_{k_1, k_2, k_3\geq 1} \frac{h^{k_1+k_2+k_3}}{k_1! k_2!k_3!} \sum_{1\leq j_1 < j_2 < j_3\leq s} b_{j_1}^{k_1} b_{j_2}^{k_2} b_{j_3}^{k_3} C(c_{j_3} h)^{k_3}C(c_{j_2} h)^{k_2}C(c_{j_1} h)^{k_1} + \cdots \\
   \label{eq:Psi(h)_expansion}
       & \quad = I  +\sum_{k\geq 1}  h^k \sum_{1\leq j_1 \leq \cdots \leq j_k \leq s} \frac{b_{j_1}\cdots b_{j_k}}{\sigma(j_1, \ldots j_k)}\, C(c_{j_k}h) \cdots C(c_{j_1}h)
 \end{aligned}
 \end{equation}
 where
\begin{eqnarray*}
      \sigma(j_1,\ldots, j_k)=1 &\mbox{ if }& j_1 < \cdots < j_k,\\
  \sigma(j_1,\ldots, j_k)=\ell ! \, \sigma(j_{\ell+1},\ldots, j_k) &\mbox{ if }& j_1=\cdots = j_{\ell} < j_{\ell+1} \leq \cdots \leq j_k.
\end{eqnarray*}
Since 
\begin{equation*}
 h^k C(c_{j_k}h) \cdots C(c_{j_1}h) = \sum_{i_1,\ldots,i_k \geq 1}h^{i_1+\cdots +i_k}  c_{j_1}^{i_1-1} \cdots c_{j_k}^{i_k-1}  \,C_{i_k} \cdots C_{i_1},
\end{equation*}
we arrive at
\begin{equation}
\label{eq:Psi(h)_multiindex}
  \e^{-h F_1}\,  \Psi(h) =  I +  \sum_{k\geq 1}
\sum_{i_1,\ldots,i_k \geq 1} h^{i_1+\cdots +i_k} \, v_{i_1,\ldots,i_k}(b_1,c_1,\ldots,b_s,c_s)
 \,C_{i_k} \cdots C_{i_1},
\end{equation}
where
\begin{equation}
\label{eq:v_multiindex}
v_{i_1,\ldots,i_k}(b_1,c_1,\ldots,b_s,c_s) 
= \sum_{1\leq j_1 \leq \cdots \leq j_k \leq s} \frac{b_{j_1}\cdots b_{j_k}}{\sigma(j_1, \ldots, j_k)}\,
c_{j_1}^{i_1-1} \cdots c_{j_k}^{i_k-1}.
\end{equation}
In this way, a consistent splitting method is at least of order $r$ if and only if 
\begin{equation}
\label{eq:ocond}
\sum_{1\leq j_1 \leq \cdots \leq j_k \leq s} \frac{b_{j_1}\cdots b_{j_k}}{\sigma(j_1, \ldots, j_k)}\,
c_{j_1}^{i_1-1} \cdots c_{j_k}^{i_k-1} =  \frac{1}{(i_1+\cdots + i_k)\cdots (i_1+i_2) i_1}
\end{equation}
 holds for each multi-index $(i_1,\ldots,i_k)$ such that $i_1+\cdots+i_k\leq r$.   For illustration, in Table \ref{tab.1} we give explicit conditions (\ref{eq:ocond}) corresponding to multi-indices with up to three indices.

However, such order conditions are not independent. The situation is very similar to that of the previous subsection:  instead of series indexed by the set of words in the alphabet $\{1,2\}$, now they are indexed by the set of words in the alphabet $\mathbb{N}=\{1,2,3,4,\ldots\}$. To distinguish the words of both sets, we will keep referring to the words in the alphabet $\mathbb{N}$ as multi-indices, and will write them as $(i_1,\ldots,i_k)$ instead of $i_1\cdots i_k$. Analogously to the previous subsection, the corresponding coefficients of the series expansions of both $ \e^{-h F_1} \Psi(h)$ and $\e^{-h F_1} \e^{h (F_1+F_2)}$ satisfy the shuffle relations
\begin{equation}
\label{eq:shuffle_relations2}
v_{i_1,\ldots,i_n} v_{i_{n+1},\ldots,i_{n+m}} = \sum_{\sigma \in \mathrm{Sh}(n,m)} v_{i_{\sigma(1)}, \dots, i_{\sigma(n+m)}}.
\end{equation}
This can be seen by showing that such coefficients can be written in both cases as iterated integrals.  From the discussion above, we already know that, for each multi-index $(i_1,\ldots,i_k)$,  the coefficients
\begin{gather*}
   \frac{1}{(i_1+\cdots + i_k)\cdots (i_1+i_2) i_1} =  \int_0^1 \int_0^{\tau_k} \cdots \int_0^{\tau_2}
 \tau_1^{i_1-1} \cdots \tau_k^{i_k-1}
  d\tau_1 \cdots d\tau_k
\end{gather*}
of the series expansion of $\e^{-h F_1} \e^{h (F_1+F_2)}$ are indeed iterated integrals.

We next show that the coefficients (\ref{eq:v_multiindex}) of the series expansion (\ref{eq:Psi(h)_multiindex}) can also be defined as iterated integrals.

We begin by showing that the infinite series expansion (\ref{eq:Psi(h)_expansion}) of $\e^{-h F_1} \, \Psi(h)$ can be represented as a truncated series plus a remainder. For that purpose, we consider (\ref{eq:Yode}) with $A(\tau)$ defined as follows:
\begin{equation}
\label{eq:A(tau)2}
A(\tau) = 
\left\{
\begin{matrix}
b_i C(h c_i) & \mbox{ if } & \tau \in [i-1,i], \ i \in \{1,\ldots,s\}, \\
0 & \mbox{ if } & \tau >s.
\end{matrix}
\right.
\end{equation}
In that case, if $\tau \in [i,i+1]$ with $i<s$,
\begin{equation*}
Y(\tau) =  
\e^{(\tau-i) h b_{i+1} C(h c_{i+1})}  \e^{h b_i C(h c_i)} \cdots \e^{h b_1 C(h c_1)} 
\end{equation*}
and in particular, $\e^{-h F_1} \, \Psi(h) = Y(s)$.  Since the solution $Y(\tau)$ of (\ref{eq:Yode}) admits the representation (\ref{eq:Yodesol}) with remainder (\ref{eq:Yoderem}), then
\begin{equation*}
\begin{split}
\e^{-h F_1} \Psi(h) &=  I + \sum_{k=1}^{n} h^k \int_0^s \int_0^{\tau_k} \cdots \int_0^{\tau_2} 
A(\tau_k) A(\tau_{k-1}) \cdots A(\tau_1) d\tau_1 \cdots d\tau_k\\
&+ h^{n+1} \mathcal{R}_{n+1}(s,h),
\end{split}
\end{equation*}
where for each $n$, 
\begin{equation*}
\mathcal{R}_{n}(s,h) = \int_0^s \int_0^{\tau_n} \cdots \int_0^{\tau_2} A(\tau_n) A(\tau_{n-1}) \cdots A(\tau_1) Y(\tau_1) d\tau_1 \cdots d\tau_n.
\end{equation*}

By comparison with (\ref{eq:Psi(h)_expansion}), we conclude that
\begin{gather*}
\sum_{1\leq j_1 \leq \cdots \leq j_k \leq s} \frac{b_{j_1}\cdots b_{j_k}}{\sigma(j_1, \ldots j_k)}\, C(c_{j_k}h) \cdots C(c_{j_1}h) \\ =  \int_0^s \int_0^{\tau_k} \cdots \int_0^{\tau_2} 
A(\tau_k) A(\tau_{k-1}) \cdots A(\tau_1) d\tau_1 \cdots d\tau_k
\end{gather*}

Now, $A(\tau)$ can be written as
\begin{equation*}
A(\tau) = d_1(\tau) C_1+ h d_2(\tau) C_2+ h^2 d_3(\tau) C_3+ \cdots
\end{equation*}
where $d_j(\tau) = b_i c_i^{j-1}$ if $\tau \in [i,i+1]$ with $i<s$. Proceeding as in previous section, one obtains 
\begin{equation}
v_{i_1,\ldots,i_k}(b_1,c_1,\ldots,b_s,c_s) 
= \int_0^{s} \int_0^{\tau_1} \cdots \int_0^{\tau_{k-1}} d_{i_1}(\tau_1)\cdots  d_{i_k}(\tau_k) d\tau_k  \cdots d\tau_1,
\end{equation}
which implies that the shuffle relations (\ref{eq:shuffle_relations2}) hold for the polynomials $v_{i_1,\ldots,i_k}$.

As in the previous subsection, a set of independent conditions that imply the order conditions (\ref{eq:ocond}) can be obtained by considering (\ref{eq:ocond})  for each Lyndon multi-index $(i_1,\ldots,i_k)$ such that $1<i_1+\cdots+i_k\leq r$. Here, we exclude the multi-index $(1)$ as in that case (\ref{eq:ocond}) coincides with the second equality in the consistency condition (\ref{eq:consistency_cond}).  For instance, the subset of Lyndon multi-indices $(i_1,\ldots,i_k)$ such that $1<i_1+\cdots+i_k\leq 5$ is
$$\{(2),(3),(4),(5),(1,2),(1,3),(1,4),(2,3),(1,1,2),(1,1,3),(1,2,2),(1,1,1,2)\}.$$

\begin{table}[h!]
\begin{center}
\begin{tabular}{|c|c|} \hline
  Multi-index  &  \hspace*{0.4cm}  Condition\\ \hline \hline
      $(i)$ & $\displaystyle \sum_{j=1}^{s} b_j \, c_j^{i-1}=\frac{1}{i}$ \\ \hline
      $(i_1,i_2)$   &  $\displaystyle  \frac{1}{2} \sum_{j=1}^{s} b_j^2 c_j^{i_1+i_2-2} + \sum_{1\leq j_1 < j_2 \leq s} b_{j_1} b_{j_2} c_{j_1}^{i_1-1} c_{j_2}^{i_2-1} = \frac{1}{(i_1+i_2) i_1}$ \\  \hline
      $(i_1,i_2,i_3)$   &  $\displaystyle  \frac{1}{6} \sum_{j=1}^{s} b_j^3 c_j^{i_1+i_2+i_3-3} +  \frac{1}{2} \sum_{1\leq j_1 < j_3 \leq s} b_{j_1}^ 2 b_{j_3} c_{j_1}^ {i_1+i_2-2} c_{j_3}^{i_3-1}$ \\
&   $\displaystyle +  \frac{1}{2} \sum_{1\leq j_1 < j_2 \leq s} b_{j_1} b_{j_2}^2 c_{j_1}^{i_1-1} c_{j_2}^{i_2+i_3-2}$ \\
 & $\displaystyle+   \sum_{1\leq j_1 < j_2 < j_3\leq s} b_{j_1} b_{j_2} b_{j_3} c_{j_1}^{i_1-1} c_{j_2}^{i_2-1}c_{j_2}^{i_3-1}    = \frac{1}{(i_1+i_2+i_3)(i_1+i_2) i_1}$ \\  \hline
\end{tabular}
\end{center}
\caption{Conditions (\ref{eq:ocond}) corresponding to multi-indices with up to three indices.\label{tab.1}}
\end{table}

For time-symmetric splitting methods, 
a set of independent order conditions will be obtained by considering (\ref{eq:ocond}) restricted to Lyndon multi-indices $(i_1,\ldots,i_k)$ with odd weight $i_1+\cdots+i_k$.
 For instance, the subset of Lyndon multi-indices $(i_1,\ldots,i_k)$ such that $1<i_1+\cdots+i_k\leq 5$ with odd weight $i_1+\cdots+i_k$
 is
$$\{(3),(5),(1,2),(1,4),(2,3),(1,1,3),(1,2,2),(1,1,1,2)\}.$$


Notice that the treatment carried out in this subsection may also be formally applied when $F_1$ is an unbounded operator. In that case, however, one
has to get rigorous estimates of the remainders to prove stability and convergence of the corresponding schemes, as done in \cite{thalhammer08hoe}, for example.

\subsection{Order conditions IV: composition methods with Lyndon multi-indices}
\label{subsec.2.5}

We now turn our attention to compositions (\ref{compm.1}) of a basic second-order time-symmetric scheme $S_h^{[2]}$ with 
appropriate coefficients $\gamma_1,\ldots,\gamma_s$ chosen to achieve higher orders. Of course, 
a set of conditions that guarantee that the scheme (\ref{compm.1}) attains a given order can be obtained by rewriting that composition in terms of basic maps, such as (\ref{eq:Psi(h)}),  and using the characterization of the order of the splitting method (\ref{eq:Psi(h)}) described in Subsection~\ref{subsec.2.4}. The corresponding parameters $a_j, b_j$ can be obtained in terms of $\gamma_1,\ldots,\gamma_s$ as follows:  $a_{1}=\gamma_{1}/2$, and 
 \begin{equation*}
a_{j+1} = \frac{\gamma_{j} + \gamma_{j+1}}{2}, \qquad\quad  b_j=\gamma_{j}, \qquad \mbox{ for } \, j=1,\ldots,s
\end{equation*}
with $\gamma_{s+1}=0$.
However, the resulting polynomial equations for any given order $r$, once written in terms of the coefficients $\gamma_1,\ldots,\gamma_s$, are no longer independent.  An alternative formulation of the order of (\ref{compm.1}) in terms of explicit independent algebraic equations in the 
coefficients $\gamma_1,\ldots,\gamma_s$ will be presented in Paragraph~\ref{sss:compmoc} below. This characterization is based on the treatment of
the more general composition (\ref{eq:compint}), which is treated next. 


\subsubsection{Order conditions of compositions of a basic method and its adjoint}
\label{sss:compintoc}

The composition (\ref{eq:compint}) is at least of order $r$ if 
$\Psi(h) - \exp(h (F_1 + F_2)) = \mathcal{O}(h^{r+1})$,
where $\Psi(h)$ is the associated Lie transformation (\ref{lt_comp1}).

To get the series expansion of $\Psi(h)$, we first consider the expansion in powers of $h$ of the Lie transformation $\mathcal{X}(h) = I + h X_1 + h^2 X_2+\cdots$ associated to the basic integrator $\chi_h$.  If $\chi_h$ is  the Lie--Trotter scheme, $\chi_h = \varphi_h^{[1]} \circ \varphi_h^{[2]}$, then 
\[
   \mathcal{X}(h) = \e^{h F_2} \e^{h F_1} = I + h (F_1 + F_2) + h^2 (\frac12 F_1^2 + F_2 F_1 + \frac12 F_2^2) + \cdots. 
\]   
To deal with the most general problem, however, from now on we only assume that
$\chi_h$ is a smooth consistent integrator, so that 
 $X_1 = F_1 + F_2$,  and  each $X_n$ can be defined so that for each smooth function $g$, $X_n g$ is a new smooth function given by (\ref{eq:X_n}).
 
Let us consider $\Psi_0(h) = I$, and for each $j\geq 1$, 
\begin{equation}
\label{eq:Psij}
\begin{split}
\Psi_{2j-1}(h) &=  \mathcal{X}(-\alpha_1 h)^{-1} \mathcal{X}(\alpha_2 h) \cdots  \mathcal{X}(-\alpha_{2j-1} h)^{-1}\\
\Psi_{2j}(h) &=  \mathcal{X}(-\alpha_1 h)^{-1} \mathcal{X}(\alpha_2 h) \cdots  \mathcal{X}(-\alpha_{2j-1} h)^{-1}  \mathcal{X}(\alpha_{2j} h),
\end{split}
\end{equation}
so that in particular, $\Psi(h) = \Psi_{2s}(h)$. Notice that
\begin{align*}
\mathcal{X}(-h)^{-1} &= I + \sum_{k \geq 1} (-1)^k (-h X_1 + h^2 X_2- h^3 X_3 + \cdots)^k \\
&= I + h X_1 + h^2 (X_1^2- X_2) + h^3 (X_1^3 - X_1 X_2 - X_2 X_1 + X_3) + \cdots,
\end{align*}
which implies that for each $k\geq 1$, there exist polynomials $w_{i_1,\ldots,i_m}(\alpha_1,\ldots,\alpha_k)$ on the coefficients $\alpha_1,\ldots,\alpha_k$ such that
\begin{equation}
\label{eq:Psik}
\Psi_k(h) = I + \sum_{n\geq 1} h^n \sum_{m\geq 1} \sum_{i_1+\cdots+i_m=n} w_{i_1,\ldots,i_m}(\alpha_1,\ldots,\alpha_k) \, X_{i_1} \cdots X_{i_m}.
\end{equation}

We next determine such polynomial coefficients $w_{i_1,\ldots,i_m}(\alpha_1,\ldots,\alpha_k)$ recursively from the relations $\Psi_{1}(h)  = \mathcal{X}(-\alpha_{2j-1} h)^{-1}$ and 
\begin{equation*}
\Psi_{2j-1}(h)  = \Psi_{2j-2}(h) \mathcal{X}(-\alpha_{2j-1} h)^{-1}, \quad
 \quad \Psi_{2j}(h)  = \Psi_{2j-1}(h)  \mathcal{X}(\alpha_{2j} h),
\end{equation*}
or equivalently, $\Psi_{1}(h) \mathcal{X}(-\alpha_{1} h)  = I$ and
\begin{equation}
\label{eq:Psi_rec}
\Psi_{2j-1}(h) \mathcal{X}(-\alpha_{2j-1} h)  = \Psi_{2j-2}(h), \quad
 \quad \Psi_{2j}(h)  = \Psi_{2j-1}(h)  \mathcal{X}(\alpha_{2j} h). 
\end{equation}
Specifically, for arbitrary coefficients $w_{i_1,\ldots,i_m}$   and $\lambda$ one has
\[
\begin{aligned}
& \left(I + \sum_{n\geq 1} h^n \sum_{m\geq 1} \sum_{i_1+\cdots+i_m=n} w_{i_1,\ldots,i_m}  \, X_{i_1} \cdots X_{i_m}\right) \mathcal{X}(\lambda h)
= I  + \sum_{i\geq 1} h^i  (w_{i} + \lambda^i )  X_n \\
& \qquad\qquad + \sum_{n\geq 1} h^n \sum_{m\geq 2} \sum_{i_1+\cdots+i_m=n} (w_{i_1,\ldots,i_m} +
\lambda^{i_m} w_{i_1,\ldots,i_{m-1}} ) \, X_{i_1} \cdots X_{i_m},
\end{aligned}
\]
so that, taking this expression into account, (\ref{eq:Psi_rec}) leads to the following identities
\begin{align*}
 w_i(\alpha_1) &=  - (-\alpha_{1})^i, \\ 
 w_i(\alpha_1,\ldots,\alpha_{2j-1}) &= w_i(\alpha_1,\ldots,\alpha_{2j-2}) - (-\alpha_{2j-1})^i, \\
 w_i(\alpha_1,\ldots,\alpha_{2j})  &= w_i(\alpha_1,\ldots,\alpha_{2j-1}) + \alpha_{2j}^i, \\
 w_{i_1,\ldots,i_m}(\alpha_1) &= - (-\alpha_{1})^{i_m}  w_{i_1,\ldots,i_{m-1}}(\alpha_1), \\
 w_{i_1,\ldots,i_m}(\alpha_1,\ldots,\alpha_{2j-1}) &= w_{i_1,\ldots,i_m}(\alpha_1,\ldots,\alpha_{2j-2}) \\
  &- (-\alpha_{2j-1})^{i_m}  w_{i_1,\ldots,i_{m-1}}(\alpha_1,\ldots,\alpha_{2j-1}), \\ 
  w_{i_1,\ldots,i_m}(\alpha_1,\ldots,\alpha_{2j}) &= w_{i_1,\ldots,i_m}(\alpha_1,\ldots,\alpha_{2j-1}) \\
  &+ \alpha_{2j}^{i_m}  w_{i_1,\ldots,i_{m-1}}(\alpha_1,\ldots,\alpha_{2j-1}).
\end{align*}
Clearly, (\ref{eq:Psik}) holds for the  coefficients  $w_{i_1,\ldots,i_m}(\alpha_1,\ldots,\alpha_k)$  determined by the relations above.  Equivalently, the functions $w_{i_1,\ldots,i_m}$ can be defined as
\begin{equation}
  \label{eq:wfun}
  \begin{split}
w_{i}(\alpha_1,\ldots,\alpha_{2\ell})&=
   \sum_{j =1}^s  \big(\alpha_{2j}^{i} - (-\alpha_{2j-1})^{i} \big), \\
w_{i_1,\ldots,i_m}(\alpha_1,\ldots,\alpha_{2\ell})&=
    \sum_{j =1}^s
\big(\alpha_{2j}^{i_m} - (-\alpha_{2j-1})^{i_m} \big) w_{i_1,\ldots,i_{m-1}}(\alpha_1,\ldots,\alpha_{2j-1}).
\end{split}
\end{equation}

Notice that  
\begin{equation*}
w_{i_1,\ldots,i_m}(\alpha_1,\ldots,\alpha_k,0,\ldots,0) = w_{i_1,\ldots,i_m}(\alpha_1,\ldots,\alpha_k),
\end{equation*}
 as expected from  (\ref{eq:Psij}),  (\ref{eq:Psik}),  and $\mathcal{X}(0)=I$. 
 
Comparing the series expansion of $\Psi(h)=\Psi_{2s}(h)$ with
\begin{equation*}
\exp(h \, (F_1 + F_2)) = \exp(h\, X_1) = I + \sum_{n\geq 1} \frac{h^n}{n!} X_1^n,
\end{equation*}
we finally conclude that the scheme (\ref{eq:compint}) is of order at least $r$ if and only if 
\begin{equation}
\label{eq:ocond3}
w_{i_1,\ldots,i_m}(\alpha_1,\ldots,\alpha_{2s})=
\left\{
\begin{array}{ll}
\frac{1}{m!} & \mbox{ if} \quad (i_1,\ldots,i_m) = \stackrel{m}{\overbrace{(1,\ldots,1)}} \\ \\
0 & \mbox{ otherwise},
\end{array}
\right.
\end{equation} 
 for each multi-index $(i_1,\ldots,i_k)$ such that $1\leq i_1+\cdots+i_k\leq r$.

 Furthermore, the order conditions (\ref{eq:ocond3}) are not all independent. For instance, it is straightforward to check from (\ref{eq:wfun}) that $w_{i_1,i_2} + w_{i_2,i_1} + w_{i_1 + i_2} = w_{i_1} w_{i_2}$ for arbitrary indices $i_1,i_2$. In particular,  $2w_{1,1} + w_{2} = w_{1}^2$, which implies that if the order conditions  (\ref{eq:ocond3}) for the multi-indices $(1)$ and $(2)$ are satisfied, then the condition for the multi-index $(1,1)$ is automatically fulfilled. Actually, such dependences are similar to the shuffle relations (\ref{eq:shuffle_relations}) that hold for the coefficients $v_{i_1,\ldots,i_m}$ considered in Subsection~\ref{subsec.2.4}.  Indeed,
 \begin{equation}
\label{eq:quasi_shuffle_relations}
w_{i_1,\ldots,i_n} w_{i_{n+1},\ldots,i_{n+m}} = \sum_{\sigma \in \mathrm{Sh}(n,m)} w_{i_{\sigma(1)}, \dots, i_{\sigma(n+m)}} + \cdots,
\end{equation}
where $\cdots$ refers to sums of products of coefficients corresponding to multi-indices with $m-1$ or fewer indices. In fact,  as shown in~\cite{chartier09aat}, 
 the dependences  (\ref{eq:quasi_shuffle_relations}) are directly related to the quasi-shuffle product $*$ on the linear span of multi-indices introduced in \cite{hoffman00qsp}, and  
 due to such dependences, it is enough to consider
   (\ref{eq:ocond3}) for Lyndon multi-indices $(i_1,\ldots,i_m)$ such that $i_1+\cdots+i_k\leq r$.  That is, the scheme (\ref{eq:compint}) is of order at least $r$ if $\alpha_1 + \cdots+\alpha_{2s}=1$ and for each Lyndon multi-index $(i_1,\ldots,i_k)$ such that $1< i_1+\cdots+i_k\leq r$,
  \begin{equation}
w_{i_1,\ldots,i_m}(\alpha_1,\ldots,\alpha_{2s})=0.
\end{equation}
In particular, a method of order 3 must satisfy, besides consistency, the conditions $w_2 = w_3 = w_{1,2}=0$.

 This provides an alternative characterization of the order conditions of general splitting schemes (\ref{eq:splitting}) in terms of the coefficients $\alpha_j$ obtained from the method parameters $a_j,b_j$ from (\ref{eq:abalpha}). Most
importantly, this also allows one to characterize the order of the scheme (\ref{compm.1}) obtained by composing Strang maps. This will be presented next in Section~\ref{sss:compmoc}.

\subsubsection{Explicit characterization of the order conditions of scheme (\ref{compm.1}) } 
\label{sss:compmoc}

In \cite{murua99ocf}, such an explicit characterization 
is obtained in terms of a set of polynomials indexed by certain sets of rooted trees decorated by the set of odd positive integers. We now describe a related formulation in terms of polynomials indexed by the set of Lyndon multi-indices with odd indices based on the formalism developed in the previous subsection for 
the more general composition (\ref{eq:compint}).

We begin by considering (\ref{compm.1}) as the particular case of (\ref{eq:compint}) with 
\begin{equation}
\label{eq:alphagamma}
\alpha_{2j-1}=\alpha_{2j}=\gamma_j/2 \quad \mbox{for} \quad  j=1,\ldots,s.
\end{equation}
For each multi-index $(i_1,\ldots,i_k)$, we define the function $u_{i_1,\ldots,i_m}$ on the set of finite sequences $(\gamma_1,\ldots,\gamma_s)$ of real numbers as follows:
\begin{equation}
u_{i_1,\ldots,i_m}(\gamma_1,\ldots,\gamma_s) = 2^{i_1+\cdots+i_m-m} w_{i_1,\ldots,i_m}(\alpha_1,\ldots,\alpha_{2s})
\end{equation}
with $\alpha_{2j-1} = \alpha_{2j}=\gamma_j/2$.
Clearly, scheme (\ref{compm.1}) is of order at least $r$ if and only if $\gamma_1 + \cdots+\gamma_s = 1$ and
\begin{equation}
\label{eq:ocond2}
u_{i_1,\ldots,i_m}(\gamma_1,\ldots,\gamma_s)=0
\end{equation} 
 for each Lyndon multi-index $(i_1,\ldots,i_k)$ such that $1<i_1+\cdots+i_k\leq r$.   However,  by definition, $u_{i_1,\ldots,i_m}(\gamma_1,\ldots,\gamma_s)\equiv 0$ if $i_m$ is even. Moreover, 
 for any multi-index $(i_1,\ldots,i_k)$ with some even index, (\ref{eq:ocond2}) holds provided that it holds for every Lyndon multi-index with fewer indices. 
 
Therefore, scheme (\ref{compm.1}) is of order at least $r$ if and only if $\gamma_1 + \cdots+\gamma_s = 1$ and (\ref{eq:ocond2}) for each Lyndon multi-index $(i_1,\ldots,i_k)$ with odd indices such that $1<i_1+\cdots+i_k\leq r$.  
For instance, the set of Lyndon multi-indices $(i_1,\ldots,i_k)$ of odd indices such that $1<i_1+\cdots+i_k\leq 7$ is
\begin{equation*}
\{(3),(5),(7),(1,3),(1,5),(1,1,3),(1,1,5),(1,3,3),(1,1,1,3),(1,1,1,1,3)\}.
\end{equation*}
 The resulting number of order conditions, denoted as $m_n$, is gathered in Table \ref{table_orcon}.

For multi-indices $(i_1,\ldots,i_m)$ with odd indices, the functions $u_{i_1,\ldots,i_m}$ can be written more explicitly as follows: 
\begin{align}
  \label{eq:ufun1}
  u_{i}(\gamma_1,\ldots,\gamma_s)&=
    \sum_{j=1}^s  \gamma_j^i, \\
     u_{i_1, i_2}(\gamma_1,\ldots,\gamma_s)&=
    \sum_{j_2=1}^s  \gamma_{j_2}^{i_2}   \sum_{j_1=1}^{j_2 *}  \gamma_{j_1}^{i_1}, \\
u_{i_1, i_2, i_3}(\gamma_1,\ldots,\gamma_s)&= \sum_{j_3=1}^s  \gamma_{j_3}^{i_3} \sum_{j_2=1}^{j_3 *}  \gamma_{j_2}^{i_2}   \sum_{j_1=1}^{j_2*}  \gamma_{j_1}^{i_1}, \\
\label{eq:ufun4}
u_{i_1, i_2, i_3, i_4}(\gamma_1,\ldots,\gamma_s)&= \sum_{j_4=1}^s  \gamma_{j_4}^{i_4} \sum_{j_3=1}^{j_4*}  \gamma_{j_3}^{i_3} \sum_{j_2=1}^{j_3 *}  \gamma_{j_2}^{i_2}   \sum_{j_1=1}^{j_2*}  \gamma_{j_1}^{i_1}, 
\end{align}
and so on. Here, as in \cite{murua99ocf}, we have used the notation
\begin{equation*}
\sum_{j=1}^{k*} A_j =  \frac{A_k}{2} + 
\sum_{j=1}^{k-1} A_j
\end{equation*}
for $A_1,\ldots,A_k \in \mathbb{R}$.

For time-symmetric integration schemes (\ref{compm.1}), 
a set of independent conditions for even order $r$ will be obtained by considering (\ref{eq:ocond2}) restricted to Lyndon multi-indices $(i_1,\ldots,i_k)$ of odd indices and  odd weight $i_1+\cdots+i_k< r$. For instance, for order $8$, one only needs to consider the following subset of Lyndon multi-indices:
\begin{equation*}
\{(3),(5),(7),(1,1,3),(1,1,5),(1,3,3),(1,1,1,1,3)\}.
\end{equation*}
 The above characterization of the order of the scheme  (\ref{compm.1}) is also true in the more general case where the Strang map $S^{[2]}_h$ is replaced by an arbitrary time-symmetric integrator of order $2\ell$.  In that case, only Lyndon multi-indices with indices from the set $\{1, 2\ell+1, 2\ell+3, 2\ell+5,\ldots \}$ have to be taken into account.

\subsection{Negative time steps}

Splitting and composition methods of order $r \ge 3$ necessarily involve some negative coefficients. This can already be observed in the simple
triple jump scheme (\ref{suzu1}), and in fact has been established as a general theorem by \cite{goldman96nno,sheng89slp,suzuki91gto}. A simple
proof can be obtained as follows \cite{blanes05otn}: given the existing relationship between the splitting method (\ref{eq:splitting}) 
and the composition (\ref{eq:compint}),
it is clear that any splitting scheme of order $r \ge 3$ has to verify the condition $w_3 = 0$, where, by virtue of (\ref{eq:wfun}),
\begin{equation} \label{nts1}
   w_3(\alpha_1,\ldots, \alpha_{2s}) = \sum_{j=1}^s (\alpha_{2j}^3 - (-\alpha_{2j-1})^3), 
\end{equation}
where the coefficients $\alpha_j$  are related to $a_j, b_j$ via (\ref{eq:abalpha}).
Since, for all $x, y \in \mathbb{R}$, it is true that $x^3 + y^3 < 0$ implies $x + y < 0$, then there must
exist some $j \in \{1, \ldots, s\}$ in the sum of (\ref{nts1}) such that
\[
  \alpha_{2j-1}^3 + \alpha_{2j}^3 < 0  \qquad \mbox{ and thus } \qquad \alpha_{2j-1} + \alpha_{2j} =
  b_j < 0.
\]
Obviously, one can also write (by taking $\alpha_{0} = 0$, $\alpha_{2s+1}=0$)
\[
   w_{3}(\alpha_1,\ldots,\alpha_{2s}) =  \sum_{i=0}^{s} (\alpha_{2j}^3 + \alpha_{2j+1}^3) = 0
\]
just by grouping terms in a different way, and thus, by repeating the argument,
 there must exist some $j \in \{0, \ldots, s\}$ such that
\[
   \alpha_{2j} + \alpha_{2j+1} = a_j < 0.
\]
This proof clearly shows the origin of the existence of backward time steps: the equation $w_3=0$ can only be
satisfied if at least one $a_j$ \emph{and} one $b_j$ are negative. 


\section{Splitting methods for special problems}
\label{sect3}

Whereas the analysis carried out in Section \ref{sect2} is completely general, there are important problems arising in applications
whose particular structure allows one to simplify the treatment and design schemes without taking into account all the order conditions.
Some of them are reviewed in this section, where we also show how to adapt splitting methods to deal with explicitly time-dependent
systems.

\subsection{RKN splitting methods}
\label{subsec3.1}

Many differential equations of practical interest are of the form 
\begin{equation}  \label{rkn.1}
    y^{\prime\prime}= g(y), 
\end{equation}
where $y \in \mathbb{R}^{d}$ and $g: \mathbb{R}^{d}\longrightarrow \mathbb{R}^{d}$. An example in point corresponds to 
 Hamiltonian systems of the form $H(q,p) = T(p) + V(q)$, where the kinetic energy $T(p)$ is quadratic in the momenta
$p$, i.e., $T(p)=\frac12 p^T M^{-1} p$ for a constant invertible symmetric 
matrix $M$, and $V(q)$ is the potential. In that case, the
corresponding Hamiltonian system can be written in the form
(\ref{rkn.1}) with $y=q$, $g(y)=- M^{-1} \nabla V(y)$.

By transforming (\ref{rkn.1}) into a first order ODE system 
(of dimension $D=2d$) in the new variables $x=(y,v)^T$, with $v = y^\prime$,  it is clear that the resulting equation
\[
   x^{\prime} = \frac{d }{dt} \left( \begin{array}{c}
   				y \\
				v
			  \end{array} \right) = 
		 \left( \begin{array}{c}
		 		v \\
				g(y)
			 \end{array} \right)
\]
can be expressed as $x^{\prime} = 	f_1(x) + f_2(x)$, with		 		
\begin{equation} \label{eq:f=fa+fbRKN}
 f_1(y,v) = (0,g(y))^T, \qquad f_2(y,v) = (v,0)^T,
\end{equation}
and splitting methods of the form (\ref{eq:splitting}) can be applied, with the exact $h$-flows $\varphi^{[1]}_h$ and $\varphi^{[2]}_h$ given by
\begin{equation}   \label{rkn.2} 
  \varphi_h^{[1]}:  \; \left( \begin{matrix} y_0 \\ v_0 \end{matrix} \right) \longmapsto 
  \left( \begin{matrix} y_0 \\ v_0 + h \, g(y_0) \end{matrix} \right);
\qquad 
  \varphi_h^{[2]}:  \; \left( \begin{matrix} y_0 \\ v_0 \end{matrix} \right) \longmapsto 
  \left( \begin{matrix} y_0+ h \, v_0 \\ v_0 \end{matrix} \right).
\end{equation}
Just as the class of Runge--Kutta methods can be conveniently adapted to (\ref{rkn.1}) to get more efficient schemes (the so-called Runge--Kutta--Nystr\"om
or RKN methods) \cite{hairer93sod}, special splitting methods can also be designed to improve the accuracy whilst reducing 
the computational cost with respect to the general composition (\ref{eq:splitting}). For analogy, they are sometimes called splitting methods of RKN type. The key point here is that
the differential operators $F_1$ and $F_2$ associated with (\ref{eq:f=fa+fbRKN}) satisfy $[F_1,[F_1,[F_1,F_2]] = 0$ identically. In other words,
$F_{1112} = 0$ in (\ref{bch2}), which introduces linear dependences among higher order terms in the expansion of
$\log(\Psi(h))$ \cite{mclachlan02sm} and therefore contributes to a reduction in the number of order conditions.
{In the notation of Subsection~\ref{subsec.2.4}, $6 C_4 = -[F_1,[F_1,[F_1,F_2]]$, and thus we have that $C_k = 0$ for $k\geq 4$. Hence, the order conditions (\ref{eq:ocond})  for multi-indices $(i_1,\ldots,i_k)$ with some index $i_j \geq 4$ need not be considered in that case. In particular, a splitting scheme is at least of order $5$ provided that  (\ref{eq:ocond}) holds for the Lyndon multi-indices
$$\{(2),(3),(1,2),(1,3),(2,3),(1,1,2),(1,1,3),(1,2,2),(1,1,1,2)\}.$$
Here, we have excluded the Lyndon multi-indices $(4),(5),(1,4)$ from the set of Lyndon multi-indices $(i_1,\ldots,i_k)$  with $i_1+\cdots+ i_k\leq 5$. For order $6$, in addition, one has to consider (\ref{eq:ocond}) for the following Lyndon multi-indices:
$$
(1,1,1,1,2),(1,1,1,3),(1,2,3) ,(1,1,2,2),(1,3,2).
$$
They are obtained by excluding from the set of 9 Lyndon multi-indices with $i_1+\cdots+ i_k=6$ the Lyndon multi-indices $(2,4), (1,1,4),(1,5),(6)$.  
Similarly, for order $7$, in addition, one has to consider (\ref{eq:ocond}) for the following Lyndon multi-indices:
\begin{gather*}
(1,1,1,1,1,2),(1,1,1,1,3),(1,1,2,3) ,(1,1,1,2,2),(1,1,3,2),\\
(2,2,3),(1,2,1,3),(1,1,2,3),(1,3,3),(1,2,2,2).
\end{gather*}
For orders higher than $7$, more reductions of the order conditions occur, in addition to those obtained by excluding Lyndon multi-indices $(i_1,\ldots,i_k)$ having some index $i_j \geq 4$, due to additional dependencies among nested commutators of $C_1,C_2,C_3$. For instance, it is straightforward to check that $[C_3,[C_2,C_3]]$ vanishes identically, which implies that, for order 8, one can also exclude the Lyndon multi-index  $(2,3,3)$.

The class of problems for which the reduction in order conditions discussed above is in fact more general than (\ref{eq:f=fa+fbRKN}). 
In particular, it includes the situation where $f_1$ depends only on $y$ and $f_2(y,v)$ is linear in $v$. 
 For Hamiltonian systems, this generalization corresponds to $H(q,p) =  V(q) + T(q,p)$, where $T(q,p)$ is quadratic in $p$. 
In that case, the flow associated with
$T(q,p)$ should be easily computed for the RKN splitting methods to be advantageous. 

Actually, such a reduction of the order conditions also holds for certain PDEs, and, in particular, for
 the time-dependent Schr\"odinger equation considered in Subsection \ref{subsec.1.3b} (with $F_1 = \hat V$ and $F_2 = \hat T$), since  
 the corresponding graded Lie algebra is isomorphic to the classes of problems discussed above~\cite{mclachlan19tla}. In fact, in 
 \cite{mclachlan19tla} it is conjectured (and checked up to order 20) that the case  (\ref{eq:f=fa+fbRKN}) (where $f_2(y,v)=v$) and the more general case (where $f_2(y,v)$ is linear in $v$ but may depend on $y$) give rise to the same reduction in the order conditions. 

 The actual number $d_r$ of order conditions for orders $r \le 11$ are given
in Table \ref{table_orcon}.  Since the order conditions up to order three are identical as in the general case, the results for negative time steps still apply.

This reduction in the number of order conditions allows one to design schemes involving a smaller number of elementary flows than in the general case, 
eventually leading to a greater efficiency. We have already illustrated a popular 4th-order method within this class (scheme RKN$_6$4 in the examples of Section \ref{sect1}.

\subsection{Methods with commutators}
\label{nested-commut}

Another possible way to improve the efficiency of splitting methods consists in incorporating into the scheme 
not only the flows of $F_1$ and $F_2$, but also the flows
of some of their commutators $[F_1,F_2]$, $[[F_1,F_2],F_1]$, etc., or convenient approximations of these flows. Of course, the strategy makes sense if the
gain in accuracy, stability or any other favorable property compensates the extra computational cost due to the presence of these additional flows. A popular 4th-order method belonging to this class
is
 \begin{equation}
  \label{eq:modpotex}
  \psi_h =
 \varphi_{a_1 h}^{[1]} \circ  \varphi_{b_1 h}^{[2]} \circ  \varphi_{a_2 h}^{[1]} \circ \varphi_{d_2 h^3}^{[112]} \circ \varphi_{a_2 h}^{[1]} \circ \varphi_{b_1 h}^{[2]} \circ \varphi_{a_1 h}^{[1]},
\end{equation}
with 
\begin{equation}
\label{eq:koseleff}
a_1 = \frac16, \quad b_1 = \frac12, \quad a_2 = \frac13, \quad d_2 = -\frac{1}{72},
\end{equation}
first proposed by \cite{koseleff93cfp} and \cite{chin97sif}. Here $\varphi^{[112]}_{h}$ denotes the $h$-flow corresponding to $F_{112} = [F_1,[F_1,F_2]]$ (or $F_{112}=2C_3$ in the notation of Subsection~\ref{subsec.2.4}).  

We next analyze the family of schemes \eqref{eq:modpotex}.
The Lie transformation $\Psi(h)$  is
\begin{equation*}
\Psi(h) = \e^{h a_1 F_1} \e^{h b_1 F_2} \e^{h a_2 F_1} \e^{h^3 2 d_2 C_3} \e^{h a_2 F_1} \e^{h b_1 F_2} \e^{h a_1 F_1}.
\end{equation*}
Under the assumption that $a_1 + a_2 = \frac12$, we have 
\begin{equation*}
\e^{-h F_1} \Psi(h) =   \e^{h b_1 C(h c_3)} \e^{h^3 d_2 D(h c_2)}  \e^{h b_1 C(h c_1)} ,
\end{equation*}
where $c_1 = a_1$, $c_2 = a_1 + a_2$, $c_3 = a_1+ 2 a_2$, $C(h)$ is given by (\ref{eq:C(h)}), and 
\begin{equation*}
D(h) = 2 \e^{ -h F_1} \,  C_3\, \e^{h F_1}  = \sum_{n=1}^{\infty} h^{n-1} (n+1)n C_{n+2},
\end{equation*}
so that
\begin{equation*}
\e^{-h F_1} \Psi(h) =   \e^{h b_1 (C_1 + h c_3 C_2 + h^2 c_3^2 C_3 )} \e^{2 h^3 d_2 C_3}  \e^{h b_1 (C_1 + h c_1 C_2 + h^2 c_1^2 C_3 )}  + \mathcal{O}(h^4).
\end{equation*}
Expanding the right-hand side in powers of $h$ and comparing the coefficients multiplying $h C_1$, $h^2 C_2$, $h^3 C_1 C_2$, and $h^3 C_3$ respectively (corresponding to the Lyndon multi-indices $(1), (2), (1,2), (3)$) with those in the expansion
\begin{equation*}
\e^{-h F_1} \e^{h (F_1+F_2)} = I + h C_1 + \frac{h^2}{2} C_2  + \frac{h^3}{3} C_3  + \frac{h^3}{6} C_1C_2 + \frac{h^3}{3} C_2C_1 + \mathcal{O}(h^4),
 \end{equation*}
 we conclude that the time-symmetric scheme (\ref{eq:modpotex}) is at least of order 4 if
 \begin{gather*}
b_1 = \frac12, \quad  b_1 (c_1 + c_3) = \frac12, \quad \frac12\, b_1^2 (c_1 + c_3) + b_1^2 c_1 = \frac16, \\
b_1 (c_1^2 + c_3^2) + 2 d_2 = \frac13.
\end{gather*}
That system of polynomial equations has a unique solution, corresponding to the choice (\ref{eq:koseleff}), as expected. The order conditions of more general products of scaled exponentials of $h F_1$, $h F_2$, and $h^3 [F_1,[F_1,F_2]]$ can be derived similarly.

Recall that in the RKN case  $C_4 = [C_3,F_1]=0$, which implies that $\varphi^{[1]}_{h}$ and $\varphi^{[112]}_{h}$ commute. Hence, the three central terms in (\ref{eq:modpotex}) can be merged
into one, the $h$-flow of the differential operator $2 a_2\, F_1 + d_2 h^2 F_{112}$, which is of the form 
\[
   \sum_j G_j(q,h) \frac{\partial}{\partial v^j}. 
\]   
   In the case of (\ref{eq:f=fa+fbRKN}),
\begin{eqnarray}
\label{eq:modpotphi_0}
  \varphi_{a_2\, h}^{[1]} \circ \varphi_{d_2\, h}^{[112]} \circ \varphi_{a_2\, h}^{[1]}(y,v) = (y, \, v+ 2 h a_2\, g(y) + h^3 d_2\, \frac{\partial g}{\partial y}(y)g(y)).
\end{eqnarray}
For Hamiltonian systems $H(p,q)= V(q) + \frac12 p^T M^{-1} p$, 
$F_{112}$ is the operator associated to the Hamiltonian function 
$ (\nabla V)^T M^{-1} \nabla V$,  depending only on $q$. Thus, (\ref{eq:modpotphi_0}) is the $h$-flow of the Hamiltonian 
\begin{equation}
  \label{eq:modpot}
  2 a_2 \, V(q) + d_2\, h^2  (\nabla V(q))^T M^{-1} \nabla V(q),
\end{equation}
which reduces to the potential $V(q)$ when $a_2=1/2$ and $d_2=0$. This explains the term ``splitting methods with modified potentials'' 
frequently used in the literature~\cite{lopezmarcos97esi,rowlands91ana,wisdom96sco}. One such method has been illustrated in practice in subsection
\ref{subsec.1.3b} (Figure \ref{fig_Schrod1}).

In addition to the reduction in the number of force evaluations, including flows associated with commutators has another advantage: since 
 the coefficients $a_i, b_i$ do not have to satisfy  all the order conditions at order $r\geq 3$, the results for negative 
time steps do not apply here, and in fact methods of order greater than two do exist within this class. In addition to (\ref{eq:modpotex}), 
other ``forward'' fourth-order methods (i.e., with all $a_j > 0$ and $b_j>0$)
involving second derivatives of the potential have been published \cite{omelyan02otc,omelyan03sai} and applied to systems where the presence of
negative coefficients leads to severe stability problems \cite{bader13sts}. Although it has been shown that achieving order six in general requires some
negative coefficients \cite{chin05sop}, it is indeed possible to construct a 6th-order processed method for cubic potentials with all $a_j, b_j$ 
being positive \cite{blanes23esm}.

Additional flows corresponding to commutators involving more operators can in principle be incorporated into the scheme. Thus, for instance, the operator 
$[F_{12},F_{112}] = -2[C_2,C_3]$ is also of the form 
\[
  \sum_j g^{[5]}_j(q,h) \frac{\partial}{\partial v^j},
\]  
   which allows us to compute its $h$-flow explicitly (see \cite{blanes01hor,blanes08sac} for more details). Again, this procedure allows one to introduce additional free parameters
into the scheme and construct more efficient integrators as long as the simultaneous evaluation of $g(y)$, $\frac{\partial g}{\partial y}(y)g(y)$, $g^{[5]}(y)$, etc. is not 
substantially more expensive than the evaluation of $g(y)$ itself.

\subsection{Near-integrable systems}
\label{near-integ}

Very often in applications one has to deal with differential equations such as
\begin{equation}     \label{qi.1}
      x^\prime = f_1(x) + \varepsilon f_2(x),
\end{equation}
where $|\varepsilon| \ll 1$ and the $h$-flows $\varphi^{[1]}_h$, $\varphi^{[2]}_h$ corresponding to $f_1$ and $ \varepsilon f_2$, respectively, are readily available. In 
classical Hamiltonian mechanics, in particular, is rather common to have a Hamiltonian function $H$ which is a small perturbation of an exactly integrable
Hamiltonian $H_1$, that is, $H = H_1 + \varepsilon H_2$, with $0 < \varepsilon \ll 1$ \cite{goldstein80cme,pars79ato}. A canonical example corresponds to
the gravitational $N$-body problem in Jacobi coordinates, already considered in subsection \ref{subsec.1.4}.

For this type of problem, splitting methods of the form (\ref{eq:splitting}), 
\[
\psi_h = \varphi^{[1]}_{a_{s+1} h}\circ \varphi^{[2]}_{b_{s}h}\circ \varphi^{[1]}_{a_{s} h}\circ \cdots\circ
 \varphi^{[1]}_{a_{2}h}\circ
 \varphi^{[2]}_{b_{1}h} \circ \varphi^{[1]}_{a_{1}h},
\] 
are especially well adapted. On the one hand, the error is at most $\mathcal{O}(\varepsilon)$ and vanishes with $\varepsilon$, since in that case the scheme
reproduces the exact solution. On the other hand, typically,  $|\varepsilon| \ll h$ (or at least $|\varepsilon| \approx h$), so that one
is mainly interested in eliminating error terms with small powers of $\varepsilon$, and its number grows as a polynomial in the order $r$, rather than
exponentially. Thus, there is only one error term of order $\varepsilon h^k$ (namely, the term $h^k \varepsilon w_{11\cdots 12}
F_{11\cdots12}$ in the expansion (\ref{bch2})), $\left\lfloor\frac12(k-1)\right\rfloor$ terms of order $\mathcal{O}(\varepsilon^2 h^k)$ and
$\left\lfloor\frac16(k-1)(k-2)\right\rfloor$ terms of order $\mathcal{O}(\varepsilon^3 h^k)$ \cite{mclachlan95cmi}. 

In the treatment of splitting methods for near-integrable systems it is convenient to introduce the notion of \emph{generalized
order}, following \cite{mclachlan95cmi}. Thus, we say that $\psi_h$ is of generalized order $(r_1, r_2, \ldots, r_m)$, with $r_1 \ge r_2 \ge \cdots \ge r_m$, if
\[
  \psi_h(x)  - \varphi_h(x) = \mathcal{O}(\varepsilon h^{r_1 +1} + \varepsilon^2 h^{r_2 +1} + \cdots + \varepsilon^m h^{r_m +1}) \quad 
  \mbox{ as } \quad (h, \varepsilon) \longrightarrow (0,0).
\]
With this notation, a method such that the local error is $\mathcal{O}(\varepsilon h^{2n+1} + \varepsilon^2 h^3)$ 
is said to be of generalized order $(2n,2)$. In this sense, the Strang scheme is of order $(2,2)$, whereas the $(10,6,4)$ integrator $\psi_h$ used in subsections \ref{subsec.1.3} and
\ref{subsec.1.4} satisfies
\[
  \psi_h(x)  - \varphi_h(x) = \mathcal{O}(\varepsilon h^{11} + \varepsilon^2 h^{7} + \varepsilon^3 h^{5}). 
\]
{The general analysis of the order conditions carried out in subsection \ref{subsec.2.4} in terms of Lyndon multi-indices readily allows us to characterize the generalized order of a given splitting scheme. We simply need to observe that, by replacing $F_2$ with $\varepsilon F_2$ in the expansions derived in subsection~\ref{subsec.2.4}, a series expansion of $\e^{-h F_1} (\Psi(h)-\e^{h(F_1 + \varepsilon F_2)})$ is obtained, where the term associated to a multi-index $(i_1,\ldots,i_k)$ with $k$ indices is affected by a $k$th power of $\varepsilon$.}
In particular, we get
Table \ref{cond.g.o}, which contains  explicitly the Lyndon multi-indices involved in several
consistent palindromic splitting methods of the given generalized order.

\begin{table}[h!]
\begin{center}
\begin{tabular}{|l|l|} \hline
   Generalized order  &  \hspace*{0.4cm} Lyndon multi-indices  \\ \hline
      $(2n,2)$  &  $(3), (5), \ldots, (2n-1)$  \\
      $(8,4)$  &  $(3)$, $(5)$, $(7)$, $(1,2)$ \\
      $(10,4)$  &  $(3)$, $(5)$, $(7)$, $(9)$, $(1,2)$ \\
      $(8,6,4)$  &  $(3)$, $(5)$, $(7)$, $(1,2)$, $(1,4)$, $(2,3)$ \\
      $(10,6,4)$  &  $(3)$, $(5)$, $(7)$, $(9)$, $(1,2)$, $(1,4)$, $(2,3)$ \\ \hline
\end{tabular}
\end{center}
\caption{Lyndon multi-indices corresponding to consistent palindromic splitting methods of a given generalized
order. }
\label{cond.g.o}
\end{table}

\subsection{Splitting methods for linear systems}
\label{sub.lisys}

In the numerical integration of the Schr\"odinger equation implemented in Subsection \ref{subsec.1.3b}
we separated the system into kinetic and potential energy and then applied several schemes based on this
splitting. It has long been recognized, however, that there exist other possibilities for splitting such a system. Given the
$N$-dimensional linear ODE
\begin{equation}   \label{ls1}
  i \frac{d }{dt} u(t) = H \, u(t),  \qquad
    u(0) = u_{0} \in \mathbb{C}^N
\end{equation}
resulting from the space discretization of equation (\ref{Schr1}), with $H$ real and symmetric, 
\cite{gray94chs,gray96sit} separate $u$ into its real and imaginary
parts, $q = \Re(u)$, $p = \Im(u)$. Then, in terms of $q,p$, equation (\ref{ls1}) leads to
\begin{equation}   \label{clas1}
  \frac{d }{dt} q = H \, p, \, \qquad  \qquad
  \frac{d }{dt} p = - H \, q,
\end{equation}
so that they can be seen as the classical evolution equations
corresponding to the Hamiltonian function 
\begin{equation} \label{ls1b}
  \hat{H}(q,p) = \frac{1}{2} p^T H p + \frac{1}{2} q^T H q 
\end{equation}  
in terms of canonical variables $q$ and $p$. The exact solution of (\ref{clas1}) is given by
\begin{equation} \label{clas2}
 \left(  \begin{array}{c}
  q(t) \\
 p(t)  \end{array} \right) = O(t H) 
 \left( \begin{array}{c}
  q_0 \\
 p_0  \end{array} \right), \quad \mbox{ where } \quad
    O(t H) = \left(
 \begin{array}{rcr}
   \cos(t H)   &  & \sin(t H)  \\
  -\sin(t H)  &   & \cos(t H)  \end{array}
  \right)
 \end{equation} 
is an orthogonal and symplectic $2N \times 2N$ matrix. As with the formal solution $u(t) = \e^{-i t H} u_0$ of eq. (\ref{ls1}),
$O(tH)$ may be very expensive to compute, so that suitable approximations might be necessary, such as those provided by
splitting methods applied to (\ref{ls1b}). In this respect, notice that if we introduce the nilpotent matrices $A$ and $B$
\begin{equation}    \label{eq.4a}
  A \equiv \left(  \begin{array}{ccc}
               0   & \, &  H  \\
               0    & \, &  0  \end{array} \right),
              \qquad\qquad
  B \equiv \left(  \begin{array}{ccc}
               0   & \, & 0  \\
              -H    & \, & 0  \end{array} \right),
\end{equation}
then it is clear that the symplectic Euler-VT method of Subsection \ref{subsec.1.3} is simply
\begin{equation} \label{ETVli}
    \left( \begin{array}{c}
          q_{n+1} \\
          p_{n+1} 
         \end{array} \right) =  \e^{h  B} \, \e^{ h A}
      \left( \begin{array}{c}
          q_{n} \\
          p_{n} 
         \end{array} \right),              
\end{equation}
whereas the St\"ormer--Verlet Algorithm \ref{alg-SV-VTV} corresponds to
\begin{equation} \label{eq.S.lin}
  x_{n+1} = \e^{\frac{h}{2} B} \, \e^{h A} \, \e^{\frac{h}{2} B} \, x_n, \qquad\quad x = (q,p)^T
\end{equation}
for a given time step $h$. Notice that all schemes based on this splitting are automatically symplectic. 

At this point, nothing prevents us from using any of the RKN splitting methods treated in subsection \ref{subsec3.1}, even with nested
commutators as in subsection \ref{nested-commut}. It turns out, however, that the particularly simple algebraic structure of the system (\ref{ls1b})
makes it possible to design more efficient schemes. Specifically, there is only one independent condition to increase the order
from $r=2k-1$ to $r=2k$, and only two to increase the order from
$r=2k$ to $r=2k+1$ for a given $k$ (see~\cite{blanes08otl} for more details). As a result, splitting methods of the form
\begin{equation} \label{spo1}
  x_{n+1} = \e^{h a_{s+1} A} \, \e^{h b_s B} \, \cdots \, \e^{h a_2 A} \, \e^{h b_1 B} \, \e^{h a_1 A} x_n
\end{equation}
of order $r$ for $r=2, 4, 6, 8, 10$ and 12 
can be obtained with $s=r$ exponentials $\e^{h b_j B}$  \cite{gray96sit}. By contrast, at least
15 and 31 exponentials are needed in general to attain orders 8 and 10, respectively.  

A couple of comments are worth making. First, this class of symplectic methods do not preserve the orthogonal character of the exact solution 
given by $O(t H)$ (or alternatively, the unitarity of (\ref{ls1}). Nevertheless, it has been shown in \cite{blanes08otl} that the average relative errors
due to the lack of preservation of orthogonality or unitarity do not grow with time, since the schemes are conjugate to orthogonal or unitary methods for
sufficiently small values of $h$. Second, although initially motivated by the time integration of the Schr\"odinger equation, methods (\ref{spo1}) can be
generalized in several ways. Thus, they have been used to construct an algorithm to approximate $\e^{-i t H} v$ 
for any real symmetric matrix $H$ and any complex vector $v$ by only carrying out matrix-vector products of the form $H v$. As shown in \cite{blanes15aea},
the algorithm is more efficient than schemes based on Chebyshev polynomials for all tolerances and values of $h$. These methods can also be adapted
for systems of the form
\[
  x' = M y, \qquad\quad y' = - N x, 
\]
with $x \in \mathbb{R}^{d_1}$, $y \in \mathbb{R}^{d_2}$, $M \in \mathbb{R}^{d_1 \times d_2}$ and $N \in \mathbb{R}^{d_2 \times d_1}$.

\subsection{Splitting methods for non-autonomous systems}
\label{sub.nas}

So far we have restricted our attention to splitting methods for autonomous differential equations. The question we analyze next
is whether the same techniques can be applied when there is an explicit time dependence in the equation
to integrate. The ideal situation would be that the methods designed for $x' = f(x)$ could
also be used (maybe with only minor modifications) when one has $x' = f(t,x)$. In addition, one would like
the schemes previously considered in this section for special problems to remain
valid in the non-autonomous case.

Let us first consider the general situation, corresponding to a system of the form
\begin{equation}    \label{naut.1}
   x' = f(t,x) = f_1(t,x) + \cdots + f_m(t,x),
   \qquad\quad  x(0) = x_0,
\end{equation}
that is, when the explicit time dependence is present in each part. Then we can take 
$t$ as a new coordinate and transform (\ref{naut.1})
 into an equivalent autonomous equation to which standard splitting algorithms can subsequently be
applied. More specifically, equation (\ref{naut.1}) is equivalent to the
enlarged system
\begin{equation}  \label{naut.2}
    \frac{d}{dt}  \left( \begin{array}{c}
          x   \\
        x_{d+1} 
            \end{array}  \right) =
 \underbrace{\left( \begin{array} {c}
         0  \\
            1
         \end{array}  \right)}_{\hat{f}_0} +
 \underbrace{\left( \begin{array} {c}
         f_1(x_{d+1},x)  \\
            0
         \end{array}  \right)}_{\hat{f}_1} + \cdots +
 \underbrace{\left( \begin{array} {c}
         f_m(x_{d+1},x) \\
            0
         \end{array}  \right)}_{\hat{f}_m},
\end{equation}
with $x_{d+1}\in\mathbb{R}$. If the resulting (autonomous) equations
\[
  y' = \hat{f}_i(y),  \qquad\quad
i=0,1,\ldots,m \qquad\qquad \mbox{ with } \quad y = (x,x_{d+1}) 
\]
can be solved, then we may use any splitting method of the form 
(\ref{eq:splitting}), since $x_{d+1}$ advances only with $\hat{f}_0$ and remains constant for the rest of the system.

It turns out that, for problems which are separable into just two parts, that is,
%
\begin{equation}    \label{naut.1b}
   x' = f(t,x) = f_1(t,x) + f_2(t,x),
   \qquad\quad  x(0) = x_0,
\end{equation}
one can do better: if $t$ is taken as a new coordinate \emph{twice}, and one writes
%
\begin{equation}  \label{naut.2b}
    \frac{d}{dt}  \left( \begin{array}{c}
          x   \\
        x_{d+1}  \\
        x_{d+2}
            \end{array}  \right) =
 \underbrace{\left( \begin{array} {c}
         f_1(x_{d+1},x)  \\
            0 \\
            1
         \end{array}  \right)}_{\hat{f}_1} +
 \underbrace{\left( \begin{array} {c}
         f_2(x_{d+2},x) \\
            1 \\
            0
         \end{array}  \right)}_{\hat{f}_2},
\end{equation}
with $x_{d+1},x_{d+2}\in\mathbb{R}$, then one can apply the same splitting schemes designed for autonomous systems separable into two pieces 
to
\[
  y' = \hat{f}_1(y),  \qquad\quad
  y' = \hat{f}_2(y), \qquad\qquad \mbox{ with } \quad y=(x,x_{d+1},x_{d+2}).
\]
This is so because  $x_{d+1}$ is constant
when integrating the first equation and $x_{d+2}$ is constant when
solving the second one. The procedure can be viewed as a
generalization of the one proposed in \cite{sanzserna96cni} for Hamiltonian systems 
$H(t,q,p)=T(t,p)+V(t,q)$: by introducing two new coordinates $q_{d+1},q_{d+2}$ and their associated momenta, $p_{d+1},p_{d+2}$, one instead deals with
the formally autonomous Hamiltonian
\[ 
  \tilde H(q_{d+1},q_{d+2},q,p_{d+1},p_{d+2},p)= \big(  T(p_{d+2},p)+p_{d+1}\big) + 
	\big( V(q_{d+1},q)- q_{d+2}  \big).
\]
In the special case of the non-autonomous second order differential equation
\begin{equation}  \label{rkn.1_na}
    y^{\prime\prime}= g(t,y), 
\end{equation}
it is convenient to split the system in the extended phase space as
\[
   \frac{d }{dt} \left( \begin{array}{c}
   				y \\
				v  \\
				y_{d+1}
			  \end{array} \right) = 
		 \left( \begin{array}{c}
		 		v \\
				0 \\
				1
			 \end{array} \right)+
		 \left( \begin{array}{c}
		 		0 \\
				g(y_{d+1},y) \\
				0
			 \end{array} \right),
\]
since then one has an autonomous system with the same algebraic structure as those considered in subsection \ref{subsec3.1}, so that RKN splitting methods
(even including commutators) can also be used.

For Hamiltonian systems $H(t,q,p)=T(p)+V(t,q)$, 
this is equivalent to introducing a new coordinate $q_{d+1}=t$ and its associated momentum, $p_{d+1}=-H$, and  considering the extended (autonomous) Hamiltonian function
\[ 
  \tilde H(q_{d+1},q,p_{d+1},p)= \big(  T(p)+p_{d+1}\big) + V(q_{d+1},q),
\]
which is still quadratic in momenta, so that symplectic RKN methods can be used. Notice that $\tilde H$ is only linear in $p_{d+1}$ and
modified potentials only involve derivatives of the potential with respect to $q$, but not with respect to $q_{d+1}$, i.e. they do not require time derivatives.
In this case the evolution for $p_{d+1}$ is irrelevant, so there is no need to compute it.

Finally, if one takes $t$ as two new coordinates in the non-autonomous near-integrable system
\begin{equation}    \label{nautNI.1b}
   x' = f(t,x) = f_1(t,x) + \varepsilon f_2(t,x),
   \qquad\quad  x(0) = x_0,
\end{equation}
then the special structure of a near-integrable system is destroyed. A partial remedy consists in separating the system as 
\begin{equation}  \label{nautNI.2b}
    \frac{d}{dt}  \left( \begin{array}{c}
          x   \\
        x_{d+1} 
            \end{array}  \right) =
 \underbrace{\left( \begin{array} {c}
         f_1(x_{d+1},x)  \\
            1
         \end{array}  \right)}_{\hat{f}_1} +
 \underbrace{ \varepsilon\left( \begin{array} {c}
         f_2(x_{d+1},x) \\
            0
         \end{array}  \right)}_{ \varepsilon \hat{f}_2},
\end{equation}
which requires to numerically solve the non-autonomous unperturbed system.

One should bear in mind that
taking time as an additional coordinate is of interest only if the time dependence in $f_i, \ i=1,\ldots,m$ is cheap to compute. Otherwise the resulting algorithm may be computationally costly, since these functions have to be evaluated $s$ times per time step. 
This drawback can be
avoided by approximating the exact solution  at each step
 by a composition of maps that in some sense incorporates average values of the vector fields
with different weights on the subinterval $[t_n, t_{n+1}]$  \cite{blanes06smf}. Specifically, the schemes read
\begin{equation}   \label{naut.3}
   \psi_{h} =
     \varphi_{h}^{[\hat{A}_{s+1}]} \circ
      \varphi_{h}^{[\hat{B}_s]} \circ \varphi_{h}^{[\hat{A}_s]}
       \circ \cdots \circ \varphi_{h}^{[\hat{A}_{2}]} \circ
       \varphi_{h}^{[\hat{B}_1]} \circ
       \varphi_{h}^{[\hat{A}_1]},
\end{equation}
where the maps $\varphi_{h}^{[\hat{A}_i]}$,
$\varphi_{h}^{[\hat{B}_i]}$ are the exact $1$-flows corresponding
to the time-independent differential equations
\begin{equation}   \label{naut.4}
   x^{\prime}  =   \hat{A}_i(x), \qquad\quad
   x^{\prime}  =   \hat{B}_i(x),  \qquad  i=1,2,\ldots
\end{equation}
respectively, with
\begin{equation}  \label{naut5}
      \hat{A}_i(x)   \equiv
      h  \sum_{j=1}^k \rho_{ij} f_1(\tau_j,x),  \qquad\quad
      \hat{B}_i(x)   \equiv
       h \sum_{j=1}^k \sigma_{ij} f_2(\tau_j,x).
\end{equation}
Here $\tau_j = t_n + c_j h$ and the (real) constants
$c_j$, $\rho_{ij}$, $\sigma_{ij}$ are chosen in such a way that $\psi_h$ provides an approximation of order $r$.
  These methods have the additional advantage that, 
when
 applied to (\ref{naut.1b}) with the time frozen, they reproduce
 the standard splitting (\ref{eq:splitting}), since $\sum_{j} \rho_{ij} = a_{i}$ and $\sum_{j} \sigma_{ij} = b_{i}$. The same technique can be applied
 to the splitting methods analyzed in subsection \ref{sub.lisys} when the linear system (\ref{ls1}) is explicitly time-dependent. In that case, the resulting scheme
 involves linear combinations of $H$ evaluated at some intermediate times \cite{blanes17sta}.
 

\section{Qualitative properties of splitting methods}
\label{sect4}

\subsection{Changes of variables and differential equations on smooth manifolds}
\label{subsec41}

Given a smooth autonomous differential equation in $\mathbb{R}^D$, 
\begin{equation}
\label{eq:x'=f(x)}
x' = f(x), 
\end{equation}
a smooth change of variables $x = \theta(\hat x)$ transforms (\ref{eq:x'=f(x)}) into a new autonomous differential equation in $\mathbb{R}^D$,
\[
\hat x' = \hat{f}(\hat x),
\]
such that their $t$-flows are related as follows: for all $x_0 \in \mathbb{R}^D$,
\[
\varphi_t^{[f]} \big(\theta(\hat x_0) \big) =  \theta \big(\varphi_t^{[\hat f]}(\hat x_0) \big), \qquad \mbox{ with } \qquad \hat x_0 = \theta^{-1}(x_0).
\]
This implies that a similar property holds for the map $\psi_h \approx \varphi_h^{[f]}$ defined by the splitting method (\ref{eq:splitting}) applied to 
(\ref{eq:x'=f(x)}) when $f(x) = \sum_{j=1}^2 f_j(x)$. That is, if the change of variables $x = \theta(\hat x)$ transforms each equation $x' = f_j(x)$ into 
$\hat x' = \hat f_j(\hat x)$, then the map $\hat \psi_h \approx \varphi_h^{[\hat f]}$ obtained by applying the splitting method with the 
same $a_j,b_j$ coefficients in the new variables $\hat x$ is related to $\psi_h$ by
\[
\psi_h \big(\theta(\hat x_0) \big) =  \theta \big(\hat \psi_h(\hat x_0) \big).
\]
In words, the following two procedures give exactly the same numerical results: (i) applying the splitting method to the ODE corresponding to the new variables $\hat x$ with the initial condition $\hat x(0)=\hat x_0\in \mathbb{R}^D$ and then transforming the result to the old variables $x$,  and (ii) applying the splitting method to the ODE formulated in the old variables $x$ with the initial condition $x(0)=x_0 :=\theta(\hat x_0) \in \mathbb{R}^D$.  This property does not hold in general for other integration schemes. For instance, in the case of Runge--Kutta methods it is true if the change of variables $\theta:\mathbb{R}^D \to \mathbb{R}^D$ 
is an affine map, but not in general.

The above property can be stated in a more abstract way by saying that splitting methods constitute a particular class of numerical integrators for differential equations defined on smooth
manifolds: since the $t$-flow of a smooth vector field on a smooth manifold $\M$ is coordinate-independent, then splitting methods applied to a differential
equation on $\M$ can also be naturally defined in a coordinate-independent way. 
More precisely,  $\psi_h:\M \to \M$ defined by (\ref{eq:splitting}) gives a natural approximation of the $h$-flow $\varphi_h^{[f_1+f_2]}$ on the smooth manifold $\M$, described only in terms of objects related to the manifold itself.
This is in contrast to other approaches for the numerical integration of differential equations on manifolds, which depend on particular choices of global or local charts, or particular embeddings of the manifold in a higher dimensional Euclidean space \cite{hairer06gni}.

In fact, the treatment of the order conditions for splitting methods carried out in section~\ref{sect2} is also valid in this more abstract setting. Suppose
we have two smooth vector fields $f_1$ and $f_2$ on $\M$, and let $F_i$ ($i=1,2$) be the  linear operators on $\cinfM$ defined as follows:
for each $g\in \cinfM$, $F_i g \in \cinfM$ is the $f_i$-directional derivative of $g$. Then, the $t$-flow $\varphi_t^{[f_i]}$ (or $\varphi_t^{[i]}$ for short) satisfies, 
for all $g\in\cinfM$ and $t \in \mathbb{R}$,
\[
\frac{d}{dt} g(\varphi_t^{[i]}) = F_i g(\varphi_t^{[i]}).
\]
This implies that the series in powers of $t$ of $g(\varphi_t^{[i]})$ can be represented as $\e^{t F_i} g$, i.e., we are in the same situation as in the
case of differential equations on $\mathbb{R}^D$, discussed in Subsection \ref{subsec.1.2}. Consequently, the series in powers of $h$ of 
$g(\psi_h(x))$ for the map $\psi_h:\M \to \M$ defined by (\ref{eq:splitting}) can be obtained by expanding $\Psi(h)g$, where $\Psi(h)$ is given by (\ref{eq:LieTransformationPsi}), 
so the analysis of the order conditions of section~\ref{sect2}
can be formally applied here as well.

\subsection{Stability}
\label{subsec:stability}

An important characteristic of  numerical integrators  is their
  \emph{stability}. Generally speaking, the
numerical solution provided by a stable method
does not tend to infinity when the exact solution is bounded.
 To analyze the (linear) stability of a given integrator, a model problem is typically chosen, so that both the numerical and exact solutions can be explicitly written out.
 In the case of a splitting method like (\ref{eq:splitting}), the model problem is the simple harmonic oscillator 
$y^{\prime\prime} + \omega^2 y = 0$, $\omega >0$ \cite{lopezmarcos96aes,mclachlan97osp}, with the standard ($x= (q,p)=(\omega y, y^{\prime})$) splitting
\begin{equation}\label{harmonic2}
    \left(\begin{array}{c}
        q^\prime \\
        p^\prime  \end{array}\right) =
    \bigg[   \underbrace{\left(\begin{array}{cc}
        0 & \omega \\
        0 & 0  \end{array}\right)}_{A} +
         \underbrace{\left(\begin{array}{cc}
        0 & 0 \\
       -\omega & 0  \end{array}\right)}_{B}
 \bigg]
     \left(\begin{array}{c}
        q \\
        p  \end{array}\right),
\end{equation}
and exact solution at time $t = h$ 
\begin{equation} \label{exactsolh}
     \left(\begin{array}{c}
        q(h) \\
        p(h)  \end{array}\right) = M_z 
     \left(\begin{array}{c}
        q(0) \\
        p(0)  \end{array}\right), \qquad 
      M_z =   \left(\begin{array}{rr}
        \cos z & \sin z \\
        - \sin z & \cos z  \end{array}\right), \qquad z = h \omega.
\end{equation}
There are at least two reasons for this choice of model. First, if a numerical method already provides unbounded numerical solutions for
a given $h$ on this system, one cannot expect good behavior for more general problems. Second, there are physically relevant problems 
that, once formulated in appropriate coordinates, are expressed as a system of uncoupled 
harmonic oscillators. Thus, a precise characterization of
the stability of a splitting method on (\ref{harmonic2})
can be useful to build accurate and stable algorithms for their numerical treatment.  The linear
system (\ref{clas1}) considered in subsection \ref{sub.lisys} belongs to this category. This is also the case of the more general equation 
\begin{equation}  \label{harmonic4}
  q^{\prime} = M^{-1} p, \qquad p^{\prime} = - N q,
\end{equation}
where $M$ and $N$ are $d \times d$ symmetric positive definite matrices.
Writing $M = L L^T$ and introducing new variables $\tilde{q} = L^T q$, then $\tilde{q}^{\prime\prime} = - L^{-1} N L^{-T} \tilde{q}$. Since $N$ is
symmetric positive definite, then $L^{-1} N L^{-T}$ is diagonalizable with positive eigenvalues. A new change of variables reduces the
system to a set of $d$ uncoupled scalar harmonic oscillators $y_i^{\prime\prime} = -\omega_i^2 y_i$, with $\omega_i^2$  the eigenvalues
of $L^{-1} N L^{-T}$ \cite{bou-rabee18gia}.

Application of the splitting method (\ref{spo1}) to (\ref{harmonic2}) results in the map
\begin{equation} \label{numerh}
     \left(\begin{array}{c}
        q_{n+1} \\
        p_{n+1}  \end{array}\right) = \tilde{M}_z 
     \left(\begin{array}{c}
        q_n \\
        p_n  \end{array}\right), \qquad 
      \tilde{M}_z =   \left(\begin{array}{rr}
        K_1(z) & K_2(z) \\
        K_3(z) & K_4(z)  \end{array}\right), 
\end{equation}
where $K_1(z)$ and $K_4(z)$  (respectively, $K_2(z)$, $K_3(z)$)  are even (respectively, odd) polynomials in $z$, 
\[
  \det \tilde{M}_z = 1, \qquad K_1(0) = K_4(0) = 1,
\]   
and, if the scheme is time-symmetric, then $K_1(z) = K_4(z)$. An essential role in the analysis is played by the stability polynomial, defined as 
\[
  p(z) = \frac{1}{2} \tr \, \tilde{M}_z = \frac{1}{2} (K_1(z) + K_4(z)).
\]  
The eigenvalues of $\tilde{M}_z$ are the zeros of 
$\lambda^2 - 2 p(z) \lambda + 1$ and they determine the stability of the given method:
if $z$ is such that $|p(z)| < 1$, then $\tilde{M}_z$ has complex conjugate eigenvalues of modulus 1 and the powers $\tilde{M}^n_z$, $n \ge 0$,
remain bounded, whereas if  $|p(z)| > 1$, then $\tilde{M}^n_z$ grows exponentially with $n$. Linear instability occurs when $+1$ or $-1$
is an eigenvalue with multiplicity 2 and $\tilde{M}_z$ is not diagonalizable.

Notice that the stability polynomial of a consistent splitting method is an even polynomial satisfying
\[
   p(z)=1- \frac{z^2}{2} + \mathcal{O}(z^4) \quad \mbox{ as } \quad z \rightarrow 0,
\]
so that for sufficiently small $z = h \omega > 0$  the scheme will be stable. The stability interval is defined as the longest interval 
$(0, z_*)$ such that $\tilde{M}_z^n$ is bounded for all the iterations $n$.

In the particular case of the St\"ormer--Verlet method (\ref{eq.S.lin}) one has
\[
  K_1(z) = K_4(z) = p(z) = 1- \frac{z^2}{2}, \qquad K_2(z) = z, \qquad K_3(z) = -z + \frac{z^3}{4},
\]
and thus $z_* = 2$. This is also true for Algorithm \ref{alg-SV-TVT}, or 
$x_{n+1} = \e^{\frac{h}{2} A} \, \e^{h B} \, \e^{\frac{h}{2} A} x_n$, since both have
the same stability polynomial $p(z)$.

It is important to stress that for systems of the form (\ref{harmonic4}) that can be reduced to a collection of $d$ scalar harmonic oscillators with
frequencies $\omega_i$, $i=1,\ldots,d$, the stability interval of St\"ormer--Verlet is restricted to $h < 2/\omega_{\max}$, where
$\omega_{\max}$ is the largest frequency of the system.

Suppose now that, given an integer $k$,  we concatenate $k$ steps of length $h/k$ of the St\"ormer--Verlet method to our model problem
(\ref{harmonic2}). The resulting scheme is stable for $0 < z/k < 2$, or alternatively, its stability interval is $(0, 2k)$. It is remarkable that this is in fact
the longest stability interval one can achieve by considering any splitting method (\ref{eq:splitting}) with $s=k$ stages. 
An elementary proof of this statement is presented in \cite{bou-rabee18gia}. In consequence, the St\"ormer--Verlet method may be applied with
longer scaled time steps $z/k$ than any other splitting method with $k$ stages. This makes it the method of choice in applications 
such as molecular dynamics, where
high accuracy is not required and one is interested in using time steps as large as possible \cite{leimkuhler96imf}.

The problem of designing splitting methods of order, say, $2r$, with extended stability intervals can be addressed by first determining the
coefficients $c_j$ in
\begin{equation} \label{exprepoli}
  p(z) = \sum_{j=0}^{r} (-1)^j \frac{z^{2j}}{(2j)!} + \sum_{j=r+1}^s c_j z^{2j},
\end{equation}  
so that $p(z)$ has the largest possible value of $z_*$. Thus, a 4th-order integrator with maximal stability interval is presented in \cite{lopezmarcos96aes}, 
whereas in \cite{mclachlan97osp} the analysis is
generalized to any order and number of stages. On the other hand, in  \cite{blanes08otl} a different strategy is proposed to determine
the coefficients $c_j$ in (\ref{exprepoli}), based on interpolatory conditions and minimization of the difference  
$(p(z) - \cos z)/z^{2r+2}$ in the stability interval. This results in  
 high order methods with a large number of stages whose stability and accuracy do not deteriorate
 for larger values of $z$. It also allows one to construct very efficient 2nd-order methods for linear systems
that outperform high-order methods for a wide range of values of the time step.

Problems of the form
\begin{equation}  \label{harrmo_per}
  q^{\prime} = M^{-1} p, \qquad p^{\prime} = - N q + f(q)
\end{equation}
derived from the Hamiltonian function
\begin{equation} \label{Hpertur}
  H(q,p) = \frac{1}{2} p^T M^{-1} p + \frac{1}{2} q^T N q + U(q),
\end{equation}
with $f(q) = - \nabla_q U(q)$, also appear frequently in applications. The simple pendulum considered in Subsection \ref{subsec.1.3} belongs to this class. Instead of the usual splitting into kinetic and potential energy,
it may be advantageous to split $H$ as $H(q,p) = H_1(q,p) + H_2(q)$, with
\begin{equation} \label{ham_pert}
  H_1(q,p) = \frac{1}{2} p^T M^{-1} p + \frac{1}{2} q^T N q, \qquad H_2(q) = U(q), 
\end{equation}
or alternatively to separate (\ref{harrmo_per}) as
\[
\begin{array}{ccl}
    q^{\prime} & = & M^{-1} p, \\
    p^{\prime} & = & - N q 
 \end{array}   \qquad\quad \mbox{ and } \qquad\quad   
 \begin{array}{ccc}
 q^{\prime} &= &  0, \\
 p^{\prime} & =  &f(q)
 \end{array}
\]
and consider Strang integrators
\begin{equation} \label{strang-rot}
  S_h^{[RKR]} = \varphi_{h/2}^{[R]} \circ \varphi_h^{[K]} \circ \varphi_{h/2}^{[R]}, \qquad\quad
   S_h^{[KRK]} = \varphi_{h/2}^{[K]} \circ \varphi_h^{[R]} \circ \varphi_{h/2}^{[K]},
\end{equation}
based on the maps
\begin{equation} \label{eqA}
\begin{aligned}
  & \varphi_t^{[R]}: \; \left( \begin{matrix} q_0 \\ p_0 \end{matrix} \right) \longmapsto 
  \e^{t A} \left( \begin{matrix} q_0 \\ p_0 \end{matrix} \right), \qquad \mbox{ with } \qquad A = \left( \begin{matrix}  0 & M^{-1} \\ -N & 0 \end{matrix} \right) \\
  & \varphi_t^{[K]}: \; \left( \begin{matrix} q_0 \\ p_0 \end{matrix} \right) \longmapsto 
  \left( \begin{matrix} q_0 \\ p_0 + t f(q_0) \end{matrix} \right). 
\end{aligned}  
\end{equation}
Integrators (\ref{strang-rot}) and other splitting methods based on sequences of \emph{rotations} $ \varphi_t^{[R]}$ and \emph{kicks}
$\varphi_t^{[K]}$ are specially suitable when $f(q)$ is a small perturbation of $-Nq$, since they provide the exact solution when the perturbation
vanishes. This happens in particular in the Hamiltonian Monte Carlo method, when one deals with target densities that are perturbations
of a Gaussian density (see subsection \ref{subsecHMC}). Section \ref{sect5} is devoted to the analysis of splitting methods for this type of system.

In view of the applications, it is relevant to analyze the stability of compositions of $\varphi_t^{[R]}$ and $\varphi_t^{[K]}$. 
As in the previous case, a sequence of linear
transformations render (\ref{harrmo_per}) into a more simplified form which is used as a model problem. Specifically, the one-dimensional oscillator
\[
  q^{\prime\prime} = -q - \varepsilon q, \qquad \varepsilon > -1
\]
is the appropriate model here \cite{bou-rabee17csf}. It turns out that the Strang integrators  (\ref{strang-rot}) are also optimal
concerning stability, in the following sense \cite{casas23ano}. Let $h_k$ be the smallest positive root of the equation
\[
  \frac{k h}{2} \sin \left( \frac{h}{k} \right) = \cos  \left( \frac{\pi}{k} \right) - \cos \left( \frac{h}{k} \right). 
\]
Then, for each fixed $h^* < h_k$, $h^* \ne \pi, 2 \pi, \ldots, (k-1)\pi$, 
the intersection
of the stability region in the $(\varepsilon, h)$ domain with the line $h = h^*$ is strictly larger for a sequence of $k$ integrators (\ref{strang-rot}) than for any other
splitting method with  $k$ stages based on rotations and kicks.  One could say, therefore, that for each given step size $h$, schemes  (\ref{strang-rot}) remain stable for larger perturbations than any other
splitting method based on rotations and kicks. The value of $\varepsilon$ where instabilities arise may be very small indeed, as shown in
\cite{casas23ano} for the particular case of $k=3$ stages.

\subsection{Modified equations}

The concept of \emph{Backward Error Analysis} (BEA), arising in several branches of numerical analysis, has also shown its effectiveness
for explaining the good behavior of splitting and composition methods in long-time integrations \cite{sanz-serna92sif,hairer94bao,reich99bea}. Generally speaking, 
 given a problem $\mathcal{P}$ with true solution $\mathcal{S}$, when
a suitable numerical solver is applied, one ends up with an approximate solution $\tilde{\mathcal{S}}$. Backward error analysis thus consists in
showing that $\tilde{\mathcal{S}}$ is indeed
the \emph{exact} solution of a problem $\tilde{\mathcal{P}}$ which is in some sense close to $\mathcal{P}$. This is in contrast to \emph{forward
error analysis}, where the aim consists in estimating an appropriate distance between $\tilde{\mathcal{S}}$ and $\mathcal{S}$.

In the domain of numerical analysis of differential equations, what lies at the heart of BEA is the idea of a \emph{modified differential equation}: given
the initial value problem  $x' = f(x), \; x(0) = x_0$ and a consistent numerical integrator $\psi_h$ producing the sequence of approximations
   $x_n$ at $t_n= n h$, $n=0,1, \ldots$, one looks for another differential equation
\begin{equation}   \label{bea.0}
 \tilde{x}^{\prime} = f_h(\tilde{x})
\end{equation}
whose vector field is defined as a formal series in powers of $h$,
\begin{equation}  \label{bea.1}
   f_h(\tilde{x}) \equiv f(\tilde{x}) + h f^{[2]}(\tilde{x}) + h^2 f^{[3]}(\tilde{x}) + \cdots
\end{equation}   
and such that $x_n = \tilde{x}(t_n)$ \cite{griffiths86ots,hairer06gni}. In this way, by analyzing the difference of the
vector fields $f(x)$ and $f_h(x)$, it is possible to extract useful information  about the qualitative behavior
of the numerical solution and the global error $e_n = x_n - x(t_n) = \tilde{x}(n h) - x(t_n)$.

In the case of a splitting method, obtaining $f_h$ is quite straightforward if one uses the BCH formula to get the formal operator associated with the whole method,
as was done in Subsection \ref{subsection:BCH} for the order conditions. Thus, for the operator $\Psi(h) = \exp(h F(h))$ associated with a scheme
(\ref{eq:splitting}) of order $r$, one has
\begin{equation}   \label{mvf2}
  F(h) = h(F_1+F_2) + \sum_{i=1}^{\infty} h^{r+i} \sum_{j=1}^{c_{r+i}} w_{i,j} E_{r+i,j},
\end{equation}
where $E_{r+i,j}$ denotes the element $j$ of the Lyndon basis of the subspace $\mathcal{L}_{r+i}(F_1,F_2)$ and $w_{i,j}$ 
are fixed real numbers determined
by the actual coefficients of the method. Since $F(h)$ lies in the same Lie algebra as $F_1$ and $F_2$, the numerical solution 
inherits the properties of the exact flow associated with this
feature (e.g., Hamiltonian, volume-preserving, etc.). Now, from (\ref{mvf2}) one can easily determine the expression of each $f^{[j]}$ in
(\ref{bea.1}).

A similar procedure can be applied to composition methods by determining the formal series $h^{-1} \log \Psi(h)$, with $\Psi(h)$ the operator
associated to the scheme.

In the case of a Hamiltonian system of the form $H(q,p) = T(p) + V(q)$, the operators $F_1$ and $F_2$ are the Lie derivatives associated
with the kinetic and potential energy, respectively, so that (\ref{mvf2}) is itself an operator associated with a modified Hamiltonian $\tilde{H}$. This is
yet another reflection of the fact that splitting methods applied to a Hamiltonian system produce maps that are symplectic.

For linear problems, the series defining the modified Hamiltonian $\tilde{H}$ associated with the numerical solution
is no longer formal, and $\tilde{H}$ can be explicitly determined in closed form. Thus, the matrix
$\tilde{M}_z$ in the map (\ref{numerh})  obtained with a time-symmetric method when $|K_1(z)| < 1$, 
can be expressed as
\[
  \tilde{M}_z = \left( \begin{array}{cc}
  			\cos \theta_z & \gamma_z \sin \theta_z \\
			- \gamma_z^{-1} \sin \theta_z & \cos \theta_z
		  \end{array} \right),
\]
where $\theta_z$ and $\gamma_z$ are real functions such that $p(z) = \cos \theta_z$ and $K_2(z) = - \gamma_z^2 K_3(z)$ (see eq. 	
(\ref{numerh}) and the subsequent discussion), and  $\theta_{-z} = - \theta_z$,  $\gamma_{-z} = - \gamma_z$. 
It is then straightforward to verify that the map $\tilde{M}_z$ is precisely the $z$-flow of the modified Hamiltonian
\cite{blanes14nif,bou-rabee18gia}
\[
  \tilde{H}(q,p) = \frac{\theta_z}{2z} \left( \gamma_z p^2 + \frac{1}{\gamma_z} q^2 \right).
\]

\subsection{Modified equations and long-term behavior}

Convergence of the series  (\ref{bea.1}) defining the modified equation, apart from the linear case, 
is the exception rather than the general rule. In consequence, an alternative strategy has to be pursued 
to get rigorous estimates concerning the long-time behavior of 
the numerical solutions. Specifically, we first gives bounds on the coefficient functions $f^{[j]}(x)$ of
the modified equation, then determine an optimal truncation index, and finally 
estimate the difference between  the numerical solution $x_n$ and the exact solution $\tilde{x}(t_n)$
of the truncated modified equation. Here we summarize only the main results, and refer the reader to \cite{hairer06gni,moan02obe} and references therein
for a more comprehensive treatment.

Suppose $f(x)$, $f_1(x)$ and $f_2(x)$ are analytic in a complex neighborhood of $x_0$ verifying $\|f(x)\| \le K$ for all $x \in B_{2\rho}(x_0)$,
where $B_{\rho}(x_0) = \{ x \in \mathbb{C}^d : \|x - x_0\| \le \rho \}$, and the same is true for the functions  $f^{[j]}(x)$ of the modified equation on $B_{\rho/2}(x_0)$. If a suitable truncation index for the formal series (\ref{bea.1}) is selected, so that one has
\begin{equation}  \label{tme1}
  \tilde{x}^{\prime} = f(\tilde{x}) + h f^{[2]}(\tilde{x}) + h^2 f^{[3]}(\tilde{x}) + \cdots + h^{N-1} f^{[N]}(\tilde{x}),
\end{equation}
with $\tilde{x}(0) = x_0$ and exact flow $\tilde{\varphi}_t^{[N]}$,
then there exist constants $h_0$ with $h \le h_0/4$ and $\gamma > 0$ such that
\begin{equation}  \label{bea.6}
  \|\psi_h(x_0) - \tilde{\varphi}_h^{[N]}(x_0) \| \le h \gamma K \e^{-h_0/h}.
\end{equation}
In other words, the difference between the numerical solution $\psi_h(x_0)$ and the exact solution $\tilde{\varphi}_h^{[N]}(x_0)$   of the 
truncated modified equation (\ref{tme1}) is exponentially small.

Based on this result it is possible to get some insight into the long-time behavior of the numerical scheme. Thus, for instance,
suppose our splitting method of order $r$ is applied to a Hamiltonian system with step size $h$. Then, the modified equation can be
derived from a (truncated) Hamiltonian  
\[
  \tilde{H}(x) = H(x) + h^r H_{r+1}(x) + \cdots + h^{N-1} H_N(x),
\]
where now $x = (q,p)$.  Denoting as before by $\tilde{\varphi}_t^{[N]}$ the flow of the truncated modified equation, it is clear that
$\tilde{H}(\tilde{\varphi}_t^{[N]}(x_0)) = \tilde{H}(x_0)$ for all $t$. 
Taking into account (\ref{bea.6}) and the bounds on the functions
appearing in the modified equation (derivatives of the $\tilde{H}$ in this case), it follows that 
\[
 \tilde{H}(x_{n+1}) - 
\tilde{H}(\tilde{\varphi}_h^{[N]}(x_n)) = \mathcal{O}(h \e^{-h_0/h})
\]
 and 
\[
   \tilde{H}(x_{n}) = \tilde{H}(x_0) + \mathcal{O}( \e^{-h_0/2h}) \quad \mbox{ for } \quad n h \le \e^{h_0/2h}.
\]
If we assume in addition that the numerical solution stays in a compact set $\mathcal{K}$, then $H_{r+1}(x) + \cdots + h^{N-r-1} H_N(x)$
is uniformly bounded on $\mathcal{K}$ independently of $h$ and $N$ \cite[p. 367]{hairer06gni} and finally
\[
   H(x_n) = H(x_0) + \mathcal{O}(h^r).
\]
Equivalently, the error in the energy corresponding to the numerical solution is of order $r$ over exponentially long time intervals when
a splitting method is applied with constant step size in a compact region of the phase space \cite{moan04otk}.

With respect to the behavior of the error in position, as shown in \cite{calvo95alt,hairer06gni}, if the Hamiltonian system
is integrable and certain conditions on the frequencies at the initial point are satisfied, then
\[
   \|(q_n,p_n) - (q(t),p(t))\| \le C t h^r, \quad \mbox{ for } \quad t = n h \le h^{-r}, \qquad C = \mathrm{const.},
\]
i.e., the global error grows at most linearly in time, whereas first integrals $I(q,p)$ \index{first integral!long time preservation} 
that only depend on the action variables are well preserved
on exponentially long-time intervals, 
\[
   \|I(q_n,p_n) - I(q_0,p_0)\| \le C h^r \qquad \mbox{ for } \; t = n h \le h^{-r}.
\]

In contrast, for a non-symplectic method (non-conjugate to a symplectic one) of order $r$ one has
\[
  H(x_n) - H(x_0) =  \mathcal{O}(n h^{r+1}) =  \mathcal{O}(t h^r),
\]
i.e., the error of the energy grows linearly, whereas
the global error in the solution typically increases quadratically with time. We have seen illustrations of this feature in Section \ref{sect1}.

It is important to remark that the modified differential equation of a numerical scheme depends
explicitly on the step size used, so that if $h$ is changed, then one has a different modified equation. 
This fact helps to explain the poor long time behavior observed in practice when
a symplectic scheme is implemented directly with a standard variable step size strategy (see e.g. \cite{calvo93tdo}).

\subsection{Processing and long-term precision}
\label{sub_processing1}

The concept of \emph{conjugacy} plays a fundamental role in the study of the long term behavior of both discrete and continuous dynamical systems. In the context of splitting methods (or more generally, numerical integration methods) for systems of ODEs one replaces the $h$-flow $\varphi_h$ of the original system with a map $\psi_h$ depending on the small parameter $h$ (the step size) such that $\psi_h \approx \varphi_h$ for $h$ small enough.  The precision of the numerical approximations $x_n = \psi_h^n(x_0) \approx \varphi_h^n(x_0) = x(n h)$ can be analyzed in one of two ways.
\begin{itemize}
\item[(i)] We use standard techniques of numerical integration of ODEs to estimate the local error $\|\psi_h(x)-\varphi_h(x)\|$ and then study how this local error is propagated~\cite{hairer93sod}, 
\item[(ii)] As described in the previous subsection, we consider a truncated modified equation (\ref{tme1}) of the numerical integration map $\psi_h$, and then study the effect of replacing the original system with the modified one, in addition to the propagation of the modified local error (\ref{bea.6}) (the local error between the map $\psi_h$ and the $h$-flow $\varphi_t^{[N]}$ of the truncated modified equation).
\end{itemize}
In both cases, a better understanding of the long-term behavior of the numerical error (of the application of a given integration scheme with constant step size) can be obtained by combining such techniques with the idea of \emph{processing} a numerical integrator~\cite{lopezmarcos96ceo}. The main idea consists in
 analyzing how close $\psi_h$ and 
$\varphi_h$  are to being conjugate to each other, and using that to estimate the long-term evolution of the errors. This is closely related to the concept of effective order~\cite{butcher69teo,butcher96tno}, and the idea of enhancing numerical integrators with correctors~\cite{wisdom96sco}. 
Essentially, the procedure is as follows. Given $\psi_h$, 
\begin{itemize}
\item we find a near-identity conjugacy map $\pi_h:\mathbb{R}^D \to \mathbb{R}^D$ (i.e., $\pi_0$ being the identity map) such that $\hat \psi_h := \pi_h^{-1} \circ \psi_h \circ \pi_h$ is as close as possible to
the $h$-flow $\varphi_h$ of the original system of ODEs;
\item we estimate the propagated error $\|\psi_h^n(x_0) - \varphi_h^n(x_0)\|$ as
\begin{equation*}
\|\psi_h^n(x_0) - \varphi_h^n(x_0)\| \leq  \|\psi_h^n(x_0) - \hat \psi_h^n(x_0)\| +  \|\hat \psi_h^n(x_0) - \varphi_h^n(x_0)\|;
\end{equation*}
\item we analyze the propagated error  $\|\hat \psi_h^n(x_0) - \varphi_h^n(x_0)\|$ of the {\em processed integrator} $\hat \psi_h$, either by standard techniques or with
the modified equation of $\hat \psi_h$,
\item and we estimate the difference between the original numerical solution $x_n = \psi_h^n(x_0)$ and the numerical approximation $\hat x_n = \hat \psi_h^n(x_0)$ that would be obtained if the processed integrator were used instead of the original integration map $\psi_h$. More precisely, using the notation  
\[
x_n = \psi_h^n(x_0), \quad
\bar x_n = \psi_h^n (\pi_h(x_0) ), \quad
\hat x_n = \hat \psi_h^n(x_0) = \pi_h^{-1}(\bar x_n ), \quad
x(n h) = \varphi_h^n(x_0),
\]
for $n\geq 0$, we have
\begin{equation*}
  \|\hat x_n - x_n\|  \leq  \|\pi_h^{-1}(\bar x_n) - \bar x_n\| + \|\bar x_n - x_n\|.
\end{equation*}

\end{itemize}
Summing up, the propagated error of the integrator $\psi_h$ can be estimated as
\begin{equation}
\label{eq:propagated_error}
\|x_n - x(n h)\| \leq  \|\pi_h^{-1}(\bar x_n) - \bar x_n\| + \|\bar x_n - x_n\| +  \|\hat  x_n - x(n h)\|.
\end{equation}
If the original method is of order $r$, that is, $\psi_h(x)-\varphi_h(x) = \mathcal{O}(h^{r+1})$ as $h\to 0$, then it makes sense to choose a conjugacy map satisfying  $\pi_h(x) = x + \mathcal{O}(h^{r})$.  Hence, provided that the sequence $\{\bar x_n\}$ stays in a compact set,  one can see that 
$\|\hat x_n - \bar x_n\| = \|\pi_h^{-1}(\bar x_n) - \bar x_n\| = \mathcal{O}(h^r)$ with a constant independent of $n$. Therefore, for sufficiently large
 time intervals, the right-hand side of (\ref{eq:propagated_error}) will be dominated either by $\|\bar x_n - x_n\|=\|\psi_h^n(\bar x_0) - \psi_h^n(x_0)\|$ (the propagation along successive iterations of $\psi_h$ of a perturbation of size $\mathcal{O}(h^r)$ in $x_0$) or by $\|\hat x_n - x(n h)\|= \|\hat \psi_h^n(x_0) - \varphi_h^n(x_0)\|$ (the sum of the propagated local errors of the processed method $\hat \psi_h$). Typically, the latter dominates  over the former if the integration interval is sufficiently large. In that case, the precision of the numerical scheme $\psi_h$ for sufficiently long-term integrations will depend on the size of the local errors $\|\pi_h^{-1} \circ \psi_h(x)\circ \pi_h-\varphi_h(x)\|$ of the processed method for an appropriately chosen conjugacy map (or processor map) $\pi_h$, rather than on the local errors $\|\psi_h(x)-\varphi_h(x)\|$ of the method itself.  
 
This is illustrated by the evolution of the error in phase space of LT and $S_2$ displayed in Figures \ref{fig_PenduloEG} and \ref{fig_NBody1EG} for the pendulum problem and the $6$-body problem, respectively. Recall that the Lie--Trotter method is conjugate to the Strang splitting (see subsection \ref{sub1.2b}).  In this case, $x_n=(q_n,p_n)^T$ is the numerical solution provided by LT, whereas  $\hat x_n=(\hat q_n, \hat p_n)^T$ corresponds to $S_2$, with the processor map $\pi_h = \varphi_{h/2}^{[T]}$. In both examples,  
$\|\pi_h^{-1}(\bar x_n) - \bar x_n\| = \mathcal{O}(h)$, which is bounded for all $n$ provided that $\|\bar p_n\|$ remains bounded. 

For the pendulum problem,  the error $\|\bar x_n - x_n\|=\|\psi_h^n(\bar x_0) - \psi_h^n(x_0)\|$ due to the propagation of the initial difference $\|(\bar q_0-q_0,\bar p_0 - p_0)\|=\mathcal{O}(h)$ does not exhibit any significant increment, because (for the considered initial value) the pendulum 
 behaves as a perturbed harmonic oscillator. On the other hand, the global error $\|\hat x_n - x(n h)\|$ of $S_2$ behaves as $\mathcal{O}(t h^2)$. Therefore, the global error $\|x_n - x(n h)\|$ of LT is dominated at the beginning of the integration interval by  $\|\hat x_n - x_n\| \approx \|\bar x_n - x_n\|$ (with no clear growth over time), until it is overcome by the linearly increasing global error of $S_2$, resulting in errors of similar size at the end of the interval for LT and $S_2$. 
 
 For the $6$-body problem,  the propagation error $\|\bar x_n - x_n\|=\|\psi_h^n(\bar x_0) - \psi_h^n(x_0)\|$ increases linearly, because now $H_1$ corresponds to a collection of Keplerian problems, where perturbations in initial states are propagated linearly. However, the slope of that linear increase is smaller than that of the propagation of the global error $\|\hat x_n - x(n h)\|$ of $S_2$. Therefore, the global error for LT is dominated by $\|\hat x_n - x_n\|$ during most of the integration interval in Figure \ref{fig_NBody1EG}, and only at the end reaches the global error of $S_2$. For longer times (not shown there),
the global errors of LT and $S_2$ will be of similar size.

So far we have focused on studying the long-term performance of a given splitting method $\psi_h$ with the help of a conjugacy map $\pi_h$. In practice, we may actually enhance the performance of a given splitting method~\cite{rowlands91ana,wisdom96sco,mclachlan96mos,laskar01hos} effectively integrating the problem with the processed integration map $\hat \psi_h = \pi_h^{-1} \circ \psi_h(x)\circ \pi_h$. 
Indeed, if output is needed only every $m$ steps, the computation of
\begin{equation}
\label{eq:hatxkm}
\hat x_{n,m} =  \pi_h^{-1} \circ \psi_h^m\circ \pi_h(\hat x_{(n-1),m}), \quad n=1, 2, 3, \ldots
\end{equation}
(with $\hat x_0 = x_0$), will not require substantially more CPU time than computing 
\begin{equation*}
x_{n,m} =  \psi_h^m(x_{(n-1),m}), \quad n=1, 2, 3, \ldots,
\end{equation*}
provided that the evaluation of $\pi_h(x)$ is computationally cheap compared to $m$ evaluations of $\psi_h(x)$.  Moreover, even if frequent output is required,  one might approximate $x(n h)$ by
 \begin{equation*}
\bar x_{n} =    \psi_h^m(\bar x_{n-1}), \quad n=1, 2, 3, \ldots,
\end{equation*}
with $\bar x_0 = \pi_h(x_0)$. For sufficiently long integrations, this will cost essentially the same as applying the original integrator $\psi_h$ in a standard way, and will be nearly as accurate as $\hat x_n$ (the full application of the processed integrator), since $\|\bar x_n - \hat x_n\| = \|\pi_h^{-1}(\bar x_n) - \bar x_n\|=\mathcal{O}(h^r)$ will be negligible compared to $\|\bar x_n - x(n h)\|$ for $n$ large enough. 

This is again illustrated in Figure \ref{fig_NBody1EG} for the $6$-body problem written in Jacobi coordinates as a perturbation of Keplerian problems. Indeed, if $x_n$ is the numerical solution labelled by LT$_{\rm{pert}}$ and $\pi_h=\varphi_{h/2}^{[H_1]}$, then $\hat x_n$ corresponds to scheme $(2,2)$ while pLT$_{\rm{pert}}$ stands for $\bar x_n$. Notice that the position error of  pLT$_{\rm{pert}}$ is very similar to that of $(2,2)$, as expected from the preceding discussion.

The enhancement of splitting integrators by processing is particularly effective for problems of the form $x' = f_1(x) + \varepsilon f_2(x)$ with $|\varepsilon| \ll 1$.
Such enhancing was first considered in \cite{wisdom96sco} in the context of $N$-body problems modeling planetary systems. Several processors 
$\pi_h$ were constructed for the Strang method, leading to processed (corrected) methods $\hat \psi_h$ of generalized order $(2k,2)$  for several $k>1$
(see also \cite{mclachlan96mos,laskar01hos}), so that their local errors are $\mathcal{O}(\varepsilon h^ {2k+1} + \varepsilon^ 2 h^3)$ for a prescribed time integration. 

In subsection \ref{apss}, we construct a different processor map that effectively reduces the local error to $\mathcal{O}(\varepsilon^ 2 h^3)$ in the case 
where $f_1$ is derived from a Hamiltonian $H_1(q,p)$ of a harmonic oscillator (or a collection of harmonic oscillators whose frequencies are integer multiples of a basic frequency) and $f_2$ comes from a Hamiltonian $\varepsilon H_2(q)$ which is a polynomial in $q$. 
This means, going back to the results displayed in Figure \ref{fig_PenduloEG} for the pendulum problem, that 
 both methods LT$_{\rm{pert}}$ and $(2,2)$ are conjugate to a more accurate integrator $\hat \psi_h$ (with global error of order $\mathcal{O}(t \varepsilon^2 h^2)$) obtained from them with an appropriate processor map $\pi_h$. 
 
 Assume that, in the previous notation, $(2,2)$ (resp. LT$_{\rm{pert}}$) corresponds to the numerical solution $x_n$, and $\hat x_n$ to the processed method. Then, the error (\ref{eq:propagated_error})  is dominated by $\|\bar x_n - x_n\|$, which does not show secular growth in this case. Eventually, for long enough integration intervals, (\ref{eq:propagated_error}) will be dominated by the linearly increasing error $\|\hat x_n - x(n h)\|$ of the more accurate processed 
 integrator\footnote{This can also be checked in Figure \ref{figss1}, where the same problem is integrated with scheme $(2,2)$ and processed Strang with a larger time step and a longer time interval.}.  The situation is very similar for the evolution of errors in position displayed for the $6$-body problem in Figure \ref{fig_NBody1EG}, the only difference being that the propagation $\|\bar x_n - x_n\|$ of differences in initial values now grows linearly. The error $\|\hat x_n - x(n h)\|$ will eventually dominate in (\ref{eq:propagated_error}) because $\|\bar x_n - x_n\|$ increases linearly with a smaller slope than $\|\hat x_n - x(n h)\|$.
 
Instead of enhancing a previously existing $r$th order integrator by processing, one may also design a processed splitting method from scratch~\cite{lopezmarcos97esi,blanes99siw}: determine the $h$-parametric maps $\psi_h$ and $\pi_h$ as compositions of basic flows $\varphi_{a_j h}^{[1]}$ and $\varphi_{b_j h}^{[2]}$ with different sequences of coefficients $a_j$ and $b_j$, such that the processed map $\hat \psi_h = \pi_h^{-1} \circ \psi_h \circ \pi_h$ is a good approximation of $\varphi_h$ for sufficiently small step sizes $h$. Typically, one requires that $\pi_h^{-1} \circ \psi_h \circ \pi_h(x) - \varphi_h(x) = \mathcal{O}(h^{r+1})$, so that the processed integrator is of order $r$. In that case, if one intends to compute the approximations (\ref{eq:hatxkm}) of $x(k m)$ (for $k=1,2,3\ldots$), there is no need for the map $\psi_h$ (referred to in this context as the kernel) to be an $r$-th order approximation of $\varphi_h$.  

In any of the situations considered above (either analyzing the performance of a given splitting integrator with the help of a conjugacy map $\pi_h$, or enhancing an existing $r$th order splitting integrator by processing, or designing an splitting processing integrator from scratch), one needs to study the {\em effective order conditions} of $\psi_h$. These are the conditions on the parameters $a_i,b_i$  that guarantee that there exists a processor map $\pi_h$ such that  $\pi_h^{-1} \circ \psi_h \circ \pi_h(x) - \varphi_h(x) = \mathcal{O}(h^{r+1})$.  A general treatment of the effective order conditions of several classes of numerical integrators including splitting methods and composition methods can be found in \cite{blanes04otn} and \cite{blanes06cmf}. That treatment is based on the series expansion (\ref{bch2}) of the formal logarithm of the Lie transformation $\Psi(h)$, and it is shown that the conditions for effective order $r$ can be written in terms of the coefficients 
$w_1$, $w_2$, $w_{12}$, $w_{122}, \ldots$ featuring in (\ref{bch2}). In addition, the formal logarithm of the Lie transformation $P(h)$ of the map  $\pi_h$ is determined as 
\begin{equation*} 
\begin{aligned}
  &  \log(P(h)) = h (p_{1} F_1 + p_{2} F_2)  + h^2 p_{12} F_{12} + h^3 (p_{122} F_{122}  + p_{112} F_{112}) \\
  & \qquad + h^4 ( p_{1222} F_{1222} + p_{1122} F_{1122} + p_{1112} F_{1112}) + \cdots + \mathcal{O}(h^{r+1}),
\end{aligned}
\end{equation*}
with the coefficients $p_{\ell_1 \cdots \ell_k}$ given as polynomials of the coefficients in (\ref{bch2}).  Alternative ways of obtaining the effective order conditions  of splitting methods can be derived following the different approaches considered in Section \ref{sect2} for analyzing the standard order conditions. This is, however, out of the scope of the present work. In any case, the analysis shows that many of the order conditions of the processed method $\hat{\psi}_h$ can be
satisfied by $\pi_h$, so that $\psi_h$ must fulfill a greatly reduced set of restrictions, also of lower complexity. As a result, it is possible to construct processed
schemes as compositions of basic maps with a reduced computational cost.


\section{Highly oscillatory problems}
\label{sect5}

In this section we consider Hamiltonian systems of the form (\ref{Hpertur}) with the splitting (\ref{ham_pert}), that is, 
  \begin{equation}
  \label{eq:hamHOS}
    H= H_1+ H_2, \quad
     H_1(q,p) = \frac{1}{2} p^T M^{-1} p + \frac{1}{2} q^T N\,  q, \quad H_2(q,p) = U(q),
 \end{equation}
when $M$ and $N$ are real symmetric matrices, and
$U:\mathbb{R}^{d} \rightarrow \mathbb{R}$ is a polynomial function. The corresponding equations of motion (\ref{harrmo_per}) can be rewritten as
      \begin{equation}
      \label{eq:odeHOS}
     x^{\prime} =  f_1(x) + f_2(x), \qquad \mbox{ with } \qquad f_1(x)=A x, \quad    f_2(x)=\left(
    \begin{matrix}
      0\\ -\nabla U(q)
    \end{matrix}
\right)
\end{equation}
in terms of $x = (q,p)^T$ and the matrix $A$ given in (\ref{eqA}). Splitting methods are advantageous for system (\ref{eq:hamHOS}) provided 
$\e^{t A} x$ can be cheaply computed for each $x\in \mathbb{R}^{D}$, $D=2d$.

For the time being, we assume that $A$ is fully diagonalizable and the eigenvalues of $A$ are 
integer multiples of $\omega \, i$ (with $i$ the imaginary unit). (The more general case where the eigenvalues of $A$ lie on the imaginary axis will be considered in Subsection \ref{ss:moregeneral}.) This implies that $\e^{t A}$ is $2\pi/\omega$-periodic in $t$. For system (\ref{eq:hamHOS}),
the present assumption is equivalent to stating that the matrix $M^{-1} N$ is fully diagonalizable with all its eigenvalues of the form $-(\omega k)^2 $, with $k \in \mathbb{Z}$. In other words, $H_1$ in (\ref{eq:hamHOS}) is just a collection of harmonic oscillators whose frequencies are integer multiples of a basic frequency $\omega$. 

Splitting methods applied to (\ref{eq:odeHOS}) can be analyzed by considering series expansions in powers of $h$ and using standard tools, in particular
the material presented in the previous sections. Thus, Subsection \ref{near-integ} is particularly relevant if the basic frequency $\omega$ is large compared to the size of the potential $U(q)$ (or more generally, the size of the components of $f_2(x)$).  Indeed, rescaling time from $t$ to
$\tau=\omega\, t$, system (\ref{eq:odeHOS}) is transformed into
\begin{equation*}
     \frac{d}{d \tau} x =  \hat f_1(x) + \varepsilon\, f_2(x), \qquad \hat f_1(x)=\hat A x,
    \end{equation*}
 where $\varepsilon=1/\omega$   and  $\hat A = \varepsilon A$, so that all the eigenvalues of $\hat A$  are 
integer multiples of the imaginary unit $i$.
This approach nevertheless has an important limitation: it is meaningful only when $\omega h=h/\varepsilon$ is sufficiently small.

Different approaches have been adopted in the literature to overcome this limitation and obtain results that remain valid when $\omega h$ is large (provided that the step size $h$ is small enough compared to the size of the perturbing potential $U(q)$ and its partial derivatives).  Among them, we mention modulated Fourier expansions (see \cite{hairer06gni} and references therein) and extended word series \cite{murua17wsf,murua16cnf}. 
Extended word series were introduced in~\cite{murua17wsf} to analyze splitting methods for a class of problems that is equivalent to  (\ref{eq:odeHOS})  under the more general assumption that all the eigenvalues of $A$  lie on the imaginary axis.
 Such expansions were further used in~\cite{murua16cnf} to analyze normal forms and formal invariants of more general classes of problems.
 The formalism of extended word series allows us to work with asymptotic expansions valid for step sizes $h$ that are sufficiently small independently of the frequencies of $\e^{t A}$.

 In the present section we provide
an elementary derivation of second-order versions in $h$ of such expansions taking the Strang splitting as a case study. In particular we provide a theoretical justification for the results presented in Subsection \ref{subsec.1.3} 
for the simple pendulum and analyze the processing technique as a way to further improve those results, which are indeed valid for the general system
(\ref{eq:hamHOS}). This is done by first constructing an asymptotic expansion of the exact solution of (\ref{eq:odeHOS}) and then comparing with the expansion
corresponding to the numerical approximation obtained by a general splitting method. 
We also get the modified equation satisfied by the numerical scheme (exact up to terms in $h^3$), with explicit formulas for the coefficients, and the corresponding modified Hamiltonian. Furthermore, we explicitly construct a processor for the
Strang splitting so that the resulting scheme leads to a better approximation to the solution of (\ref{eq:odeHOS}) (in particular, with a better preservation of the energy $H$ for large time intervals). 
Finally, we indicate how the preceding results can be generalized to the more general case where the eigenvalues of $A$ lie on the imaginary axis 
 (so that, in general, $\e^{t A}$ is quasi-periodic in $t$). 
 
For the analysis it is convenient to reformulate the problem (\ref{eq:odeHOS}) into the new variables $y(t)$ given through $x(t) = \e^{t A} y(t)$, so that
now $y^{\prime} =  \e^{-t A} f_2(\e^{t A} y)$. 
Let us write
\begin{equation}
\label{eq:Fourier1}
 \e^{-t\, A}  f_2(\e^{t\, A} x) = \sum_{k \in  \mathbb{Z}} \e^{i k \omega t}\, g_k(x),
\end{equation}
i.e., the right-hand side of \eqref{eq:Fourier1} is the Fourier series expansion of $\e^{-t\, A}  f_2(\e^{t\, A} x)$.  The map $f_2$ being real implies that $g_k:\mathbb{R}^D \to \C^D$ is such that each component of $g_{-k}(x)$ is the complex conjugate of the corresponding component of $g_k(x)$.

Since we have assumed that $U(q)$ is a polynomial, 
each component of $f_2(x)$ is also a polynomial in the variables $x$, which guarantees that there is a finite number of non-zero terms in (\ref{eq:Fourier1}). We will use the notation
\begin{equation}
\label{eq:I}
\mathcal{I} = \{k \in  \mathbb{Z}\ : \ g_k \neq 0\}.
\end{equation}

\paragraph{Remarks:}
\begin{itemize}
\item The assumption that  $U(q)$ and each component of $f_2(x)$ are polynomials in the variables $x$ might seem too restrictive. However,  for more general assumptions (for instance, real-analyticity), one can always replace $U(q)$ and $f_2(x)$ with sufficiently accurate polynomial approximations. Furthermore, the material in the present section is formally valid for more general assumptions on $U(q)$ and $f_2(x)$ if one allows the set (\ref{eq:I}) to be infinite. In that case, (\ref{eq:Fourier1}) will be an infinite series, several of the formulae derived here will also involve infinite series, and appropriate assumptions should be made so as to guarantee convergence. 
\item The assumptions that $A$ is a real matrix and that $f_2:\mathbb{R}^D\to\mathbb{R}^D$ are not essential. One could consider a complex matrix $A$ and $f_2:\C^D\to\C^D$, with 
no changes in the formulae that follows.
\end{itemize}

By substitution of $t=0$ into \eqref{eq:Fourier1}, one gets 
$f_2(x) = \sum_{k \in \mathcal{I}} g_k(x)$.
Thus, the solution of the initial value problem defined by (\ref{eq:odeHOS}) and $x(0) = x_0 \in \mathbb{R}^{D}$ can be expressed as 
 $x(t) = \e^{t A} y(t)$ where $y(t)$ is the solution of 
 \begin{equation}  \label{eq:yode}
    \frac{d }{dt}y = \sum_{k \in \mathcal{I}} \e^{i k \omega  t} g_k(y), \qquad  y(0) = x_0.
\end{equation}
From the definition of the Fourier coefficients $g_k$ in (\ref{eq:Fourier1}), one can prove that
\begin{equation}
\label{eq:eigen2}
\e^{-t A} g_k(\e^{t A} x) = \e^{i k \omega t}\, g_k(x),  \quad \mbox{for} \quad k \in  \mathcal{I},
\end{equation}
and, by applying the operator $\frac{d}{dt}\big|_{t=0}$ on both sides, this is equivalent to \cite{murua16cnf}
\begin{equation}
\label{eq:eigen}
(f_1, g_k) = i\, k\, \omega \, g_k,  \quad \mbox{for} \quad k \in \mathcal{I}.
\end{equation}
Here,  $(\cdot,\cdot)$ represents the usual Lie--Poisson bracket already defined in Section \ref{sect1}: $(f_1, g_k) (x)= g_k'(x) f_1(x) - f_1'(x) g_k(x)$.

\subsection{Expansion of the exact solution}

We next obtain an approximate representation of the flow $\varphi_h^{[f_1 + f_2]}$ of (\ref{eq:odeHOS}) with initial condition $x_0 \in \mathbb{R}^D$
that is valid for sufficiently small values of $|h|$ independently of the basic frequency $\omega$. This is done by getting an expansion
of the solution $y(t)$ of (\ref{eq:yode}) valid for $|t|\leq h$. 
To begin with, we apply the substitution
$y(t) = x_0 + \mathcal{O}(t)$ on the right-hand side of the integral form 
\begin{equation*}
y(t) = x_0 + \sum_{k \in \mathcal{I}} \int_{0}^{t} \e^{i k \omega s} g_k(y(s))\, ds
\end{equation*}
of the initial value problem. This gives
\begin{equation}
\label{eq:y(t)expansion}
y(t) = x_0 + \sum_{k \in \mathcal{I}} \int_{0}^{t} \e^{i k \omega s} g_k(x_0)\, ds + \mathcal{O}(t^2).
\end{equation}
Furthermore, $y(h)$ admits the estimate
\begin{equation}
\label{eq:y(h)expansion0}
y(h) = x_0 + \sum_{\ell \in \mathcal{I}} \int_{0}^{h} \e^{i \ell \omega t}  \big(g_{\ell}(x_0) + g'_{\ell}(x_0) (y(t)-x_0)\big)\, dt + \mathcal{O}(h^3),
\end{equation}
where $g'_{\ell}(x_0)$ is the value at $x=x_0$ of the Jacobian matrix of $g_{\ell}(x)$.

Substitution of (\ref{eq:y(t)expansion}) into the right-hand side of (\ref{eq:y(h)expansion0}) finally leads to
\begin{align*}
y(h) &= x_0  + \sum_{k \in  \mathcal{I}}\left( \int_0^h \e^{i k \omega t} dt \right) g_{k}(x_0) \\
& +  \sum_{k, \ell \in \mathcal{I} } \left( \int_0^h \int_0^t \e^{i (\ell t+ k s) \omega} \, ds \, dt \right) \,  g'_{\ell}(x_0) g_{k}(x_0)+ \mathcal{O}(h^3).
\end{align*}
Equivalently,
\begin{equation}
\label{eq:y(h)expansion}
y(h) = x_0  + h \sum_{k \in  \mathcal{I}} \alpha_k(h) g_{k}(x_0)  + 
h^2 \sum_{k, \ell \in \mathcal{I} } \alpha_{k \ell}(h) g'_{\ell}(x_0) g_{k}(x_0)+ \mathcal{O}(h^3),
\end{equation}
 where the coefficients are defined as follows.
 \begin{itemize}
 \item For $k \in \mathcal{I}$,
 \begin{equation}
 \label{eq:alpha_k}
\alpha_{k}(h) = \int_0^1  \e^{i k \omega h \tau} \, d\tau = 
\left\{
\begin{array}{cl}
1 & \mbox{if } k=0, \\
\displaystyle
\frac{\e^{i  k  \omega  h}-1}{i k \omega  h} & \mbox{otherwise,}
\end{array}
\right. 
\end{equation}
\item For $k,\ell \in \mathcal{I}$,
\begin{equation}
\label{eq:alpha_kl}
\alpha_{k \ell}(h) =  \int_0^1  \int_0^{\tau_2} \e^{i  h \omega (k  \tau_1 +  \ell  \tau_2)} \, d\tau_1 \, d\tau_2. 
\end{equation}
Specifically, $\alpha_{0 0}(h) = 1/2$, and the following recursions hold: 
\begin{equation}
\alpha_{0 k}(h) = \frac{\e^{i k \omega h} - \alpha_k(h)}{i k \omega  h}, \quad 
\alpha_{k \ell}(h) = \frac{\alpha_{k + \ell}(h) - \alpha_{\ell}(h)}{i k \omega h}, \quad \mbox{ for } \quad k,\ell \in \mathcal{I} \backslash \{0\}.
\end{equation}

\item The constant in the $\mathcal{O}(h^3)$ term depends on upper bounds of the norm of $g_k$ and its partial derivatives, but is independent of the basic frequency $\omega$.
 \end{itemize}

From this we conclude that 
\begin{equation}
\label{eq:x(t)}
\begin{split}
\varphi_h(x_0) &= \e^{h A} \left( x_0  + h \sum_{k \in  \mathcal{I}} \alpha_k(h) g_{k}(x_0) + 
h^2 \sum_{k, \ell \in  \mathcal{I}} \alpha_{k \ell}(h) g'_{\ell}(x_0) g_k(x_0)+ \mathcal{O}(h^3)
\right),
\end{split}
\end{equation}
where we have dropped the upper index in $\varphi_h$ for clarity. 
Observe  that 
\begin{equation*}
|\alpha_{k}(h)| = |\mathrm{sinc}(k h \omega/2)| \leq  1, \qquad 
|\alpha_{k \ell}(h)| \leq   \frac{1}{2}.
\end{equation*}

\subsection{Expansion of the discrete solution given by Strang splitting}

We next proceed to construct an analogous expansion for the (discrete) solution furnished by the Strang splitting given in (\ref{strang-rot}), namely
\[
    S_h^{[RKR]} = \varphi_{h/2}^{[R]} \circ \varphi_h^{[K]} \circ \varphi_{h/2}^{[R]},
\]
based on the maps (\ref{eqA}). In fact, it is straightforward to check that the approximation $\tilde x(h)=S_h^{[RKR]}(x_0)$ satisfies $\tilde x(h) = \e^ {h A} \tilde y(h)$, where $\tilde y(t)$ is the solution of
  \begin{equation*}
   \frac{d }{dt} \tilde y = \e^{-\frac{h}{2}A} g(\e^{\frac{h}{2}A} \tilde y(t))= \sum_{k \in \mathcal{I}} \e^{i\, k \omega h/2}\, g_k(\tilde y(t)), \qquad \tilde y(0) = x_0.
  \end{equation*}
In consequence,  proceeding as in the previous subsection we get
\begin{equation*}
\tilde y(h) = x_0  + h \sum_{k \in  \mathcal{I}} \tilde \alpha_k(h) g_{k}(x_0) + 
h^2 \sum_{k, \ell \in  \mathcal{I}} \tilde \alpha_{k \ell}(h) g'_{\ell}(x_0) g_k(x_0)+ \mathcal{O}(h^3),
\end{equation*}
where for $k, \ell  \in \mathcal{I}$,
\begin{equation}
\label{eq:alphaStrang}
\tilde \alpha_{k}(h) = \e^{i k \omega h/2}, \qquad 
\tilde \alpha_{k \ell}(h) = \frac{1}{2} \e^{i (k + \ell) \omega  h/2},
\end{equation}
and finally
\begin{equation}
\label{eq:splittingexpansion}
S_h(x_0)= \e^{h A} \left( x_0  + h \sum_{k \in  \mathcal{I}} \tilde \alpha_k(h) g_{k}(x_0) + 
h^2 \sum_{k, \ell \in  \mathcal{I}} \tilde \alpha_{k \ell}(h) g'_{\ell}(x_0) g_k(x_0)+ \mathcal{O}(h^3)
\right),
\end{equation}
where we have also dropped the upper index in $S_h$. Notice that, from (\ref{eq:alphaStrang}),
\begin{equation*}
|\tilde \alpha_k(h)| \leq 1, \quad |\tilde \alpha_{k \ell}(h)| \leq \frac{1}{2}
\end{equation*}
and 
\begin{equation}
\label{eq:sympl_cond_alpha}
\tilde \alpha_{k\ell}(h) + \tilde \alpha_{\ell k}(h) = \tilde \alpha_k(h)  \tilde \alpha_{\ell}(h), \quad k,\ell \in \mathcal{I}. 
\end{equation}

\subsection{Composition formulae}

If one is interested in extending the previous analysis to more general splitting methods of the form (\ref{eq:splitting}), then a composition rule
concatenating the expansions corresponding to different basic methods is clearly needed. This can be obtained as follows. 

Suppose the map $\psi_{h}:\mathbb{R}^D\rightarrow \mathbb{R}^D$ can be expanded as 
 \begin{equation}
 \label{eq:chi}
\psi_h(x) =
x+ h \sum_{k \in \mathcal{I}} \kappa_{k}\,   g_{k}(x)+  
h^2 \sum_{k,\ell \in \mathcal{I}}
\kappa_{k \ell}\,   g_{\ell}'(x) g_{k}(x) + \mathcal{O}(h^3),
  \end{equation}
for some coefficients $\kappa_k, \kappa_{k,\ell}$. Then, from (\ref{eq:eigen2}), one has  for all $s\in \mathbb{R}$ that
  \begin{equation}
  \label{eq:swittching_rule}
 \begin{split} 
\psi_{h}(\e^{s A} x) &= \e^{s A} \left( x+ h \sum_{k \in \mathcal{I}} \e^{i\,k \omega s}\kappa_{k}\,   g_{k}(x) \right. \\
 &\left. \hspace{3em} + \, h^2 \sum_{k,\ell \in \mathcal{I}}
\e^{i\,(k+\ell) \omega s} \kappa_{k \ell}\,   g_{\ell}'(x) g_{k}(x) + \mathcal{O}(h^3)\right).
\end{split}
  \end{equation}
If in addition the map $\hat \psi_{h}:\mathbb{R}^D\rightarrow \mathbb{R}^D$ can be expanded as 
 \begin{equation*}
\hat \psi_h(x) =
x+ h \sum_{k \in \mathcal{I}} \hat  \kappa_{k}\,   g_{k}(x)+  
 h^2 \sum_{k,\ell \in \mathcal{I}}
\hat \kappa_{k \ell}\,   g_{\ell}'(x) g_{k}(x) + \mathcal{O}(h^3),
  \end{equation*}
then one has the following expression for the composition $\psi_h \circ \hat \psi_h$:
\begin{align*}
\psi_h(\hat \psi_h(x)) &=  x+ h \sum_{k \in \mathcal{I}} (\kappa_k + \hat  \kappa_{k})\,   g_{k}(x) \\
&+  
 h^2 \sum_{k,\ell \in \mathcal{I}}
(\kappa_{k \ell} + \hat \kappa_k \, \kappa_{\ell} + \hat \kappa_{k \ell})\,   g_{\ell}'(x) g_{k}(x) + \mathcal{O}(h^3).
\end{align*}
 This is also valid if the coefficients $\kappa_{k}$, $\hat \kappa_{k}$, $\kappa_{k\ell}$, and $\hat \kappa_{k\ell}$ depend on $h$, although in that case 
the constant in the $\mathcal{O}(h^3)$ term will also depend on the bounds of the coefficients.
   
By combining the previous results one finally arrives at the following composition rule:
\begin{equation}
\label{eq:comp}
\e^{s A} \psi_h(\e^{\hat s A} \hat \psi_h(x)) =  \e^{(s+\hat s) A} \Big(x+ h \sum_{k \in \mathcal{I}} \gamma_k \,   g_k(x) +  
 h^2 \sum_{k,\ell \in \mathcal{I}}
\gamma_{k \ell}\,   g_{\ell}'(x) g_{k}(x) \\
+ \mathcal{O}(h^3)\Big),
\end{equation}
where
\begin{equation}
\label{eq:gammacomp}
\begin{split}
\gamma_k &= \e^{i\, \hat s \omega k} \, \kappa_k + \hat  \kappa_{k}, \\
\gamma_{k \ell} &= \e^{i\, \hat s \omega (k+\ell)} \, \kappa_{k\ell} + \e^{i\, \hat s \omega \ell} \, \hat \kappa_k \, \kappa_{\ell} + \hat \kappa_{k\ell}.
\end{split}
\end{equation}

\subsection{Expansions for arbitrary splitting methods}

We now have all the required ingredients to extend the expansion (\ref{eq:splittingexpansion}) for Strang to a more general splitting method of the form 
(\ref{eq:splitting}) based on kicks and rotations, i.e., on the maps (\ref{eqA}). Specifically, if $\psi_h$ denotes such a splitting, then 
$\psi_h(x_0)$ can be expanded as
\begin{equation}
\label{eq:splittingexpansion2}
\e^{h\,(a_1 + \cdots + a_{s+1})  A} \left( x_0  + h \sum_{k \in  \mathcal{I}} \tilde \alpha_k(h) g_{k}(x_0) + 
h^2 \sum_{k, \ell \in  \mathcal{I}} \tilde \alpha_{k \ell}(h) g'_{\ell}(x_0) g_k(x_0)+ \mathcal{O}(h^3)
\right)
\end{equation}
 with the following coefficients:
\begin{itemize}
\item  for $k \in  \mathcal{I}$,
 \begin{equation}
 \label{eq:alphasplitting1}
\tilde \alpha_{k}(h) =  \sum_{j=1}^{s} b_j  \, \e^{i\, k c_j \omega  h }, \qquad \mbox{ with } \quad c_j = a_1+\cdots + a_j, \quad \mbox{for} \quad j=1,\ldots,s;
\end{equation}
\item for $k, \ell \in \mathcal{I}$,
\begin{equation}
 \label{eq:alphasplitting2}
\tilde \alpha_{k \ell}(h) = \sum_{1\leq j <  m \leq s} b_{j}  b_{m} \,  \e^{i (c_j k + c_m \ell) \omega h} +
\sum_{1\leq j \leq  r} \frac{1}{2} \,  b_{j}^2 \, \e^{i c_{j} (k + \ell) \omega h},
\end{equation}
\end{itemize}
so that they again satisfy the relation (\ref{eq:sympl_cond_alpha}), since it is preserved under compositions (\ref{eq:comp})--(\ref{eq:gammacomp}).

We can now estimate the local error of a consistent splitting method by taking into account (\ref{eq:splittingexpansion2})
and (\ref{eq:x(t)}) as
\begin{equation*}
\psi_h(x) - \varphi_h(x) = \e^{h A} \left( x_0  + h \sum_{k \in  \mathcal{I}} \delta_k(h) g_{k}(x_0) + 
h^2 \sum_{k, \ell \in  \mathcal{I}} \delta_{k \ell}(h) g'_{\ell}(x_0) g_k(x_0)+ \mathcal{O}(h^3)
\right),
\end{equation*}
where
\begin{equation}
\label{eq:deltak}
\begin{split}
\delta_{k}(h) &=  \sum_{j=1}^{s} b_j \, \e^{i k c_j \omega  h } - \int_0^1  \e^{i k \omega h \tau} \, d\tau, \\
\delta_{k \ell}(h) &= \sum_{1\leq j <  m \leq s} b_{j}  b_{m} \,  \e^{i (c_j k + c_m \ell) \omega h} +
\sum_{1\leq j \leq  r} \frac{1}{2} \,  b_{j}^2 \, \e^{i c_{j} (k + \ell) \omega h} \\
&-  \int_0^1  \int_0^{\tau_2} \e^{i  h \omega (k  \tau_1 +  \ell  \tau_2)} \, d\tau_1 \, d\tau_2. 
\end{split}
\end{equation}

It is worth remarking that if one expands the exponentials in (\ref{eq:deltak}) in series of powers of $h$ up to a certain degree $r$, then one recovers the order conditions obtained in Subsection~\ref{subsec.2.4} for multi-indices with one and two indices. In other words, using the terminology introduced in Subsection~\ref{near-integ} to analyze perturbed problems of the form (\ref{qi.1}), one obtains the conditions for a splitting method to be of generalized order $(r_1,r_2,1)$ (or $(r_1,r_2,2)$ in the case of time-symmetric splitting methods). 
Notice, however, that replacing these exponentials with such truncated series expansions will only give useful information about the size of the local error coefficients $\delta_{k}(h)$ and $\delta_{k\ell}(h)$ provided that the scaled step sizes $\{|k|  \omega h\ : \ k \in \mathcal{I}\}$ are sufficiently small.

\subsection{Modified ODE for the discrete flow furnished by a splitting method}
\label{ss:mode}

At this point it is useful to formulate a modified differential equation whose $h$-flow is closer in a certain sense to the map $\psi_h$ corresponding 
to a general splitting method admitting an expansion of the form (\ref{eq:splittingexpansion2}). More precisely, the idea is to construct a modified ODE
whose $h$-flow is $\mathcal{O}(h^3)$ close to $\psi_h$.

From the composition formula (\ref{eq:comp}), the $n$-th iterate of the map $\psi_h$ admits an expansion of the form
\begin{equation*}
 \psi_h^n(x_0)= e^{h n A} \left( x_0  + h \sum_{k \in  \mathcal{I}} \gamma_k(n,h) g_{k}(x_0) + 
h^2 \sum_{k, \ell \in  \mathcal{I}} \gamma_{k \ell}(n,h) g'_{\ell}(x_0) g_k(x_0)+ \mathcal{O}(h^3)
\right),
\end{equation*}
where $\gamma_k(0,h)=0$, $\gamma_{k\ell}(0,h)=0$, $\gamma_k(1,h)=\tilde \alpha_k(h)$, $\gamma_{k \ell}(1,h)=\tilde \alpha_{k\ell}(h)$.
Thus, to find a suitable modified ODE, it makes sense to assume that its corresponding $t$-flow $\tilde \varphi_{t}$ 
can be expanded for all $t \in \mathbb{R}$ as 
\begin{equation} \label{expvp1}
\tilde \varphi_t(x_0)= e^{t A} \left( x_0  + h \sum_{k \in  \mathcal{I}} \gamma_k(t/h,h) g_{k}(x_0) + 
h^2 \sum_{k, \ell \in  \mathcal{I}} \gamma_{k \ell}(t/h,h) g'_{\ell}(x_0) g_k(x_0)+ \cdots
\right).
\end{equation}
Then, the right-hand side of the corresponding ODE must be of the form
\begin{equation*}
\left. \frac{d}{dt}\tilde \varphi_{t}(x)\right|_{t=0}= A\, x + 
 \sum_{k \in  \mathcal{I}} \beta_k(h) g_{k}(x) + 
h \sum_{k, \ell \in  \mathcal{I}} \beta_{k \ell}(h) g'_{\ell}(x) g_k(x)+ \cdots,
\end{equation*}
where
\begin{equation}
\label{eq:beta0}
\beta_k(h) = \left. \frac{d}{d\tau}  \gamma_k(\tau,h) \right|_{\tau=0}, \qquad
\beta_{k \ell}(h) = \left. \frac{d}{d\tau}  \gamma_{k\ell}(\tau,h) \right|_{\tau=0}.
\end{equation}
We are then bound to study the family of perturbed ODEs
      \begin{equation}
      \label{eq:odebeta}
     \frac{d}{d t} x               
     = A\, x+  \sum_{k \in  \mathcal{I}} \beta_k(h)\, g_k(x) 
     +  h \sum_{k, \ell \in  \mathcal{I}} \beta_{k \ell}(h)\, g'_{\ell}(x) g_k(x),
    \end{equation}
with initial condition $x(0) = x_0$. The coefficients $\beta_k(h)$ and $\beta_{k \ell}(h)$ are, for fixed $h$, arbitrary complex numbers. Notice that this
equation generalizes the original ODE (\ref{eq:odeHOS}), which corresponds to the case where  $\beta_k(h)=1$  and $\beta_{k \ell}(h)=0$. 
 
 Assuming that the $t$-flow $\tilde \varphi_{t}$ of (\ref{eq:odebeta}) can be expanded as in (\ref{expvp1}), the group 
property of the flow (i.e, $\tilde \varphi_{t+s} = \tilde \varphi_t \circ \tilde \varphi_s$), together with the composition formula (\ref{eq:comp}), lead to
\begin{equation} \label{ecgamma}
\begin{aligned}
  & \gamma_k(\tau + \sigma,h) = \e^{i \omega h \tau k}\, \gamma_k(\sigma,h) + \gamma_k(\tau,h), \\
  & \gamma_{k\ell}(\tau + \sigma,h) = \e^{i \omega h \tau (k+\ell)}\, \gamma_{k\ell}(\sigma,h)
+ \e^{i \omega h \tau \ell}\, \gamma_k(\tau,h)\, \gamma_{\ell}(\sigma,h) + \gamma_{k\ell}(\tau,h).
\end{aligned}
\end{equation}
Now, applying the operator $\left. \frac{d}{d\sigma} \right|_{\sigma=0}$ to both sides of (\ref{ecgamma}) and using (\ref{eq:beta0}), we obtain the following system of ODEs for the coefficients $\gamma_k(\tau,h)$, $\gamma_{k\ell}(\tau,h)$:
\begin{align*}
\frac{d}{d\tau}  \gamma_k(\tau,h) &= \e^{i \omega h \tau k}\, \beta_k(h), \\
\frac{d}{d\tau} \gamma_{k\ell}(\tau,h) &= \e^{i \omega h \tau (k+\ell)} \, \beta_{k\ell}(h)
+ \e^{i \omega h \tau \ell}\, \beta_{\ell}(h) \, \gamma_k(\tau,h).
\end{align*}
This, together with the initial conditions $\gamma_k(0,h)=0$, $\gamma_{k\ell}(0,h)=0$, allows us to express
$\gamma_k(\tau,h)$, $\gamma_{k\ell}(\tau,h)$ in terms of the coefficients $\beta_k(h)$, $\beta_{k\ell}(h)$ of the perturbed ODE (\ref{eq:odebeta}):
\[
\begin{aligned}
 & \gamma_k(\tau,h) = \beta_k(h) \int_0^{\tau}  \e^{i \omega h \sigma k}\, d\sigma = \tau\, \beta_{k}(h)\, \alpha_{k}(\tau h),\\
 & \gamma_{k\ell}(\tau,h) = \tau\, \beta_{k\ell}(h)\, \alpha_{k+\ell}(\tau h) + \tau^2\,  \,  \beta_{\ell}(h)  \, \alpha_{k \ell}(\tau h).
\end{aligned}
\]
Specifying these equations to the case $\tau = 1$, and replacing $\gamma_k(1,h)$ and $\gamma_{k\ell}(1,h)$ 
 with $\tilde\alpha_k(h)$ and $\tilde \alpha_{k\ell}(h)$ (the coefficients corresponding to $\psi_h$), respectively, we have
\begin{equation} \label{fftc1} 
 \begin{aligned}
 & \tilde \alpha_k(h) = \beta_k(h) \, \alpha_k(h), \\
 & \tilde \alpha_{k\ell}(h) = \beta_{k\ell}(h)\, \alpha_{k+\ell}(h) 
+ \beta_k(h) \,  \beta_{\ell}(h) \,  \alpha_{k\ell}(h),
\end{aligned}
\end{equation}
where the expression of $\alpha_k(h)$, for all $k\in \Z$, is given by (\ref{eq:alpha_k}). Hence, $\alpha_0(h) = 1$ and, for any
$k\neq 0$,  $\alpha_k(h) \neq 0$ if and only if  $k h \neq 2\pi j$ for all $j \in \Z$. Therefore, if we assume that 
 $h \in \mathbb{R}$ is such that, for all $k, \ell \in  \mathcal{I}\backslash\{0\}$,
\begin{equation}
\label{eq:non-resonant-h}
\frac{k \omega h}{2\pi} \not \in \Z\backslash\{0\} \qquad \mbox{ and } \qquad \frac{(k+\ell) \omega h}{2\pi} \not \in \Z\backslash\{0\},
\end{equation}
(so that $\alpha_k(h) \ne 0$ and $\alpha_{k+\ell}(h) \ne 0$), then equations (\ref{fftc1}) can be solved in the $\beta_k(h)$ and $\beta_{k+\ell}(h)$
coefficients. Thus, for each $k,\ell \in  \mathcal{I}$, one has 
\begin{equation} \label{eq:beta}
\beta_{k}(h) =
\frac{ \tilde \alpha_{k}(h)}{\alpha_{k}(h)} ,  \qquad\quad
\beta_{k \ell}(h) =
\frac{ \tilde \alpha_{k \ell}(h) -  \alpha_{k \ell}(h) \, \beta_{k}(h) \beta_{\ell}(h) }{\alpha_{k+\ell}(h)}
\end{equation}
and the $h$-flow of the corresponding modified equation (\ref{eq:odebeta}) agrees up to terms of order $\mathcal{O}(h^3)$ with the expansion (\ref{eq:splittingexpansion2}) of the splitting method $\psi_h$. Moreover, the identity (\ref{eq:sympl_cond_alpha}) satisfied by the coefficients in the
expansions of both the splitting method $\psi_h$ and the exact solution, eqs. (\ref{eq:alpha_k})--(\ref{eq:alpha_kl}), implies 
\begin{equation}
\label{eq:sympl_cond_beta}
\beta_{k\ell}(h) + \beta_{\ell k}(h) = 0, \qquad k,\ell \in \mathcal{I},
\end{equation}
so that the modified equation (\ref{eq:odebeta}) can also be expressed as
\begin{equation} \label{eq.mod.eq}
\begin{aligned} 
     \frac{d}{d t} x               
     &= A\, x+  \sum_{k \in  \mathcal{I}} \beta_k(h)\, g_k(x) 
		     +  \frac{h}{2} \sum_{k, \ell \in  \mathcal{I}} \beta_{k \ell}(h)\, (g'_{\ell}(x) g_k(x)-g'_k(x) g_{\ell}(x)),\\
 &= A\, x+  \sum_{k \in  \mathcal{I}} \beta_k(h)\, g_k(x) 
		     +  \frac{h}{2} \sum_{k, \ell \in  \mathcal{I}} \beta_{k \ell}(h)\, (g_k,g_{\ell})(x).
\end{aligned}
\end{equation}
In fact, equation (\ref{eq.mod.eq}) is itself Hamiltonian, in the sense that there exists a Hamiltonian function $\tilde{H}(x;h)$ such that
\begin{equation} \label{eqHmod1}
    \frac{d}{d t} x  = J \, \nabla \tilde H(x; h),
\end{equation}
where $J$ is the canonical symplectic matrix (\ref{canonJ}). This can be seen as follows.
 First, $H_1$ in the Hamiltonian (\ref{eq:hamHOS}) can be written as
\begin{equation}
\label{eq:ham2}
 H_1(x) = \frac{1}{2} x^T Q x, \qquad \mbox{ with } \qquad Q = \left(
  \begin{matrix}
    N & 0 \\ 0 & M^{-1}
  \end{matrix}
\right),
\end{equation}
whereas $H_2(\e^{t A} x)$ has a Fourier expansion of the form
\begin{equation}  \label{expan_H}
H_2(\e^{t A} x) = \sum_{k \in \mathcal{I}} \e^{i k \omega t} G_k(x), \quad \mbox{ with } \quad g_k(x) = J \, \nabla G_k(x),
\end{equation}
and furthermore, for all $k \in \mathcal{I}$,
\begin{equation}
\label{eq:Gdec}
G_k(\e^{t\, A} x) = \e^{i k \omega t} G_k(x),
\end{equation}
so that (\ref{eq:eigen}) is equivalent to
\begin{equation}
\label{eq:eigenHam}
\{H_1, G_k\} = i k \omega \, G_k,  \quad \mbox{for} \quad k \in  \mathcal{I}.
\end{equation}
Here $\{A,B\}$ stands for the Poisson bracket of  $A, B\in \mathcal{C}^1(\mathbb{R}^D)$ defined as follows: for each $x \in \mathbb{R}^{D}$, 
\begin{equation*}
\{A,B\}(x) = (\nabla A(x))^T J\,  \nabla B(x).
\end{equation*}
It is then clear that the modified ODE (\ref{eq.mod.eq}) can be written like (\ref{eqHmod1}) with
    \begin{equation}
    \label{eq:mHam}
\tilde H(x; h) =  \frac12\, x^T Q x + \sum_{k \in  \mathcal{I}} \beta_k(h)\, G_k(x) 
   +  \frac{h}{2} \sum_{k, \ell \in  \mathcal{I}} \beta_{k \ell}(h)\, \{G_k,G_{\ell}\}(x).
\end{equation}
We should stress that both (\ref{eq:odebeta}) and the modified Hamiltonian $\tilde H(x; h)$
are well defined as long as the non-resonance assumptions (\ref{eq:non-resonant-h}) hold.

\subsection{Splitting methods with processing}

In the spirit of Subsection \ref{sub_processing1}, we now consider  a
 {\em processed} splitting integrator 
\begin{equation*}
\hat \psi_{h} = \pi_{h}^{-1}\circ \psi_{h} \circ \pi_{h},
\end{equation*}
where $\psi_h$ is a composition of type (\ref{eq:splitting}) based on kicks and rotations 
and $\pi_h:\mathbb{R}^D \to \mathbb{R}^D$ is a near-to-identity map with an expansion of the form
\begin{equation}
\label{eq:chiexpansion}
 \pi_h(x) =
x+ h \sum_{k \in \mathcal{I}} \kappa_{k}(h)\,   g_{k}(x)+  
h^2 \sum_{k,\ell \in \mathcal{I}}
\kappa_{k \ell}(h)\,   g_{\ell}'(x) g_{k}(x) + \mathcal{O}(h^3).
\end{equation}

In contrast to Subsection \ref{sub_processing1}, where the processor map $\pi_h$ is analyzed by its power series expansion in $h$, 
the coefficients $\kappa_{k}(h)$ and $\kappa_{k\ell}(h)$ featuring in  (\ref{eq:chiexpansion}) will depend on $h$ and the frequency $\omega$ in a non-polynomial way. This will allow us to find appropriate processors valid for step sizes $h$ that are not necessarily small compared with $\omega$.

Application of the composition rule (\ref{eq:comp})--(\ref{eq:gammacomp}) to both sides of the identity $\pi_{h} \circ \hat \psi_{h} = \psi_{h} \circ \pi_{h}$
implies that $\hat \psi_h$ can be expanded as
\begin{equation}
\label{eq:hatpsiexpansion}
\hat \psi_h(x) = \e^{h A} \left( x  + h \sum_{k \in  \mathcal{I}} \hat \alpha_k(h) g_{k}(x) + 
h^2 \sum_{k, \ell \in  \mathcal{I}} \hat \alpha_{k \ell}(h) g'_{\ell}(x) g_k(x)+ \mathcal{O}(h^3)
\right),
\end{equation}
where
\begin{equation} \label{eq:halphak}
\begin{aligned}
\hat \alpha_k(h) &= (1- \e^{i k \omega h}) \, \kappa_k(h) + \tilde \alpha_k(h), \\
\hat \alpha_{k \ell}(h) &= (1- \e^{i (k+\ell) \omega h}) \, \kappa_{k \ell}(h) 
   + \kappa_k(h)\,  \tilde \alpha_{\ell}(h) 
   - \e^{i \ell \omega h} \kappa_{\ell}(h)\,  \hat \alpha_k(h) 
   +\tilde \alpha_{k \ell}(h). 
\end{aligned}
\end{equation}
The expansion of the processed scheme $\hat{\psi}_h$ coincides with that of the exact flow (\ref{eq:x(t)}) if $\hat \alpha_k(h)=\alpha_k(h)$ for each
$k \in \mathcal{I}\backslash\{0\}$, and this is possible only when the following non-resonance condition holds:
\begin{equation}
\label{eq:non-resonant-h1}
\frac{k \omega h}{2\pi} \not \in \Z,
\end{equation}			
in which case
\begin{equation}
\label{eq:kappa_k_splitting}
  \kappa_{k}(h)=\frac{\tilde \alpha_{k}(h)-\alpha_{k}(h)}{\e^{i  k \omega h}-1}.
\end{equation}
Observe that, for $k=0$, $\hat \alpha_k(h)=\alpha_k(h)$ regardless of the chosen value of $\kappa_0(h)$, since $\tilde \alpha_0(h)=\alpha_0(h)=1$. For simplicity, it makes sense to choose the processor map in such a way that $\kappa_0(h)=0$.

In consequence, if $h$ satisfies the non-resonance condition (\ref{eq:non-resonant-h1}) for all $k \in \mathcal{I}\backslash\{0\}$, and the coefficients $\kappa_k(h)$ for $k\neq 0$ are chosen as (\ref{eq:kappa_k_splitting}), then the local error of the processed scheme reads
\begin{equation*}
\hat \psi_h(x)-\varphi_h(x)= \e^{h A} \left( h^2\sum_{k,\ell \in \mathcal{I}}
\big(\hat \alpha_{k \ell}(h) -\alpha_{k \ell}(h) \big)
g'_{\ell}(x)\, g_{k}(x) + \mathcal{O}(h^3)\right)
\end{equation*}
and the $h$-flow of $\hat{\psi}_h$ is also Hamiltonian, with the modified Hamiltonian function
    \begin{equation*}
\hat H(x; h) =  \frac12\, x^T Q x + \sum_{k \in  \mathcal{I}} \hat \beta_k(h)\, G_k(x) 
   +  \frac{h}{2} \sum_{k, \ell \in  \mathcal{I}} \hat \beta_{k \ell}(h)\, \{G_k,G_{\ell}\}(x),
\end{equation*}
where
\begin{equation} \label{eq:betahat}
\hat \beta_{k}(h) =1 ,  \qquad\quad
\hat \beta_{k \ell}(h) =
\frac{ \hat \alpha_{k \ell}(h) -  \alpha_{k \ell}(h) }{\alpha_{k+\ell}(h)}.
\end{equation}
Equivalently,
    \begin{equation}
    \label{eq:mHamhat}
\hat H(x; h) =  H(x)
   +  \frac{h}{2} \sum_{k, \ell \in  \mathcal{I}} \hat \beta_{k \ell}(h)\, \{G_k,G_{\ell}\}(x),
\end{equation}
which is an $\mathcal{O}(h)$ perturbation of the original Hamiltonian, and thus one may expect that the value of the $H(x)$ (typically, the energy of the original system) will be well approximated for relatively large time intervals. A word of caution is in order here: by construction, the difference between the processed map $\hat \psi_h$ and the $h$-flow of the modified Hamiltonian (\ref{eq:mHamhat}) is formally of order $\mathcal{O}(h^3)$, but the constant in 
$\mathcal{O}(h^3)$ depends on the size of the modulus of the coefficients
\begin{equation*}
\kappa_k(h), \quad \kappa_{k,\ell}(h), \quad \hat \beta_{k \ell}(h), \quad \mbox{for} \quad k,\ell \in \mathcal{I}.
\end{equation*}
For near-resonant step sizes $h$, that is, for step sizes such that $\e^{i k \omega h}-1$ is small for some index $k$ belonging to $\mathcal{I}$ or $\mathcal{I} +\mathcal{I}$, the size of some of these critical coefficients may become large. In such cases, one cannot expect that the processed map $\hat \psi_h$ will be close from the $h$-flow of the modified Hamiltonian (\ref{eq:mHamhat}).
 
 \subsection{A processed Strang scheme}
 \label{apss} 
  In the particular case in which $\psi_h$ is the Strang splitting $S_h^{[RKR]}$, the coefficients $\kappa_k(h)$ in the corresponding expansion
(\ref{eq:chiexpansion}) verifying (\ref{eq:kappa_k_splitting}) read  
\begin{equation}
\label{eq:kappa_k}
  \kappa_{k}(h)
  = 	\frac{1}{i k \omega h} \left(1 - \mathrm{sinc}(k \omega h/2)^{-1}\right), \qquad k \in \mathcal{I}\backslash\{0\},
\end{equation}
and $\kappa_{-k}(h)=-\kappa_k(h)$ for all $k\neq 0$.

Assuming that the potential function $U(q)$ in (\ref{eq:hamHOS}) is a polynomial of degree $m$, so that  
  \begin{equation*}
\mathcal{I} \subset \{-m,\ldots,-1,0,1,\ldots m\},
\end{equation*}
we next construct a fully explicit processor $\pi_h$ for $S_h^{[RKR]}$. We define it as a composition of basic flows as follows:
  \begin{equation} 
  \label{eq:proc_splitting}
\pi_h = \varphi^{[R]}_{\alpha}\circ \varphi^{[K]}_{b_{2m}(h)}\circ \varphi^{[R]}_{\alpha}\circ  \varphi^{[K]}_{b_{2m-1}(h)} \circ \cdots \circ
\varphi^{[K]}_{b_{2}(h)} \circ  \varphi^{[R]}_{\alpha}\circ
 \varphi^{[K]}_{b_{1}(h)} \circ \varphi^{[R]}_{\alpha},
\end{equation}
where $\alpha=\frac{2\pi}{2m+1}$, and the coefficients $b_j(h)$, $j=1,\ldots, 2m$, depend on $h$. Since $S_h^{[RKR]}$ is time-symmetric, then
it makes sense to construct the processor $\pi_h$ such that $\hat \psi_h = \pi_h^{-1} \circ \psi_h \circ \pi_h$ is also time-symmetric. This can
be achieved by requiring that $\pi_{-h} = \pi_{h}$, or equivalently, by requiring that $b_j(-h)=b_j(h)$ for $j=1,\ldots, 2m$. This condition, together with
\[
  \kappa_0(h)=0, \qquad \tilde \alpha_k(h)=\alpha_k(h) \quad \mbox{ for all } \quad k \in \{-m,\ldots,-1,0,1,\ldots m\}
\]
uniquely determines the $b_j(h)$ coefficients as
  \begin{equation}
 \label{eq:bsol}
 b_j(h)  = -b_{2m-j+1}(h) = \frac{2}{2m+1} \sum_{k=1}^m  \frac{1}{k} \left(\mathrm{sinc}(k \omega h/2)^{-1} - 1 \right)\,  \sin\left(\frac{2 k j \pi}{2 m+1}\right),
\end{equation}
 for $j=1,\ldots,m$. 
This can be seen as follows: successive application of the composition rule (\ref{eq:gammacomp}) shows that (\ref{eq:proc_splitting}) admits the expansion 
(\ref{eq:chiexpansion}) with
\begin{align}
 \label{eq:bsystem}
h \, \kappa_{k}(h) &=  \sum_{j=1}^{2m} b_j(h) \, \e^{2 i k j \pi/(2m+1)}, \quad k=-m,\ldots,-1,0,1,\ldots m,
\end{align}
and for $k,\ell \in \mathcal{I}$, 
 \begin{align}
 \label{eq:kappakl}
h^2 \kappa_{k \ell}(h) &= \sum_{1\leq j <  n \leq 2m} b_{j}(h)  \, b_{n}(h) \,  \e^{2 \pi (j k + n \ell)/(2m+1)} \\
\nonumber
&+
\sum_{1\leq j \leq  2m} \frac{1}{2} \,  b_{j}(h)^2 \, \e^{i2 \pi  j (k + \ell)/(2m+1)}.
\end{align}
One can check that $\kappa_0(h)=0$ and the symmetry condition $\kappa_{-k}(h)=-\kappa_k(h)$ ($k\neq 0$) that holds for (\ref{eq:kappa_k}) imply that $b_j(h)  = -b_{2m-j+1}(h)$ for $j=1,\ldots,2m$.  Now, (\ref{eq:bsystem}) means that application of the inverse discrete Fourier transform to the vector $(0,b_1,b_2,\ldots,b_m, -b_m,\ldots,-b_1)$ gives the vector \
\begin{equation*}
 \frac{h}{2m+1}(0,\kappa_1,\ldots,\kappa_m,-\kappa_m,\ldots,-\kappa_1).
\end{equation*}
Equivalently, the former is obtained by applying the discrete Fourier transform to the latter. Rearranging terms, one finally arrives at (\ref{eq:bsol}).

The modified Hamiltonian (\ref{eq:mHamhat}) can be obtained in the following way: first, 
the coefficients (\ref{eq:kappakl}) can be used to
compute $\hat \alpha_{k\ell}(h)$ from (\ref{eq:halphak}),
 and then the coefficients $\hat \beta_{k\ell}(h)$ are determined from (\ref{eq:betahat}).

\paragraph{Example: simple pendulum.}  
  
As an illustrative example we again consider the simple pendulum of Subsection \ref{subsec.1.3}, described by the Hamiltonian function
(\ref{pendulum1}). As we saw there, for initial conditions in a neighborhood of the stable equilibrium $(0,0)$ it is advantageous to decompose $H$
as in (\ref{pendulum2}), so that it constitutes a particular example of system (\ref{eq:hamHOS}):
\[
  H = H_1 + H_2, \qquad H_1(q,p) = \frac{1}{2} p^2 + \frac{1}{2} q^2, \qquad H_2(q,p) = U(q) = 1 - \frac{1}{2} q^2 - \cos q.
\]
Although $U(q)$ is not a polynomial, and therefore the set $\mathcal{I}$ in (\ref{eq:I}) is infinite, we can truncate the Fourier expansion (\ref{expan_H}) and work instead with
$\mathcal{I} = \{ -m, \ldots, -1, 0, 1, \ldots, m \}$ for a given $m$ to construct a processed Strang scheme based on kicks and rotations as proposed earlier.
Specifically, we take $m=4$, a step size $h=5/6$, then determine the corresponding coefficients $b_j(h)$ and form the integrator 
$\hat{\psi}_h = \pi_h^{-1} \circ S_h^{[RKR]} \circ \pi_h$, with $\pi_h$ given by (\ref{eq:proc_splitting}) (with $m=4$), 
whereas $S_h^{[RKR]}$ corresponds to the map (\ref{soint}).

Figure \ref{figss1} shows the relative error in energy (left) and phase space (right) over the time interval $[0,500]$ 
corresponding to the solution initiated at $(q_0,p_0) = (1/10,0)$. As usual, $S_2$ denotes the 
St\"ormer--Verlet method applied to (\ref{pendulum1}), (2,2) is the Strang splitting $S_h^{[RKR]}$, and  P$(2,2)$  is the
processed $S_h^{[RKR]}$ scheme with step size $h=5/6$. 
Since the integrator $(2,2)$ is conjugate to P$(2,2)$,  the error of the former is eventually dominated by the error of the latter, as expected from the discussion in
subsection \ref{sub_processing1}.   Observe that the step size used here is larger than those taken in Figure \ref{fig_Pendulo}.

 \begin{figure}[htb] \label{figss1}
\centering
\includegraphics[scale=0.52]{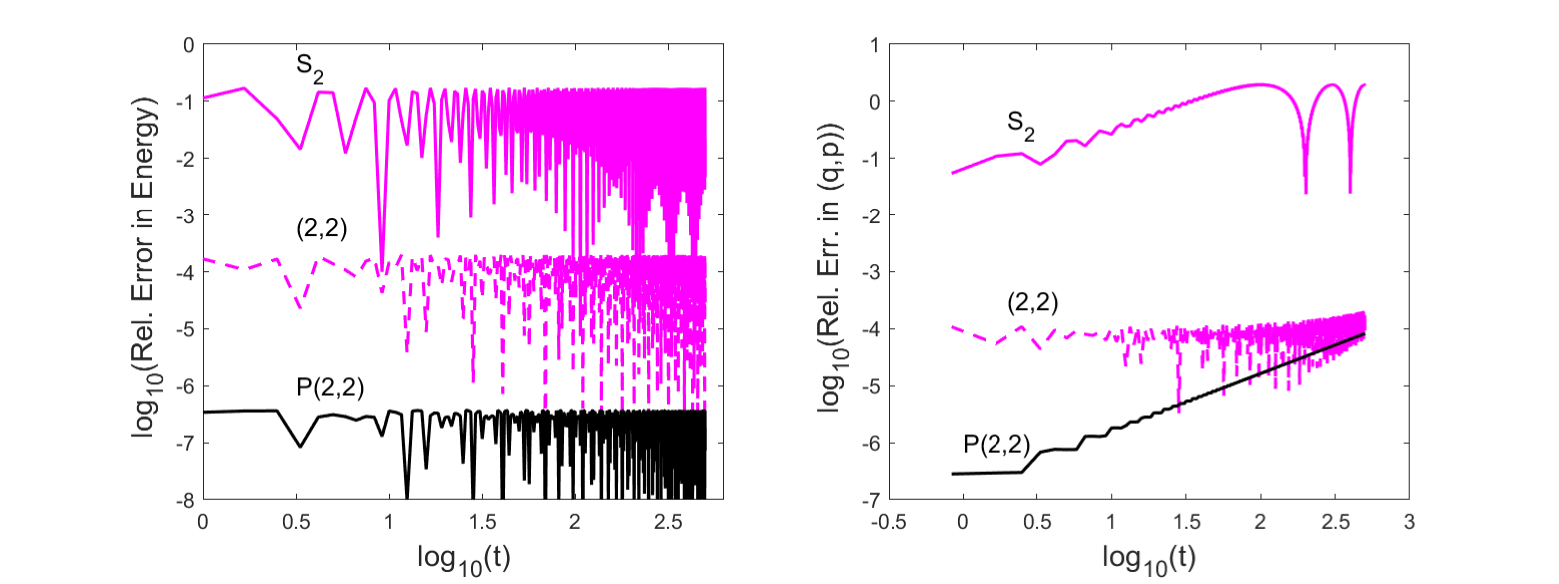} 
 \caption{Pendulum. Evolution of relative errors in energy (left) and in phase space (right) obtained with St\"ormer--Verlet, $S_2$, scheme $S_h^{[RKR]}$,
 denoted (2,2), and the
 processed (2,2), P$(2,2)$ with initial state $(q_0,p_0)=(1/10,0)$ and step size $h=5/6$.}
 \end{figure} 

Next, we consider the Strang method $S_h^{[RKR]}$ and its processed version P$(2,2)$ for the same initial condition and a time interval $[0,500]$
and show the maximum energy error in this interval for step sizes in the range $[0, 3 \pi]$ (Figure \ref{figss2}).
The spikes in the curve of the error of Strang correspond to the step-sizes violating the non-resonance condition (\ref{eq:non-resonant-h1}) for $\omega=1$ and $k=2,4$. Indeed, the potential $U(q)$ can be well approximated near the origin by $q^ 4/24$, which implies that $\mathcal{I}=\{-4,-2,0,2,4\}$. The curve of the error of processed Strang also has spikes for such resonant step sizes, and additionally, for the step sizes $h$ such that $3 h/(2\pi) \in \mathbb{Z}$, which is due to the fact that the processor map $\pi_h$ has been designed to work well for problems for which $\mathcal{I}=\{-4,-3,-2,-1,0,1,2,3,4\}$. Away from such resonant step sizes, the energy error of processed Strang is considerably smaller than that of unprocessed Strang.

 \begin{figure}[h] \label{figss2}
\centering
\includegraphics[scale=0.5]{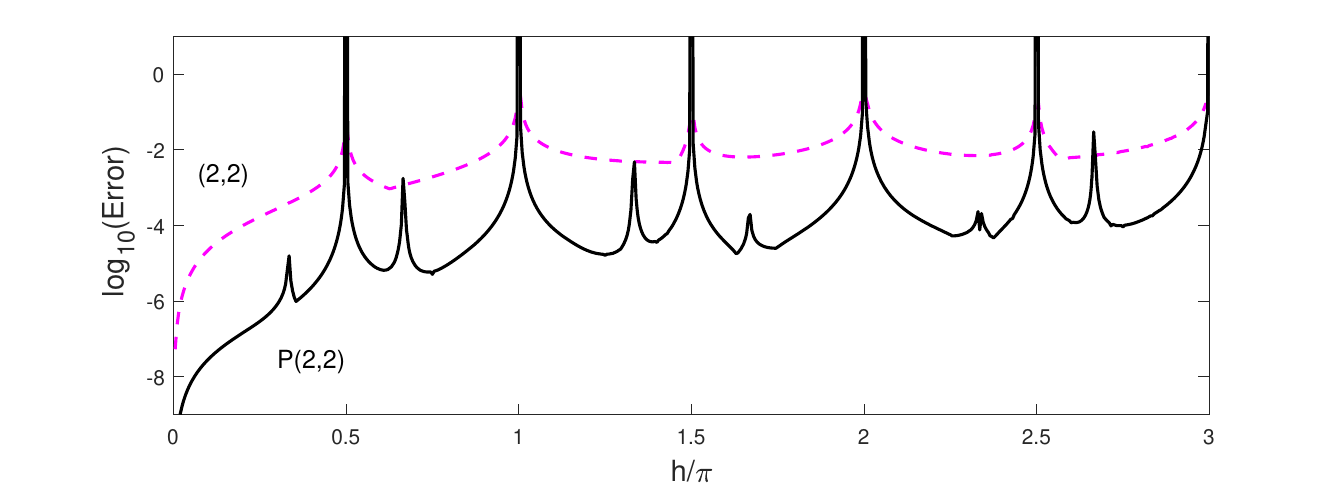} \\
 \caption{Pendulum. Maximum relative energy error in the time interval $[0,500]$, initial state $(q_0,p_0)=(1/10,0)$, and step sizes in the range $[0,3\pi]$ committed 
 by the Strang splitting $S_h^{[RKR]}$ and its processed version P$(2,2)$.}
 \end{figure}

 \subsection{More general assumptions}
 \label{ss:moregeneral}
 
 The results given in the present section are also valid with small modifications in the more general setting of systems (\ref{eq:odeHOS}) such that the eigenvalues of $A$ lie on the imaginary axis. In the most general case, where such eigenvalues are not integer multiples of $\omega\, i$, the exponential 
 $\e^ {t A}$ is quasi-periodic with a finite number of basic frequencies $(\omega_1,\ldots,\omega_r) \in \mathbb{R}^ r$. Under such conditions,  
 $\e^{-t A}  f_2(\e^{t\, A} x)$ admits a multi-variable Fourier expansion 
\begin{equation}
\label{eq:Fourierg}
 \sum_{k \in  \mathbb{Z}^r} e^{i\, \langle k, \omega \rangle\, t}\, g_k(x),
\end{equation}
where $k$ now denotes a $r$-tuple of integers $k=(k_1,\ldots,k_r)\in \mathbb{Z}^r$, $\omega=(\omega_1,\ldots,\omega_r)$ is the vector of basic frequencies, and $\langle k,\omega\rangle = k_1 \omega_1 + \cdots + k_r \omega_r$.  The set of indices $\mathcal{I}$ is then defined as
\begin{equation*}
\mathcal{I} = \{k \in  \mathbb{Z}^r\ : \ g_k \neq 0\}.
\end{equation*}
 Under these more general assumptions and notation, all previous formulae are valid if each occurrence of $k \omega$ with $k\in \mathcal{I}$ is replaced by $\langle k,\omega\rangle$, the only exception being the explicit construction of the processing map $\pi_h$ carried out in subsection \ref{apss}. It is worth remarking that, compared to the periodic case, in the quasi-periodic one there are typically more resonant step sizes, that is, step sizes $h$ such that
 \begin{equation*}
\frac{\langle k, \omega \rangle h}{2\pi} \in \Z\backslash\{0\} \quad \mbox{ or } \quad \frac{\langle (k+\ell), \omega\rangle h}{2\pi} \in \Z\backslash\{0\}
\end{equation*}
for some $k,\ell \in \mathcal{I}$.


\section{Splitting methods for PDEs}
 \label{sect6}
 
 \subsection{Splitting, LOD and ADI methods}
 
 Splitting methods can also be applied to partial differential equations (PDEs), in which case equation (\ref{ivp.1}) has to be viewed
 as the abstract system associated with the PDE initial value problem in autonomous form, and $f$ as a spatial partial differential operator. For clarity,
in this section  we write 
 \begin{equation} \label{pde.1}
   u_t(x,t) = f(x,u(x,t)), \qquad u(x,0) = u_0(x)
\end{equation}
to distinguish between the unknown $u(x,t)$ defined in a certain function space and the spatial variable $x \in \mathbb{R}^d$, but for notational purposes
it is convenient to drop the dependence of $u$ and $f$ on $x$. To introduce the basic concepts and methods it is not necessary at this stage
to specify the dimension $d$,
the number of components of $u$, the relevant function space and the boundary conditions. These will be considered when analyzing particular applications. Moreover, for 
simplicity in the presentation, only the autonomous case will be treated. If $f$  depends explicitly on time, then one can take $t$ as a new
coordinate, as in done 
in subsection \ref{sub.nas}.

Very often, the operator $f$ contains contributions coming from very different physical
sources, so that one may decompose it into two (or more) parts, 
and use different schemes to solve each sub-problem approximately. For instance,
in a
reaction-diffusion system, $f(u) =  \nabla \cdot (D \nabla u) + g(u)$, where $D$ and $g$ may also depend on $x$, it makes sense to split the diffusion
from the reaction terms, i.e.,
\[
  f(u) = f_1(u) + f_2(u), \qquad \mbox{ with } \qquad f_1(u) = \nabla \cdot (D \nabla u), \qquad f_2(u) = g(u).
\]
Algorithm  \ref{alg-LT} can of course be used in this setting, but one still has to specify how to solve each initial value sub-problem in practice with appropriate
boundary conditions. A simple
possibility consists in applying the backward Euler scheme, thus resulting in the so-called \emph{Marchuk--Yanenko} operator-splitting scheme
 \begin{equation} \label{m-y}
 \begin{aligned}
  & u_{n+1/2} = u_n + h f_1(u_{n+1/2}) \\
  & u_{n+1} = u_{n+1/2} + h f_2(u_{n+1}), \qquad n = 0, 1, 2, \ldots
 \end{aligned}
\end{equation}
where, for consistency with the rest of the paper, we have denoted $h \equiv  \Delta t$, the time step size. In spite of its low order of consistency 
(order 1)\footnote{Here `order' should be 
understood as the order of consistency with respect to the solution of
the ODE problem on a fixed spatial grid, not with respect to the underlying PDE solution \cite{hundsdorfer03nso}.} and the fact that the intermediate stage
$u_{n+1/2}$ is not a consistent approximation to the exact solution, its simplicity and robustness make it a useful alternative way to deal with complicated
problems and even non-smooth operators \cite{glowinski16smi}. It is also appropriate for parabolic problems, since it incorporates the damping properties
of the backward Euler method \cite{hundsdorfer03nso}.

If instead of using the backward Euler scheme to integrate each sub-problem in Algorithm  \ref{alg-LT} one applies the second-order implicit trapezoidal rule, 
it results in Yanenko's Crank--Nicolson method \cite{hundsdorfer03nso,marchuk90saa}:
\begin{equation} \label{y-cn}
 \begin{aligned}
  & u_{n+1/2} = u_n + \frac{h}{2} f_1(u_{n})  + \frac{h}{2} f_1(u_{n+1/2}) \\
  & u_{n+1} = u_{n+1/2} + \frac{h}{2} f_2(u_{n+1/2}) +  \frac{h}{2}  f_2(u_{n+1}), \qquad n = 0, 1, 2, \ldots
 \end{aligned}
\end{equation}
In the end, however, it is also of first order of consistency.

Another widely popular class of spitting methods in the domain of PDEs is the Peaceman--Rachford scheme and its variants. Although
initially designed for the numerical solution of elliptic and parabolic equations \cite{peaceman55tns,douglas56otn}, they also apply to more general situations. 
The procedure goes as follows. Given an approximation $u_n$ for the solution of (\ref{pde.1}) at $t = t_n$, the approximation $u_{n+1}$ is computed using
the backward (respectively, forward) Euler scheme with respect to $f_1$ (resp., $f_2$) on the  sub-interval $[t_n, t_{n+\frac{1}{2}}]$. Then, the roles of $f_1$
and $f_2$ are interchanged on the sub-interval $[t_{n+\frac{1}{2}}, t_{n+1}]$. In other words, the \emph{Peaceman--Rachford} scheme corresponds to
the sequence
\begin{equation} \label{p-r}
 \begin{aligned}
  & u_{n+1/2} = u_n + \frac{h}{2} f_1(u_{n+1/2}) + \frac{h}{2}  f_2(u_n) \\
  & u_{n+1} = u_{n+1/2} + \frac{h}{2} f_1(u_{n+1/2}) + \frac{h}{2} f_2(u_{n+1}), \qquad n = 0, 1, 2, \ldots,
 \end{aligned}
\end{equation}
which is of second order of consistency.
Notice that, in contrast to methods (\ref{m-y}) and (\ref{y-cn}), both $f_1$ and $f_2$ appear in each of the two stages, and thus the intermediate value
$u_{n+1/2}$ provides a consistent approximation at $t= t_{n+1/2}$. On the other hand, it does not possess a natural formulation where $f$
is split into more than two operators. In general, 
$f_1$ and/or $f_2$ can be nonlinear, unbounded and even multi-valued. For a more detailed
treatment, the reader is addressed to \cite{glowinski16sfa} and references therein.

A classical alternative to scheme (\ref{p-r}), of first order, is the \emph{Douglas--Rachford} method \cite{douglas56otn}, which instead reads 
\begin{equation} \label{d-r}
 \begin{aligned}
  & \hat{u}_{n+1} = u_n + h f_1(\hat{u}_{n+1}) + h f_2(u_n) \\
  & u_{n+1} = u_n + h f_1(\hat{u}_{n+1}) + h f_2(u_{n+1}), \qquad n = 0, 1, 2, \ldots
 \end{aligned}
\end{equation}
and can be generalized to decompositions of $f$ involving more than two operators. Notice that the roles of $f_1$ and $f_2$ in (\ref{d-r}) are not symmetric,
in contrast to the Peaceman--Rachford method. On the basis of many numerical experiments, \cite{glowinski16sfa} conclude that the scheme
(\ref{d-r}) is more robust and faster than (\ref{p-r}) for problems where one of the operators is non-smooth, in particular when one is interested in approximating
steady state solutions.

Let us now analyze the particular case when $f$ in (\ref{pde.1}) is a linear spatial differential operator. Assuming that an appropriate semidiscretization of 
(\ref{pde.1}) in the space variable $x$ has been carried out, one ends up with the system
\begin{equation} \label{pde.2}
  \frac{dU}{dt} = F_1 U + F_2 U,
\end{equation}  
 where $F_1, F_2 \in \mathbb{C}^{M \times M}$, $ F_1 F_2 \ne F_2 F_1$ in general, and $U \in \mathbb{C}^M$ approximates $u$ on the space
 grid points $x_1, \ldots, x_M$ (see below). Then, a step
of the Marchuk--Yanenko scheme (\ref{m-y}) reads
\begin{equation} \label{m-y2}
  U_{n+1} = (I - h F_2)^{-1} \, (I - h F_1)^{-1} U_n,
\end{equation}
where $U_n \approx (u(x_1,t_n), \ldots, u(x_M, t_n))^T$. Notice that (\ref{m-y2}) corresponds to applying the $[0/1]$ Pad\'e approximant to the exponentials in the Lie--Trotter scheme
$U_{n+1} = \e^{h F_2} \, \e^{h F_1} U_n$, whence the first order of the approximation is obtained at once. 
On the other hand, applying (\ref{y-cn}) results in the sequence
\begin{equation} \label{y-cn2}
  U_{n+1} = (I - \frac{h}{2} F_2)^{-1} \, (I +\frac{h}{2} F_2) \, (I - \frac{h}{2} F_1)^{-1} \, (I + \frac{h}{2} F_1) \, U_n.
\end{equation}  
In other words, it corresponds to the application of the $[1/1]$ Pad\'e approximant to the exponentials in the Lie--Trotter scheme. It is then clear that, although
a second order Crank--Nicolson method is carried out for each exponential, their combination (\ref{y-cn2}) is of first order, since in the end we are using only a
variant of the Lie--Trotter scheme. 

By the same token,  a straightforward computation shows that the Peaceman--Rachford scheme (\ref{p-r}) applied to the linear equation (\ref{pde.2}) can be written as
\begin{equation} \label{p-r2}
  U_{n+1} = (I - \frac{h}{2} F_2)^{-1} \, (I + \frac{h}{2} F_1) \, (I - \frac{h}{2} F_1)^{-1} \, (I + \frac{h}{2} F_2) \, U_n.
\end{equation}
As a matter of fact, all these algorithms can be formulated by applying properly chosen compositions of the implicit and explicit Euler methods. 
This observation may eventually lead to the construction of methods of higher order and/or improved behavior, but in the same family.
Suppose $f$ in (\ref{pde.1}) is of the form $f = f_1 + f_2$ and that the solution of each equation $u_t=f_k(u)$, $k=1,2$, is numerically approximated by the maps
\[
\begin{aligned}
  & u_{n+1}=\psi_{h}^{ke}(u_n) \equiv u_n+ h f_k(u_n), \qquad \quad \mbox{ explicit Euler } \\
  & u_{n+1}=\psi_{h}^{ki}(u_n) \equiv u_n+ h f_k(u_{n+1}), \qquad\;\,  \mbox{ implicit Euler } 
\end{aligned}  
\]
so that, by combining all variants, we form the following first order schemes (and their corresponding adjoints):
\begin{eqnarray} \label{options1o}
  &&  \phi_{h}^{1}=\psi_{h}^{1e} \circ \psi_{h}^{2e}, \qquad 
\qquad 
    \phi_{h}^{1*}=\psi_{h}^{2i} \circ \psi_{h}^{1i}, \nonumber \\ 
  &&  \phi_{h}^{2}=\psi_{h}^{1e} \circ \psi_{h}^{2i}, \qquad 
\qquad 
    \phi_{h}^{2*}=\psi_{h}^{2e} \circ \psi_{h}^{1i}, \\ 
  &&  \phi_{h}^{3}=\psi_{h}^{1i} \circ \psi_{h}^{2e}, \qquad 
\qquad 
    \phi_{h}^{3*}=\psi_{h}^{2i} \circ \psi_{h}^{1e}, \nonumber \\ 
  &&  \phi_{h}^{4}=\psi_{h}^{1i} \circ \psi_{h}^{2i}, \qquad 
\qquad 
    \phi_{h}^{4*}=\psi_{h}^{2e} \circ \psi_{h}^{1e}.  \nonumber 
\end{eqnarray}
Then, one can conclude the following:
\begin{itemize}
	\item The Marchuk--Yanenko operator-splitting scheme \eqref{m-y} can be expressed as $u_{n+1}= \phi_h^{1*}(u_n)$. Therefore, 
	compositions 
	\[
	  u_{n+1}= \phi_{{h}/2}^{1*}\circ\phi_{{h}/2}^{1}(u_n)
		\qquad \mbox{ and} \qquad
	  u_{n+1}= \phi_{{h}/2}^{1}\circ\phi_{{h}/2}^{1*}(u_n)
	\]
	yield time-symmetric second order approximations.
	\item Yanenko's Crank--Nicolson method \eqref{y-cn} corresponds to the composition
	\begin{equation} \label{ycn}
	  u_{n+1}=  \psi_{{h}/2}^{2i}\circ\psi_{{h}/2}^{2e}\circ\psi_{{h}/2}^{1i}\circ\psi_{{h}/2}^{1e}(u_n),
	\end{equation}
	which is not time-symmetric and therefore only of first order. Notice, however, that (\ref{ycn}) can also be expressed as 
	\[
	  u_{n+1} =\psi_{{h}/2}^{2i}\circ\left( \phi_{{h}/2}^{2*}\circ\phi_{{h}/2}^{2} \right)\circ (\psi_{{h}/2}^{2i})^{-1} (u_n),
	\]
	so it is conjugate to a time-symmetric second order method. This feature could account for its observed good behavior in practice.
	\item The Peaceman--Rachford scheme \eqref{p-r} is just the symmetric composition
	\[
	 u_{n+1}=  \phi_{{h}/2}^{3*}\circ \phi_{{h}/2}^{3} (u_n),
	\]
	and is therefore of second order.
	It is also worth mentioning that whereas all symmetric second order compositions 
	\[
	\phi_{{h}/2}^{\ell*}\circ \phi_{{h}/2}^{\ell}\qquad \mbox{ and } \qquad 
	\phi_{{h}/2}^{\ell}\circ \phi_{{h}/2}^{\ell*}, \qquad \ell=1,\ldots,4
	\]
provide consistent approximations at the midpoint, the steady state solution is captured only when $\ell=3$: 
 if $f(w)=0$, then 
 \[
 \phi_{{h}/2}^{3*}\circ \phi_{{h}/2}^{3}(w)=\phi_{{h}/2}^{3}\circ \phi_{{h}/2}^{3*}(w)=w.
 \]
\end{itemize}
The Douglas--Rachford method \eqref{d-r} can be alternatively formulated as a splitting method in an extended space as follows. 
Consider the enlarged system 
\begin{equation} \label{mts2_1}
 w_t= \left( \begin{array}{c}
          \hat{u}_t \\
          v_t 
         \end{array} \right)  =  
   \left( \begin{array}{c}
          f_1(\hat{u})+f_2(v)\\
          0 
         \end{array} \right)+ 
            \left( \begin{array}{c}
          0  \\
          f_1(\hat{u})+f_2(v)\ 
         \end{array} \right)
				= g_1(w) + g_2(w),
\end{equation}
with initial condition $w(0)= \left( \hat{u}(0),v(0) \right)=\left(u_0,u_0\right)$ and solution $w(t) = (\hat{u}(t), v(t)) = (u(t), u(t))$.  Take $w_n=(u_n,u_n)$ and form
the composition
\[
	\hat w_{n+1}=  \phi_{h}^{1*} (w_n) = \psi_{h}^{2i}\circ \psi_{h}^{1i}(w_n), \qquad n \ge 0,
\]
providing $\hat w_{n+1}=(\hat{u}_{n+1},v_{n+1})$. Then scheme  \eqref{d-r} is recovered by taking $u_{n+1}=v_{n+1}$, and considering
$w_{n+1}=({u}_{n+1},u_{n+1})$ as the starting point for the next iteration.

\paragraph{Example.}
The two-dimensional heat equation with source term and Dirichlet boundary conditions on the unit square may serve as illustration of these methods
\cite{hundsdorfer03nso}. 
Specifically, the system reads
\begin{equation} \label{h-e}
 \begin{aligned}
   & u_t = u_{xx} + u_{yy} + g(x,y,t) \qquad \mbox{ on } \, \Omega = (0, 1) \times (0, 1) \\
   & u(x,y,t) = u_{\Gamma}(x,y,t) \qquad \mbox{ on } \Gamma =\partial\Omega \\
   & u(x,y,0) = u_0(x,y) \qquad \mbox{ on } \, \Omega.
 \end{aligned}
\end{equation}
Suppose we take a Cartesian grid in $\Omega$ based on $M+1$ equally spaced intervals  in the $x$ and $y$ directions, so that $\Delta x=\Delta y=\frac1{M+1}$
and apply finite differences to approximate the space derivatives. Then
we end up with $M^2$ interior points and the aim is to get approximations of $u$ at these points, i.e., to determine $u_{i,j}(t)\simeq u(x_i,y_j,t)$ for $i,j=1,2,\ldots,M$.
The problem can be conveniently formulated in terms of the `supervector' $U=({u}_1^T,\ldots,{u}_M^T)^T$, with ${u}_i=(u_{i,1},\ldots,u_{i,M})^T$, and $M^2$ components
\[
  U_{\ell}(t)=u_{i,j}(t), \qquad \ell= j + M (i-1). 
\]
Using standard second-order finite differences for $\partial_{xx}$ and $\partial_{yy}$ we get a linear system
of the form $U' = F \, U + s(t)$, where $s(t)$ contains sources and boundary data, and $F = F_1 + F_2$, with
\[
   F_1 = I_M \otimes B_M, \qquad F_2 = B_M \otimes I_M, \qquad
  B_M = \frac1{(\Delta x)^2}\left(\begin{array}{cccc}
  -2 & 1 & \cdots & 0 \\
  1 & -2 & \ddots & \vdots\\
   \vdots & \ddots &\ddots &1  \\
  0 &  & 1 & -2
 \end{array}  \right).
\]
Here $I_M$ is the identity matrix, $B_M$ the differentiation matrix and $\otimes$ denotes the tensor product. 
Notice that $F_1$ acts in the $x$-direction and $F_2$ in the $y$-direction, so 
both $F_1$ and $F_2$ are essentially one-dimensional.
Since $F_1$ is a tridiagonal matrix, the linear system resulting from the application of the previous implicit schemes can be solved in an efficient
way, whereas $F_2$ is equivalent to a tridiagonal matrix, so the same considerations apply \cite{iserles96afc}.

First consider the homogeneous equation with zero Dirichlet boundary conditions. The solution for one time step is
\[
  U(t_n+h) = \e^{{h}(F_1+F_2)} \, U(t_n),
\]
whereas the application of the Peaceman--Rachford scheme leads to the approximation (\ref{y-cn2})
\[
   U_{n+1} = (\phi^{3*}_{h/2}\circ\phi^{3}_{h/2}) \, U_n = \e^{{h}(F_1+F_2)} \, U_n +{\cal O}({h}^3). 
\]
When a dimensional splitting is done, as in this case,
applying schemes (\ref{m-y2}) and (\ref{y-cn2}) corresponds essentially to carrying out computations in only one dimension.
This is the reason why (\ref{m-y}) and
(\ref{y-cn}) are called \emph{locally one-dimensional (LOD)} methods. Analogously, given the alternate use of $F_1$ and $F_2$ in this setting, the name
\emph{alternating direction implicit (ADI)} is usually attached to methods (\ref{p-r}) and (\ref{d-r}).

Regarding these ADI/LOD methods as composition schemes allows one to get approximations in the non-homogeneous case too, with only minor changes,
while keeping any favorable properties (if no order reduction occurs due to the Dirichlet boundary conditions). To illustrate this point,
we again apply the explicit and implicit Euler methods to the equation $U' = s(t)$ corresponding to the non-homogeneous term:
\[
  U_{n+1}= \psi_{h}^{s_e}(U_n) = U_n+{h} \, s(t_n), 
	\qquad 	\qquad 
  U_{n+1}=\psi_{h}^{s_i}(U_n)=U_n+{h} \, s(t_{n+1}),
\]
and consider the first order scheme $\hat \phi_{h}^3=\psi_{{h}/2}^{1i}\circ \psi_{{h}/2}^{2e}\circ \psi_{h}^{s_e}$. Then
\[
\begin{aligned}
 & U_{n+1}  = \left( \hat \phi^{3*}_{{h}/2}\circ\hat \phi^{3}_{{h}/2} \right) \, U_n  = 
    \left( \psi_{{h}/2}^{s_i}\circ \psi_{{h}/2}^{2i}\circ \psi_{{h}/2}^{1e}\circ 
 \psi_{{h}/2}^{1i}\circ \psi_{{h}/2}^{2e}\circ \psi_{{h}/2}^{s_e} \right) \, U_n \\
  & = \frac{h}{2} s(t_{n+1}) +
\left(I-\frac{h}2 F_2\right)^{-1}
\left(I+\frac{h}2 F_1\right)
\left(I-\frac{h}2 F_1\right)^{-1}
\left(I+\frac{h}2 F_2\right) \Big(U_n+ \frac{h}2 s(t_{n})\Big)
\end{aligned}
\]
produces a symmetric second-order scheme.

\noindent
$\Box$

\

It is not only reaction-diffusion problems that can be treated in this way. For instance, systems of hyperbolic conservation laws in three dimensions, such as
\[
  u_t + \nabla \cdot f(u) = 0, \qquad\quad  u(x,y,z,0) = u_0(x,y,z)
\]
can be numerically approximated with dimensional splitting by applying a specially tailored numerical scheme to each scalar conservation law
$u_t + f(u)_x = 0$, etc.  \cite{holden10smf}.  

The one-dimensional convection-diffusion problem
\[
  u_t + f(u)_x = A(u)_{xx}, \qquad u(x,0) = u_0(x)
\]
with a scalar non-decreasing function $A(\cdot)$, $A(0) = 0$, possesses a rich set of phenomena depending on the interplay of the different nonlinearities. In
this case one can split the problem into a convective and a diffusion part and apply Algorithm \ref{alg-LT}. This formally results in the so-called Godunov
split \cite{holden10smf}. Specifically, denoting the solution operator corresponding to the scalar conservation law $u_t + f(u)_x = 0$
as $\mathcal{U}_{h}^{[1]}$, and  the one corresponding to the (weak) solution of the nonlinear heat equation $u_t = A(u)_{xx}$ by
$\mathcal{U}_{h}^{[2]}$, then the scheme reads
\[
  u(x, t_n+ h) \approx u_{n+1} = \mathcal{U}_{h}^{[2]} (\mathcal{U}_{h}^{[1]} (u_n)).
\]
Of course, one can also use the Strang splitting
\[
  u(x, t_n+ h) \approx u_{n+1} = \mathcal{U}_{h/2}^{[2]} (\mathcal{U}_{h}^{[1]} (\mathcal{U}_{h/2}^{[2]}  (u_n))).
\]
To get a numerical approximation, each of the two operators must be approximated. This can be done, for example, by a front-tracking method for $\mathcal{U}_{h}^{[1]}$
and by a standard implicit finite-difference method for the parabolic operator $\mathcal{U}_{h}^{[2]}$. A convergence analysis of such schemes has been
carried out in \cite{holden10smf}.   

Although  the analysis of splitting methods can be done by 
power series expansions and the formalism of Lie operators, there are
fundamental differences with respect to the ODE case. 
Nonlinear PDEs in general possess solutions that
exhibit complex behavior in small regions of space and time, such as sharp transitions and discontinuities, and thus
they lack the usual smoothness required for the analysis. Moreover, 
even if the exact solution of the original problem is smooth, it may well happen that the composition
defining the splitting method provides non-smooth approximations. Therefore, it is necessary to develop 
an appropriate mathematical framework to
analyze the convergence of the numerical solution to the correct solution of the original problem, and this has to be done very often on
a case by case basis (see  e.g. \cite{holden10smf} and references therein). Thus, in particular, the first- and second-order convergence of
the Godunov and Strang splitting methods on the Korteweg--de Vries equation $u_t - u u_x + u_{xxx} = 0$ has been proved 
in \cite{holden11osf} if the initial data are sufficiently regular, 
whereas the result has been extended in \cite{holden13osf} to equations of the form $u_t = A u + u u_x$, when $A$ is a linear
differential operator such that the full equation is well-posed. More recently, convergence results and error estimates have also been obtained 
for initial conditions with low regularity, namely $u_0 \in H^{s}$ with $0 < s \le 3/2$, where $H^s$ denotes the Sobolev space \cite{rousset21cee}.

Another source of difficulties related with the application of splitting methods to PDEs is the treatment of boundary conditions. In this respect,
one should take into account that the boundary conditions are defined for the whole operator $f$ in (\ref{pde.1}), and they do not necessarily 
hold for the subproblems defined by each part $f_1$ and $f_2$. Therefore, one cannot
expect the numerical solution obtained by a splitting method to belong to the domain of $f$. This results in severe order reductions in reaction-diffusion
problems when Dirichlet or Neumann boundary conditions are considered \cite{hansen09hos,hundsdorfer03nso}. In particular, the order
reduction for the Strang splitting is one in the infinity norm. Similar order reductions for advection-reaction problems have also been reported
\cite{hundsdorfer95ano}.

Several procedures have been considered in the literature to avoid this order reduction in the case of reaction-diffusion problems. One possibility, proposed in
\cite{einkemmer15oor,einkemmer16oor} consists in introducing a smooth correction function in such a way that the new reaction flow is compatible with the prescribed boundary conditions. For time-invariant Dirichlet boundary conditions, this correction can be computed only once at the beginning of the simulation, but
for time-dependent Dirichlet, Neumann, or Robin boundary conditions, the correction is time-dependent and has to be computed at each time step. Various
techniques to deal with this problem are explored in \cite{einkemmer18ebc}. An alternative approach requiring additional calculations with
grid values on the boundaries, and not on grid values on the total domain, is proposed in \cite{alonsomallo17aor,alonsomallo19aor}.

\subsection{IMEX methods}

IMEX schemes are suitable combinations of implicit and explicit schemes and constitute a popular technique for approximating the solution of PDEs
that involve terms of different nature \cite{ascher95iem,ascher97ier,hundsdorfer03nso}. Thus, for convection-diffusion or reaction-diffusion problems 
where the convection or reaction terms are moderately stiff, it might be appropriate to use an explicit scheme for these parts and
 an implicit scheme for the diffusion term. We next analyze the connections of some popular IMEX methods with splitting and composition methods.

Suppose we have the semi-discrete system $u_t = f_1(u) + f_2(u)$, where $f_1$ is (still) a diffusion term and $f_2$ a nonlinear term suitable for explicit integration.
Then the simple composition (see (\ref{options1o}))
\[
  u_{n+1}= \phi_h^3 (u_n) = \psi_{h}^{1i}\circ \psi_{h}^{2e}(u_n)=u_n+{h} \, (f_1(u_{n+1})+f_2(u_n))
\]
corresponds to the linear one-step IMEX scheme of \cite{ascher97ier}.  

Letting $\phi_{h}^{[2,k]}$, $k \ge 2$, denote an explicit $k$th-order (Runge--Kutta or multistep) method for the equation $u_t = f_2(u)$, it turns out that many
IMEX methods from the literature have the structure
\begin{equation}  \label{IMEX2}
  u_{n+1}= \psi_{{h}/2}^{1i}\circ\phi_{h}^{[2,k]} \circ \psi_{{h}/2}^{1e}(u_n),
\end{equation}
thus yielding second-order approximations.

Let us consider, for instance, 
the popular Crank--Nicolson--leapfrog (IMEX-CNLF) method \cite{hundsdorfer03nso}
\begin{equation}  \label{IMEX-CNLF}
  u_{n+1}= u_{n-1}+2{h}f_2(u_{n})+{h}(f_1(u_{n+1})+f_1(u_{n-1})),
\end{equation}
to be initiated with, for example, $u_{1}= u_0+{h}\big( f_1(u_0)+f_2(u_{0}) \big)$.
It is equivalent to the one-step method
\begin{equation}  \label{IMEX-CNLFb}
  u_{n+1}= u_{n}+{h^*}f_2(u_{n+1/2})+\frac{h^*}{2}(f_1(u_{n+1})+f_1(u_{n})),
\end{equation}
with $h^*=2h$, which requires us to compute the approximation at the midpoint and so advances every half time step.
 If we take the explicit second order midpoint rule as $\phi_{h}^{[2,2]}$ in \eqref{IMEX2}, then we get the following sequence of maps:
\[
  \begin{array}{l}
	  \psi_{{h}/2}^{1e}: \quad U= u_n+\frac{h}2f_1(u_{n}) \\
		\phi_{h}^{[2,2]}: \ 
				\left\{   \begin{array}{l}
										V= U+\frac{h}2f_2(U) \\
										\hat V=  U+hf_2(V)
									\end{array}		  \right\} \\
		\psi_{{h}/2}^{1i}: \quad	u_{n+1}= \hat V+\frac{h}2f_1(u_{n+1})	
	\end{array}
\]
or, equivalently
\begin{equation}  \label{IMEX-CNLF2}
  \begin{aligned}
	& U= u_n+\frac{h}2f_1(u_{n}) \\
	& V= U+\frac{h}2f_2(U) \\
  & u_{n+1}= u_n+{h}f_2(V)+\frac{h}2(f_1(u_{n})+f_1(u_{n+1})).
 \end{aligned}
\end{equation}
This scheme is in fact quite similar to (\ref{IMEX-CNLFb}) for $h=h^*$ since $V$ is a first order approximation to $u_{n+1/2}$, the solution at the midpoint. Notice that by replacing $f_2(u_{n+1/2})$ in \eqref{IMEX-CNLFb} with  $f_2(V)$, where $V$ depends explicitly on $u_{n}$, it allows one 
to advance from $u_{n+1}$ to $u_{n+2}$ without evaluating the solutions at $u_{n+3/2}$ and therefore halving the computational cost to solve the implicit equations involved.  

Although the method is not symmetric (due to the lack of symmetry of the explicit scheme $\phi_{h}^{[2,2]}$), by instead using an explicit method
of order $k>2$, the overall scheme will be time-symmetric up to this order $k$ and therefore one can apply extrapolation to get a method of order $k$ in an
efficient way.

Another IMEX Runge-Kutta method that combines the implicit and explicit trapezoidal methods and shows a fairly good performance on examples is \cite{hundsdorfer03nso}:
  \begin{equation}
  \begin{aligned}
	& U  =  u_n+{h} \big( f_1(u_n)+f_2(u_{n}) \big)  \\
    & u_{n+1} =  u_n+\frac{h}2 \big(f_1(u_n)+f_2(u_{n}) \big)+\frac{h}2 \big(f_1(u_{n+1}) + f_2(U)\big).  \label{eq.IMEX2b}
 \end{aligned}
 \end{equation}

If we take the system 
\begin{equation} \label{eq.IMEX3}
 w'= \left( \begin{array}{c}
          u \\
          U 
         \end{array} \right)'  =  
   \left( \begin{array}{c}
          f_1(u)+f_2(U)\\
          0 
         \end{array} \right)+ 
            \left( \begin{array}{c}
          0  \\
          f_1(u)+f_2(u)\ 
         \end{array} \right)
				= g_1(w) + g_2(w),
\end{equation}
and form, as before, the composition
\begin{equation}  \label{IMEX3b}
  w_{n+1}= \psi_{{h}/2}^{1i}\circ\phi_{h}^{[2]} \circ \psi_{{h}/2}^{1e}(w_n),
\end{equation}
where $\phi_{h}^{[2]}$ denotes the exact solution of $w' = g_2(w)$, and $w_n=(U_n,u_n)=(u_n,u_n)$, then we get
  \begin{eqnarray}
	u_{1/2} & =&  u_n+\frac{h}2 \big(f_1(u_n)+f_2(u_{n}) \big) \nonumber \\
	U_{n+1} & =&  u_n+{h} \big(f_1(u_{1/2})+f_2(u_{1/2}) \big) \label{eq.IMEX3c} \\
  u_{n+1}& =&  u_n+\frac{h}2 \big(f_1(u_n)+f_2(u_{n}) \big)+\frac{h}2 \big(f_1(u_{n+1})+f_2(U_{n+1})\big). \nonumber 
 \end{eqnarray}
Scheme \eqref{eq.IMEX2b} is recovered by replacing $u_{1/2}$ with $u_n$ in the expression of $U_{n+1}$. 
However, if in \eqref{eq.IMEX3c}
we do not restart the value of $U_n$ at each step, the method will be symmetric and therefore one can apply extrapolation to increase its order.

Higher order IMEX methods have been built involving implicit multistep or Runge--Kutta methods. Although splitting methods are no longer appropriate to
advance the diffusion term due to the presence of negative coefficients for orders higher than two, one can incorporate higher derivatives or complex coefficients
and form new higher order splitting IMEX methods (or consider extrapolation from a basic symmetric second order method).

\subsection{Schr\"odinger equations}

\paragraph{$T$--$V$ splitting.}

Until now, in the context of PDEs, we have dealt with splitting methods of orders one and two. 
There are, however, relevant problems where high order splitting
methods can be and have been safely used, and where a rigorous convergence analysis can be established. This is the case, in particular, of the 
time-dependent Schr\"odinger equation, already considered in Section \ref{subsec.1.3b}. 
The numerical experiments presented there clearly indicate that the Strang splitting based on kinetic and potential energy
 in combination with a pseudo-spectral space discretization,  i.e., method (\ref{ssTV}), provides approximations of order two in the time step, in accordance
 with
 \begin{equation} \label{schro1}
   \e^{\frac{\tau}{2} V} \, \e^{\tau T} \, \e^{\frac{\tau}{2} V} = \e^{\tau (T+V)} + \mathcal{O}(\tau^3 (\|T\| + \|V\|)^3)
\end{equation}
where $\tau=-i h$. One should be aware, however, that this error estimate only makes sense for bounded $T$ and $V$. In fact, when the norm of
$T$ or $V$ is very large (as is usually the case when the number of space discretization points $M$ is large), then (\ref{schro1}) is of no practical use, and
thus other estimates are necessary to explain the observed good behavior.    

Error bounds for the Strang splitting are actually derived in \cite{jahnke00ebf}. They clearly show that the method is indeed of order two
when applied to pseudo-spectral discretizations  of the time-dependent Schr\"odinger equation under some regularity conditions and periodic boundary conditions. 
Specifically, assume that the potential $V(x)$ is $\mathcal{C}^5$-smooth, periodic and
bounded, $\|V \psi\| \le \beta  \|\psi\|$, $\beta > 0$. Then,  if $u(x,t)$ denotes the trigonometric interpolation
polynomial of solution of the pseudo-spectral method and $u_n(x)$ is the corresponding trigonometric interpolation polynomial built from the numerical approximations
obtained from the Strang splitting (\ref{ssTV}) at time $t=t_n = n h$, for the local and global errors one gets the following bounds:
\begin{equation} \label{boS}
 \begin{aligned}
  &  \|u_1 - u(\cdot, \tau)\|_{L^2} \le C_1 |\tau|^3 \|u_0\|_{H^2} \\
  &  \|u_n - u(\cdot, t_n)\|_{L^2} \le C_2 |\tau|^2 \|u_0\|_{H^2},
 \end{aligned} 
\end{equation}
respectively. Here $\| \cdot\|_{H^2}$ denotes the usual Sobolev norm and the constants $C_1$, $C_2$ are independent of the initial data $u_0$ and
the discretization parameters $M$, $n$ and $\tau$, with $0 \le t_n \le t_f$ for some finite $t_f$. The case when $V$ is time-dependent and bounded for any $t$
has been recently treated in \cite{valle23fos}, where new schemes are proposed and analyzed.

The previous results for the Strang splitting have been extended in \cite{thalhammer08hoe} to splitting methods of the general form
\begin{equation} \label{mts1}
  u_{n+1} =  \prod_{j=1}^s \e^{b_j \tau V} \, \e^{a_j \tau T} \, u_n = \e^{b_s \tau V} \, \e^{a_s \tau T} \, \cdots \,    \e^{b_1 \tau V} \, \e^{a_1 \tau T} \, u_n
\end{equation}
whose coefficients $a_j$, $b_j$ satisfy the order conditions up to order $r$. In that case
\begin{equation} \label{boMT}
  \|u_n - u(\cdot, t_n)\|_{L^2} \le C \|u(\cdot, 0) - u_0\|_{L^2} + C |\tau|^r \|u(\cdot,0)\|_{H^{r}}, \qquad 0 \le t_n \le t_f
\end{equation}
is valid with some constant $C$  depending on $t_f$, but not on $n$ and $h$. This error bound implies, in particular, that the splitting methods 
of section \ref{sect8} retain their order of convergence  when applied to the Schr\"odinger equation with periodic boundary conditions, provided that 
the  data are sufficiently differentiable (see also \cite{hansen09esf}). Otherwise, an order reduction may occur.  

One can also take advantage of the property (\ref{rknS}) and include the commutator $[V,[T,V]]$ in the composition (\ref{mts1}), as in RKN splitting methods,
so that one ends up with 
 \begin{equation} \label{mts1mod}
  u_{n+1} =  \prod_{j=1}^s \e^{b_j \tau V + c_j \tau^3 [V,[T,V]]} \, \e^{a_j \tau T} \, u_n. 
\end{equation}
Even in this case the resulting schemes retain their order of convergence if the solution is sufficiently regular, as shown in \cite{kieri15sco}. It is then  also possible
to apply the RKN schemes presented in section \ref{sect8}.

\paragraph{Symplectic splitting.}

In section \ref{sub.lisys} we reviewed the symplectic structure involved in the (semi-discretized) Schr\"odinger equation and illustrated how the Strang splitting
(and, for that matter, all methods presented in Section \ref{sect8}) can be applied if $H$ is a real and symmetric matrix. They are formulated as products of
exponentials of the nilpotent matrices $A$ and $B$ given in (\ref{eq.4a}), 
\begin{equation} \label{mts2}
  \left( \begin{array}{c}
          q_{n+1} \\
          p_{n+1} 
         \end{array} \right)  =  \prod_{j=1}^s \e^{b_j h B} \, \e^{a_j h A} \, 
   \left( \begin{array}{c}
          q_{n} \\
          p_{n} 
         \end{array} \right)         
         = \exp \left[ h \left(\begin{array}{cc}
         		0 & H \\
		       -H & 0
		      \end{array} \right) \right] 
            \left( \begin{array}{c}
          q_{n} \\
          p_{n} 
         \end{array} \right)         
          + \mathcal{O}(h^{r+1})
\end{equation}
and orders $r = 2, 4, 6, 8, 10$ and $12$ have been achieved with only $s=r$ exponentials $\e^{a_j h A}$ and $\e^{b_j h B}$ \cite{gray96sit,zhu96nmw,liu05oos}. 

The processing technique has also been used to construct splitting schemes
with two different goals in mind: to attain maximal stability and maximal accuracy. They have the general structure 
$P^{-1}(h H) K(h H)P(h H)$, 
where $K$ (the kernel) is built as a composition (\ref{mts2}) with a large number of stages $s$ and $P$ (the processor) is
taken as a polynomial.
Although these methods are neither unitary nor unconditionally stable, they are symplectic and conjugate 
to unitary schemes. In consequence, neither the average error in energy nor the norm of the solution increases with time. Specifically, in  
\cite{blanes06sso,blanes08otl} kernels with up to 
$19$, $32$ and $38$ stages have been proposed, either to construct methods of orders $r=10$, $16$ and $20$, or to bring highly
accurate \emph{second} order methods with an enlarged stability domain. 

This approach to approximating $\e^{\tau H} u$  is closely related to other polynomial approximations of the form
\begin{equation} \label{poly1}
  \e^{\tau H} u \approx P_m(h H) u,
\end{equation}
where $P_m(y)$ is a polynomial in $y$ approximating the exponential $\e^{-i y}$. Different choices for such $P_m(y)$ are available: truncated Taylor
or Chebyshev series expansions of $\e^{-i y}$ for an appropriate real interval of $y$, or a Lanczos approximation, where the polynomial is determined
by a Galerkin approximation on the Krylov space  spanned by $u, H u, \ldots, H^{m-1} v$ \cite{lubich08fqt}.

Given a prescribed error tolerance, some appropriate estimates of the upper and lower bounds of the eigenvalues of the matrix $H$, $E_{\min}$ and $E_{\max}$, and a time integration interval, $[t_0,t_f]$, in the Chebyshev approach one chooses, according to some known error bounds, the lowest-degree polynomial that provides the solution with such accuracy. The coefficients of the polynomial are determined for each case and the action of the polynomial on a vector is computed recursively using the Clenshaw algorithm \cite{lubich08fqt}. On the other hand, in the Taylor approach one has to adjust the maximum degree allowed with the time step $h$ to reach the desired accuracy with the minimum number of matrix-vector products. As a result, in general, Chebyshev turns out to be between two and three times faster than Taylor, depending on the final time at which the output is desired.  Chebyshev and Lanczos approximations have quite similar error bounds (see \cite{lubich08fqt}) and their relative performance depends on the particular problem considered.

It is worth remarking that, whereas in the approach 
(\ref{poly1}) the approximation of $\e^{\tau H} u$ is constructed by computing products of the form $H u$, where $u \in \mathbb{C}^M$, with symplectic
splitting methods of the form (\ref{mts2}) one proceeds by 
%
%
successively computing 
\emph{real} matrix-vector products $H q$ and $H p$ with different weights. With splitting methods, the real and imaginary parts of $\e^{\tau H} u\equiv \e^{-ih H} (q+ip)$    
are approximated in a different way, with a considerably reduced computational cost. 

In \cite{blanes15aea}, several optimized
symplectic splitting methods are constructed and an algorithm is presented that automatically selects the most efficient one for a prescribed error
tolerance under the same conditions as when using the Chebyshev method.
The resulting algorithm is between 1.4 and 2 times faster than the Chebyshev method for the same accuracy, with reduced energy and unitarity errors for large values of $h$. The computation of the coefficients of the schemes
is largely based on the stability and error analysis of splitting methods carried out in \cite{blanes08otl,blanes11eao}.

In contrast with $T$--$V$ splitting methods, which preserve unitarity by construction and are thus unconditionally stable, the previous
polynomial approximations suffer from a step size restriction. Given $h$, $\Delta x$ and $m$, the degree of the polynomial, these methods must satisfy
the restriction 
\[
  \displaystyle \frac{h}{m\Delta x^2} \leq C,
\] 
and so the time interval that one can advance per matrix-vector product is proportional to $\Delta x^2$ or, equivalently, the number of matrix-vector products to reach the final time  is inversely proportional to $\Delta x^2$.

\paragraph{Time-dependent potentials.}

We have assumed so far that the potential in the Schr\"odinger equation (\ref{Schr0}) does not depend explicitly on time or, if it does, only varies slowly with time, so that in each sub-interval
$[t_n, t_{n+1}]$ the corresponding matrix $V$ is obtained from the average of $V(x,t)$ on this interval. In general, however, one has to deal
with situations in which this approximation is no longer valid. In that case, and in contrast with other approaches based on Chebyshev or Lanczos approximations,
splitting methods can still be applied with some appropriate modifications. 

In a similar way to classical Hamiltonian problems, one may take $t$ in the potential as an additional coordinate, $x_{d+1}=t$, introduce its canonical momentum, $p_{d+1}=-i \, \partial_{x_{d+1}}$, and deal with the system in the extended phase space, $\tilde H=(T+p_{d+1})+V(x,x_{d+1})$. Since the action of the operator $p_{d+1}$ just corresponds to a shift in the variable $x_{d+1}$,
i.e., $\e^{ah\partial_t}V(t,x)=V(t+ah,x)$, and moreover
\[
\e^{ah\partial_t}\e^{-ibh V(t,x)}\e^{-ah\partial_t}=\e^{-ibh \e^{ah\partial_t}V(t,x)}=\e^{-ibh V(t+ah,x)},
\]
composition \eqref{mts1} applied to the corresponding non-autonomous problem now simply reads (notice that $T$ and $p_{d+1}$ commute)
\begin{equation} \label{mts1t}
  u_{n+1} =  \prod_{j=1}^s \e^{b_j \tau V(t_n+c_j h)} \, \e^{a_j \tau T} \, u_n , 
	\qquad \mbox{ with } \;\; c_j = \sum_{k=1}^j a_k, \qquad j \in \{1,\ldots,s\},
\end{equation}
and analogously if the scheme includes modified potentials (see also e.g. \cite{chin02gsa} and references therein).

For a generic time-dependent Hamiltonian $H(t)$, the 
 Schr\"odinger equation $i \partial_t \psi = H(t) \psi$ can be recast as a non-autonomous evolution equation of the form
\begin{equation} \label{tdp1}
  u'(t) = A(t) u(t), \qquad u(t_0) = u_0, \qquad t \in [t_0,t_f],
\end{equation}
defined by a family of time-dependent linear operators $(A(t))_{t \in [t_0,t_f]}$, which, assuming a spatial discretization has been carried out, are generally complex
matrices of large dimension and large norm.  

It turns out that standard $r$th-order splitting methods defined by coefficients $(a_{j}, b_{j})_{j=1}^s$ can be applied in this setting simply by adding the
trivial relation $\frac{d}{dt} t = 1$ 
to equation (\ref{tdp1}). This results in the scheme
\begin{equation} \label{tdp2}
  u_{n+1} = \prod_{j=1}^s \e^{h b_j A(t_n + c_j h)} u_n, \qquad \mbox{ with } \;\; c_j = \sum_{k=1}^j a_k, \qquad j \in \{1,\ldots,s\}
\end{equation}  
of the same formal order $r$ as the method originally designed for autonomous problems. 

A different approach is based on the use of the Magnus expansion \cite{magnus54ote} to get a formal solution representation of 
(\ref{tdp1}) as the exponential of an infinite series:
\begin{equation} \label{tdp2b}
  u(t_n + h) = \e^{\Omega(h)} u_n, \qquad\quad  \Omega(h) = \sum_{m=1}^{\infty} \Omega_m(h),
\end{equation}
where each term $\Omega_m$ involves multiple integrals of nested matrix-commutators \cite{blanes09tme}. 
By appropriately truncating this series and approximating the integrals by
quadratures, efficient integrator schemes can be constructed \cite{iserles00lgm}.  Thus, for instance, taking the two-stage Gauss--Legendre quadrature rule with
nodes $c_{1,2} = \frac{1}{2} \mp \frac{\sqrt{3}}{6}$, results in the scheme
\begin{equation} \label{tdp3}
  \begin{aligned}
   & \Omega^{[4]}(h) = \frac{1}{2} h (A_1 + A_2) + \frac{\sqrt{3}}{12} h^2 [A_2, A_1], \\
   & u_{n+1} = \e^{\Omega^{[4]}(h)} u_n,
  \end{aligned}
\end{equation}
with $A_j = A(t_n + c_j h)$, $j=1,2$. For the Schr\"odinger equation with $A(t) = -i (T + V(t))$ and a smooth time-dependent potential $V(t)$, it is shown in
\cite{hochbruck03omi} that  Magnus integrators retain their full order of convergence (without bounds on $T$ in the error bound) for sufficiently regular solutions,
uniformly with respect to the space discretization.   

There are, however, several issues related to Magnus integrators due to the presence of iterated
commutators. Thus, computing the action of iterated commutators on vectors can be 
very costly due to the number of matrix-vector products required. This is particularly relevant when considering problems in two and three space
dimensions \cite{bader16emf}. In addition, the evolution equations defining high-order Magnus integrators in general involve differential
operators of different nature from the original problem \cite{blanes21ao}.

A different class of exponential integrators that circumvent these difficulties whilst still retaining the favorable properties of Magnus integrators is formed
by the so-called commutator-free quasi-Magnus (CFQM) methods: the basic idea is to replace the single exponential in (\ref{tdp2b}) with a composition
of several exponentials involving linear combinations of the values of the operator $A$ at certain nodes, $c_k$, of a quadrature rule:
\begin{equation*}
\begin{gathered}
u_{n+1} = \e^{h \, B_{nJ}} \cdots \, \e^{h \, B_{n1}} \, u_n \; \approx \; u(t_{n+1}) = \e^{\Omega(h)} \, u(t_n)\,, \\
c_k \in [0,1]\,, \qquad A_{nk} = A(t_n + c_k h)\,, \qquad k \in \{1, \dots, K\}\,, \\
B_{nj} = a_{j1} \, A_{n1} + \cdots + a_{jK} \, A_{nK}\,, \qquad j \in \{1, \dots, J\}\,.
\end{gathered}
\end{equation*}
Particular examples of CFQM exponential integrators are the exponential midpoint rule (order 2)
\begin{equation}
\label{eq:CFQMOrder2JK1}
\begin{gathered}
J = K = 1 \,, \quad c_1 = \tfrac{1}{2}\,, \quad a_{11} = b_1 = 1\,, \\
u_{n+1} = \e^{h \, A(t_n + \frac{1}{2} h)}\, u_n
\end{gathered}
\end{equation}
and the fourth-order scheme
\begin{equation}
\label{eq:CFQMOrder4JK2}
\begin{gathered}
 J = K = 2\,, \quad \alpha = \tfrac{\sqrt{3}}{6}\,, \quad c_1 = \tfrac{1}{2} - \alpha\,, \quad c_2 = \tfrac{1}{2} + \alpha\,, \\
a_{11} = a_{22} = \tfrac{1}{4} + \alpha\,, \quad a_{12} = a_{21} = \tfrac{1}{4} - \alpha\,,
 \\
B_j(h) = a_{j1} \, A(t_n + c_1 h) + a_{j2} \, A(t_n + c_2 h)\,, \quad j \in \{1, 2\}\,, \\
u_{n+1} = \e^{h \, B_2(h)} \, \e^{h \, B_1(h)}\, u_n.
\end{gathered}
\end{equation}
A detailed treatment and specific schemes up to order six can be found in \cite{blanes17hoc,blanes18cao} and up to order eight in \cite{alvermann11hoc}, 
whereas it is proved in \cite{blanes21ao}
that CFQM methods applied to
the Schr\"odinger equation with Hamiltonian $H(t) = -\frac{1}{2} \Delta + V(t)$ are
unconditionally stable in the underlying Hilbert space and retain full order of convergence under low regularity requirements on the initial state.

\paragraph{Semiclassical regime.}

The so-called semi-classical Schr\"odinger equation 
\begin{equation} \label{scse1}
  i \varepsilon \partial_t \psi(x,t) = \left( -\frac{\varepsilon^2}{2} \Delta + V(x) \right) \psi(x,t)
\end{equation}
(in atomic units) with a small parameter $\varepsilon \ll 1$, arises in particular when applying the time-dependent Born--Oppenheimer approximation for the motion
of nuclei as driven by the potential energy surface of the electrons \cite{lubich08fqt}. In that case $\varepsilon^2$ represents the mass ratio of
nuclei and electrons. Recall that
eq. (\ref{scse1}) has highly oscillatory solutions with wavelength $\sim \varepsilon$, so grid-based numerical schemes require a resolution of this
order in both space and time, which is computationally very expensive. One of the challenges, therefore, is 
to construct numerical methods that are robust in the limit $\varepsilon \rightarrow 0$.

Several options have been
proposed and analyzed in detail (see, for instance, the recent review \cite{lasser20cqd}). Among others, we can recount the following: 
\begin{itemize} 
 \item Split the equation into the usual kinetic and potential energy parts and apply
the Strang splitting in time in combination with trigonometric spectral methods \cite{bao02ots}. Although the resulting scheme is unconditionally stable, time-reversible, and
preserves the position density, it requires very fine resolution, in both space and time, for small $\varepsilon$ \cite{jin11mac}.
\item Use a Gaussian wave packet as an approximation for the wave function $\psi(x,t)$ depending on certain parameters and apply a variational
splitting to get approximate solutions for the differential equations they satisfy \cite{faou06api}. 
The resulting algorithm is symplectic, time-reversible and preserves the unit $L^2$-norm of the wave packets.
 \item 
Another variant of this approach instead consists in taking Hagedorn wave packets. They provide a spectral approximation in space with a time-dependent set of basis functions giving the exact solution of
the Schr\"odinger equation with the potential locally approximated by a quadratic function.
The potential $V(x)$ is split into the quadratic term $U(q(t),x)$ in the Taylor expansion of $V$,  around the time-dependent classical position $q(t)$ and the remainder
 \cite{faou09csq}. The overall algorithm has a number of conservation and limit properties, as listed in \cite{faou06api}, for instance. In addition, highly efficient splitting methods for perturbed problems can also be applied \cite{blanes20hoe}.
 \item A different class of exponential splittings is proposed in \cite{bader14eaf} for the one-dimensional case. Essentially, the formal solution of the space-discretized
 equation $u' = i (\varepsilon T - \varepsilon^{-1} V) u$ is approximated as
 \[
   \e^{i h (\varepsilon T - \varepsilon^{-1} V)} \approx \e^{R_0} \e^{R_1} \cdots \e^{R_s} \e^{T_{s+1}} \e^{R_s} \cdots \e^{R_1} \e^{R_0}
 \]
 with error $\mathcal{O}(\varepsilon^{2s+2})$. Here $R_0 = \mathcal{O}(\varepsilon^0)$, $R_k = \mathcal{O}(\varepsilon^{2k-2})$, $k \ge 1$, and
 $T_{s+1} = \mathcal{O}(\varepsilon^{2s})$. In this approach the number of exponentials grows linearly with $s$ and the exponentials can be computed efficiently, although
 the terms $R_k$ and $T_{s+1}$ contain nested commutators. 
   
\end{itemize} 

\paragraph{Nonlinear Schr\"odinger equations.}

Introducing nonlinear effects in the Schr\"odinger equation allows one to model some relevant physical phenomena taking place
in nonlinear optics, quantum superfluids,  plasmas, water waves, etc. (see e.g. \cite{sulem99tns} and references therein). Consider in particular a
Bose--Einstein condensate (BEC), the ground state of a system of interacting bosons very close to
zero temperature. It was first predicted by Einstein in 1925 and experimentally realized in 1995 \cite{anderson95oob}. Mathematically, a BEC of 
an atomic species trapped in an external  potential
$V(x)$ is modeled by the (normalized) Gross-Pitaevskii equation (GPE)
\begin{equation} \label{g-p1}
     i  \partial_t \psi (x,t) = \left(%
                                                    -\frac{1}{2} \Delta + V(x) +   \sigma |\psi (x,t)|^2%
                                                \right) \psi(x,t),
\end{equation}
with asymptotic boundary conditions $\psi(x,t) \rightarrow 0$ as $|x| \rightarrow \infty$.
Here the parameter $\sigma$ originates from the mean-field
interaction between the particles: repulsive forces lead to  $\sigma > 0$, whereas $\sigma < 0$ represents attractive forces. Equation
(\ref{g-p1}) has been the subject of many different studies, including the existence of solutions and its numerical treatment. Concerning the first aspect,
we refer to \cite{cazenave03sse,carles08sca} and references therein. 

With respect to the numerical integration of the GPE equation, a combination of spectral discretization in space with splitting methods in time constitutes
a natural option, and in fact has been explored in detail in the literature (see e.g. \cite{bao03nso,bao03nso2,bao04tds,thalhammer09hot}). If
 (\ref{g-p1}) is expressed as 
\begin{equation}\label{SepAB}
    i \partial_t \psi(x,t) = (A + B(x,\psi)) \psi(x,t), \qquad \psi(x,0) = \psi_0(x),
\end{equation}
with
\begin{equation} \label{FourierSplit}
    A= -\frac{1}{2} \Delta ,   \qquad\quad   B(x,\psi) = V(x) + \sigma |\psi |^2,
\end{equation}
it is clear that the solution of the initial value problem $i \psi_t = A \psi$, $\psi(x,0) = \psi_0(x)$ is
\[
  \psi(x,t) = \e^{-i t A} \, \psi_0(x)
\]
and (an approximation to) $\psi(x,t)$ is obtained by representing the initial value with respect to the (truncated) Fourier basis functions.

On the other hand, given a real function $G$, the solution of 
\begin{equation}  \label{NLS-pot}
    i \partial_t \psi (x,t) = G \big(x,|\psi (x,t)| \big) \, \psi (x,t), 
\end{equation}
leaves the norm invariant, i.e., $|\psi(x,t)|=|\psi(x,0)|$, and therefore
\begin{equation} \label{solNLS-pot}
  \psi(x,t) = \e^{-itG(x,|\psi(x,0)|)} \, \psi(x,0).
\end{equation}
In consequence, the initial value problem $i \partial_t \psi = B \psi$, $\psi(x,0) = \psi_0(x)$ is also solvable. Splitting methods can also be applied in the more general situation
when the potential $V$ is explicitly time-dependent, as explained earlier for the linear Schr\"odinger equation.

\subsection{Parabolic evolution equations}

Let us now consider the evolution equation
\begin{equation} \label{pee1}
  u'(t) = L u(t) = A u(t) + B u(t), \qquad t \ge 0, \qquad  u(0) = u_0,
\end{equation}
where the linear, possibly unbounded, operators $A$, $B$ and $L$ generate $C_0$ semigroups over an infinite dimensional Banach space $X$. We recall
that if $L$ is the infinitesimal generator of the $C_0$ semigroup $T(t)$ on $X$ and $u_0 \in \mathcal{D}(L)$, the domain of $L$ (which is dense in $X$), 
then $u(t) = T(t) u_0$ is a classical solution of (\ref{pee1}) \cite{partington04loa}. Since in the special case of a bounded linear operator $L$ the solution is
given by the familiar expression $u(t) = \e^{t L} u_0$ \cite{engel06asc}, the semigroup $T(t)$ is also denoted by the symbol $\e^{t L}$ (see 
\cite{engel06asc,pazy83sol,yosida71fan} for an introduction to the theory of $C_0$ semigroups). 

A prototypical example is the linear heat equation with potential
\begin{equation} \label{pee2}
   u_t(t,x) = \frac{1}{2}\Delta u(t,x) - V(x) u(t,x).
\end{equation}
Here $V(x)\geq 0$, 
$t \ge 0$ and $x \in \mathbb{R}^d$ (or $x \in \mathbb{T}^d$). In that case $(A u)(x) =  \frac{1}{2}\Delta u(x)$, $(B u)(x) = -V(x) u(x)$ and $A$ generates only
a $C_0$ semigroup. This can be seen by considering the equation $u_t = \frac{1}{2}\Delta u$ on $0 < x < 1$ with Dirichlet boundary conditions: the $k$th Fourier
mode of the solution is $c_k \e^{-\frac{1}{2}(k \pi)^2 t}$, which is in general not well defined for $t < 0$.

In \cite{hansen09esf}, the following result has been
established. Assuming that $\|\e^{t A}\| \le \e^{\omega t}$, $\|\e^{t B}\| \le \e^{\omega t}$ for the same value of $\omega \ge 0$ and all $t \ge 0$, and that
for any operator $E_{r+1}$ obtained as the product of exactly $r+1$ factors chosen amongst $A$ and $B$, there is a constant $\tilde{C} > 0$ such that
\[
  \|E_{r+1} \e^{t(A+B)} u_0 \| \le \tilde{C}.
\]
Then a splitting method $\Psi(h) = \prod_{j=1}^s \e^{b_j h B} \, \e^{a_j h A}$, with all $a_j \ge 0$ and $b_j \ge 0$ and
of classical order $r$ retains its order when applied to (\ref{pee1}):
\begin{equation} \label{ppe3}
  \| (\Psi(h)^n - \e^{n h L}) u_0  \| \le C h^r, \qquad \mbox{ for } \qquad n \, h \le t_f,
\end{equation}
where the constant $C$ is independent of $n$ and $h$ on the bounded time interval $[0, t_f]$. 

In practice, however, the positivity requirement on the coefficients restricts the splitting method to be at most order two. If, in addition
$[B,[A,B]]$ is 
a bounded operator
and $\|\e^{t [B,[A,B]]}\| \le \e^{\omega t}$, 
then the same result (\ref{ppe3}) also holds
for splitting methods involving double commutators, 
\[
 \Psi(h) = \prod_{j=1}^s \e^{b_j h B + c_j h^3 [B,[A,B]]}\, \e^{a_j h A}
\]
  and
positive $a_j$, $b_j$ \cite{kieri15sco}. Efficient schemes within this class up to order 4 specially tailored to the problem at hand have recently been proposed in
\cite{blanes23esm}.

\paragraph{Example: imaginary time propagation.}  
 
The imaginary-time evolution method is a well-known approach to computing the ground state (and its corresponding eigenvalue) of a quantum system
with Hamiltonian $H =  -\frac{1}{2} \Delta +  V$ \cite{auer01afo,lehtovaara07sot,bader13sts}. 
Essentially, under the time transformation $t = -i s$, the 
time-dependent Schr\"odinger equation (\ref{Schr0}) is transformed into ($\hbar = 1$)
\begin{equation} \label{itp1}
 \partial_s \psi(x,s) = \frac{1}{2} \Delta \psi(x,s) - V(x) \psi(x,s), \qquad \psi(x,0) = \psi_0(x),
\end{equation}
i.e., a linear heat equation of the form (\ref{pee2}). If we denote
the (real) eigenvalues of $H$ as $E_j$, with $E_0 < E_1 < \cdots$, and the corresponding eigenfunctions as
$\phi_j$, $j=0,1,2, \ldots$, the initial wave function $\psi_0(x)$ can be 
expanded in the orthonormal basis $\{\phi_j\}$,
\[
  \psi_0(x) = \sum_{j \ge 0} c_j\;\phi_j(x), \qquad
  c_j=\left\langle\phi_j\,|\,\psi(\cdot,0)\right\rangle,
\]
where $\langle\cdot\,|\,\cdot\rangle$ is the usual $L^2$-scalar product. Then the solution of (\ref{itp1}) can be written as
 \begin{equation} \label{eq:evolution_imagtime}
    \psi(x,s) = \e^{-s H} \psi(x,0) = \sum_{j \ge 0} \e^{-s E_j}\,c_j\;\phi_j(x).
\end{equation}
Notice that, for sufficiently large $s$ one gets $\psi(x,s) \rightarrow \e^{-s E_0} c_0 \phi_0$, since
the other exponentials decay more rapidly. In other words, any given wave function at $s = 0$ for which $c_0\neq 0$, converges towards the ground state solution 
when $s \rightarrow \infty$. Once an accurate approximation to $\phi_0$ is obtained, the associated eigenvalue $E_0$ is easily obtained by computing $E_0=\langle\phi_0\,|\,H \phi_0\rangle$. Other functions $\phi_j$ can also be approximated, for example by propagating different wave functions
simultaneously in time \cite{aichinger05afc}. 

To illustrate the technique we next consider the same example as in Subsection \ref{subsec.1.3b}, but now integrating in imaginary time up to $s_f=1$. We take the same initial conditions (which is slightly closer to one of the minima of the potential) and normalize the solution at the final time. Figure~\ref{fig_Schrod6IT} shows in the left panel the potential, initial conditions and the normalized wave function at the final time when the spatial interval, $x\in[-13,13]$ is divided into $M=256$ and $M=512$ parts, leading to the same visual results. In the right panel we show the $L_2$-norm error in the normalized wave function at the final time versus the number of FFTs required by the same methods as previously, both for  $M=256$ (solid lines) and $M=512$ (dashed lines). We notice that one scheme suffers from step size restriction because it involves negative coefficients, and this is inversely proportional to $\Delta x^2$. However, the schemes with positive coefficients are insensitive to the spatial mesh.

\begin{figure}[!htb]
\includegraphics[scale=0.6]{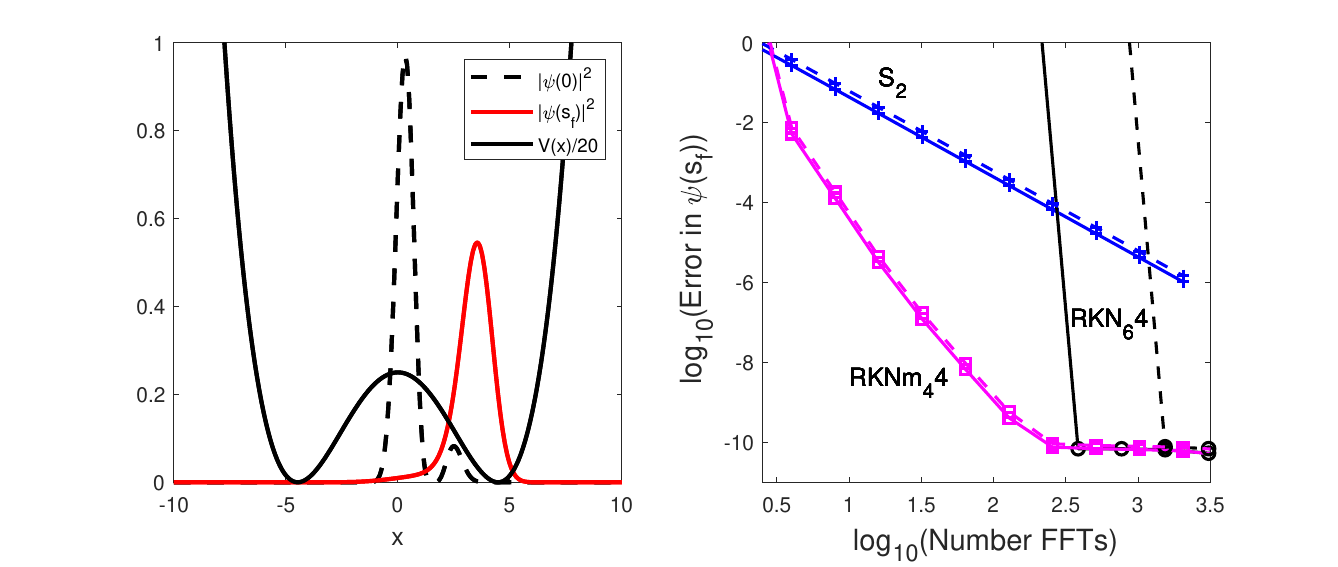}
\caption{Imaginary time propagation of the Schr\"odinger equation. Left: Double-well potential, initial and (normalized) final wave function obtained with a mesh with  $N=256$ and repeated with  $N=510$ (the two results overlap visually). Right: discrete $L_2$-norm error in the normalized solution at the final time, $s_f=1$,  vs. number of FFTs for different values of the time step obtained with the same methods as in Figure~\ref{fig_Schrod1}: solid lines for $M=256$ and dashed lines for  $M=512$
discretization points of the space interval.
}
\label{fig_Schrod6IT}
\end{figure}
\noindent
$\Box$
 
 \

One possible way to circumvent this order barrier for splitting methods when applied to problem (\ref{pee1}) consists in considering schemes
with \emph{complex} coefficients $a_i$, $b_i$ having \emph{positive} real part. 
As shown in 
\cite{castella09smw,hansen09hos}, any splitting method 
\[
 \Psi(h) = \prod_{j=1}^s \e^{b_j h B} \, \e^{a_j h A}
\] 
 within this class still retains its classical order if the previous assumptions 
 are conveniently modified. Specifically, one has to extend 
the notion of a $C_0$ semigroup $T(t)$  to the sector $\Sigma_{\theta}$ in the complex plane, 
for some angle $0 < \theta < \pi/2$, 
\[
   \Sigma_{\theta} = \{ t \in \mathbb{C} : |\arg(t)| < \theta \},
\]
so that $T(t)$ is analytic in $t$ for all $t \in   \Sigma_{\theta}$ \cite{engel06asc,pazy83sol}. Thus, if 
\begin{itemize}
  \item  $L$, $A$ and $B$ generate analytic semigroups on $X$ (now a complex Banach space) in the sector $\Sigma_{\theta}$, $0 < \theta < \pi/2$, with
\[
   \|\e^{t A}\| \le \e^{\omega |t|}, \qquad \|\e^{t B}\| \le \e^{\omega |t|}
\]
for $\omega \ge 0$ and all $t \in  \Sigma_{\theta}$;
 \item  $ \|E_{r+1} \e^{t(A+B)} u_0 \| \le C$, 
 where $E_{r+1}$ is a composition of the operators $A$ and $B$ that consist of exactly $r+1$ factors,
\end{itemize} 
then the splitting method of classical order $r$
with coefficients $a_i, b_i \in \Sigma_{\theta} \subset \mathbb{C}$ satisfies
\begin{equation} \label{ppe4}
    \| (\Psi(h)^n - \e^{n h L}) u_0 \| \le C h^r, \qquad 0 \le  n h < t_f.
\end{equation}     
Since $\Psi(h) u_0$ is complex-valued, this approach cannot be applied in principle when the operators $A$ and $B$ are real, i.e., for problems
defined in a real Banach space $X$. The most straightforward remedy consists in projecting the numerical solution after each time step on the
real axis, i.e., computing the approximations $u_n \approx u(t_n)$, as $u_n = \mathrm{Re}(\Psi(h)u_{n-1})$. In that case one still has the bound (\ref{ppe4})
for the resulting integration scheme \cite{hansen09hos}.

Although this is only valid in the linear case, similar results are observed in practice in the nonlinear heat equation 
$\partial_t u = \Delta u + F(u)$ with periodic boundary conditions, at least when $F$ is analytic \cite{castella09smw,blanes13oho}, and so they constitute
a strong motivation to study splitting methods with complex coefficients in general. This is precisely
the subject of the following section.



\section{Splitting methods with complex coefficients}
\label{sect7}

It has been known for a long time that, besides real solutions, the order conditions arising from splitting and composition methods (see section \ref{sect2})
also admit complex solutions \cite{bandrauk91ies,suzuki90fdo,suzuki91gto,suzuki95hep}. In fact, some of the resulting 
methods were explored in the context of Hamiltonian \cite{chambers03siw} and quantum mechanics  \cite{bandrauk06cis,prosen06hon} before being used 
more recently to
overcome the order barrier in parabolic partial differential equations \cite{castella09smw,hansen09hos}.

The use of the complex plane to solve problems formulated 
in the real line has proved to be very fruitful in many branches of mathematics, as illustrated by 
Painlev\'e's famous dictum\footnote{\textit{Il apparut que, entre deux v\' erit\'es du domaine r\'eel, le chemin le plus facile et le plus court passe bien souvent par le domaine
complexe} \cite{painleve00adt}.}. By applying the same logic in the particular case of splitting methods, 
it is reasonable to ask what benefits (if any) might result from carrying out 
the integration along paths in the complex plane. 

We have already seen at the end of Section \ref{sect6} that a strategy to develop effective methods of order higher than two for systems evolving in a semigroup,
such as the heat equation\footnote{This in fact constitutes Problem 10 in the list of open problems posed in \cite{mclachlan02sm}.}, 
consists precisely in applying splitting schemes with complex coefficients having positive real part. Moreover, the large number of complex solutions for
the order conditions might offer more flexibility in the final choice of coefficients, and perhaps lead to schemes with smaller truncation errors and new symmetries.

On the other hand, when dealing with problems formulated in the real line, the use of complex arithmetic introduces an additional computational cost  
with respect to methods with purely real coefficients. From a theoretical point of view, the vector field appearing in 
the differential equation has to be analytic at least in a domain
containing the path where the actual integration is carried out in the complex plane. Otherwise, order reductions are to be expected unless the implementation
is not conveniently adapted.

In this section we summarize some of the issues involved in the construction and analysis of splitting methods with complex coefficients, and also review some
of their most salient properties with regard to preservation of qualitative properties, in both classical Hamiltonian and quantum systems. 

For brevity, henceforth we denote a complex number $a$ with positive real part by writing $a \in \mathbb{C}_+$, so we are only interested in methods whose
coefficients $a_j, b_j \in \mathbb{C}_+$.

\subsection{Compositions}

Most of the existing splitting methods with complex coefficients have been constructed by applying the composition technique of subsection
\ref{subsec21}. Here, for completeness, we review the process, starting with the Lie--Trotter scheme $\chi_h = \varphi_h^{[2]} \circ \varphi_h^{[1]}$
as the basic method and composing it with different weights.  The simplest situation corresponds of course to 
\begin{equation} \label{cc1}
 \phi_h^{[2]} \equiv \chi_{\gamma_{1,2} h} \circ  \chi_{\gamma_{1,1} h}.
 \end{equation}
In accordance with Subsection \ref{ss2cltm}, $ \phi_h^{[2]}$ is of order 2 if the coefficients satisfy equation (\ref{or2_lt}), i.e., 
\[
  \gamma_{1,1} = \frac{1}{2} \pm \frac{1}{2} i, \qquad\quad  \gamma_{1,2} = \overline{\gamma}_{1,1} = \frac{1}{2} \mp \frac{1}{2} i.
\]  
If one tries to construct a method of order 3, a composition of at least four maps $\chi_h$
are needed, since  in that case there are 4 order conditions. Although there are 8 solutions (+cc), only 4 have positive
real part. In fact, two of these solutions result from composing $\phi_h^{[2]}$, namely
\begin{equation} \label{cc3}
   \phi_h^{[3]}  =  \phi_{\gamma_{2,2} h}^{[2]} \circ  \phi_{\gamma_{2,1} h}^{[2]} 
\end{equation}
by requiring that $\gamma_{2,1} + \gamma_{2,2} = 1, \; \gamma_{2,1}^3 + \gamma_{2,2}^3 = 0$. Notice that in
(\ref{cc3}) $\gamma_{2,2} =  \overline{\gamma}_{2,1}$.
This procedure can in principle be repeated by considering the recurrence
\begin{equation} \label{ccr}
  \phi_h^{[r+1]} = \phi_{\gamma_{r,m_r}h}^{[r]} \circ \cdots  \circ  \phi_{\gamma_{r,1}h}^{[r]}, \qquad r=1,2,\ldots
\end{equation}  
to construct higher order methods with \emph{any} first-order integrator $\phi_h^{[1]}$ (not necessarily the Lie--Trotter scheme). It turns out, however, that
$\phi_h^{[4]}$ and higher order schemes obtained
from (\ref{ccr}) have at least one coefficient with
negative real part \cite{blanes13oho}. Thus, if one is interested in schemes having only coefficients in $\mathbb{C_+}$, $r=3$
indeed constitutes an order barrier for this type of composition \cite{hansen09hos}.

Let us now consider the Strang splitting as the basic method. Then, 
composition (\ref{compm.1}) with $s=2$ already provides a method of order 3:
\begin{equation} \label{ccmo3}
  S_h^{[3]} = S_{\gamma_2 h}^{[2]} \circ S_{\gamma_1 h}^{[2]}, \qquad \mbox{ with } \qquad 
   \gamma_1 = \frac{1}{2} \pm i \frac{\sqrt{3}}{6}, \qquad \gamma_2 = \overline{\gamma}_1.
\end{equation}
Again, higher order methods can be obtained by applying the recursive procedure
 \begin{equation} \label{ccr2}
 S_{h}^{[r+1]} = S_{\gamma_{r,2} h}^{[r]} \circ S_{\gamma_{r,1} h}^{[r]}, \qquad r=2,3,\ldots,
 \end{equation}  
 with $\gamma_{r,1} + \gamma_{r,2} = 1$, $\gamma_{r,1}^{r+1} + \gamma_{r,2}^{r+1} = 0$ (see eq. (\ref{eq:ufun1})). 
 The solution with the smallest phase is given by
\begin{align} \label{eq:gam}
  \gamma_{k,1} =  \frac{1}{2} \pm \frac{i}{2} \tan \left( \frac{\pi}{2(r+1)} \right), \qquad\quad  \gamma_{r,2} = \overline{\gamma}_{r,1}, 
\end{align}
but only schemes up to order 6 with coefficients in $\mathbb{C}_+$ are possible with this approach \cite{hansen09hos}. 

As we have already seen, the triple jump (\ref{eq:triple_jump}) allows one to raise the order by two. Thus, starting from $S_h^{[2]}$, the composition
\begin{equation} \label{trij.3}
     S_{h}^{[4]} = S_{\gamma_3 h}^{[2]} \circ S_{\gamma_2 h}^{[2]} \circ  S_{\gamma_1 h}^{[2]},
\end{equation}
has, apart from the real solution (\ref{suzu1}), 
\begin{equation} \label{or4p}
     \gamma_1=\gamma_3  = \frac{1}{2 - 2^{1/3} \e^{2 i \ell \pi/3}}, \qquad
	\gamma_2 = 1- 2 \alpha_1, \qquad  \ell=1,2,
\end{equation}
corresponding to a time-symmetric composition, and 
\begin{equation} \label{or4sc}
	\gamma_1=\overline{\gamma}_3  =  \frac{1}4 \pm i \, \frac14\sqrt{\frac53}
	,  	\qquad\quad \gamma_2 = \frac12,  
\end{equation}
possessing the same symmetry as (\ref{cc3}) and (\ref{ccmo3}). 
Palindromic methods up to order 8 with coefficients having positive real part are possible
by applying the triple jump composition, whereas if one instead considers the \emph{quadruple jump},
\[
   S_{\gamma_{k,1} h}^{[2k]} \circ S_{\gamma_{k,2} h}^{[2k]} \circ S_{\gamma_{k,2} h}^{[2k]} \circ S_{\gamma_{k,1} h}^{[2k]},
\]
one can achieve order 14 with all coefficients in $\mathbb{C}_+$. These order barriers have been rigorously proved in
\cite{blanes13oho}. 

Compositions (\ref{ccmo3}) and (\ref{trij.3}) can be reformulated as splitting methods when $S_h^{[2]}$ is the Strang splitting. Thus,
\begin{equation} \label{sym_conj0}
   S_h^{[3]} = \varphi^{[1]}_{\overline{a}_1 h} \circ \varphi^{[2]}_{\overline{b}_1 h} \circ \varphi^{[1]}_{a_2 h} \circ \varphi^{[2]}_{b_1 h} 
   \circ \varphi^{[1]}_{a_1 h},
\end{equation}   
with 
\[
 a_1 = \frac{1}{4} \pm i \frac{\sqrt{3}}{12}, \qquad b_1 = \frac{1}{2} \pm i \frac{\sqrt{3}}{2},  \qquad a_2 = \frac{1}{2}, 
\] 
 whereas the palindromic version of (\ref{trij.3}) reads
\begin{equation} \label{pa4p}
   S_{h,P}^{[4]} = \varphi^{[1]}_{a_1 h} \circ \varphi^{[2]}_{b_1 h} \circ \varphi^{[1]}_{a_2 h} \circ \varphi^{[2]}_{b_2 h} \circ \varphi^{[1]}_{a_2 h} \circ 
  \varphi^{[2]}_{b_1 h} \circ \varphi^{[1]}_{a_1 h}, 
\end{equation}
with 
\[
 a_1 = \frac{\gamma_1}{2}, \quad b_1 = \gamma_1, \quad a_2 = \frac{1}{2}-a_1, \quad b_2 = 1-2 b_1, 
 \]
 respectively. Finally, solution (\ref{or4sc}) leads to
\begin{equation} \label{sc4sc}   
    S_{h,C}^{[4]} = \varphi^{[1]}_{\overline{a}_1 h} \circ \varphi^{[2]}_{\overline{b}_1 h} \circ \varphi^{[1]}_{\overline{a}_2 h} 
    \circ \varphi^{[2]}_{b_2 h} \circ \varphi^{[1]}_{a_2 h} \circ 
  \varphi^{[2]}_{b_1 h} \circ \varphi^{[1]}_{a_1 h}, 
\end{equation}
with
\[
  a_1 = \frac{b_1}{2}, \qquad b_1 = \frac{1}4 \pm i \, \frac14\sqrt{\frac53}, \qquad a_2 = \frac{1}{4}(2 b_1 + 1), \qquad b_2 = \frac{1}{2}.
\]
Before proceeding further, it is useful to illustrate how the different symmetries of these third- and fourth-order schemes manifest in practice on a simple
example.

\paragraph{Example: two-level system.}
We consider the time evolution of a two-level quantum system, described by \cite{messiah99qme}
\begin{equation} \label{eq.unit.1}
    i \, \frac{dU}{dt} = H U = ( \sigma_1 + \sigma_2) \, U, \quad \mbox{ with } \quad 
     \sigma_1 = \left( \begin{array}{rr}
 			0 & 1 \\
			1 & 0
		\end{array} \right), \;\;
 \sigma_2 = \left( \begin{array}{rr}
 			0 & -i \\
			i & 0
		\end{array} \right)
\end{equation}
and $U(0)=I$.
Clearly, $U(t)$ is a $2 \times 2$ unitary matrix and $\det U(t) = 1$. We test both the preservation of unitarity in the numerical approximations
and the computational efficiency of the previous 3th- and 4th-order schemes. To this end, we take $\varphi_h^{[j]} = \e^{\tau \sigma_j}$, $j=1,2$, with $\tau = -i h$, 
and compare with the exact solution $\e^{\tau(\sigma_1 + \sigma_2)}$.
Notice that, for problems formulated in complex variables, the use of complex coefficients does not generally increase the overall cost
of the algorithm. 

We fix $t_f=1000$ as the final time, and adjust $h$ (and therefore 
the number of steps $n$) so that all the methods require the same computational effort. Specifically, $n=6000$ ($h=1/6$) for scheme  
(\ref{sym_conj0}) and $n=4000$ ($h=1/4$) for the 4th-order methods (\ref{pa4p}) and (\ref{sc4sc}). Finally 
we compute $|\|U_{\mathrm{app}}(jh)\|-1|, \ j=1,2,\ldots, n$, where $U_{\mathrm{app}}(jh)$ denotes the approximate solution after $j$ steps, which
is shown in Figure \ref{figu71cc} (left). Dashed (red) line corresponds to scheme (\ref{pa4p}), solid (black) line is produced
by (\ref{sc4sc}), and finally the dashed-dotted line is the result obtained with the 3th-order scheme (\ref{sym_conj0}). Notice the different behavior
of different classes of integrators with complex coefficients: whereas
error in the
preservation of unitarity does not grow with time for methods (\ref{sc4sc}) and (\ref{sym_conj0}), it certainly does for the palindromic scheme (\ref{pa4p}).

Figure \ref{figu71cc} (right) shows an efficiency diagram of these integrators, together with the classical triple jump with real coefficients (dotted line). It is produced
by integrating eq. (\ref{eq.unit.1}) with different values of $h$ and computing
the error of the approximation (in the 2-norm) at the final time $t_f = 10$ as a function of the computational cost (estimated as the number of exponentials 
involved in the whole integration). We observe the surprisingly good performance of scheme (\ref{sym_conj0}) and the superior efficiency of the 4th-order
methods with complex coefficients in comparison with their real counterpart. This
is related to the difference in size of the main error term in the expansion of the modified operator associated with each 4th-order
method, which is typically
associated with the coefficient multiplying ${\cal E}_5=|\sum_{j=1}^3 \gamma_j^5|$. Thus, $\mathcal{E}_5$ 
is about 200 times smaller for the schemes (\ref{or4p}) and
(\ref{or4sc}) than for the triple jump with real coefficients \cite{blanes22asm}.

\begin{figure}[htb]
\centering
  \includegraphics[scale=0.42]{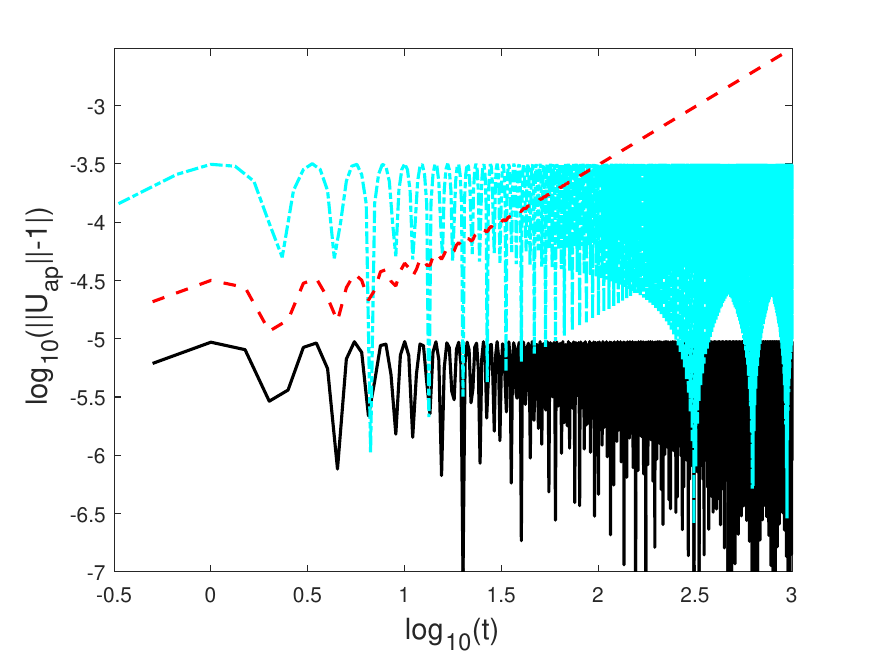}
  \includegraphics[scale=0.42]{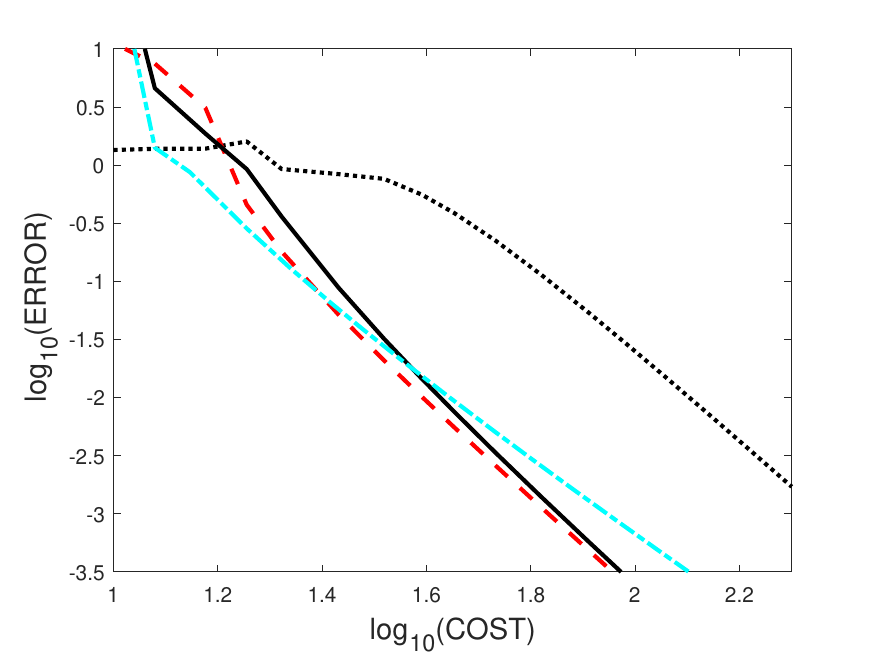}
\caption{Left: evolution of error in unitarity for system (\ref{eq.unit.1}) integrated with 3th- and 4th-order splitting methods with complex coefficients. Right: 
efficiency diagram computed for a final time $t_f = 10$. Different patterns in the coefficients lead to different qualitative behavior.}
\label{figu71cc}
\end{figure}

\noindent
$\Box$
 
 \

\subsection{Symmetric-conjugate methods and unitary problems}

The previous example shows that the time-symmetric character of a method with complex coefficients does not necessarily guarantee good preservation
properties, whereas the situation is different for schemes (\ref{sc4sc}) and (\ref{sym_conj0}). They belong to the general class of methods
\begin{equation} \label{sym_conj1}
  \psi_h = \varphi^{[1]}_{a_{s+1} h}\circ \varphi^{[2]}_{b_{s}h}\circ \varphi^{[1]}_{a_{s} h}\circ \cdots\circ
 \varphi^{[1]}_{a_{2}h}\circ
 \varphi^{[2]}_{b_{1}h} \circ \varphi^{[1]}_{a_{1}h}
\end{equation}
whose coefficients verify
\begin{equation} \label{sym_conj0a}
  a_{s+2-j} = \overline{a}_j, \qquad b_{s+1-j} = \overline{b}_j, \qquad j=1,2,\ldots
\end{equation}  
and can be properly called \emph{symmetric-conjugate}. They possess the following distinctive feature: assume our differential equation 
$x' = f(x)$ is reversible with respect to complex conjugation, in the sense that $\overline{f(x)} = - f(\overline{x})$ for all $x \in \mathbb{C}$, and the same
holds for each piece $f_j$ in $f = f_1 + f_2$. Then scheme (\ref{sym_conj1})-(\ref{sym_conj0a}) verifies
\begin{equation} \label{rever1}
  \overline{\psi}_h = \psi_h^{-1},
\end{equation}
so that the map $\psi_h$ is also reversible with respect to complex conjugation \cite[section V.1]{hairer06gni}. Notice that (\ref{rever1}) is not satisfied by
palindromic compositions involving complex coefficients.

Property (\ref{rever1}) has some major consequences. Suppose one is dealing with the linear problem
\begin{equation} \label{scli1}
   i \frac{du}{dt} = H u, \qquad u(0) = u_0, \qquad\qquad u \in \mathbb{C}^N,
\end{equation}    
where $H$ is a $N \times N$ real matrix of the form $H = H_1 + H_2$, with $H_1$, $H_2$ also real, so that the symmetric-conjugate
method (\ref{sym_conj1})-(\ref{sym_conj0a}) reads
\begin{equation} \label{scli2}
   \psi_{\tau} = \e^{\overline{a}_1 \tau H_1} \, \e^{\overline{b}_1 \tau H_2} \, \e^{\overline{a}_2 \tau H_1} \, \cdots \, \e^{a_2 \tau H_1} \, \e^{b_1 \tau H_2}
   \, \e^{a_1 \tau H_1}
\end{equation}   
with $\tau = -i h$. It has been proved in \cite{bernier23scs} that if
all the eigenvalues of $H$ are real and simple, then, for sufficiently small $h$, there exist real matrices $D_{\tau}$ (diagonal) and $P_{\tau}$ (invertible) such that
$\psi_{\tau}^n = P_{\tau} \, \e^{n \tau D_{\tau}} P_{\tau}^{-1}$, all the eigenvalues of the linear map $\psi_{\tau}$ have modulus 1 and the norm $|u|^2$, 
and the energy $\overline{u}^T H u$ are almost
preserved for long times.  In other words, any symmetric-conjugate splitting method $\psi_{\tau}$ applied to (\ref{scli1})
is similar to a unitary method for sufficiently small values of the step size $h$. 

One should recall that the time-dependent Sch\"odinger equation, once discretized in space, leads to an equation of the
form (\ref{scli1}), with both $H_1$ and $H_2$ real and symmetric matrices: $H_1 = T$ is associated with the second-order periodic spectral
differentiation matrix and $H_2 = V$ is the (diagonal) matrix corresponding to the potential evaluated at the grid points. It is then admissible to use
symmetric-conjugate methods in this setting, since they also guarantee preservation of invariants for sufficiently small values of $h$. 
We should notice, however, that both positive and negative imaginary parts are present in the scheme if 
$a_j \in \mathbb{C}$ (since the method is consistent), and this
may lead to severe instabilities due to the unboundedness of the Laplace operator \cite{castella09smw,hansen09hos}.
It then seems appropriate to consider in this setting only symmetric-conjugate splitting methods with $0 < a_j < 1$ and $b_j \in \mathbb{C}_+$. 
Different schemes of orders 3 up to 6 within this family have been proposed and tested in \cite{bernier23scs}, also showing promising results in terms
of efficiency.

If, on the other hand, the equation to integrate is of the form 
\begin{equation} \label{reli1}
    \frac{du}{dt} = H u = (H_1 + H_2) u, \qquad u(0) = u_0, \qquad u \in \mathbb{R}^N
\end{equation}    
with $H_1$, $H_2$ real symmetric matrices, then a symmetric-conjugate splitting method $\psi_h$ of order $r$ satisfies
\[
  \left( \overline{\psi}_h \right)^T = \psi_h,
\]
for all values of $h$. In consequence, there exists a family of
unitary matrices $U_h$  and a family of real diagonal matrices $D_h$ depending smoothly on $h$ such that
\[
  (\psi_h)^n = U_h \, \e^{n h D_h} \,   (\overline{U}_h)^T,
\]
where $\mathrm{Im}(U_h) =  \mathcal{O}(h^r)$ and $D_h$ is a perturbation of order $r$ of the matrix $D_0$ diagonalizing $H$. If $u_n = (\psi_h)^n u_0$,
then the quantity $\| \mathrm{Im}(u_n) \|/\|u_n\|$ remains bounded along the numerical trajectory, since the error  in the imaginary part is only due
to the transformation $U_h$. Again, this is not true for palindromic schemes with complex coefficients. 

We can then conclude that symmetric-conjugate splitting methods can be safely used to integrate equation (\ref{reli1}) when $H_1$ and $H_2$ are
 real symmetric, as for palindromic (or time-symmetric) schemes with real coefficients, with one important difference: whereas in
 splitting methods with real coefficients at least one $a_j$ and at least one $b_j$ are negative when the order $r \ge 3$, 
symmetric-conjugate methods of order $r \ge 3$ do exist with coefficients having positive real part.

Several parabolic PDEs lead to equation (\ref{reli1}) after space discretization. This is the case, in particular, for the time-dependent
Schr\"odinger equation in imaginary time \eqref{itp1}.

\subsection{Projecting on the real axis}

In the general case of an equation  $x' = f(x) = f_1(x) + f_2(x)$ with $f_1$, $f_2$ real, the usual practice consists 
in projecting the numerical evolution after each time step to its real part, since one is interested in getting real approximations to the exact
solution. Thus, if $\psi_h$ is a splitting method of order $r$, by following this approach, we are actually applying the scheme
\begin{equation} \label{ra1}
  R_h = \frac{1}{2} \big( \psi_h + \overline{\psi}_h \big),
\end{equation}
which is also at least of order $r$, in agreement with the comments at the end of section \ref{sect6} \cite{hansen09hos}. In fact,  
the order of $R_h$ is actually $r+1$ if $\psi_h$ is symmetric-conjugate and $r$ is odd. This
was already noticed in \cite{chambers03siw} for the 3rd-order method (\ref{sym_conj0}) and proved in general in \cite{blanes22osc}. 

This feature of symmetric-conjugate methods has some effects, especially when constructing high order schemes to be used by projecting on the real axis. 
Although apparently
a symmetric-conjugate composition requires solving the same number of order conditions to achieve order $r$ as a general composition (and more than
a palindromic one for orders higher than four), additional reductions take place for the projected method (\ref{ra1}) \cite{blanes22osc}, and in fact the resulting 
$R_h$ requires less stages. Thus, in particular, it is possible to construct a composition of the form
\begin{equation} \label{ra2}
  \psi_h = S_{\overline{\alpha}_1 h}^{[2]} \circ S_{\overline{\alpha}_2 h}^{[2]} \circ \cdots \circ
    S_{\alpha_2 h}^{[2]} \circ S_{\alpha_1 h}^{[2]},
\end{equation}
with $S_{h}^{[2]}$ a time-symmetric 2nd-order method, with only 5 stages so that its projected part $R_h$ is of order 6, whereas
schemes based on palindromic compositions require at least 7 stages. The reduction is more notable for higher orders: thus, a composition (\ref{ra2}) of
order 5 involving 9 appropriately chosen stages leads to a projected method $R_h$ of order 8. By contrast, 15 stages are required by palindromic
compositions.

One of the salient features of splitting and composition methods is that they preserve by construction qualitative features of the exact solution, as 
seen in Section \ref{sect4}. 
These favorable properties are of 
course lost when the special linear combination (\ref{ra1}) is considered, 
but nevertheless the resulting scheme $R_h$ still preserves them up to an order much higher than
the order of the method itself. More specifically, suppose $S_{h}^{[2]}$ is a time-symmetric 2nd-order and symplectic method applied to
a Hamiltonian system. Then, one has the following result \cite{blanes22osc}:
\begin{itemize}
 \item If the symmetric-conjugate composition $\psi_h$ given in (\ref{ra2}) is of odd order $r=2k-1$, $k \ge 2$, then the method $R_h$ (of order $2k$) 
 preserves time-symmetry up to order $4k-1$, i.e.,
 \[
    R_h \circ R_{-h} = \mathrm{id} + \mathcal{O}(h^{4k}),
\]
and symplecticity up to the same order,    
\[
  (R_h^{\prime})^T \, J \, R_h^{\prime} = J + \mathcal{O}(h^{4k}),
\]
where $J$ denotes the canonical symplectic matrix and $R_h^{\prime}$ is the Jacobian.
\item If $\psi_h$ is  symmetric-conjugate and of even order $r=2k$, then $R_h$ (which is also of order $2k$) preserves time-symmetry and symplecticity up to order $4k+3$.
\item If $\psi_h$ is palindromic of order $r=2k$, then the resulting $R_h$ (also of order $2k$)  preserves time-symmetry and symplecticity up to order $4k+1$.
\end{itemize}



\section{A collection of splitting methods}
\label{sect8}

Once the order conditions for a splitting or composition method of a given order $r$ have been
explicitly obtained, the next step in the construction of particular schemes consists of course in
solving these polynomial equations to determine the coefficients. To begin, one considers
compositions with as many parameters as equations and try to use a computer algebra system to
determine all real solutions. Nevertheless, solving the order conditions in this way is only possible for
moderate values of $r$, and thus one has to turn to numerical techniques. Since the number of real
solutions usually increases with the order, the problem is how to select the particular solution
expected to give the best performance when the integrator is applied to practical problems. This
is typically done by minimizing some objective function, depending on the particular class of schemes
considered. Thus, in the case of splitting methods of the form (\ref{eq:splitting}), one tries to minimize
the leading error term in the asymptotic expansion of $\log(\Psi(h))$ in (\ref{bch2}), namely
$\sum_{j=1}^{c_{r+1}} \beta_j F_{{r+1},j}$, where $F_{{r+1},j}$ denote the elements of the Lyndon
basis associated with Lyndon words with $r+1$ letters. The objective function is then
\[
  E_{r+1} = \left( \sum_{j=1}^{c_{r+1}} |\beta_j|^2 \right)^{1/2}.
\]  
One has to take into account, however, that $E_{r+1}$ will change if another basis of Lie brackets is considered. We are also assuming that all these brackets contribute in a similar way, something that is not guaranteed to take place in general. 
It makes sense, then, to introduce other quantities as possible
estimators of the error committed. In particular, it has been noticed that large coefficients $a_j$, $b_j$ in the splitting method 
usually lead to large truncation errors, since higher order terms in $\Psi(h)$ depend on increasingly higher powers of these coefficients. 
For this reason, it is also convenient to keep track of the quantities
\begin{equation} \label{sizeco}
  \Delta \equiv \sum_j \left( |a_j| + |b_j| \right) \qquad \mbox{ and } \qquad \delta \equiv \max_{j} \, (|a_j|, |b_j|)
\end{equation}
and eventually discard solutions with large values of $\Delta$ and/or $\delta$.  In the case of
compositions (\ref{eq:compint}) of a basic first-order method and its adjoint, and compositions of Strang maps, (\ref{compm.1}), 
a frequently used criterion is to choose the solution which minimizes 
$\sum_{j=1}^{2s} |\alpha_j|$, and $\sum_{j=1}^{s} |\gamma_j|$, respectively. Other possibilities include minimizing only those terms in
the truncation error that are not removable by a processor, in accordance with the analysis carried out in Sections \ref{sect4} and \ref{sect5}.

Including additional maps (stages) in the composition provides additional parameters that may lead to smaller values of the chosen
objective function. It is important at this point to remark
that the efficiency of a method is
measured by taking into account the computational cost required to achieve a
given accuracy. Thus, in case one has
several methods of order $r$ with different computational cost
(usually measured as the number of stages or evaluations of the functions involved), the most efficient
method does not necessarily correspond to the cheapest method: the 
extra cost of some methods can be compensated by an improvement in
the accuracy obtained. In fact, this is what usually happens in practice, although solving the polynomial
equations with additional stages (and free parameters) is by no means a trivial task. Continuation techniques
have been shown to be very useful in this context (see e.g. \cite{blanes13nfo,alberdi19aab}).

The use of the processing technique also allows one to construct methods of a given order typically requiring a reduced computational effort,
so that the overall efficiency of the resulting schemes is enlarged if the output is not frequently computed.

With all these considerations in mind, our purpose in this section is to present 
a comprehensive overview of (most of) the existing methods with real coefficients by classifying them
into different families and giving the appropriate references. The corresponding coefficients can be found at the website
\url{https://www.gicas.uji.es/SplittingMethods.html}.

\subsection{Symmetric compositions of time-symmetric 2nd-order schemes}

Perhaps the first method in this family corresponds to the 4th-order method obtained from the triple jump composition (\ref{eq:triple_jump})
when $k=2$. It was known and internally used in the accelerator physics community during the 1980s \cite{forest90fos} but was independently discovered 
afterwards in several settings: in \cite{candy91asi}, in \cite{campostrini90aco} as an algorithm for hybrid Monte Carlo simulations,  and also in 
\cite{creutz89hoh} as the more general composition
\begin{equation} \label{eq:S_h-recP}
   S_{h}^{[2k+2]} = \left(S_{\gamma_1 h}^{[2k]}\right)^p 
	\circ 
	S_{(1- 2 p\gamma_1) h}^{[2k]} 
	\circ 
	\left(S_{\gamma_1 h}^{[2k]}\right)^p,  \quad \mbox{ with } \quad
     \gamma_1= \frac{1}{2p - (2p)^{1/(2k+1)}}.
\end{equation}
This recursion was also obtained in \cite{suzuki90fdo,suzuki91gto}, whereas the case $p=1$ became immensely popular following \cite{yoshida90coh}.
It was in a certain sense generalized in  \cite{mclachlan02foh}, where a rule of thumb was provided to select the optimal value of $p$ for different orders, also valid
for processed methods.

As stated in Subsection \ref{subsec21}, other choices for the coefficients $\gamma_j$ in the more general composition 
\begin{equation} \label{ss_again}
   \psi_h = S_{\gamma_s h}^{[2]} \circ S^{[2]}_{\gamma_{s-1} h} \circ \cdots \circ S^{[2]}_{\gamma_1 h}
\end{equation}
lead to more efficient schemes when $r \ge 6$. It turns out that virtually all published methods of this form correspond to time-symmetric
compositions, i.e, $\gamma_{s+1-j}=\gamma_j$, and are therefore of even order. There are at least two reasons for that: (i) the task of constructing new
schemes is simplified, since the number of order conditions to be solved is much reduced, and (ii) the resulting methods also require a smaller number of stages
$s$, since the number of order conditions at even orders (automatically satisfied by a time-symmetric composition) is greater than the total number of order
conditions divided by 2, for sufficiently large $r$. For instance, constructing a non-symmetric 8th-order scheme requires at least $s=16$ (and therefore
to solve a system of 16 polynomial equations), whereas $s=15$ is the minimum number for a symmetric composition, so that only 8 polynomial equations have
to be solved.

\begin{table}
  \centering
\begin{tabular}{|l|l|l|l|} \hline
\multicolumn{4}{|c|}{\bfseries  $\psi_h = S_{\gamma_s h}^{[2]} \circ S^{[2]}_{\gamma_{s-1} h} \circ \cdots \circ S^{[2]}_{\gamma_1 h}$}
\\ \hline
 $r=4$ & $r=6$  & $r=8$ & $r=10$ \\ \hline
\textbf{3}-CR90 & \textbf{7}-YOS90 & {\bf 15}-YOS90,SUZ93, &  {\bf 31}-SUZ93,KR97, \\
\framebox[3mm]{{\bf 5}}-SUZ90,    & {\bf 9}-MCL95b, &\quad \; MCL95b,KR97 &  \quad \; SS05 \\
\quad MCL95b & \quad KR97 &	{\bf 17}-MCL95b, KR97 \hspace*{-0.2cm} &  
		{\bf 33}-KR97,TSI99,   \\
          & {\bf 11-\framebox[4.5mm]{13}}-SS05 \hspace*{-0.2cm}  &  {\bf 19-\framebox[4.5mm]{21}}-SS05 &  \quad  HLW02,SS05 	 \\
\qquad  &  &  {\bf 24}-CSS93 & {\bf \framebox[4.5mm]{35}}-SS05   \\
\hline \hline
   {\bf P:3-17}-MCL02  \hspace*{-0.2cm}& 
		{\bf P:\framebox[4.5mm]{11}-13} 
  	& {\bf P:\framebox[4.5mm]{13}-19}-BCM06 & {\bf P:19-\framebox[4.5mm]{23}}-BCM06 \\
   &	\ \ -BCM06  &  & \\ \hline
\end{tabular}
  \caption{Symmetric compositions of 2nd-order symmetric methods of order $r=4,6,8,10$ published in the literature. We indicate the
  number of stages \textbf{s} and the pertinent reference. Processed methods are preceded by
  \textbf{P}. The recommended methods, in accordance with the experiments carried out in Appendix \ref{apenda}, are \fbox{framed}.}\label{tableSS}
\end{table}

In Table~\ref{tableSS} we present the most relevant schemes of this type found in the literature.
At each order, $r$, we label each method by the number of
stages ${\bf s}$ and an acronym indicating the author(s) and year when it was first published. 
 We also include processed methods, referred to as {\bf P}:${\bf s}$,
where $\mathbf{s}$ is the number of stages of the
kernel. The following list provides additional information about the collected schemes:
\begin{itemize}
	\item CR90: collective name given to the triple jump composition, (\ref{eq:S_h-recP}) with $k=p=1$.
	\item SUZ90: recursion (\ref{eq:S_h-recP}) with $k=1$, $p=2$, proposed in \cite{suzuki90fdo}.
	\item YOS90: \cite{yoshida90coh} gives three 6th-order compositions with $s=7$ and five 8th-order methods with $s=15$.
	\item CSS93: \cite{calvo93hos} obtain a 7th-order non-symmetric scheme which can be written as composition \eqref{ss_again} with $s=12$. After symmetrization, it provides a 24-stage 8th-order method.
	\item SUZ93: \cite{suzuki93hod} provides one 6th-order method with $s=14$, six 8th-order methods with $s=15$ and four 10th-order schemes with $s=31$. 
	\item MCL95b: \cite{mclachlan95otn} gives several optimized schemes: one 4th-order method with $s=5$, one 6th-order method with $s=9$, and two 8th-order methods with $s=15$ and $s=17$. 
	\item KR97: \cite{kahan97ccf} also present several optimized schemes: two 6th-order methods with $s=9$, three 8th-order methods (one with $s=15$ and two with $s=17$), and five 10th-order methods (two with $s=31$ and three with $s=33$).
	\item TSI99: \cite{tsitouras99ato} constructs one optimized 10th-order method with $s=33$.
	\item MCL02: construction of kernel with recurrence (\ref{eq:S_h-recP}) for several values of $p$ in  \cite{mclachlan02foh}.
	\item HLW02: \cite[first edition published in 2002]{hairer06gni} provide one optimized 10th-order method with $s=33$, the most efficient at the time.
	\item SS05: \cite{sofroniou05dos} carry out an exhaustive search of 6th-order schemes with $s=11,13$, 8th-order methods with $s=19,21$, and 10th-order methods with $s=31,33,35$.
	\item BCM06: \cite{blanes06cmf} give several processed methods: two 6th-order kernels with $s=11$ and $s=13$, two 8th-order kernels with $s=13$ and $s=19$, and  two 10th-order kernels with $s=19$ and $s=23$.
\end{itemize}

\subsection{Splitting into two parts / composition of a basic first-order method and its adjoint}

Although the order conditions for the two types of scheme, 
 \begin{equation} \label{split_again}
\psi_h = \varphi^{[1]}_{a_{s+1} h}\circ \varphi^{[2]}_{b_{s}h}\circ \varphi^{[1]}_{a_{s} h}\circ \cdots\circ
 \varphi^{[1]}_{a_{2}h}\circ
 \varphi^{[2]}_{b_{1}h} \circ \varphi^{[1]}_{a_{1}h},
\end{equation}
and 
\begin{equation} \label{compo_again}
\psi_h =  \chi_{\alpha_{2s} h} \circ \chi^*_{\alpha_{2s-1}h} \circ \cdots \circ \chi_{\alpha_{2}h} \circ \chi^*_{\alpha_{1}h},
\end{equation}
are equivalent by virtue of the relationship (\ref{eq:abalpha}),
the optimization
procedures to get the most efficient schemes may
differ. In consequence, a particular method optimized for
systems that are separable into two parts is not necessarily the best scheme when written as
(\ref{compo_again}), although in practice their performances
are closely related.

Since schemes of order $r \ge 6$  require more stages than taking the composition (\ref{ss_again}) with  
 $S_{h}^{[2]}=\chi_{h/2} \circ \chi^*_{h/2}$ or $S_{h}^{[2]}=\varphi^{[1]}_{h/2}\circ \varphi^{[2]}_{h} \circ \varphi^{[1]}_{h/2}$, only methods with $r \le 4$ seem
 promising. Note, however, that a 6th-order symmetric composition (\ref{ss_again}) with $s=7$ has only 3 real solutions for the $\gamma_j$, whereas
symmetric versions of (\ref{split_again}) and (\ref{compo_again}) require at least $s=9$ stages to solve the 9 order conditions. This enlarged number of equations
might then provide some solution leading to smaller error terms.

\begin{table}
  \centering
\begin{tabular}{|l|l|l|l|} \hline
\multicolumn{4}{|c|}{\bfseries 
$\psi_h =  \chi_{\alpha_{2s} h} \circ \chi^*_{\alpha_{2s-1}h} \circ \cdots \circ \chi_{\alpha_{2}h} \circ \chi^*_{\alpha_{1}h}$}
\\
\multicolumn{4}{|c|}{\bfseries $\psi_h = \varphi^{[1]}_{a_{s+1} h}\circ \varphi^{[2]}_{b_{s}h}\circ \varphi^{[1]}_{a_{s} h}\circ \cdots\circ
 \varphi^{[1]}_{a_{2}h}\circ
 \varphi^{[2]}_{b_{1}h} \circ \varphi^{[1]}_{a_{1}h}$ }
\\
\hline
  $r=2$ & $r=3$ & $r=4$ & $r=6$    \\ \hline
 \framebox[3mm]{{\bf 2}}-MCL95b, & {\bf 3}-RUT83, & 
  {\bf 4-5}-MCL95b 
  &        {\bf 9}-FOR92   \\
	\quad OMF03, BCS14    & \ \ \ SUZ92   & 
 \framebox[3mm]{{\bf 6}}-BM02 
	&  \framebox[4.5mm]{{\bf 10}}-BM02 \\
 {\bf 3}-BCS14 
	&    &   &   \\
						\hline  \hline
&  & {\bf P:3,\framebox[3mm]{4}}-BCR99  & {\bf P:5}-BCR99    \\
 &  & {\bf P:\framebox[3mm]{6}}-BCM06  & {\bf P:\framebox[3mm]{9}}-BCM06 \\ \hline
\end{tabular}
  \caption{Symmetric composition schemes of the form  (\ref{split_again}) (appropriate when the ODE is
  split in two parts) and
(\ref{compo_again}) (composition of a first-order method and its adjoint) of order $r=2,3,4,6$ with real coefficients. 
The notation is the same as in Table~\ref{tableSS}.
}\label{tableAB}
\end{table}

The schemes presented in Table~\ref{tableAB} have been specifically designed to deal with these problems, and cannot be obtained as particular cases of 
the composition (\ref{ss_again}). In more detail, they correspond to the following.
\begin{itemize}
	\item RUT83: \cite{ruth83aci} gives the first 3-stage third order non-symmetric method for systems separable into two parts.
	\item SUZ92: \cite{suzuki92gnh} presents a family of 3-stage third order methods for the composition \eqref{compo_again} with $s=3$ and $\alpha_6=0$.
	\item FOR92: \cite{forest92sol} gives a 6th-order method for separable systems with $s=9$.
	\item MCL95b: \cite{mclachlan95otn} gives an optimized 2-stage 2nd-order method and two optimized 4th-order schemes of the form
	(\ref{split_again}) with $s=4$ and $s=5$. 
	\item BCR99: \cite{blanes99siw} present several processed methods for separable systems: 4th-order methods with $s=3$ and $s=4$, one 5th-order method with a non-symmetric kernel and $s=4$, and finally one 6th-order scheme with $s=5$.
	\item BM02: \cite{blanes02psp} obtain one 4th-order method with $s=6$ and one 6th-order method with $s=10$, both optimized for separable systems.
	\item OMF03:  \cite{omelyan03sai} rediscover the optimized 2-stage 2nd-order scheme of \cite{mclachlan95otn} as an efficient method with positive coefficients.
	\item BCM06: \cite{blanes06cmf} provide one 4th-order processed method with $s=6$ and two 6th-order processed schemes with $s=9$ and $s=10$.
	\item BCS14: \cite{blanes14nif} propose efficient symmetric 2nd-order methods with $s=2$ and $s=3$, designed for hybrid Monte Carlo simulations.
\end{itemize}

\subsection{Runge--Kutta--Nystr\"om methods}

Given the relevance of the problem $y^{\prime\prime} = g(y)$, and the simplifications arising when splitting methods are applied to this system (see
subsection \ref{subsec3.1}), many RKN splitting schemes have been obtained in the literature. 
Since $F_1$ and $F_2$ play different roles here, it is
convenient to classify them according to the sequence of coefficients as
\begin{equation}  \label{ab-ba}
\begin{aligned}
 \mbox{ AB: }
 &  & \varphi^{[A]}_{a_{s}h}\circ \varphi^{[B]}_{b_{s} h}\circ
 \cdots\circ
 \varphi^{[A]}_{a_{1}h}\circ \varphi^{[B]}_{b_{1}h}, \\
 \mbox{ BA: } 
&  & \varphi^{[B]}_{b_{s} h} \circ \varphi^{[A]}_{a_{s}h}\circ
 \cdots\circ
 \varphi^{[B]}_{b_{1}h}\circ \varphi^{[A]}_{a_{1}h},
\end{aligned} 
\end{equation}
although in fact both types are conjugate to each other, thus providing similar efficiency.
A word of caution is necessary: \emph{in all the RKN
splitting methods published in the literature the operator $F_1$ as defined in subsection \ref{subsec3.1} is in fact associated with the map $\varphi^{[B]}_{t}$, and $F_2$ is associated with $\varphi^{[A]}_{t}$}}, so that $[F_1,[F_1,[F_1,F_2]]]=0$.
Thus, in particular, in classical and quantum mechanics, $A$ corresponds to the kinetic energy and $B$ to the potential energy. This should be taken
into account when implementing the schemes collected here.

 It is also convenient to take profit of the FSAL (First Same As
Last) property, and thus one may consider the non-equivalent
compositions
\begin{equation}  \label{aba-bab}
\begin{aligned}
 \mbox{ ABA: } &  & \varphi^{[A]}_{a_{s+1}h}\circ \varphi^{[B]}_{b_{s} h}\circ
 \varphi^{[A]}_{a_{s}h}\circ
 \cdots\circ
 \varphi^{[B]}_{b_{1}h}\circ
 \varphi^{[A]}_{a_{1}h} \\
  \mbox{ BAB: }  &  & \varphi^{[B]}_{b_{s+1} h}\circ
 \varphi^{[A]}_{a_{s}h}\circ \varphi^{[B]}_{b_{s} h}\circ
 \cdots\circ
 \varphi^{[A]}_{a_{1}h} \circ \varphi^{[B]}_{b_{1}h},
\end{aligned} 
\end{equation}
as well as their time-symmetric versions ($a_{s+2-i}=a_i,b_{s+1-i}=b_i$ for ABA and 
$b_{s+2-i}=b_i,a_{s+1-i}=a_i$ for BAB, respectively).

\begin{table}
  \centering
\begin{tabular}{|l|l|l|l|} \hline
\multicolumn{4}{|c|}{RKN splitting methods}\\ \hline
 $r=4$ & $r=5$ & $r=6$  & $r=8$  \\ \hline
 {\bf 4N}$_{\mathrm{AB}}$-MA92    & {\bf 5N}$_{\mathrm{ABA}}$-OS94 
  &  {\bf 7S}$_{\mathrm{ABA}}$-FOR92, 	& {\bf 17S}$_{\mathrm{ABA}}$-OL94  
	\\
 {\bf 4N}$_{\mathrm{BAB}}$-CSS93     & {\bf \framebox[3mm]{6}N}$_{\mathrm{AB}}$-MA92  
		& \ \ OS94		&   {\bf 17,18,\framebox[4.5mm]{19}S}$_{\mathrm{ABA},}$  \\
 {\bf 4,5S}$_{\mathrm{ABA}}$-MCL95b     & {\bf 6N}$_{\mathrm{AB}}$-CHO00
		&   {\bf 7S}$_{\mathrm{BAB}}$-FOR92 	&  $_{\mathrm{BAB}}$-BCE22  \\
 {\bf \framebox[3mm]{6}S}$_{\mathrm{ABA},\mathrm{BAB}}$-BM02  &  &  {\bf 11S}$_{\mathrm{BAB}}$-BM02	& \\
  {\bf 4,5S}$_{\mathrm{ABA},\mathrm{BAB}}$-OMF03  &  &  {\bf \framebox[4.5mm]{14}S}$_{\mathrm{ABA}}$-BM02  &  \\
		\hline\hline
 {\bf P:\framebox[3mm]{2}N}$_{\mathrm{AB}}$-BCR99  &   & {\bf P:4-6S}$_{\mathrm{ABA,BAB}}$
 & {\bf P:9S}$_{\mathrm{ABA}}$   \\
  &     & \quad -BCR01a 	& \quad -BCR01a  	\\ 
	&   &  {\bf P:\framebox[3mm]{7}S}$_{\mathrm{BAB}}$-BCR01b	&  {\bf P:\framebox[4.5mm]{11}S}$_{\mathrm{BAB}}$	\\ 
	&   &  		& \quad -BCR01b			\\  \hline
\end{tabular}
  \caption{RKN splitting integrators of order $r=4,5,6,8$. Since the role of the flows
$\varphi_h^{[A]}$ and $\varphi_h^{[B]}$ is not interchangeable
here, we distinguish symmetric {\bf S} and non-symmetric {\bf N}
compositions with a subscript $\mathrm{AB,ABA,BAB}$. As usual, processed
methods are preceded by {\bf P} and the recommended schemes are \fbox{framed}.}\label{tableRKN}
\end{table}

In
Table~\ref{tableRKN} we present the most representative methods  within
this class. We add {\bf S} or {\bf N} to distinguish symmetric from non-symmetric schemes
and the subscript AB, ABA and BAB to denote the different alternatives (\ref{ab-ba}), (\ref{aba-bab}). 
Processed methods have also been included. Explicitly, they correspond to the following.
\begin{itemize}
	\item FOR92: \cite{forest92sol} gives a 6th-order method with $s=7$.
	\item MA92: \cite{mclachlan92tao} present non-symmetric 4th- and 5th-order methods with $s=4$ and $s=6$, respectively.
	\item CSS93: \cite{calvo93tdo} obtain one 4th-order non-symmetric BAB method with $s=4$. 
	\item OS94: \cite{okunbor94crk} obtain four 5th-order methods with $s=5$ and sixteen 6th-order symmetric methods with $s=7$.
         \item OL94: \cite{okunbor94eoe} obtain one 8th-order ABA method with $s=17$.
	\item MCL95b: \cite{mclachlan95otn} gives optimized RKN methods 	of order 4 with $s=4$ and $s=5$, and of order 6 with $s=7$.
	\item BCR99: \cite{blanes99siw} present a processed 4th-order method with a non-symmetric kernel and $s=2$.
	\item CHO00 \cite{chou00ofs} get a non-symmetric scheme of order 5 with $s=6$.
	\item BCR01a: \cite{blanes01hor} present one processed 6th-order method with $s=6$
	\item BCR01b: \cite{blanes01nfo} give one 6th- and one 8th-order processed method with $s=7$ and $s=11$, respectively.
	\item BM02: \cite{blanes02psp} obtain one BAB 4th-order method with $s=6$ and two 6th-order methods: one BAB with $s=11$ and one ABA with $s=14$.
	\item OMF03: \cite{omelyan03sai} propose 4th-order methods of type ABA and BAB with $s=4$ and $s=5$.
	\item BCE22: \cite{blanes22rkn} obtain optimized 8th-order methods with $s=17,18$ and $19$ of type ABA and BAB.
\end{itemize}

As we have seen in Subsection \ref{nested-commut},
the particular structure of problem $y^{\prime\prime} = g(y)$ allows one to include  the flows corresponding to
nested commutators in the previous compositions (\ref{ab-ba})-(\ref{aba-bab}). This may on one hand lead to
more efficient schemes, and on the other hand allows one to achieve order 4 with only positive coefficients. 

\begin{table}
  \centering
\begin{tabular}{|l|l|l|} \hline
\multicolumn{3}{|c|}{RKN splitting methods with commutators}\\ \hline
 $r=4$ & $r=6$  & $r=8$  \\ \hline
{\bf 2S}$_{\mathrm{ABA,BAB}}$
& {\bf 4,\framebox[3mm]{5}S}$_{\mathrm{ABA,BAB}}$
 & {\bf 11S}$_{\mathrm{ABA,BAB}}$-OMF02  \\
    	-KOS93,CHI97 	& -OMF02 		&  \\
       	 {\bf 4S}$_{\mathrm{ABA,BAB}}$-SUZ95 		&   &  \\
          	 {\bf 3,\framebox[3mm]{4},5S}$_{\mathrm{ABA,BAB}}$		&     &       \\
          \ \ 	-CHI06,OMF02				&     &       \\
 \hline\hline
  {\bf P:\framebox[3mm]{1}S}$_{\mathrm{BAB}}$-TI84,ROW91,		&  {\bf P:\framebox[3mm]{3}S}$_{\mathrm{ABA,BAB}}$
			& {\bf P:4S}$_{\mathrm{ABA}}$-BCR01a			\\
    WHT96,BCR99 		&-BCR99  & {\bf P:5S}$_{\mathrm{BAB}}$			\\ 
 {\bf P:2S}$_{\mathrm{BAB}}$-LSS97	&  & \quad -BCR01a,BCR01b	\\ \hline
\end{tabular}
  \caption{RKN splitting methods of order $r=4,6,8$ with modified potentials. Schemes are coded as in Table
  \ref{tableRKN}.}\label{tableRKNm}
\end{table}

In Table~\ref{tableRKNm} we present some relevant methods within this class, both processed and non-processed. They correspond to the following:

\begin{itemize}
	\item TI84: \cite{takahashi84mco} obtain the kernel of a 4th-order processed method with $s=1$. 
	\item ROW91: \cite{rowlands91ana} rediscovers the kernel of the 4th-order processed method with $s=1$. 
	\item KOS93: \cite{koseleff93cfp} obtains the first 4th-order BAB method with positive coefficients with $s=2$ (but with a misprint in one of the coefficients).
	\item SUZ95: \cite{suzuki95hep} gives 4th-order ABA and BAB methods with $s=3$ 
	\item WHT96: \cite{wisdom96sco} obtain one 4th-order processed method with $s=1$.
	\item LSS97: \cite{lopezmarcos97esi} obtain 4th-order processed methods with $s=2$.
	\item CHI97: \cite{chin97sif} rediscovers the 4th-order method of \cite{koseleff93cfp} (but now with the correct coefficients) and shows its efficiency in practice.
	\item BCR99: \cite{blanes99siw} obtain the following processed methods: one 4th-order method with $s=1$, one non-symmetric 5th-order method with $s=2$ and one 6th-order method with $s=3$.
	\item BCR01a: \cite{blanes01hor} present one 6th-order method with $s=2$ and higher order commutators and two 6th-order methods with $s=3$.
	\item BCR01b: \cite{blanes01nfo} propose one optimized 6th-order method with $s=3$ and one 8th-order method with $s=5$.
  \item OMF03:  \cite{omelyan03sai} present a detailed treatment of 4th-order methods with $2\leq s\leq5$ and 6th-order methods with $s=4,5$ both ABA and BAB and different number of modified potentials.
     \item CHI06: \cite{chin06cco} obtains families of 4th-order schemes with coefficients written analytically in terms of free parameters for $s>2$.
\end{itemize}

\subsection{Methods for near-integrable systems}
\label{mfnis}

As we have seen in Subsection \ref{near-integ}, splitting methods are outstanding for near-integrable systems of the form $x' = f_1(x) + \varepsilon f_2(x)$: 
the presence of
two parameters, $h$ (the step size) and $\varepsilon$ (the size of the perturbation), allows one to reduce the number of order
conditions and still construct highly efficient schemes. There are problems where, in addition, one can use the same techniques as for RKN methods and also incorporate the flows of nested commutators into the algorithm. This is the case, in particular, for Hamiltonian systems
of the form $H = H_1 + \varepsilon H_2$, when $H_1$ is the sum of 2-body Kepler problems and $H_2$ depends only on coordinates. If we denote 
by $F_1$ the operator associated with $H_1$, and by $F_2$
the operator associated to the perturbation $H_2$, then $[F_2,[F_2,[F_2,F_1]]] =  0$, and the analysis done in  Subsection \ref{subsec3.1}
can be applied here (with the obvious interchange $F_1 \leftrightarrow F_2$).

\begin{table}
  \centering
\begin{tabular}{|l|l|l|} \hline
\multicolumn{3}{|c|}{Splitting methods of generalized order}\\ \hline
 $(n,2)$ & $(n,4)$ & $(n,6,4)$    \\ \hline
{\bf \framebox[3mm]{1}}(2,2){\bf S}-WH91 
 & {\bf 4}(6,4){\bf S}$_{\mathrm{BAB}}$-MCL95a 
  & {\bf \framebox[3mm]{7}}(8,6,4){\bf S}$_{\mathrm{ABA}}$ \\
 {\bf \framebox[3mm]{n}}(2n,2){\bf S}$_{\mathrm{ABA,BAB}}$ & {\bf \framebox[3mm]{5}}(8,4){\bf S}$_{\mathrm{ABA,BAB}}$-MCL95a \hspace*{-0.2cm}
 & \quad -BCF13   \\ 
 \ \ -MCL95a,LR01 & {\bf 7}(10,4){\bf S}$_{\mathrm{ABA}}$-BCF13 & 
 {\bf \framebox[3mm]{8}}(10,6,4){\bf S}$_{\mathrm{ABA}}$ \\
  	&  
	  \quad   &    \quad -BCF13  \\ \hline\hline
 {\bf P:1}(17,2)-WHT96 \hspace*{-0.2cm}
 &  &  {\bf P:\framebox[3mm]{3}}(7,6,4){\bf S}$_{\mathrm{ABA}}$-BCR00 
 \\ \hline
  & {\bf P:\framebox[3mm]{2}}(6,4){\bf S}$_{\mathrm{AB}}$-BCR00 
	& {\bf P:\framebox[3mm]{3}}(7,6,5){\bf S}$_{\mathrm{AB}}$-BCR00 
		\\ \hline
  & {\bf P:1}(6,4){\bf S}$^*_{\mathrm{ABA}}$-BCR00 
	& {\bf P:2}(7,6,5){\bf S}$^*_{\mathrm{AB}}$-BCR00 
		\\
  & {\bf P:n}(n,4){\bf S}$^*_{\mathrm{ABA}}$-LR01
	&   \\ \hline
\end{tabular}
  \caption{Splitting methods for near-integrable systems of the form $x' = f_1 + \varepsilon f_2$.
For processed methods we also include methods applicable when
$[F_2,[F_2,[F_2,F_1]]] =  0$  (the last two rows of processed methods) and schemes with modified flows (labelled $^*$).}
\label{tableNI}
\end{table}

In Table \ref{tableNI} we separate, as usual, non-processed from
processed schemes (preceded by {\bf P}). We also
include methods when 
$[F_2,[F_2,[F_2,F_1]]] =  0$
 (second row) and schemes
with nested commutators (denoted by a label $^*$) in the last two rows of processed methods. In more detail,
\begin{itemize}
	\item WH91: \cite{wisdom91smf} present the first $(2,2)$ method.
	\item MCL95a: \cite{mclachlan95cmi} gives families of symmetric $(2s,2)$ schemes of type ABA and BAB with $s\leq 5$ and positive coefficients, 
	one BAB (6,4) method with $s=4$ and two (8,4) ABA and BAB schemes with $s=5$.
	\item WHT96: \cite{wisdom96sco} present one processed $(k,2)$ method with $k=16$ and $s=1$.
	\item BCR00: \cite{blanes00psm} give several processed methods: (i) one non-symmetric (6,4,3) method with $s=2$, one (7,6,4) method with $s=3$, one (7,6,5,4) method with $s=4$, and (ii) several RKN methods, one (6,4) with $s=2$ and one (7,6,5) with $s=3$ both non-symmetric. In addition, with modified potentials, 
	they give one (6,4) with $s=1$ and a non-symmetric (7,6,5) with $s=2$.
	\item LR01: \cite{laskar01hos} carry out a systematic study of $(2s,2)$ methods, and obtain new schemes up to $s=10$ with positive 
coefficients. The order of these methods is increased to  $(2s,4)$ by including an appropriate modified potential.
	\item BCF13: \cite{blanes13nfo} present several schemes of type ABA: one (10,4) with $s=7$, one (8,6,4) with $s=7$ and one (10,6,4) with $s=8$.
\end{itemize}



\section{Splitting everywhere: some relevant applications}
\label{sect9}

\subsection{Splitting and extrapolation}

Extrapolation methods constitute a class of efficient high order schemes for solving initial value problems when the local error of the 
basic integrator has an asymptotic expansion containing only even powers of the step size $h$. In this subsection we review how some of the
previous techniques used in the construction and analysis of splitting and composition methods can also be applied to 
simplify the analysis of extrapolation.

\subsubsection{The Gragg/GBS method}

Given the initial value problem
$x' = f(x)$, $x(t_0) = x_0$,
Gragg proposed the following algorithm to produce the quantity $T_{h^*}(t)$. Denoting 
$t=t_0+2nh^*, \ t_i=t_0+ih^*$, it reads
\begin{equation} \label{eq.1.Gragg}
\begin{aligned}
		& x_{1} = x_0+h^*f(x_0)  \\
		& x_{i+1} = x_{i-1}+2h^*f(x_i), \qquad i=1,2,\ldots,2n   \\
		& T_{h^*}(t) =\frac14(x_{2n-1}+2x_{2n}+x_{2n+1}).  
\end{aligned}
\end{equation}
Subsequently, he proved that $T_{h^*}(t)$ possesses an asymptotic expansion
in even powers of $h^*$ \cite{gragg65oea}, so that
it can then be used for Richardson extrapolation.  The original proof was simplified by Stetter when he realized that (\ref{eq.1.Gragg}) can be interpreted as a one-step algorithm \cite{stetter70sts}.
Later on, \cite{hairer93sod} showed that this scheme
is consistent with the differential equation in the extended phase space
\begin{equation} \label{eq.3.StrangAut2}
	\begin{array}{ll}
u'= f(v), \qquad \qquad &u(t_0)=u_0 = y_0\\
v'= f(u), \qquad \qquad &v(t_0)=v_0 = y_0,
\end{array}
\end{equation}
whose exact solution is $u(t)=v(t)=x(t)$. In fact, to prove that $T_{h^*}(t)$ only contains even powers of $h^*$ is trivial by noticing that system
(\ref{eq.3.StrangAut2}) can be expressed as
\begin{equation}\label{eq.3.ODE}
		\frac{d}{dt} \left(
		\begin{array}{c} u \\  v   \end{array} \right)=
		 \left(\begin{array}{c} 0  \\ f(u)  \end{array} \right)	
		+
			\left(\begin{array}{c} f(v)  \\ 0 \end{array} \right) 
=
 f_1(u) + f_2(v), \qquad
\left(
		\begin{array}{c} u(t_0)  \\ v(t_0)   \end{array} \right)=	
\left(
		\begin{array}{c} x_0  \\ x_0  \end{array} \right),
\end{equation}
and that the Strang splitting with step size $h = 2 h^*$ gets
\[
\begin{aligned}
  & v_{1/2}=v_{i-1}+\frac{h}2f(u_{i-1}) \\
  &  u_{i}=u_{i-1}+hf(v_{1/2}), \qquad\qquad i=1,2,\ldots,n \\
 & v_{i}=v_{1/2}+\frac{h}2f(u_{i}).
\end{aligned}
\]
Then, clearly, $T_{h^*}=\frac12 (u_n+v_n)$, and the result is obtained by noticing that the Strang splitting is time-symmetric.

On the other hand, application of the Gragg method (\ref{eq.1.Gragg}) to the second order equation
\begin{equation}\label{eq.1.ODE2}
		y''=g(y),  \qquad  y(t_0)=y_0, \quad y'(t_0)=y'_0,
\end{equation}
produces the equivalent formula \cite{hairer93sod}
\begin{equation} \label{eq.2.StormerV}
		 y_{i+1}-2y_i +y_{i-1} = h^2g(y_i),  \qquad i=1,2,3,\ldots, 
\end{equation}
again with $h = 2 h^*$. As is well known, 
by considering the equivalent first-order system
\begin{equation} \label{eq.1.ODE2b}
		\frac{d}{dt} \left(
		\begin{array}{c} y \\ y' \end{array} \right)=
		 \left(\begin{array}{c} y' \\ 0 \end{array} \right)+
		 \left(\begin{array}{c} 0 \\ g(y) \end{array} \right)
		=f_1(y')+f_2(y),
\end{equation}
and then applying the Strang splitting $\varphi_{h/2}^{[2]} \circ \varphi_h^{[1]} \circ	\varphi_{h/2}^{[2]}$:
\[
\begin{aligned}
 & y_{n+1}=y_{n}+hy_{n}'+\frac{h^2}2g(y_n) \\
 & y'_{n+1}=y'_n+\frac{h}2(g(y_n)+g(y_{n+1})),
\end{aligned}
\]
reproduces 
St\"ormer's rule (\ref{eq.2.StormerV}) after simplification.
 In fact, one can also apply 
the other variant of the Strang splitting, i.e.,
$\varphi_{h/2}^{[1]} \circ \varphi_h^{[2]} \circ	\varphi_{h/2}^{[1]}$, to (\ref{eq.1.ODE2b}), thus yielding
\begin{equation} \label{eq.2.StormerV2}
\begin{aligned}
 & y'_{n+1}=y'_n+hg(y_n+\frac{h}2y_n') \\
 & y_{n+1}=y_{n}+\frac{h}2(y_{n}'+y_{n+1}').
\end{aligned} 
\end{equation}
This (slightly cheaper) scheme could also be used as the basic method for extrapolation, as for the St\"ormer rule, but since $f_1(y')$ and $f_2(y)$
may have different qualitative properties, the resulting schemes would also be different.

\subsubsection{Multi-product expansions}

Linear combinations of splitting methods
\[
   \sum_{k} c_k \, \big(\varphi_{b_{k,i_k} h}^{[2]} \circ \varphi_{a_{k,i_k} h}^{[1]} \circ \cdots \circ \varphi_{b_{k,1} h}^{[2]} \circ \varphi_{a_{k,1} h}^{[1]} \big)
\]
offer much freedom in the choice of parameters $\{c_k, a_{k,i_j}, b_{k,i_j} \}$ to approximate the flow of $x' = f_1(x) + f_2(x)$, although they no longer
provide approximations with the same preservation properties as the original system. As shown in
\cite{sheng89slp}, if  all  $\{a_{k,i_j}, b_{k,i_j} \}$ were to be positive, then any individual composition has to be at most of second order, and some
coefficients $c_k$ must be negative. The choice of the Strang splitting $S_h^{[2]}$ (or more generally a time-symmetric 2nd-order integrator) as the basic scheme in the composition is particularly relevant due to the structure of the truncation error. Suppose $\{k_1, k_2, \ldots, k_m \}$ denotes a given set of $m$ integer numbers, and form the linear
combination
\begin{equation} \label{mp.1}
  \phi_h \equiv \sum_{i=1}^m c_i \left( S^{[2]}_{h/k_i} \right)^{k_i}.
\end{equation}  
Then, $\phi_h$ furnishes an approximation of order $2m$ if the coefficients satisfy the linear equations
\[
  G_0 = \sum_{i=1}^m c_i = 1, \qquad G_{2+2j} = \sum_{i=1}^m \frac{c_i}{k_i^{2 + 2j}} = 0, \qquad j=0,1, \ldots, m-2,
\]
with solutions
\[
  c_i = \prod_{j=1 (\ne i)}^m \frac{k_i^2}{k_i^2 - k_j^2}.
\]  
Schemes of the form (\ref{mp.1}) are referred to as multi-product expansions \cite{chin11mpo} and allow one to construct in an easy way 
efficient high-order integrators.  
In general, if one starts with a basic time-symmetric geometric method of order $2n$, the corresponding linear combination (\ref{mp.1}) provides a
scheme of order $2n + 2 \ell$, $\ell = 1, \ldots, n$, which, quite remarkably,
preserves the geometric properties of the system (e.g., symplecticity for Hamiltonian problems) 
up to order $4n+1$ or higher \cite{blanes99eos,chan00eos}. This can be shown by appropriately rewriting  the differential operator associated with (\ref{mp.1})
as a product of exponentials of operators \cite{blanes16aci}.

\subsection{Crouch--Grossman methods}

The nonlinear differential equation
\begin{equation}\label{eq.Crouch-Grossman0}
	Y' = A(Y) \, Y, \qquad Y(0)=Y_0 \in \mathcal{G},
\end{equation}
where $A\in\mathbb{R}^{d\times d}$ and $\mathcal{G}$ is a matrix Lie group, 
appears in relevant physical fields such as rigid
body mechanics, in the calculation of Lyapunov exponents ($\mathcal{G}
\equiv \mathrm{SO}(d)$) and other problems arising in 
Hamiltonian dynamics ($\mathcal{G} \equiv \mathrm{Sp}(d)$). In fact, it can be shown
that every differential equation evolving on a matrix Lie group
$\mathcal{G}$ can be written in the form (\ref{eq.Crouch-Grossman0}) \cite{iserles00lgm}. 

A class of methods providing by construction approximations in $\mathcal{G}$ was proposed in \cite{crouch93nio}
by appropriately modifying Runge--Kutta schemes as follows (see also \cite[sect. IV.8]{hairer06gni}).

Let $b_i,a_{i,j}, \ (i,j=1,\ldots,s)$ be real numbers. Then, an explicit $s$-stage Crouch--Grossman method is given by 
\begin{equation}\label{eq.Crouch-Grossman}
  \begin{array}{cclcl}
        Y^{(1)}  & =  & Y_n, &  & K_1=A(Y^{(1)})  \\
	& \vdots & & &\\
	Y^{(j)} & = & \e^{ha_{j,j-1}K_{j-1}}\cdots \,
	 \e^{ha_{j,1}K_{1}} \, Y_n, &   &K_j=A(Y^{(j)}), \qquad 2 \le j \le s  \\
	Y_{n+1} & = & \e^{hb_{s}K_{s}}\cdots	\, \e^{hb_{1}K_{1}} Y_n &  & 
  \end{array}	
\end{equation}
As usual, the order conditions to be satisfied by the coefficients $b_i,a_{i,j}, \ (i,j=1,\ldots,s)$ can be found by comparing the Taylor series expansions of the
exact and numerical solutions. It has been shown that the order conditions for classical Runge--Kutta methods form 
a subset of those for the Crouch--Grossman methods \cite{owren99rkm}. Notice that the $s$-stage scheme (\ref{eq.Crouch-Grossman}) involves a total of
$s(s+1)/2$ matrix exponentials.

In \cite{crouch93nio} there are several 3-stage methods of order 3, whereas in \cite{owren99rkm} a 4th-order scheme with $s=5$ stages 
(or equivalently, 15 exponentials) is presented. Achieving higher order within this approach is by no means simple, and the resulting methods require a large
number of exponentials. It turns out, however, that splitting can also be used here to construct schemes of any order with a reduced number of matrix
exponentials. To illustrate the technique, we consider the separable system in the extended phase space
\begin{equation} \label{CG-extended}
		\frac{d}{dt} 
		\left(
		\begin{array}{c} U \\ V  \end{array} \right)=
		\left(
		\begin{array}{c} A(V) \, U  \\ 0 \end{array} \right)+
		\left(
		\begin{array}{c}  0 \\ A(U) \, V  \end{array} \right),
\end{equation}
to which we apply the Strang splitting. This leads to
\begin{equation} \label{eq.3.LieGroup2}
\begin{array}{rcl}
 U_{1/2}&=&e^{{\frac{h}2}A(V_n)} \, U_n \\
 V_{n+1}&=&e^{{h}A(U_{1/2})} \, V_n\\
 U_{n+1}&=&e^{{\frac{h}2}A(V_{n+1})} \, U_{1/2}.
\end{array}
\end{equation}
Notice that (\ref{eq.3.LieGroup2}) is simply a particular example of \eqref{eq.Crouch-Grossman} with $s=3$, i.e., 
\begin{equation} \label{eq.Crouch-Grossman1}
  \begin{array}{cclcl}
        Y^{(1)} & = & Y_n, & & K_1=A(Y^{(1)})  \\
	Y^{(2)} &=& \e^{ha_{2,1}K_{1}} \, Y_n, & & K_2=A(Y^{(2)})  \\
	Y^{(3)} &=& \e^{ha_{3,2}K_{2}} \, \e^{ha_{3,1}K_{1}} \, Y_n,  & & K_3=A(Y^{(3)})  \\
	Y_{n+1} &=& \e^{hb_3 K_{3}}\,\e^{hb_2 K_{2}} \,\e^{hb_1 K_{1}} \, Y_n & & 
 \end{array}	
\end{equation}
with
\[
  a_{2,1}=\frac12, \quad 
  a_{3,2}=1, \quad
  a_{3,1}=0, \quad 
  b_{1}=\frac12, \quad 
  b_{2}=0, \quad 
  b_{3}=\frac12.
\]
Since $a_{3,1}=b_{2}=0$ and $a_{2,1}=b_{1}$, then the whole scheme can be computed with only three exponentials, instead of 6 by applying the
recursion \eqref{eq.Crouch-Grossman}. Higher order methods can be obtained in the same way by applying appropriate splitting schemes to the enlarged
system (\ref{CG-extended}). Thus, in particular, it is possible to get a Crouch--Grossman scheme of order 4 with 7 exponentials.

\subsubsection{Positivity-preserving splitting methods}

Differential equations $x' = f(x)$ modeling population dynamics must satisfy 
positivity and frequently mass preservation, that is,
\begin{eqnarray*}
 \mbox{(i)} &  & x_i(0) \ge 0  \quad \Rightarrow \quad x_i(t) \ge 0 \;\; \forall \, t, \; i=1,\ldots,D \quad \mbox{(\emph{positivity preservation})} \\
 \mbox{(ii)} &  & {\bf 1}^T{x}(t)={\bf 1}^T{x}(0), \quad \mbox{ with } \quad {\bf 1}=(1,\ldots,1)^T \qquad\quad \mbox{(\emph{mass conservation})}
 \end{eqnarray*}
Both of them are automatically fulfilled if the differential equation can be written as $x' = A(x) x$, where $A$ is a matrix such that
\begin{itemize}
 \item if $x_i \ge 0$, $i=1,\ldots, D$, then $A(x)_{k,\ell}\geq0$ when $k\neq\ell$,  and $A_{k,k}\leq0$, for $k,\ell=1,\ldots,D$ (to guarantee
 positivity)
 \item $\sum_{k=1}^D A_{k,\ell}=0$ for $\ell=1,\ldots,D$ (to ensure mass conservation).
\end{itemize} 
 In such situations, it is of course expedient to use a numerical integrator that also satisfies these features at the discrete level when carrying out
 simulations. Unfortunately, there is also an order barrier here for two of the most used families of schemes, namely Runge--Kutta and multistep methods:
 if they unconditionally preserve positivity, then 
they cannot be better than first order \cite{bolley78cdl}.

It turns out, however, that splitting methods are able to overcome this order barrier. To see this point, again consider the problem in the extended
phase space
\begin{eqnarray} \label{ppm}
 {y}' &\!\!\!=\!\!\!& A({z}){y}, \qquad {y}(0)={y}_0={x}_0, \nonumber \\
 {z}' &\!\!\!=\!\!\!& A({y}){z}, \qquad {z}(0)={z}_0={x}_0,
\end{eqnarray}
where  ${x}(t)={y}(t)={z}(t)$. Now, as in the previous subsection, we apply the Strang splitting to the system (\ref{ppm}), i.e.,
\begin{eqnarray}
 {y}_{1/2} &\!\!\!=\!\!\!& \e^{\frac12 hA({z}_0)}{y}_0, \nonumber \\
 {z}_1 &\!\!\!=\!\!\!& \e^{hA({y}_{1/2})}{z}_0,   \label{eq:SplittingMethod} \\
 {y}_{1} &\!\!\!=\!\!\!& \e^{\frac12 hA({z}_1)}{y}_{1/2}.   \nonumber 
\end{eqnarray}
Then, both ${y}_1$ and ${z}_1$ are 2nd-order approximations to the exact solution $x(h)$ preserving mass and positivity. This can be easily checked
by realizing that ${\bf 1}^T A=0$ and so ${\bf 1}^T \e^Ax={\bf 1}^\top x$. On the other hand, for $h > 0$, if $x_i \ge 0$ then $(\e^{hA(x)})_{jk} \ge 0$
for all $i,j,k$.

If $h < 0$, positivity is no longer preserved with this approach (see \cite{blanes22ppm} for a more detailed treatment of this topic).

\subsection{Exact splitting}

There are equations for which a clever decomposition of the vector field in fact provides the exact solution. The resulting factorization is called 
\emph{exact splitting}
and, in the case of partial differential equations, allows one to get very accurate approximations once each subsystem is solved by pseudo-spectral
methods and pointwise multiplication. Early applications of this idea to the time integration of the Gross--Pitaevskii equation are contained in
\cite{chin05foa,bader11fmf}, and a convenient generalization to a large class of PDEs, called inhomogeneous quadratic differential equations, has been recently proposed in \cite{bernier21esm2}. These PDEs are of the form
\begin{equation} \label{eq:exs1}
  u_t(t,x) = -p^{w} u(t,x), \qquad u(0,x) = u_0(x), \qquad\qquad t \ge 0, \; x \in \mathbb{R}^d,
\end{equation}  
where $d \ge 1$, $u_0 \in L^2(\mathbb{R}^d)$ and $p^w$ is the quadratic differential operator 
\[
  p^w = (x, \,  -i \nabla_x) \, Q \left( \begin{array}{c}
  					x \\
					-i \nabla_x
				\end{array} \right) + Y^T 	\left( \begin{array}{c}
  					x \\
					-i \nabla_x
				\end{array} \right)  + c.
\]
Here $Q$ is a constant $2d \times 2d$ symmetric matrix with complex coefficients, $Y \in \mathbb{C}^{2d}$  and $c \in \mathbb{C}$ is constant.
If $(-p^w)$ is such that the real part of the polynomial $p(X) = X^T Q X + Y^T X + c$ 
is bounded from below on $\mathbb{R}^{2d}$, then it generates a strongly continuous semigroup, denoted 
$\e^{-t p^w} u_0$. It was shown in \cite{bernier21esm2} that $\e^{-t p^w}$ can be factorized as 
a product of operators that can be exactly evaluated. Notice that
some relevant equations are of the form (\ref{eq:exs1}), so this technique allows one to construct numerical methods that are spectral in space and
exact in time using only a small number of fast Fourier transforms. This has been illustrated in 
\cite{bernier21esm} in the case of kinetic and nonlinear Schr\"odinger equations.

As an illustrative example of this technique, let us consider the Schr\"odinger equation with a quadratic potential
\[
  i\, u_t = \frac12 (-\partial_{xx} + x^2) \, u   \equiv \frac{1}{2}(P^2+X^2) \, u. 
\]
Then, for $|h|<\pi$, one has
\[
u(h,x)=\e^{-ih\left( P^2/2+X^2/2\right)} \, u(0,x) =
  \e^{-if(h)\left(X^2/2 \right)} \, \e^{-ig(h)\left( P^2/2\right)} \, 
  \e^{-if(h)\left(X^2/2 \right)} \, u_0(x), 
\]
with 
\[
   f(h)= \frac{1-\cos(h)}{\sin(h)}, \qquad \qquad g(h)=\sin(h),
\]
so the action of each exponential can be easily evaluated.

\subsection{Gravitational $N$-body problem}

The long-time numerical integrations of the whole Solar System carried out in \cite{sussman92ceo,wisdom91smf},
suggesting its chaotic nature, represented a major boost for the development of new and efficient symplectic algorithms. The integrator used in these simulations, referred to in the dynamical astronomy literature as the Wisdom--Holman map, is simply the Strang splitting method applied to
 the Hamiltonian system 
\begin{equation} \label{grav.1}
 \hat H(\hat{q},\hat{p}) = H_1(\hat{q},\hat{p}) + H_2(\hat{q}),
\end{equation}
where $H_1$ and $H_2$ are given by (\ref{jacobi.1})) in terms of Jacobi coordinates, as described in subsection~\ref{subsec.1.4}. Here 
$H_1$ represents the $N-1$ independent Keplerian motions of the planets (which can be efficiently solved) and $H_2$ accounts for the
gravitational interactions among the planets. Moreover, as shown in subsection~\ref{subsec.1.4},  $H_2(\hat q) = \mathcal{O}(\varepsilon)$, with $\varepsilon\equiv \frac{1}{m_0} \max_{1\leq i \leq N-1} m_i$. For the Solar System, $\varepsilon \approx10^{-3}$. In consequence, 
splitting methods for near-integrable systems, such as those presented
in subsections \ref{near-integ} and \ref{mfnis}, are particularly appropriate. In fact, the Wisdom--Holman map, in the terminology of
Subsection~\ref{near-integ}), corresponds to a method of generalized order $(2,2)$.

Very often, splitting methods are used  when very long-time integrations are involved 
and getting high accuracy is not of particular concern,
but rather when having good behavior of the error over the whole integration period is the crucial point. For instance, in the situation analyzed in
\cite{sussman92ceo}, the important issue was to decide whether the motion is regular or chaotic, and not to compute accurately the state of the Solar System
after a long time.

The experiment we describe next illustrates that even when short-time integrations are involved and very high accuracy is required, splitting methods can outperform
standard integrators once the particular structure of the problem is incorporated into their very formulation. 

\begin{figure}[htb]
\centering
\includegraphics[scale=0.47]{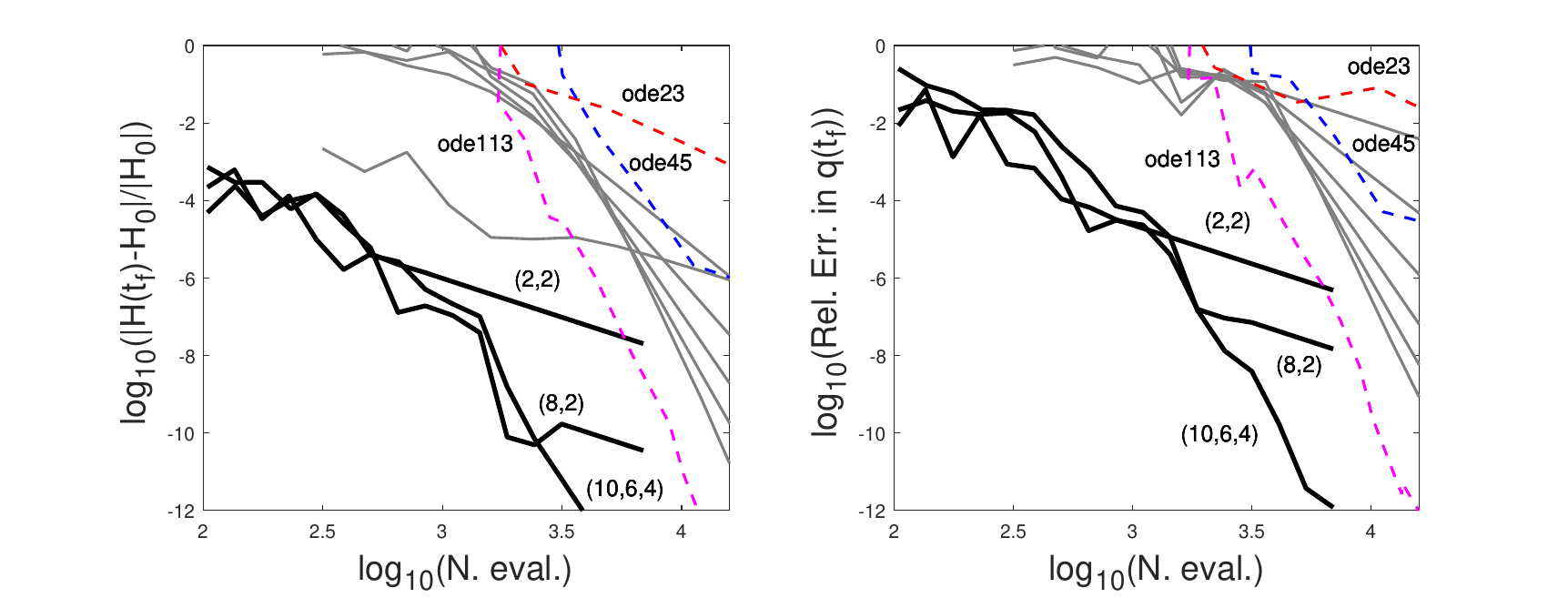} \\
\caption{Outer Solar System. Efficiency diagrams showing the relative error in energy  (left panel) and in positions $q$ (right panel) at the final time $t_f = 2 \cdot 10^5$ days vs. the number of force evaluations, obtained with different schemes: extrapolation of order $r=2,4,6,8,10,12$ (grey curves), MATLAB routines \texttt{ode23}, \texttt{ode45} and \texttt{ode113},
and splitting methods of generalized order (2,2), (8,2) and (10,6,4).}
\label{fig_NBody_Efficiency}
\end{figure}

Specifically, we again deal with the simplified model of the outer Solar System addressed in Subsection \ref{subsec.1.4}, with the same initial conditions and
an integration interval of $2 \times 10^5$ days, and compare splitting methods specially tailored for near-integrable systems with some popular standard
integrators. Thus, Figure \ref{fig_NBody_Efficiency} shows efficiency diagrams for the error in energy (left) and in positions (right) at the final time vs. the number of force evaluations, obtained with three different
splitting methods of generalized orders (2,2), (8,2) and (10,6,4) in comparison with extrapolation methods of orders 2-12 (grey lines) and the MATLAB routines
\texttt{ode23}, \texttt{ode45} and \texttt{ode113}.

As stated before, extrapolation methods are particularly appropriate for second-order differential equations when high accuracy is required
 \cite{hairer93sod}. They take the symmetric second order scheme \eqref{eq.2.StormerV} as the basic method in the linear combination
\eqref{mp.1} with the harmonic sequence $k_i=i, \ i=1,2,\ldots,m$, to get a method of order $r=2m$ at the cost of $m(m+1)/2+1$ 
evaluations of the force\footnote{One evaluation can be saved if the basic scheme \eqref{eq.2.StormerV2} is taken instead. This is the actual implementation
we use in our experiments here.}. On the other hand, \texttt{ode23} and \texttt{ode45} solve the ODE with variable time step using embedded 
Runge--Kutta methods of orders 2-3 and 4-5, respectively, whereas \texttt{ode113} does it with variable time step and variable order 
using multistep methods up to order 13. 
Due to the smoothness of the problem and the low cost of evaluating the highest order, the latter procedure provides better results than extrapolation methods of
high order, as shown in the figure. Nevertheless, none of them can compete with the specially designed splitting methods for this particular problem at
all accuracies. This is so even when a very short integration interval is considered. In fact, this relative superiority will only increase for longer integration intervals
due to their different error propagation mechanism.

\subsection{Molecular simulations}

In the simulation of the dynamics of large molecules two classes of algorithms are commonly used: (i) molecular dynamics and (ii) those based on stochastic differential equations and the Monte Carlo method.

In classical molecular dynamics, the motion of the atoms in the molecule is determined by integrating Newton's second law,
\[
   M q^{\prime\prime} = -\nabla_q U(q),
\]
where $q$ is a vector containing all positions (in Cartesian coordinates), $M$ is a diagonal matrix whose elements are the atomic masses and $U(q)$ is
the (empirical) potential function modeling the inter-particle interactions. Popular models include the Lennard--Jones and Morse potentials, although
many other choices are being used for different molecules (see e.g. \cite{leimkuhler15md} for more details). In any case, the number of atoms may be
exceedingly large (up to 100000) and the resulting system is highly nonlinear, exhibiting sensitive dependence on perturbations. In addition, the initial
velocities are typically assigned randomly, so  there is no point in trying to obtain accurate trajectories. Thus, although in principle 
Runge--Kutta--Nystr\"om splitting methods can be used in this setting, the method of choice in practice is the St\"ormer--Verlet scheme, given its
excellent stability properties and low computational cost \cite{leimkuhler96imf}. 
Other variants include the use of multiple time-stepping (and in particular the mollified impulse
method \cite{garcia-archilla99lts}), and symplectic schemes specially designed for Hamiltonian systems with constraints, such as SHAKE 
\cite{ryckaert77nio} and RATTLE \cite{andersen83rav}. 
   
To take into account the influence of a medium (say, a solvent or air) on a molecular system, typically as external random impacts, a common practice
consists in extending the molecular dynamics approach by incorporating a stochastic component in the equations of motion. The resulting system is
referred to as Langevin dynamics and  is described by the stochastic differential equation
\begin{equation} \label{lange.1}
  \begin{aligned}
    & dq = M^{-1} p dt \\
    & dp = - \nabla_q U(q) dt - \gamma p dt + \sqrt{2 \gamma k_B T} M^{1/2} dW
  \end{aligned}
\end{equation}
in terms of coordinates $q \in \mathbb{R}^d$ and momenta $p \in \mathbb{R}^d$. 
Here $\gamma > 0$ is the friction coefficient or collision frequency, $k_B$ is the Boltzmann constant, $T$ is the temperature and $W$ is a Wiener process.
We may regard (\ref{lange.1}) as modeling a system of particles immersed in a fluid bath consisting of many particles giving rise to much weaker
interactions than those modeled by the potential $U$ \cite{leimkuhler15md}.

As with deterministic equations, splitting methods can be constructed to deal with Langevin dynamics. One possible approach consists in decomposing
(\ref{lange.1}) into three parts as follows:
\begin{equation} \label{lange.2}
  \left( \begin{array}{c}
 			dq \\
			dp
	  \end{array} \right) = 		
\underbrace{ \left( \begin{array}{c}
 			M^{-1} p \\
			0
	  \end{array} \right) dt }_{A} + 
	  \underbrace{ \left( \begin{array}{c}
 			0 \\
			-\nabla_q U(q)
	  \end{array} \right) dt }_{B} + 
	  \underbrace{ \left( \begin{array}{c}
 			0 \\
			- \gamma p dt + \sigma M^{1/2} dW
	  \end{array} \right) }_{O}, 
\end{equation}
with $\sigma \equiv  \sqrt{2 \gamma k_B T}$. Each of the parts can be exactly solved: the solution of the first is a \emph{drift} in position,
\[
  \varphi_h^{[A]}(q,p) = (q + h M^{-1} p, p),
\]  
the second corresponds to a \emph{kick} in momentum,
\[
  \varphi_h^{[B]}(q,p) = (q,  p - h \nabla_q U(q)),
\]  
and the third piece defines an Ornstein--Uhlenbeck process in $p$ which can be exactly sampled:
\[
  \varphi_h^{[O]}(q,p) = (q, \e^{-\gamma h} p + \sqrt{k_B T (1 - \e^{-2 \gamma h})} M^{1/2} R),
\]
where $R$ is a vector of independent and identically distributed normal random numbers. In \cite{leimkuhler13rco,leimkuhler13rae}, several splitting methods
constructed along these lines are designed and tested on numerical experiments. They are denoted by the acronym resulting from concatenating
the previous symbols $A$, $B$ and $O$. Thus, BAOAB corresponds to
\[
   \varphi_{h/2}^{[B]} \circ   \varphi_{h/2}^{[A]} \circ \varphi_{h}^{[O]} \circ \varphi_{h/2}^{[A]} \circ \varphi_{h/2}^{[B]}, 
\]
and so on. The analysis carried out in \cite{leimkuhler13rco,leimkuhler13rae,leimkuhler16tco} for the large friction limit, and in \cite{alamo16atf} for the general case,
shows that BAOAB offers an improved behavior with respect to other members of this family. Other Langevin integrators are given in 
\cite{leimkuhler15md}.

Another approach to molecular simulations, rather than approximating Hamiltonian trajectories, consists in studying the paths originating from the collection
of all initial conditions within a given set. This perspective allows one to apply statistical mechanics for calculating averages.

Specifically, if $H(q,p)$ denotes the Hamiltonian function of the system in thermal equilibrium at temperature $T$, 
then the probability measure $\mu$ in $\mathbb{R}^{2d}$
with density
\begin{equation} \label{density.1}
  Z^{-1} \e^{-\beta H(q,p)}, \qquad \mbox{ where } \qquad Z = \int_{\mathbb{R}^{2d}} \e^{-\beta H(q,p)} dq dp < \infty
\end{equation}
is preserved by the flow of $H(q,p)$. Here $\beta = 1/(k_B T)$. The measure $\mu$ is called the Boltzmann--Gibbs distribution and $Z$ is the
partition function. Intuitively, $\mu$ provides the distribution of $(q,p)$ over an ensemble of many copies
of the given system when it is embedded in a much larger system acting as a ``heat bath'' at constant temperature. In other words, 
$Z^{-1} \e^{-\beta H(q,p)} dq dp$ represents the fraction of copies with momenta between $p$ and $p +dp$, and configuration between $q$ and
$q+ dq$ \cite{bou-rabee18gia}. Then, the average energy is given by
\[
  \tilde{E} = Z^{-1}  \int_{\mathbb{R}^{2d}} H(q,p) \, \e^{-\beta H(q,p)} dq dp.
\]
In general, $\tilde{E}$ cannot be obtained analytically, and even the use of numerical quadratures is unfeasible if the dimension $d$ is large. 
A common practice to approximate this class of integrals then consists in applying Monte Carlo methods, and
more specifically, Markov chain Monte Carlo methods \cite{brooks11homc}.

\subsection{Hamiltonian Monte Carlo}
\label{subsecHMC}

For Hamiltonian systems of the form
\begin{equation} \label{Ham-HMC}
  H(q,p) = \frac{1}{2} p^T M^{-1} p + U(q),
\end{equation}
we can factorize the (non-normalized) density (\ref{density.1}) as
\[
  \e^{-\beta H(q,p)} = \e^{-  \frac{1}{2} \beta  p^T M^{-1} p} \, \e^{-\beta U(q)},
\]
so that $q$ and $p$ are stochastically independent \cite{bou-rabee18gia}: the marginal distribution of the configuration variables $q$ has probability density
$\propto \e^{-\beta U(q)}$, whereas the momenta $p$ obey a Gaussian distribution with zero mean and covariance matrix $M$. Therefore, samples from
the $p$-marginal are easily obtained, and thus one may concentrate on generating samples with probability density
\begin{equation} \label{density.2}
  Z_q^{-1} \e^{- \beta U(q)}, \qquad \mbox{ where } \qquad Z_q = \int_{\mathbb{R}^d} \e^{-\beta U(q)} dq.
\end{equation}
This can be carried out with the so-called Hamiltonian (or hybrid) Monte Carlo method (HMC). In fact, HMC was proposed in the landmark paper
\cite{duane87hmc} not in the context of molecular simulation, but in lattice quantum chromodynamics. Later on, it was used in data science and
statistics \cite{liu08mcs,neal11muh}. More generally, HMC can be used to sample from any continuous probability
distribution on $\mathbb{R}^d$ for which the density function can be evaluated (perhaps up to an unknown normalizing constant): given a target
distribution $\Pi(q)$, if $U(q)$ denotes the negative logarithm of the (not necessarily normalized) probability function of the target, then it is clear that
$\Pi(q)$ is given by (\ref{density.2}) (with $\beta = 1$). HMC then generates samples $(q_i, p_i) \in \mathbb{R}^{2d}$ from the Boltzmann--Gibbs 
distribution corresponding to $H$, i.e., 
\begin{equation} \label{BG}
 P(q,p) = (2 \pi)^{-d/2} \, |\det M|^{-1/2} \e^{-  \frac{1}{2}   p^T M^{-1} p} \, Z_q^{-1} \, \e^{- U(q)},
\end{equation}
by means of a Markov chain so that $P(q,p)$ is an invariant of this chain. The corresponding marginal $q_i \in \mathbb{R}^d$ chain then leaves
invariant the target distribution $\Pi(q)$, and this allows one to estimate the multidimensional integral of a certain function $F$
 with respect to $\Pi$ by averaging $F$ at the points $(q_i)$ of the Markov chain.
 
 Since the HMC algorithm has been thoroughly reviewed in the excellent survey \cite{bou-rabee18gia}, 
 we only summarize here those aspects most closely related to splitting methods. The basic HMC procedure is described in Table \ref{table-HMC},
 although more elaborate possibilities exist. Notice that, when generating the proposal, the integrator to be used has to be both reversible 
 and volume-preserving. 
 Given the structure
 of $H$, it is clear that palindromic splitting methods constitute the natural choice for this task.

  \begin{table}
\hrule
\bigskip
{{}

Given $q^{(0)}\in\mathbb{R}^d$, $m_{\max} \geq 1$, set $m=0$.

\begin{enumerate}
\item (Momentum refreshment.) Draw $p^{(m)}\sim \mathcal{N}(0,M)$.

\item (Integration leg.) Compute $(q^*,p^*)$  ($q^*$ is the proposal) by  integrating, by means of
a reversible, volume-preserving integrator with step size \(h\), the equations of motion derived from the Hamiltonian  \eqref{Ham-HMC} over 
an interval  \(0\leq t\leq Nh\). The initial condition is \((q^{(m)},p^{(m)})\).

\item (Accept/reject.) Calculate \[a^{(m)} = \min\big(1, \exp(H(q^{(m)},p^{(m)})-H(q^*,p^*))\big)\]
and draw $u^{(m)} \sim \mathcal{U}(0,1)$. If $a^{(m)}>u^{(m)}$, set $q^{(m+1)}=q^*$ (acceptance); otherwise set $q^{(m+1)} = q^{(m)}$ (rejection).

\item Set $m= m+1$. If $m = m_{\max}$ stop; otherwise go to step 1.
\end{enumerate}
\hrule
\caption{HMC algorithm.  The function \(H =(1/2) p^TM^{-1}p+U(q)\) is the Hamiltonian.
 The algorithm generates a Markov chain $q^{(0)}\mapsto q^{(1)} \mapsto \dots \mapsto q^{(m_{\max})}$ reversible with respect to the target probability distribution \(\propto \exp(-U(q))\).}
 \label{table-HMC}
}
\end{table}

One of the most salient features of the algorithm is that it is able to generate proposal moves that, while being far away from the current state of the Markov
chain, has a high probability of acceptance, thus reducing the correlation between samples. 

As the main contribution to the computational cost of HMC resides in the numerical integration of the equations of motion, it is of
paramount interest to use methods requiring as few evaluations of the force $-\nabla_q U(q)$ as possible, with small energy errors (to avoid rejections)
and able to use large step sizes $h$. Taking these considerations  into account, it is hardly surprising that the St\"ormer--Verlet algorithm is the method of choice,
especially in low dimensions. If $d$ increases, however, there are more favorable alternatives. In particular, the  
3-stage method 
\[
 \varphi_{a_1h}^{[T]}\circ
 \varphi_{b_1h}^{[V]}\circ
 \varphi_{a_2h}^{[T]}\circ
 \varphi_{b_2h}^{[V]}\circ
 \varphi_{a_2h}^{[T]}\circ
 \varphi_{b_1h}^{[V]}\circ
 \varphi_{a_1h}^{[T]}
\]
with
\[
  a_1=0.11888010966548, \;\; b_1=0.29619504261126, \;\;
	a_2=\frac12-a_1, \;\; b_2=1-2b_1,
\]
provides better results at the same computational cost as the standard St\"ormer--Verlet method \cite{blanes14nif}. 
The rationale behind this and other multi-stage splitting methods presented in  \cite{blanes14nif} is that the coefficients
$a_j$, $b_j$ are chosen to minimize the energy error in the proposal for relatively large values of $h$ (and not, as in the usual numerical integration
domain, to increase the accuracy in the limit $h \rightarrow 0$). An extension of this idea is the so-called Adaptive Integration Approach (AIA)
\cite{fernandez16ams,akhmatskaya17asi}: here the user
chooses the value of $h$ to be used, and then the AIA algorithm itself finds the coefficients of the method within a given family of $s$-stage splitting integrators
providing the best acceptance rate. The efficiency can be further enhanced by using a conveniently modified version of processing so that time-reversible 
kernels provide time-reversible integrations \cite{blanes21sps}.

In many situations of interest in statistics, the target density is a perturbation of a Gaussian density, so the corresponding Hamiltonian is given by
eq. (\ref{eq:hamHOS}), i.e., 
\[
    H(q,p) = \frac{1}{2} p^T M^{-1} p + \frac{1}{2} q^T N q + U(q),
\]
where $N$ is a constant symmetric, positive definitive $d \times d$ matrix whose spectral radius typically grows with the dimension $d$. 
In that case it is advantageous to precondition the dynamics by choosing the mass matrix (which is free in this setting) as
$M = N$, since now all the $d$ frequencies of the preconditioned system are 1, and the only restriction on the step size used by
St\"ormer--Verlet is $h < 2$, independently of $h$ \cite{bou-rabee18gia}. Furthermore, the considerations on stability exposed in 
subsection \ref{subsec:stability} indicate that
the Strang integrators (\ref{strang-rot}) based on the exact solution of the quadratic part of $H$ and the perturbation $U(q)$ constitute the best option. In fact,
the numerical experiments collected in \cite{casas23shm} 
show that these methods, together with preconditioning, dramatically reduce the computational cost in all
test problems and all observables considered, in comparison with the standard St\"ormer--Verlet scheme. It is worth remarking that all the theory developed in 
Section \ref{sect5} is valid here, and in particular the explicit processor built in subsection \ref{apss} for the Strang splitting is expected to reduce the error in $H$
and thus it might contribute to increasing the acceptance probability.

\subsection{Quantum statistical mechanics}

The description of a quantum system in thermal equilibrium at temperature $T$ with Hamiltonian $H$
is based on the thermal density matrix
\[
   \rho = \e^{-\beta H},
\]
where, again, $\beta = 1/(k_B T)$ \cite{feynman72st}. In fact, most of the properties of the system can be obtained
from $\rho$, but now as an average on the different quantum states. 
Thus, the equilibrium value of an operator $\hat{O}$ corresponding to a physical observable $O$ for
a system of $N$ quantum particles in a volume $V$,
is given by
\begin{equation}   \label{mc.1}
  \langle \hat{O} \rangle = Z^{-1} \, \mathrm{Tr} \left( \e^{-\beta H} \hat{O} \right)=
  Z^{-1}  \sum_n \langle n | \e^{-\beta H} \hat{O} | n \rangle,
\end{equation}
where the partition function $Z$ now reads
\begin{equation}   \label{mc.2}
   Z = \mathrm{Tr} \, (\e^{-\beta H}) = \sum_{n} \langle n | \e^{-\beta H} | n \rangle,
\end{equation}
 and the states 
$|n \rangle$ form a complete, orthonormal basis set \cite{ceperley95pii,landau05agt}. 
Since the eigenvalues of the
Hamiltonian $H$ are not generally known, one tries to evaluate the traces in
(\ref{mc.1}) and (\ref{mc.2}) without diagonalizing the Hamiltonian. This can be done
with the Feynman path-integral approach. To proceed, one considers
the position representation where the particle is labelled. Then the density matrix is given by
$\rho(R, R'; \beta) \equiv \langle R | \e^{-\beta H} | R' \rangle$,
where $R \equiv  \{   r_1, \ldots, r_N \}$, $r_i$ is the
position of the $i$th particle and the elements of $ \rho(R, R'; \beta) $ are positive
and can be interpreted as probabilities. The partition function is then
\begin{equation}   \label{mc.3}
   Z = \int d R \, \langle R | \e^{-\beta H} | R \rangle  =
    \int d R \, \rho(R, R; \beta).
\end{equation}   
Since $\e^{-\beta H} = \left( \e^{-\varepsilon H} \right)^M$,  
with $\varepsilon = \beta /M$ for any positive integer $M$, 
the density matrix can be expressed as
\begin{eqnarray}   \label{mc.4}
  \rho( R_0, R_M; \beta) & = & \int \cdots \int d R_1 d R_2 \cdots
  d R_{M-1} \, \rho(R_0, R_1; \varepsilon) \,  \nonumber \\
    &  &  \qquad  \times \, \, \rho(R_1, R_2; \varepsilon)  \cdots 
    \rho(R_{M-1}, R_M; \varepsilon).
\end{eqnarray}
The action for a given link $k$ is defined as 
\[
S_k \equiv S(R_{k-1},R_k; \varepsilon) = - \ln \big( \rho(R_{k-1},R_k; \varepsilon) \big), 
\]
so
(\ref{mc.4}) becomes
\begin{equation} \label{mc.5}
   \rho( R_0, R_M; \beta)  =  \int \cdots \int d R_1 d R_2 \cdots d R_{M-1} \, \exp \left( - \sum_{k=1}^M S_k \right),
\end{equation}   
and the goal is then to construct a sufficiently accurate approximation
to $\rho$ while minimizing the number of integrals involved in (\ref{mc.5}) (the number
of \textit{beads} $M$).

On the other hand, the
partition function (\ref{mc.3}) can be written as
\begin{eqnarray}   \label{mc.4.b}
     Z  & = &  \int \cdots \int d R_0 d R_1 d R_2 \cdots
  d R_{M-1} \, \rho(R_0, R_1; \varepsilon) \,  \nonumber \\
    &  &  \qquad  \times \, \, \rho(R_1, R_2; \varepsilon)  \cdots 
    \rho(R_{M-1}, R_0; \varepsilon),
\end{eqnarray}
where the first $| R_0 \rangle$ and the last $| R_M \rangle$ elements are
identified as required by the trace operations.  

In practical applications, one must generally use approximations to $\rho$. To this end, notice that, typically
\begin{equation}   \label{ham1}
  H = T + V  = -\frac{1}{2} \sum_{i=1}^N \Delta_i + V.
\end{equation}
 It then makes sense to approximate 
$\e^{-\varepsilon H}$ by symmetric products of $\e^{-\varepsilon T}$ and $\e^{-\varepsilon V}$.
The simplest approximation is given, of course, by the Lie--Trotter scheme, 
known in this setting as the \emph{primitive action}, in which case
\[
  \rho(R_0, R_2; \varepsilon) \approx \int d R_1 \langle
     R_0 | \e^{-\varepsilon T} |  R_1 \rangle \, \langle
     R_1 | \e^{-\varepsilon V} | R_2 \rangle.
\]
The operator $V$ is diagonal in the position representation, whereas the kinetic matrix
can be evaluated by using the eigenfunction expansion of $T$ \cite{ceperley95pii}. It is then possible
to arrive at the discrete path-integral expression for the density matrix 
\begin{eqnarray}   \label{pri-ac.2}
  \rho( R_0, R_M; \beta) & = & \int \cdots \int d R_1 d R_2 \cdots
  d R_{M-1} \,   \left( \frac{1}{2 \pi \varepsilon \hbar^2} \right)^{3N M/2}  \\
  &  & \qquad  \times \, \exp \left( -  \sum_{j=1}^M \left( \frac{1}{2 \varepsilon \hbar^2}
  \| R_{j-1} - R_j\|^2 + \varepsilon V(R_j) \right) \right).  \nonumber
\end{eqnarray}
providing an approximation of order $\mathcal{O}(\varepsilon^2)$, since the Lie--Trotter method is of effective order two. 
A significant improvement with respect to the primitive approximation is achieved by considering the so-called
Takahashi--Imada action \cite{takahashi84mco}, i.e., the approximation
\[
  \e^{-\varepsilon (T + V)} \approx \e^{-\frac{\varepsilon}{2} T} \, \e^{-\varepsilon V - \frac{\varepsilon^3}{24} [V,[T,V]]} \, \e^{-\frac{\varepsilon}{2} T}
\]  
of effective order 4 in $\varepsilon$.

A typical approach to compute the multidimensional integrals appearing in (\ref{mc.3}) and (\ref{mc.4.b}) is to apply a Monte Carlo sampling
according to the probability density $\pi \propto \exp(- \sum_{k=1}^M S_k)$, where $Z$ normalizes $\pi$. 
In this respect, notice that for the Hamiltonian
(\ref{ham1}) one has  $\langle R_{k-1} | \e^{-a_i \varepsilon T} |  R_k \rangle \propto \exp( - \| R_{k-1} - R_k\|^2 / (2 a_i \varepsilon))$, so 
the coefficient $a_j$ cannot be negative for a probabilistic based simulation. Here again the order barrier for splitting methods having
positive coefficients is important: no splitting method of order higher than two can be used for doing quantum statistical calculations, unless nested
commutators enter into their formulation. In this context, the 4th-order scheme (\ref{eq:modpotex}) is widely used, as well as the more general 
2-parameter family of 4th-order methods 
\begin{equation}  \label{pf.3}
  \e^{a_1 \varepsilon T} \e^{\varepsilon  W_{b_1,c_1}} \e^{a_2 \varepsilon T} 
       \e^{\varepsilon  W_{b_2,c_2}} \e^{a_2 \varepsilon T} 
     \e^{\varepsilon  W_{b_1,c_1}}  \e^{a_1 \varepsilon T}, 
\end{equation}
with modified potential
\begin{equation}   \label{eq:modpotphi}
        W_{b_i,c_i} = b_i V + c_i \varepsilon^2 [V,[T,V]].
\end{equation}
In this case, the presence of two free parameters makes it possible to minimize some of the 6th-order error terms and thus yield more efficient
schemes. Numerical simulations 
carried out in \cite{sakkos09hoc}  show that the 
required number $M$ in (\ref{mc.5}) to reproduce the exact energy of the system at low
temperatures is much smaller 
with scheme (\ref{pf.3}) than 
with the Lie--Trotter and Takahashi--Imada methods (see \cite{chin23aop} for a recent review).

\subsection{Vlasov--Poisson equations}

When a gas is brought to a very high temperature, electrons leave the atoms, thus leading to an overall mixture of electrically charged particles (ions and
electrons) usually called \emph{plasma}. There is a hierarchy of models to describe plasmas, ranging from those based on $N$ particles evolving with the
laws of classical relativistic mechanics and forces due to external and self-consistent electromagnetic fields, to kinetic equations and 
fluid models, more appropriate when the
plasma is in thermodynamic equilibrium.

In kinetic models each species $s$ in the plasma is characterized by a distribution function $f_s(x,v,t)$, so that $f_s dx dv$ is the average number of particles
of species $s$ with position and velocity in a box of volume $dx dv$ centered at $(x,v)$. In the limit where the collective effects are dominant over collisions
between particles, the kinetic equation describing the system (in the non-relativistic regime) is the Vlasov equation \cite{vlasov61mpt}
\begin{equation} \label{Vla1}
  \frac{\partial f_s}{\partial t} + v \cdot \nabla_x f_s + \frac{q_s}{m_s} (E + v \times B) \cdot \nabla_v f_s = 0,
\end{equation}  
where $q_s$ and $m_s$ denote the charge and mass of the particles of species $s$, and $E$ and $B$ stand for the electric and magnetic field, respectively.
 Equation (\ref{Vla1}) just expresses the fact that the distribution function
$f_s$ is conserved along the trajectories of the particles,
and is typically coupled with the Maxwell equations
to take into account the self-consistent electromagnetic field generated by the particles \cite{sinitsyn11kbv}. 

The so-called Vlasov--Poisson equation describes a plasma with only one atomic species (alternatively, electrons and positive ions) in the mean
electric field derived from the potential $\phi$ created by the particles. 
Since $m_e \ll m_i$, then the effect of the ions can be treated as a uniform neutralizing background. Taking all constants
equal to 1 and denoting the distribution by $f$, the relevant system of equations describing the electron dynamics reads
\begin{equation} \label{Vla2}
\begin{aligned}
  & \frac{\partial f}{\partial t} + v \cdot \nabla_x f - \nabla_x \phi(f) \cdot \nabla_v f = 0 \\
  & \Delta_x \phi(f)(x) = - \left( \int_{\mathbb{R}^d} f(x,v) dv  - \frac{1}{(2\pi)^d} \int_{\mathbb{T}^d \times \mathbb{R}^d} f(x,v) dx dv \right)
\end{aligned}
\end{equation}
in the domain $(x,v) \in \mathbb{T}^d \times \mathbb{R}^d$, $d=1,2,3$. The energy associated with equations (\ref{Vla2}) is
\begin{eqnarray} \label{Vla3}
  \mathcal{H}(f) & = & \int_{\mathbb{T}^d \times \mathbb{R}^d} \frac{|v|^2}{2} f(x,v) dx dv + \int_{\mathbb{T}^d} \frac{1}{2} |E(f)(x)|^2 dx \\ \nonumber
  & = & \mathcal{T}(f) + \mathcal{U}(f)
\end{eqnarray}  
and is preserved along the solution. In fact, the system has infinitely many invariants, but not all of them can be preserved by numerical schemes. 
Taking into account that the solution of the equations associated with $\mathcal{T}$ and $\mathcal{U}$ can be solved exactly (up to a phase space
discretization), it is then natural to consider splitting methods for the time integration. Specifically, denoting by $\varphi_t^{\mathcal{T}}(f)$ the solution of
$\partial_t f + v \cdot \nabla_x f = 0$, i.e.,
\[
  f(t,x,v) = f(0, x-tv,v),
\]
and by $\varphi_t^{\mathcal{U}}(f)$  the solution 
of $\partial_t f - \nabla_x \phi(f) \cdot \nabla_v f = 0$, which reads
\[
 f(t,x,v) = f(0, x, v -t E(f(0))),
 \]
 where $E(f(0))$ is the value of the electric field at time $t=0$, the numerical integrators are of the form
 \[
   \varphi_{b_{s+1} h}^{\mathcal{U}}(f) \circ \varphi_{a_s h}^{\mathcal{T}}(f) \circ \varphi_{b_{s} h}^{\mathcal{U}}(f) \circ \cdots \circ
    \varphi_{b_{2} h}^{\mathcal{U}}(f) \circ \varphi_{a_1 h}^{\mathcal{T}}(f) \circ \varphi_{b_{1} h}^{\mathcal{U}}(f).
\]
The convergence of these schemes can be established by requiring the appropriate smoothness of $f$ and using the Hamiltonian structure
of the system. Moreover, the functionals $\mathcal{T}$ and $\mathcal{U}$ in the decomposition (\ref{Vla3}) satisfy the formal relation
\[
  [[[\mathcal{T},\mathcal{U}], \mathcal{U}], \mathcal{U}](f) = 0 \qquad \mbox{ for all } \;\; f,
\]
where $[\cdot, \cdot]$ is the Poisson bracket associated with the infinite dimensional Poisson structure of the system (see \cite{casas17hoh}
 for a detailed treatment),
so that Runge--Kutta--Nystr\"om splitting methods can be safely used in this setting. In addition, when $d=1$, 
\[
  [[\mathcal{T},\mathcal{U}], \mathcal{U}](f) =  2 m(f) \mathcal{U}(f), \qquad \mbox{ where } \qquad m(f) = \frac{1}{2\pi} \int_{\mathbb{T} \times \mathbb{R}}
   f(x,v) dx dv.
\]
Since $m(f)$ is a constant of the motion, this introduces additional simplifications. In particular, methods involving nested commutators only require the
evaluation of   $\mathcal{U}$ with appropriate coefficients, and thus methods up to order 6 can be designed involving a reduced number of maps
that, when combined with semi-Lagrangian techniques in phase space, provides high efficiency \cite{casas17hoh}


\subsection{Quantum simulation of quantum systems}

As noted by Feynman in his celebrated contribution \cite{feynman82spw}, simulating the full time evolution of arbitrary quantum systems on a classical computer
requires exponential amounts of computational resources: the states of the system are wave functions that belong to a vector space whose dimension
grows exponentially with the size of the system, so that merely to record the state of the system is already cumbersome. For this reason, he
conjectured the possibility of using a specific quantum system (a quantum computer) to simulate the behavior of arbitrary quantum systems whose
dynamics are determined by local interactions. This conjecture was later shown to be correct \cite{lloyd96uqs}. Today, quantum computers
can simulate a variety of systems arising in quantum chemistry, quantum field theory, many-body physics, etc. (see e.g. 
\cite{berry15shd,childs19nol} and the recent review \cite{miessen23qaf}). 

Simulating the time evolution of a quantum system with Hamiltonian $H$ requires approximating $\e^{-i t H}$ ($\hbar = 1$). If 
the Hamiltonian is the sum over many local interactions, then one can also use splitting methods for this purpose. 
Specifically, suppose that for a system composed of $n$ variables, $H$ can be decomposed as
\begin{equation}  \label{qs1}
   H = \sum_{j=1}^{\ell} H_j,
\end{equation}
where each $H_j$ acts on a space of dimension $m_j$ encompassing at most $k$ of the variables, and $\ell$ is a polynomial in $n$. 
For instance, the well known Hubbard, Ising and Heisenberg models belong to this class \cite{nielsen10qca}. The important point is that, whereas $\e^{-i t H}$
is difficult to compute, each $\e^{-i t H_j}$ acts on a much smaller subsystem and can be straightforward to evaluate by quantum circuits. In fact,
the explicit quantum simulation algorithm proposed in \cite{lloyd96uqs} is based on the Lie--Trotter method
\[
  \chi_t = \e^{-i t H_{\ell}} \cdots \e^{-i t H_2} \, \e^{-i t H_1},
\]
whereas subsequent proposals include the Strang splitting
\[
  S_t^{[2]} =   \e^{-i \frac{t}{2} H_1} \, \e^{-i \frac{t}{2} H_2} \, \cdots \e^{-i t H_{\ell}} \, \cdots \e^{-i \frac{t}{2} H_2} \, \e^{-i \frac{t}{2} H_1},
\]
and especially the quintuple jump recursion (\ref{eq:quintuple_jump})
\begin{equation} \label{suzu-qs}
  S_t^{[2k+2]} = \left( S_{\gamma_{2k} t}^{[2k]} \right)^2 \circ S_{(1- 4 \gamma_{2k}) t}^{[2k]} \circ
     \left( S_{\gamma_{2k} t}^{[2k]} \right)^2
\end{equation}
with $\gamma_{2k} = 1/(4 - 4^{1/(2k+1)})$.
These
schemes are known in this setting as \emph{product formulas} \cite{childs19nol,chen22epf} and the procedure is called \emph{Trotterization}. 

Product formulas provide approximations $U_{\mathrm{app}}$ to the exact evolution $\e^{-i t H}$ and the goal is, given a time $t$ and a maximum simulation
error $\varepsilon$, to find an algorithm (a quantum circuit) $U_{\mathrm{app}}$ such that $\| U_{\mathrm{app}} - \e^{-i t H} \| < \varepsilon$. In consequence,
it is of paramount importance to analyze and eventually provide tight bounds for the error $\varepsilon$ committed by product formulas / splitting methods
when applied to Hamiltonian systems of the form (\ref{qs1}). Thus, in \cite{childs21tot}, it is shown that 
\[
  \big\| S_t^{[2k]} - \e^{-i t H} \big\| = \mathcal{O} \left( (B_H)^{2k+1} \right), \qquad \mbox{ with } \qquad B_H \equiv \sum_{j=1}^{\ell} \|H_j\| \, t 
\]
if $H_j$ are Hermitian.  As usual, if $t$ is large, then the whole interval is divided into $N$ steps and (\ref{suzu-qs}) is applied within each step. In that case,
$\big\| \left(S_{t/N}^{[2k]} \right)^N - \e^{-i t H} \big\| = \mathcal{O}(\varepsilon)$ provided that
\[
  N = \mathcal{O} \left(   \frac{(B_H)^{1+1/(2k)}}{\varepsilon^{1/(2k)}}   \right).
\]  
As we have already pointed out, whereas the recursion (\ref{suzu-qs}) (as well as the triple jump) constitutes a systematic way to achieve high order approximations, it is not necessarily the most efficient, both in terms of errors and the number of exponentials involved. It makes sense, then, to consider
some other methods in this setting, 
such as those collected in Section \ref{sect8}. This requires, in particular, a detailed analysis to achieve more stringent bounds for the corresponding 
errors than those obtained in \cite{childs21tot}.

Hamiltonian systems of the form $H = -\frac{1}{2} \Delta + V(x)$
can also be simulated on quantum computers by using an appropriate representation of the states and the previous product formulas with the quantum
Fourier transform \cite{nielsen10qca}.

\subsection{Other topics}

Space and time constraints prevent us from including in this review additional relevant applications where splitting methods have shown their merits, as well as
other closely related important issues. Let us briefly mention some of them.

\begin{itemize}

 \item Except for a brief incursion into Langevin dynamics, we have restricted our treatment to deterministic problems, although splitting methods have
been widely applied to both ordinary and partial stochastic differential equations. Recent references in this area 
include \cite{brehier23sif,brehier23ssf,foster22hos}. 

\item The evolution of a particle of mass $m$ and charge $q$ in a given electromagnetic field is modeled by the Lorentz equation 
$m x^{\prime\prime} = q( E + x^{\prime} \times B)$. By introducing the velocity $v = x^{\prime}$ as a new variable, the resulting system of first-order
ODEs can be split into three explicitly solvable parts, so that composition methods can be applied preserving volume in phase space
\cite{he15vpa,he16hov}. The treatment can also be generalized to relativistic charged particles \cite{zhang15vpa}.
 

 \item In addition to the Schr\"odinger equation, both linear and nonlinear, splitting methods have been applied to 
 other relevant partial differential equations arising in quantum physics. This is the case, in particular, for the Dirac equation 
 \cite{bao20sro} and the Klein--Gordon equation
 \cite{bao22iue}. 
 
  \item One important theoretical aspect not treated here refers to the convergence analysis of splitting methods applied to the nonlinear
 Schr\"odinger equation and other semi-linear Hamiltonian PDEs, as well as the use of the Birkhoff normal forms. This fascinating topic 
 is the subject of much attention in the recent literature (see e.g. \cite{faou12gni,faou11hio,bernier22bnf}) and probably deserves a review in its own right.  It might be of interest to explore whether the formalism of extended word series, successfully applied  in~\cite{murua16cnf} to construct formal invariants and normal forms of general classes of finite dimensional Hamiltonian systems could be adapted and applied for semi-linear Hamiltonian PDEs.   

\item Here we have only considered the application of splitting and composition methods with constant step size $h$. This is essential
 if one is interested in preserving qualitative properties of the system. As we have seen in section \ref{sect4}, the
 modified equation corresponding to the numerical scheme depends explicitly on $h$, so that
 if $h$ is changed then so is the modified equation, and 
the preserving properties that geometric integrators possess when $h$ is constant are no longer guaranteed. 
In other problems, such as those defined by PDEs with no particular structure, this is not a matter of concern. There are systems which
 do possess a geometric structure that is advantageous to preserve by the numerical scheme and where the use of an adaptive step size is
 of the utmost importance to get efficient approximations. A relevant example is the gravitational $N$-body problem when there are close encounters
 between some of the bodies. In that case, one may apply splitting methods with variable step size by using some specifically
 designed transformations involving the time variable, in such a way that in the new variables
 the resulting time step is constant (see, e.g., 
\cite{mikkola97psm,calvo98vsi,blanes05agi,blanes12eas} and \cite[sect. VIII.2]{hairer06gni}).

\end{itemize}

\begin{center}
\textit{Forse altro canter\`a con miglior plectro.}
\end{center}

\


\subsection*{Acknowledgments}
The authors wish to thank
J. Bernier, J. Boronat, A. Iserles, A. P\'erez, J.M. Sanz-Serna, L. Shaw, and M. Thalhammer for reading parts of this manuscript and for their comments and suggestions. They would also like to express their gratitude to their long-suffering families for their patience, understanding and wholehearted support.
This work has been
funded by Ministerio de Ciencia e Innovaci\'on (Spain) through projects PID2022-136585NB-C21 and PID2022-136585NB-C22, \\
MCIN/AEI/10.13039/501100011033/FEDER, UE, and also by Generalitat Valenciana (Spain) through project CIAICO/2021/180. AM has also received funding from the
Department of Education of the Basque Government through the Consolidated Research Group MATHMODE (ITI456-22).

\

\appendix

\section{Testing splitting methods on matrices}
\label{apenda}

In Section~\ref{sect8} we collected more than 80 different splitting and composition methods that are available in the literature. 
Although they have been classified into different categories according to the structure of the system to be integrated (RKN splitting methods, near-integrable
systems, etc.) and the type of scheme (composition of a basic 2nd-order time-symmetric method, splitting for a system separated into two parts, etc.), it is not
obvious in advance which particular method is the most suitable for carrying out a given integration of a specific problem. There are many factors
involved in the final choice: type of problem, the qualitative properties one is interested in preserving, computational cost, desired accuracy, time integration 
interval and even 
initial condition.

For these reasons, it might be illustrative to fix a list of examples and apply to them the most representative numerical integrators gathered in Section
\ref{sect8}, just to get some clues about their relative performances. To this end, we take the linear system 
\begin{equation} \label{apen1.1}
    \frac{dX}{dt} = F \, X, \qquad X(0)=I, 
\end{equation}
where $F$ is a constant real  $d \times d$ matrix and $I$ is the identity, so that the exact solution $X(t) = \e^{t F}$ can be easily
computed with any desired accuracy. Obviously, none of the methods we test here should be used, in general, to approximate the exponential of a matrix or its action on a vector (unless additional information on the structure of the matrices is known which could make splitting an efficient technique). Our purpose 
 is rather
to compare the relative performance of the different schemes.

Let $\Phi_h^N$, with $t_f=Nh$, denote the matrix obtained with a given method applied $N$ times with step size $h$ that approximates $\e^{t F}$
at $t=t_f$. We then compute the following relative errors:
\[
 {\cal E}_1=\frac{\|\e^{tF}-\Phi_h^N\|}{\|\e^{tF}\|}, \qquad\quad
 {\cal E}_2=\frac{|\tr(\e^{tF})-\tr(\Phi_h^N)|}{|\tr(\e^{tF})|}.
\] 
The first, ${\cal E}_1$, measures the accuracy of the method, whereas ${\cal E}_2$ provides an estimate of the accuracy of the scheme in the case where it is 
used as a kernel with an ideal processor. In all cases, $\| \cdot \|$ refers to the 2-norm of the matrices considered.

The numerical experiments reported here have been carried out with MATLAB, and the function \texttt{randn( )} has been used to construct the matrix elements 
randomly 
 from a normal distribution, initiated with the seed \texttt{rng(1)} for reproducibility. Additional material related with these simulations can be found at our
 website accompanying this paper.

\paragraph{Symmetric compositions of time-symmetric 2nd-order schemes.}

\begin{figure}[!htb]
\includegraphics[width=1.1\textwidth]{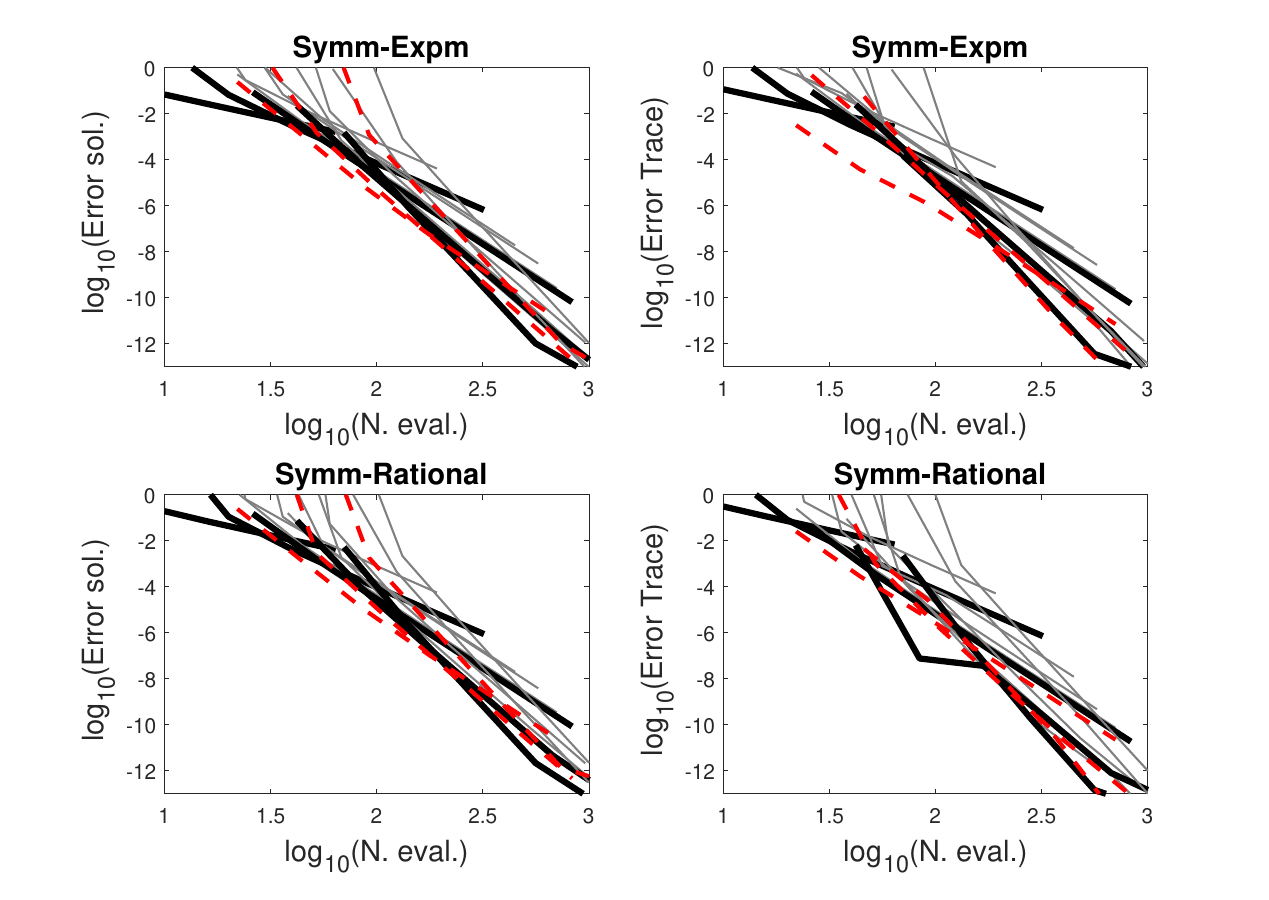}
\caption{Errors ${\cal E}_1$ (left) and ${\cal E}_2$ (right) vs. the number of evaluations of the basic scheme for $S_h^{[2,1]}$  (top) and $S_h^{[2,2]}$ (bottom).
The non-processed schemes selected in Table~\ref{tableSS} are showed in thick solid lines and the selected processed methods are drawn as thick red dashed lines.}
\label{figSS}
\end{figure}

We take $F=A+B+C$, with $A,B,C$ matrices of dimension $50 \times 50$ constructed as previously indicated, i.e., 
\begin{verbatim}
d=50; rng(1);  A=randn(d);   A=A/norm(A); B=randn(d); ...
\end{verbatim}
To illustrate the role that  the basic scheme may play in the overall performance of the composition, we take the following two 2nd-order time-symmetric methods:
\begin{equation} \label{ch.1}
  S_h^{[2,1]}=	\e^{hA/2}	\e^{hB/2}	\e^{hC}	\e^{hB/2}	\e^{hA/2}
\end{equation}
and 
\begin{equation} \label{ch.2}
  S_h^{[2,2]}=	(I-\frac{h}2A)^{-1}(I-\frac{h}2B)^{-1}(I-\frac{h}2C)^{-1}
	      (I+\frac{h}2C)(I+\frac{h}2B)(I+\frac{h}2A).
\end{equation}
The two differ in their computational cost and error terms. For the 
most relevant methods presented in Table~\ref{tableSS}, we plot the errors ${\cal E}_1$ and ${\cal E}_2$ as a function of the total number of evaluations of the basic
scheme. Figure \ref{figSS} shows the results obtained for $t_f=10$ when the basic scheme is $S_h^{[2,1]}$ (top) and $S_h^{[2,2]}$ (bottom).
Compositions with $s=5,13,19$ and $35$ stages of order $r=4, 6, 8$ and $10$, respectively, which on average show good performance in most tested examples, correspond to the thick black lines, whereas the remaining non-processed methods are shown as thin lines. On the other hand, the recommended
processed methods in Table \ref{tableSS} with kernels having $s=13, 19$ and $23$ stages of order $6, 8$ and 10, respectively, are represented by thick 
red dashed lines. 

Looking at these figures, it is hardly possible to recommend a particular scheme as the one leading to the best results for the example considered. 
We may only conclude that, for some intervals of accuracy, there are certain methods exhibiting the best performance.. Moreover, this performance also
depends on the particular basic scheme chosen as $S_h^{[2]}$ in the compositions and the nature of the composition itself: whether it is intended
to be used as the kernel of a processed method (with an optimal processor). In many situations, the relative performance 
displayed in the left graphs should be close to that of the right panels, when long-time integrations are considered and the contribution from the processable error terms can be neglected.

\paragraph{Splitting into two parts / composition of a basic method and its adjoint.}

\begin{figure}[!htb]
\centering
\includegraphics[width=1.1\textwidth]{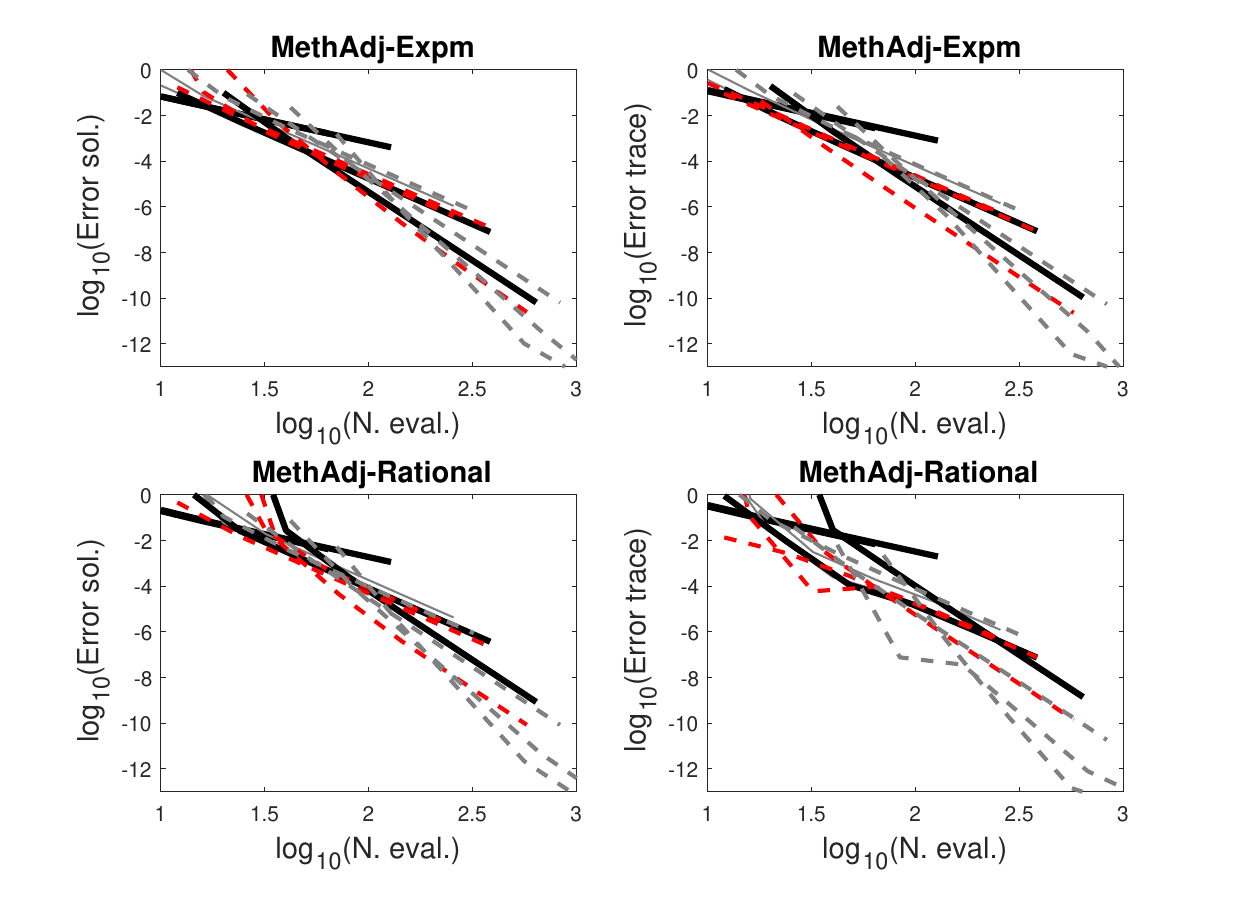}
\caption{Errors ${\cal E}_1$ (left) and ${\cal E}_2$ (right) vs. the number of evaluations of the basic scheme for $\chi_h^{[1,1]}$  (top) and $\chi_h^{[1,2]}$ (bottom). Now the schemes are compositions of $\chi_h^{[1,j]}$ and $(\chi_h^{[1,j]})^*$.}
\label{figMetAdj}
\end{figure}

We again consider $F=A+B+C$, with $A,B,C$ the same matrices as in the previous example, but instead we take as basic methods the first-order
compositions
\[
	 \chi_h^{[1,1]}=	\e^{hA}	\e^{hB}	\e^{hC}, \qquad
	 \big(\chi_h^{[1,1]} \big)^*=	\e^{hC}	\e^{hB}	\e^{hA},
	\]
and
	\[
	 \chi_h^{[1,2]}=(I-hA)^{-1}(I-hB)^{-1}(I-hC)^{-1},  \qquad
	 \big(\chi_h^{[1,2]} \big)^*=	(I+hC)(I+hB)(I+hA).
	\]
We plot the errors $\mathcal{E}_1$ and $\mathcal{E}_2$ versus the total number of evaluations of the basic scheme for most of the methods presented in Table~\ref{tableAB}. Figure \ref{figMetAdj} shows the results obtained for $t_f=10$. As before, the thick black solid lines correspond to the selected non-processed schemes and the thick red dashed lines correspond to the selected processed ones. The dashed grey lines correspond to the most efficient non-processed symmetric--symmetric schemes of orders 6, 8 and 10. 

\begin{figure}[htb]
\centering
\includegraphics[width=1.1\textwidth]{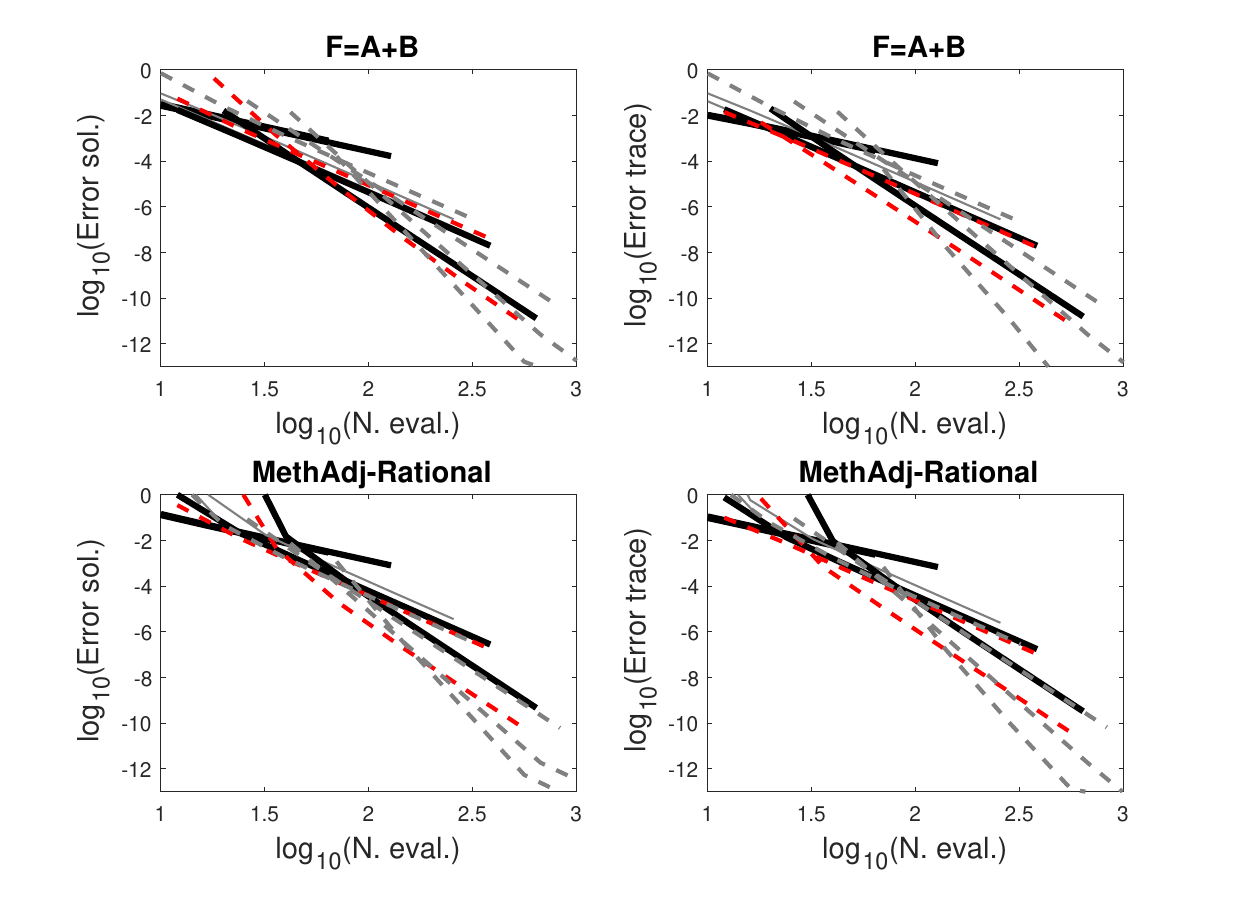}
\caption{The same as in Figure \ref{figMetAdj}, but now for the splitting $F=A+B$.}
\label{figAB}
\end{figure}

We repeat the same numerical experiments, but now taking $C=0$, while keeping $A,B$ the same matrices as previously. Hence, the problem has to be seen as separable into two parts. Notice that some methods have been optimized for this particular case. Figure \ref{figAB} shows the results obtained, where the superiority of the methods for this class of problems is clear for low-to-medium accuracy.

\paragraph{Runge--Kutta--Nystr\"om methods.}

To test the RKN splitting methods gathered in Tables \ref{tableRKN} and \ref{tableRKNm}, we take 
$F=A+B$, with 
	\[
	 A=\left( 
	\begin{array}{cc} 
	O_d & O_d \\ A_1 & O_d
	\end{array} \right), \qquad
	 B=\left( 
	\begin{array}{cc} 
	B_1 & B_2 \\ B_3 & B_4
	\end{array} \right), 
	\]
where $A_1,B_i, \ i=1,2,3,4$ are matrices of dimension $d=50$ with elements chosen as in the previous cases, with the same seed, and $O_d$ is the null matrix. 
With this choice we have
	\[
	 [A,[A,B]]=\left( 
	\begin{array}{cc} 
	O_d & O_d \\ 2A_1B_2A_1 & O_d
	\end{array} \right)
	\]
	and $ [A,[A,[A,B]]]=0.$ For this particular choice, the computation of $\e^{hB}$ dominates the total cost of the method and the cost of evaluating $\e^{hA}$ and $\e^{hA+h^3[A,[A,B]]}$ can be neglected. Obviously, we can add an artificial cost to these exponentials as a test for different problems where this term can be more expensive to evaluate.
	
\begin{figure}[!htb]
\centering
\includegraphics[width=1.05\textwidth]{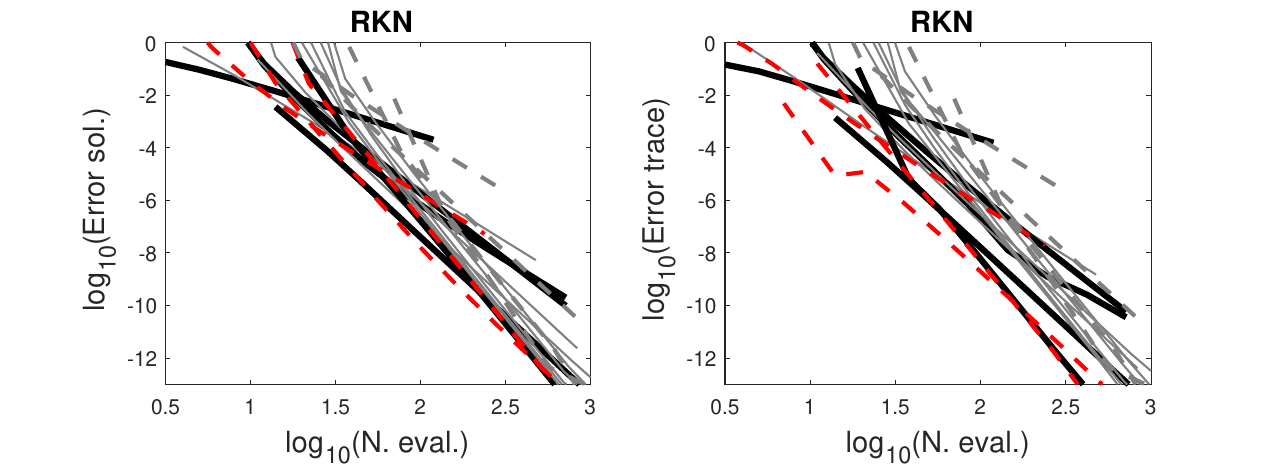}
\includegraphics[width=1.05\textwidth]{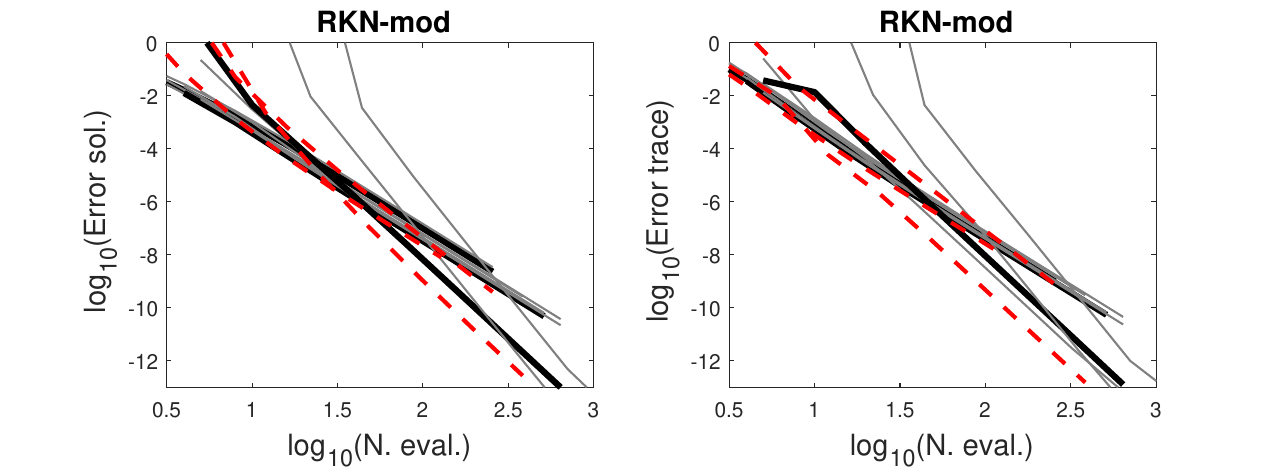}
\caption{RKN methods without (top) and with (bottom) modified potentials.}
\label{figRKN}
\end{figure}

	Figure~\ref{figRKN} (top) shows the results for most of the methods from Table~\ref{tableRKN} for $t_f=10$ where, as previously, the thick black lines correspond to the selected non-processed schemes, the thin lines correspond to the remaining non-processed schemes, the thick red dashed lines are the processed ones and the 
	dashed grey lines correspond to the selected high order symmetric--symmetric methods. Bottom graphs show,
for the same matrices, the results for the methods from Table~\ref{tableRKNm} with modified potentials.

\paragraph{Methods for near-integrable systems.}

\begin{figure}[!htb]
\centering
\includegraphics[width=1.1\textwidth]{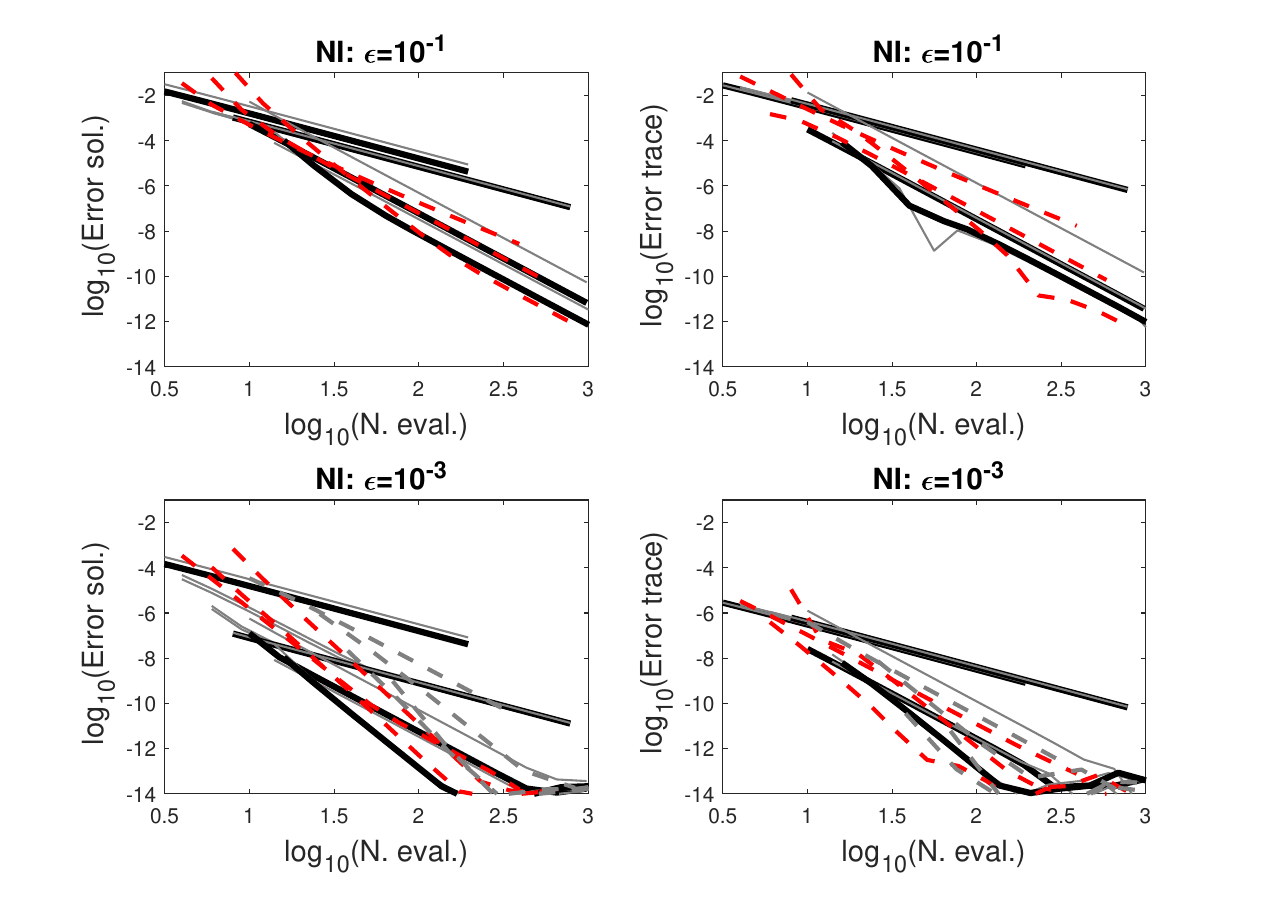}
\caption{Methods for near-integrable systems.}
\label{figNI}
\end{figure}

Let $F=A+\varepsilon B$ with $A,B$ the same matrices as for the separable problem in two parts and $\varepsilon$ a small parameter which corresponds to the relative norm of the matrices. We analyze the performance of methods tailored for perturbed problems for two choices of the small parameter: $\varepsilon=10^{-1}$ and $\varepsilon=10^{-3}$.	Figure~\ref{figNI} shows the results for the methods from Table~\ref{tableNI} for $t_f=10$. For small values of $\varepsilon$ none of the previous splitting methods are competitive against the most efficient ones from this family.

\addcontentsline{toc}{section}{References}
\bibliographystyle{actaagsm}

\

\label{lastpage}
\end{document}